\documentstyle{book}

\makeindex

\def\prop#1#2{\vspace{2ex} \noindent\textsc{#1.} \emph{#2}\par\vspace{2ex}}
\def\dkz{\noindent\textsc{Proof. }}
\def\qed{\hfill $\dashv$\vspace{2ex}}

\def\pl{\!+\!}
\def\mn{\!-\!}

\def\cirk{\,{\raisebox{.3ex}{\tiny $\circ$}}\,}

\def\mj{\textbf{1}}

\def\HDS{\vrule width0pt height2.3ex depth1.05ex\displaystyle}
\def\f#1#2{{{\HDS #1}\over{\HDS #2}}}
\def\afrac#1{{\phantom{\HDS #1}\atop{\HDS #1}}}
\def\lpravilo#1{ \makebox[-.5em][r]{\mbox{#1}} {\mbox{\hspace{1.5em}}}}

\begin{document}

\begin{titlepage}
\begin{center}

{\huge Graphs for Juncture}

\vspace{10ex}

{\sc Kosta Do\v sen} and {\sc Zoran Petri\' c}

\vspace{90ex}

{August 2012}
\end{center}
\end{titlepage}

\clearpage \pagestyle{empty} \makebox[1em]{} \clearpage

\pagestyle{myheadings}\markboth{}{}

\pagenumbering{roman} \thispagestyle{plain} \setcounter{page}{5}

\noindent \textbf{\textit{Abstract}}\label{Abstract}

\vspace{2ex}

\noindent An alternative foundation for 2-categories is explored by studying graph-theoreti\-cal\-ly a partial operation on 2-cells named juncture, which can replace vertical and horizontal composition. Juncture is a generalized vertical composition of 2-cells that need not involve the whole target and the whole source; it may involve them only partly, provided the result is again a 2-cell. Since commuting diagrams of arrows of ordinary categories may be conceived as invertible 2-cells, this study concerns ordinary category theory too. The operation of juncture has a connection with proof theory, where it corresponds to a kind of cut rule on sequents, and it is related also to an operation on which the notion of operad can be based. The main achievement of the work is a detailed description of the specific planarity involved in juncture and graphs of 2-cells, comparable to the usual combinatorial characterizations of planarity in graph theory. This work points out to an alternative foundation for bicategories, i.e.\ weak 2-categories, and more generally weak $n$-categories.

\vspace{6ex}

\noindent \textbf{\textit{Subject Classification}} (MSC2010)\label{Subject Classification}

\vspace{2ex}

\noindent{\sc Graph Theory}:

\vspace{1ex}

\noindent 05C10 (Planar graphs), 05C20 (Directed graphs), 05C62 (Graph representations), 05C76 (Graph operations)

\vspace{2ex}

\noindent{\sc Category Theory}:

\vspace{1ex}

\noindent 18A10 (Graphs, diagram schemes), 18D05 (2-categories)

\vspace{6ex}

\noindent
\textbf{\textit{Acknowledgements}}\label{Acknowledgements}

\vspace{2ex}

\noindent This work was supported by a project of the Ministry of
Science of Serbia. We are grateful to John Power for encouraging
comments.

\baselineskip=1.2\baselineskip

\newpage

\clearpage \pagestyle{empty} \makebox[1em]{} \clearpage

\pagestyle{myheadings}\markboth{Contents}{Contents}

\thispagestyle{plain}

\begin{center}
{\bf \large CONTENTS}
\end{center}

\vspace{2ex}

\begin{list}
\labelwidth{} \item[]
\emph{Abstract}\hfill\pageref{Abstract}

\vspace{-9pt}

\item[]
\emph{Subject Classification}\hfill\pageref{Subject Classification}

\vspace{-9pt}

\item[] \emph{Acknowledgements}\hfill\pageref{Acknowledgements}

\item[{\sc Chapter} 1.] {\bf Introduction}\hfill\pageref{1}

\vspace{-9pt}

\begin{list}
\labelwidth{} \item[\S 1.1.] Aim and scope\hfill\pageref{1.1}

\vspace{-5pt}

\item[\S 1.2.] D-graphs\hfill\pageref{1.2}

\vspace{-5pt}

\item[\S 1.3.] Cocycles and juncture\hfill\pageref{1.3}

\vspace{-5pt}

\item[\S 1.4.] Edge-graphs\hfill\pageref{1.4}

\vspace{-5pt}

\item[\S 1.5.] The system S$\Box$\hfill\pageref{1.5}

\vspace{-5pt}

\item[\S 1.6.] The completeness of S$\Box$\hfill\pageref{1.6}

\vspace{-5pt}

\item[\S 1.7.] Compatible lists\hfill\pageref{1.7}

\vspace{-5pt}

\item[\S 1.8.] P$'$-graphs\hfill\pageref{1.8}

\vspace{-5pt}

\item[\S 1.9.]  P$''$-graphs\hfill\pageref{1.9}

\vspace{-5pt}

\item[\S 1.10.]  P$'''$-graphs\hfill\pageref{1.10}

\end{list}

\vspace{-1.1ex}

\item[{\sc Chapter} 2.] {\bf P$'$-Graphs and P$'''$-Graphs}\hfill\pageref{2}

\vspace{-9pt}

\begin{list}
\labelwidth{} \item[\S 2.1.] Interlacing and parallelism
\hfill\pageref{2.1}

\vspace{-5pt}

\item[\S 2.2.] P$'$-graphs and grounding\hfill\pageref{2.2}

\vspace{-5pt}

\item[\S 2.3.] P$'$-graphs are P$'''$-graphs \hfill\pageref{2.3}

\end{list}

\vspace{-1.1ex}

\item[{\sc Chapter} 3.] {\bf Grounding and
Pivots}\hfill\pageref{3}

\vspace{-9pt}

\begin{list}
\labelwidth{} \item[\S 3.1.] Grounding and juncture
\hfill\pageref{3.1}

\vspace{-5pt}

\item[\S 3.2.] Pivots and their ordering \hfill\pageref{3.2}

\vspace{-5pt}

\item[\S 3.3.] Further results for the Pivot Theorem
\hfill\pageref{3.3}

\vspace{-5pt}

\item[\S 3.4.] The Pivot Theorem \hfill\pageref{3.4}

\end{list}

\vspace{-1.1ex}

\item[{\sc Chapter} 4.] {\bf P$''$-graphs and
P$'$-graphs}\hfill\pageref{4}

\vspace{-9pt}

\begin{list}
\labelwidth{} \item[\S 4.1.] Petals \hfill\pageref{4.1}

\vspace{-5pt}

\item[\S 4.2.] P-moves \hfill\pageref{4.2}

\vspace{-5pt}

\item[\S 4.3.] Completeness of P-moves \hfill\pageref{4.3}

\vspace{-5pt}

\item[\S 4.4.] P$''$-graphs are P$'$-graphs \hfill\pageref{4.4}

\end{list}

\vspace{-1.1ex}

\begin{samepage}

\item[{\sc Chapter} 5.] {\bf P$'''$-graphs and
P$''$-graphs}\hfill\pageref{5}

\vspace{-9pt}

\begin{list}
\labelwidth{} \item[\S 5.1.] B$_m$-moves \hfill\pageref{5.1}

\vspace{-5pt}

\item[\S 5.2.] Completeness of B$_m$-moves \hfill\pageref{5.2}

\vspace{-5pt}

\item[\S 5.3.] P$'''$-graphs are P$''$-graphs \hfill\pageref{5.3}

\end{list}

\vspace{-1.1ex}

\end{samepage}

\item[{\sc Chapter} 6.] {\bf The Systems {\rm S1} and {\rm S2}}\hfill\pageref{6}

\vspace{-9pt}

\begin{list}
\labelwidth{} \item[\S 6.1.] The system S$\Box_P$\hfill\pageref{6.1}

\vspace{-5pt}

\item[\S 6.2.] The system S1\hfill\pageref{6.2}

\vspace{-5pt}

\item[\S 6.3.] The system S2\hfill\pageref{6.3}

\vspace{-5pt}

\item[\S 6.4.] The equivalence of S1 and S2\hfill\pageref{6.4}

\vspace{-5pt}

\item[\S 6.5.] The completeness of S1\hfill\pageref{6.5}

\vspace{-5pt}

\item[\S 6.6.] M-graphs\hfill\pageref{6.6}

\vspace{-5pt}

\item[\S 6.7.] The completeness of S2\hfill\pageref{6.7}
\end{list}

\vspace{-1.1ex}

\item[{\sc Chapter} 7.] {\bf Disk D-Graphs and P-Graphs}\hfill\pageref{7}

\vspace{-9pt}

\begin{list}
\labelwidth{} \item[\S 7.1.] Disk D-graphs\hfill\pageref{7.1}

\vspace{-5pt}

\item[\S 7.2.] P-graphs are realizable as disk D-graphs\hfill\pageref{7.2}

\vspace{-5pt}

\item[\S 7.3.] Disk D-graphs are P-graphs\hfill\pageref{7.3}

\vspace{-5pt}

\item[\S 7.4.] D1$'$-graphs\hfill\pageref{7.4}

\vspace{-5pt}

\item[\S 7.5.] Realizing D1$'$-graphs\hfill\pageref{7.5}

\vspace{-5pt}

\item[\S 7.6.] Duality\hfill\pageref{7.6}
\end{list}

\vspace{-1.1ex}

\item[] {\bf Bibliography}\hfill\pageref{Bibliography}

\vspace{-2ex}

\item[] {\bf Index}\hfill\pageref{Index}

\end{list}

\newpage

\setcounter{chapter}{0} \setcounter{page}{0}
\pagenumbering{arabic}

\baselineskip=0.9\baselineskip

\chapter{\huge\bf Introduction}\label{1}
\pagestyle{myheadings}\markboth{CHAPTER 1. \quad
INTRODUCTION}{right-head}

\section{\large\bf Aim and scope}\label{1.1}
\markright{\S 1.1. \quad Aim and scope}

Our aim is to explore a matter related to an alternative
foundation for 2-categories (see \cite{KS74} and \cite{ML98},
Section XII.3, for the standard notion of 2-category). We study
graph-theoreti\-cal\-ly a partial operation on 2-cells we call
\emph{juncture}\index{juncture}, which can replace vertical and
horizontal composition. Juncture corresponds to the gluing of
diagrams of 2-cells along their borders so that the result is
again a diagram of 2-cells. We do not study everything needed for a theory of 2-categories,
but only these matters related to horizontal and vertical
compositions.

Commuting diagrams of arrows of ordinary categories may be conceived as invertible 2-cells, and the gluing of such commuting diagrams along their borders so as to make other such commuting diagrams is what juncture is about. So our study of juncture concerns ordinary category theory too. It is a contribution to the theory of diagrams of ordinary categories (see the end of \S 7.5 for some further comments concerning that matter).
The operation of juncture, for which we will try to show that it is worth studying from the point of view of graph theory, has also a connection with proof theory, where it corresponds to a kind of cut rule on sequents (see the last paragraph of this section).

Juncture, which is definable in terms of vertical composition, horizontal composition
and identity 2-cells, permits us to define vertical composition,
but with its help we can define also horizontal composition only in
the presence of identity 2-cells (see Chapter~6). Juncture is a generalized vertical
composition of 2-cells, where the target of the first 2-cell
may coincide only partly with the source of the second 2-cell, provided the result is again a 2-cell, as, for example,~in
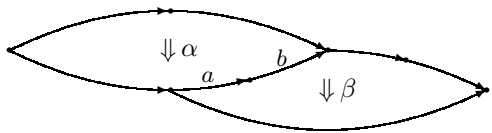
\begin{center}
\begin{picture}(180,50)(0,15)

\put(0,45){\circle*{2}} \put(61,30){\circle*{2}}
\put(61,60){\circle*{2}} \put(91,34){\circle*{2}}
\put(120.5,45){\circle*{2}}

\put(150,41.3){\circle*{2}} \put(180.5,30){\circle*{2}}

\qbezier(0,45)(60,15)(120,45) \qbezier(0,45)(60,75)(120,45)
\qbezier(60,30)(120,0)(180,30) \qbezier(120,45)(150,45)(180,30)

\put(58,30){\vector(1,0){2}} \put(88,33.5){\vector(4,1){2}}
\put(118,44){\vector(2,1){2}} \put(118,46){\vector(2,-1){2}}
\put(58,60){\vector(1,0){2}} \put(148,41.8){\vector(4,-1){2}}
\put(178.5,31){\vector(2,-1){2}} \put(178.5,29.2){\vector(2,1){2}}

\put(60,45){\makebox(0,0){$\Downarrow$}}
\put(65,45){\makebox(0,0)[l]{$\alpha$}}

\put(120,30){\makebox(0,0){$\Downarrow$}}
\put(125,30){\makebox(0,0)[l]{$\beta$}}

\put(75,33){\small\makebox(0,0)[b]{$a$}}
\put(103,39){\small\makebox(0,0)[b]{$b$}}

\end{picture}
\end{center}
The associativity of vertical and horizontal composition
and the \emph{intermuting} of these two operations (see the
equation $(\otimes\cirk)$ in \S 6.3) are now replaced by two kinds
of associativity of juncture (see the equations of $S\Box$ in \S
1.5), for which we will prove completeness (see \S 1.6 and Chapter~6).

An essential ingredient of juncture is that its correct
application, where the result is a 2-cell, is based on conditions that are respected in the
example above, but are violated in, for example,
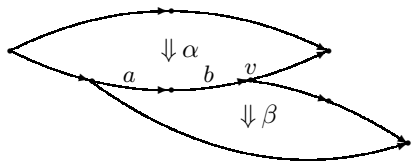
\begin{center}
\begin{picture}(150,60)(0,5)

\put(0,45){\circle*{2}} \put(61,30){\circle*{2}}
\put(61,60){\circle*{2}} \put(91,34){\circle*{2}}
\put(120.5,45){\circle*{2}} \put(31,33.5){\circle*{2}}

\put(120.5,26){\circle*{2}} \put(150.5,10){\circle*{2}}

\qbezier(0,45)(60,15)(120,45) \qbezier(0,45)(60,75)(120,45)

\qbezier(29.5,34)(90,-10)(150,10)
\qbezier(90,33.7)(120,30)(150,10)

\put(58,30){\vector(1,0){2}} \put(88,33.5){\vector(4,1){2}}
\put(28,34.5){\vector(4,-1){2}} \put(118,44){\vector(2,1){2}}
\put(118,46){\vector(2,-1){2}} \put(58,60){\vector(1,0){2}}

\put(118.5,26.7){\vector(3,-1){2}}
\put(148.5,11){\vector(3,-2){2}} \put(148.5,9.3){\vector(3,1){2}}

\put(60,45){\makebox(0,0){$\Downarrow$}}
\put(65,45){\makebox(0,0)[l]{$\alpha$}}

\put(90,20){\makebox(0,0){$\Downarrow$}}
\put(95,20){\makebox(0,0)[l]{$\beta$}}

\put(45,33){\small\makebox(0,0)[b]{$a$}}
\put(75,33){\small\makebox(0,0)[b]{$b$}}
\put(91,36){\small\makebox(0,0)[b]{$v$}}

\end{picture}
\end{center}
In this last picture the bifurcation at the point $v$
makes it impossible to say what 1-cells are sources and targets of
the result, and so the result is not a 2-cell.

We will formulate these conditions by treating diagrams of 2-cells
as planar graphs of a specific sort, and then by considering dual graphs
of these graphs. By this, and by further modifying the dual graphs
(see \S 7.6 for details), the juncture in the first example above
becomes an operation, which we call juncture too, that transforms
the two graphs on the left into the graph on the right
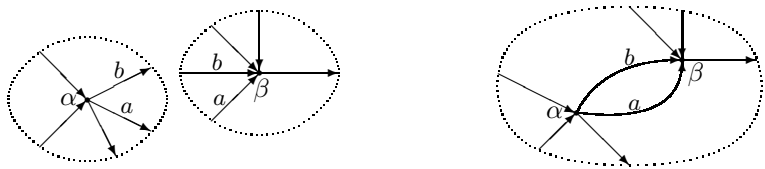
\begin{center}
\begin{picture}(320,60)(0,-5)

\put(45,20){\circle*{2}}

\qbezier[30](30,0)(0,20)(30,40) \qbezier[15](30,40)(45,47)(60,40)
\qbezier[15](30,0)(45,-7)(60,0) \qbezier[30](60,0)(90,20)(60,40)

\put(45,20){\line(-1,1){18}} \put(42.5,17.5){\vector(1,1){2}}
\put(45,20){\line(-1,-1){18}} \put(42.5,22.5){\vector(1,-1){2}}

\put(45,20){\vector(2,1){24}} \put(45,20){\vector(2,-1){24}}
\put(45,20){\vector(1,-2){10.7}}

\put(41,20){\makebox(0,0)[r]{$\alpha$}}
\put(60,15){\small\makebox(0,0)[b]{$a$}}
\put(57,28){\small\makebox(0,0)[b]{$b$}}

\put(110,30){\circle*{2}}

\qbezier[30](95,10)(65,30)(95,50)
\qbezier[15](95,50)(110,57)(125,50)
\qbezier[15](95,10)(110,3)(125,10)
\qbezier[30](125,10)(155,30)(125,50)

\put(110,30){\line(-1,0){30}} \put(107.5,30){\vector(1,0){2}}
\put(110,30){\line(-1,-1){18}} \put(107.5,27.5){\vector(1,1){2}}
\put(110,30){\line(-1,1){18}} \put(107.5,32.5){\vector(1,-1){2}}
\put(110,30){\line(0,1){23.5}} \put(110,33){\vector(0,-1){2}}

\put(110,30){\vector(1,0){30}}

\put(114,28){\makebox(0,0)[tr]{$\beta$}}
\put(95,18){\small\makebox(0,0)[b]{$a$}}
\put(94,31){\small\makebox(0,0)[b]{$b$}}


\put(230,15){\circle*{2}} \put(270,35){\circle*{2}}

\qbezier[30](200,25)(200,-5)(250,-5)
\qbezier[30](250,-5)(300,-5)(300,25)
\qbezier[30](300,25)(300,55)(250,55)
\qbezier[30](250,55)(200,55)(200,25)

\put(230,15){\line(-2,1){29.5}} \put(227.5,16){\vector(2,-1){2}}
\put(230,15){\line(-1,-1){14}} \put(227.5,12.5){\vector(1,1){2}}
\put(230,15){\vector(1,-1){20}}

\qbezier(230,15)(240,35)(270,35) \put(267.5,35){\vector(1,0){2}}
\qbezier(230,15)(270,10)(270,35) \put(270,32.5){\vector(0,1){2}}
\put(270,35){\line(-1,1){20}} \put(267.5,37.5){\vector(1,-1){2}}
\put(270,35){\line(0,1){18.5}} \put(270,37.5){\vector(0,-1){2}}
\put(270,35){\vector(1,0){28}}

\put(225,15){\makebox(0,0)[r]{$\alpha$}}
\put(272,34){\makebox(0,0)[tl]{$\beta$}}
\put(252,16){\small\makebox(0,0)[b]{$a$}}
\put(250,33){\small\makebox(0,0)[b]{$b$}}

\end{picture}
\end{center}
while in the second example, which violates conditions for
correctness, we have
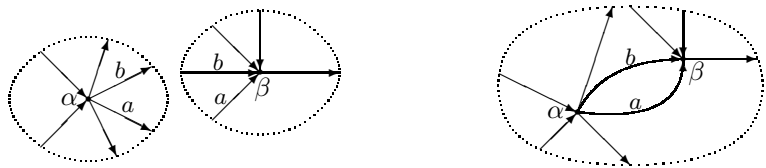
\begin{center}
\begin{picture}(320,60)(0,-5)

\put(45,20){\circle*{2}}

\qbezier[30](30,0)(0,20)(30,40) \qbezier[15](30,40)(45,47)(60,40)
\qbezier[15](30,0)(45,-7)(60,0) \qbezier[30](60,0)(90,20)(60,40)

\put(45,20){\line(-1,1){18}} \put(42.5,17.5){\vector(1,1){2}}
\put(45,20){\line(-1,-1){18}} \put(42.5,22.5){\vector(1,-1){2}}

\put(45,20){\vector(2,1){24}} \put(45,20){\vector(2,-1){24}}
\put(45,20){\vector(1,-2){10.7}} \put(45,20){\vector(1,3){7.5}}

\put(41,20){\makebox(0,0)[r]{$\alpha$}}
\put(60,15){\small\makebox(0,0)[b]{$a$}}
\put(57,28){\small\makebox(0,0)[b]{$b$}}

\put(110,30){\circle*{2}}

\qbezier[30](95,10)(65,30)(95,50)
\qbezier[15](95,50)(110,57)(125,50)
\qbezier[15](95,10)(110,3)(125,10)
\qbezier[30](125,10)(155,30)(125,50)

\put(110,30){\line(-1,0){30}} \put(107.5,30){\vector(1,0){2}}
\put(110,30){\line(-1,-1){18}} \put(107.5,27.5){\vector(1,1){2}}
\put(110,30){\line(-1,1){18}} \put(107.5,32.5){\vector(1,-1){2}}
\put(110,30){\line(0,1){23.5}} \put(110,33){\vector(0,-1){2}}

\put(110,30){\vector(1,0){30}}

\put(114,28){\makebox(0,0)[tr]{$\beta$}}
\put(95,18){\small\makebox(0,0)[b]{$a$}}
\put(94,31){\small\makebox(0,0)[b]{$b$}}


\put(230,15){\circle*{2}} \put(270,35){\circle*{2}}

\qbezier[30](200,25)(200,-5)(250,-5)
\qbezier[30](250,-5)(300,-5)(300,25)
\qbezier[30](300,25)(300,55)(250,55)
\qbezier[30](250,55)(200,55)(200,25)

\put(230,15){\line(-2,1){29.5}} \put(227.5,16){\vector(2,-1){2}}
\put(230,15){\line(-1,-1){14}} \put(227.5,12.5){\vector(1,1){2}}
\put(230,15){\vector(1,-1){20}} \put(230,15){\vector(1,3){13.3}}

\qbezier(230,15)(240,35)(270,35) \put(267.5,35){\vector(1,0){2}}
\qbezier(230,15)(270,10)(270,35) \put(270,32.5){\vector(0,1){2}}
\put(270,35){\line(-1,1){20}} \put(267.5,37.5){\vector(1,-1){2}}
\put(270,35){\line(0,1){18.5}} \put(270,37.5){\vector(0,-1){2}}
\put(270,35){\vector(1,0){28}}

\put(225,15){\makebox(0,0)[r]{$\alpha$}}
\put(272,34){\makebox(0,0)[tl]{$\beta$}}
\put(252,16){\small\makebox(0,0)[b]{$a$}}
\put(250,33){\small\makebox(0,0)[b]{$b$}}

\end{picture}
\end{center}
For the graph on the right in this last picture, the dotted
circle surrounding it cannot be divided into two semicircles, one with
outgoing arrows and the other with incoming arrows. This division
of the surrounding circle is what the exclusion of the
bifurcation mentioned above corresponds~to.

The modified dual graphs we have just introduced are analogous to the \emph{string diagrams} of \cite{JS91} (Chapter~1), \cite{S95} (Sections 4-5) and \cite{S96} (Section~4), while graphs that correspond to the diagrams of 2-cells like that in the first picture are the \emph{pasting schemes} of \cite{P90} (see \S 7.3 for the definition of this notion; see also \S 6.6  and \S 6.7). In the definition of string diagram, as in the definition of pasting scheme, planarity is assumed.

Juncture for our modified dual graphs is
applicable to a wider class of graphs than these modified dual
graphs; we call the graphs in this wider class D1-graphs (see \S
6.5). We do not assume for D1-graphs the special kind of
planarity, which consists in these graphs being realizable within
a disk as in all the pictures above except the last (where the
dotted circle could not be divided in an appropriate manner into
two semicircles). Let us call this special planarity \emph{disk
planarity}\index{disk planarity}.

Planarity need not be taken as a difficult notion from a
geometrical point of view, but from a combinatorial, i.e.\
properly graph-theoretical, point of view, it is not simple, and
our goal is to replace the assumption of disk planarity by purely
combinatorial assumptions. In other words, our goal is to
characterize disk planarity, i.e.\ the disk planar realizability
of D1-graphs, in combinatorial terms. This is a goal analogous to
that achieved by Kuratowski's and other characterizations of
planarity in graph theory (see \cite{H69}, Chapter 11; as a
byproduct of our characterization of planarity in this work, we
obtained in \cite{DP12} another characterization of planarity for
ordinary graphs, akin to Kuratowski's). Our reason for dualizing
the graphs of 2-cells is this characterization of disk planarity,
which otherwise we could not give.

We will find it more practical for our characterization of disk
planarity of D1-graphs to concentrate on juncture in the absence
of identity 2-cells, which yields the notion of \emph{D-graph}.
The D-graphs that have a disk planar realization will be called
\emph{P-graphs}\index{P-graph}. This notion is extended later (in Chapter~6) to
the notion of \emph{P1-graph} (see \S 6.5), a disk realizable
D1-graph, which has what is needed for identity 2-cells, and where
what corresponds to horizontal composition is definable.

We define P-graphs in an inductive manner involving juncture with
the notion of P$'$-graph (see \S 1.8), and non-inductively, again
involving juncture, with the notion of P$'''$-graph (see \S 1.10).
The notion of P$'''$-graph gives in the most accomplished form our
combinatorial characterization of disk planarity for D-graphs,
which provides the gist of what we need for D1-graphs. The
inductively defined notion of P$''$-graph is intermediary, and
serves as a tool to prove the equivalence of the notions of
P$'$-graph and P$'''$-graph.

The proof of this equivalence, which will occupy us
in most of our work (see Chapters~2-5), is interesting not only because of the final result it yields, but because of the light it sheds on the articulation of the notion of P-graph. We believe that the notions and techniques this proof relies on are of an intrinsic interest too. The length and the difficulties of this proof may come as a surprise, because our three definitions of P-graph are not that different. If however there is no proof much simpler than the one we found, then, judging by the distance our proof covers, these notions are indeed wide apart.

The last chapter of our work (Chapter~7) is about geometrical
realizations of P-graphs and P1-graphs. Having both the P$'$ and
P$'''$ version of the notion of P-graph will help us for that
matter.

Another result of our work is a criterion for a graph to be
realizable as a graph associated with a diagram of 2-cells, a
criterion not based on dualizing (see \S 7.5). For that we rely on
our combinatorial characterization of disk planarity.

Juncture is related to the operation of cut on sequents that one
encounters in proof theory. The aspects of juncture as they occur
in proof theory were treated in \cite{DP11}, the results of which
are related to the definition of planar polycategory---a notion
that generalizes the notions of multicategory and operad. As the notions of polycategory and multicategory, the notion of operad may be based on an operation like cut (see \cite{DP10}). The operation of juncture treated in the present work is more general than all these operations related to cut.
Proof-theoretically, it allows for cuts via finite non-empty
sequences of formulae, and with the commas on the left and right
of the turnstile being of the same nature; i.e., they are both
understood conjunctively, or both disjunctively. Moreover, what we have in this work is about sequents where we do not assume the structural rule of permutation; i.e., neither of the commas corresponds to a commutative operation. Gentzen's cut,
the plural (multiple-conclusion) form, which is treated in
\cite{DP11}, proceeds via sequences that have just a single
formula, and the two commas are of different nature.
Graph-theoretically, matters are more complicated with this more
general cut, i.e.\ with juncture, than with Gentzen's cut.

Our work points out to an alternative foundation for bicategories, i.e.\ weak 2-categories, and, further, to an alternative foundation for weak $n$-categories, a matter much debated these days. We hope that our approach may shed new light on this matter. Our equations for juncture would be replaced by cells of a higher level, cells that are isomorphisms.

Although the motivation for this work comes from category theory,
we do not deal much with this theory, and do not presuppose the
reader has any extended knowledge of it, except for the sake of
motivation. We deal in this work with matters of graph theory, but in
that theory we define everything we need, and do not presuppose
anything in particular. For the remaining mathematical disciplines
touched in our work, like geometry, topology and logic, we
presuppose only elementary matters.

\section{\large\bf D-graphs}\label{1.2}
\markright{\S 1.2. \quad D-graphs}

We introduce first the notion of graph that is common in category
theory (see \cite{ML98}, Sections I.2 and II.7). This notion,
under the label \emph{graph}, tout court, may be found in
\cite{B73} (Section 1.1), and, under the labels \emph{directed
graph} and \emph{digraph}, in \cite{BM76} (Section 10.1),
\cite{W01} (Section 1.4) and \cite{D10} (Section 1.10). What we
call graphs are not the pseudographs of \cite{H69} (Chapter~2),
which are not directed. In the style of \cite{H69}, we could call
the graphs of this work \emph{directed pseudographs}.

So a \emph{graph}\index{graph} in this work is given by two
functions $W,E\!:A\rightarrow V$, where the elements of the set $A$
are called \emph{edges}\index{edge} and those of the set $V$
\emph{vertices}\index{vertex} (in category theory, they would be
respectively \emph{arrows}, or \emph{morphisms}, and
\emph{objects}). The names of the functions $W$\index{W@$W$} and
$E$\index{E@$E$} come from \emph{West} and \emph{East}, which
accords with how we will draw pictures for particular kinds of
these graphs, from left to right (in category theory, $W$ and $E$
would be respectively the \emph{source} and \emph{target}, or
\emph{domain} and \emph{codomain}, functions). In these pictures,
an edge $a$ is represented by an arrow going from the point
representing the vertex $W(a)$ to the point representing the
vertex $E(a)$. The names we use for $W$ and $E$ (rather than
something derived from \emph{left} and \emph{right}, or the
categorial terminology) become natural when we deal with
geometrical realizations in Chapter~7, where \emph{North} and
\emph{South} appear too (see also \S 6.6; for such reasons we used
already an analogous terminology in~\cite{DP11}).

We use $X$\index{X@$X$} as a variable for $W$ and $E$. We write $\bar{W}$\index{W bar@$\bar{W}$} and
$\bar{E}$\index{E bar@$\bar{E}$} for $E$ and $W$ respectively.\index{X bar@$\bar{X}$}

We say sometimes that an edge $a$ of a graph is an edge
\emph{from} $W(a)$ \emph{to} $E(a)$, and we say that $W(a)$ and
$E(a)$ are \emph{incident}\index{incident} with~$a$.

A \emph{graph morphism}\index{graph morphism} from the graph
$W_1,E_1\!:A_1\rightarrow V_1$ to the graph
$W_2,E_2\!:A_2\rightarrow V_2$ (which is analogous to a
\emph{functor} between categories) is a pair of functions
$F_A\!:A_1\rightarrow A_2$ and $F_V\!:V_1\rightarrow V_2$ such
that for every edge $a$ in $A_1$ and every $X$ in $\{W,E\}$ we
have $F_V(X_1(a))=X_2(F_A(a))$. This means that $F_A(a)$ is an
edge from $F_V(W_1(a))$ to $F_V(E_1(a))$.

A graph morphism is an \emph{isomorphism}\index{isomorphism}\index{graph isomorphism} when both $F_A$ and $F_V$ are bijections.

Let a graph $W,E\!:A\rightarrow V$ be \emph{distinguished}\index{distinguished graph} when
$A$ and $V$ are disjoint. (This condition of disjointness for
graphs does not seem to be often mentioned in textbooks of graph
theory---exceptions are \cite{BM76} and \cite{D10}---but it is
presumably tacitly assumed by many authors.) A non-distinguished
graph is given, for example, by
\begin{tabbing}
\hspace{1.6em}\=$A=\{a,b\}$,\hspace{2em}\=$V=\{u,v,w,a\}$, \begin{picture}(40,0)(-70,28)

\put(0,30){\circle*{2}} \put(40,30){\circle*{2}} \put(20,10)

\put(0,30){\vector(1,0){40}} \put(20,30){\vector(0,-1){20}}

\put(-3,30){\small\makebox(0,0)[r]{$u$}}
\put(43,30){\small\makebox(0,0)[l]{$v$}}
\put(20,7){\small\makebox(0,0)[t]{$w$}}
\put(20,33){\small\makebox(0,0)[b]{$a$}}

\put(25,19){\small\makebox(0,0)[b]{$b$}}

\end{picture}\\[.3ex]
\>$W(a)=u$,\>$E(a)=v$,\\
\>$W(b)=a$,\>$E(b)=w$.
\end{tabbing}

More natural examples of non-distinguished graphs are
found in category theory, where sometimes objects, i.e.\ vertices,
are identified with identity arrows on these objects; so all
vertices are edges.

It is trivial to show that every graph is isomorphic to a
distinguished graph. Just replace either the set of edges or the
set of vertices by a new set in one-to-one correspondence with the
original one. Every plane graph (see \S 7.1) is distinguished, and
because of that the picture of the non-distinguished graph we had
above as an example is not very natural.

From now on we assume that \emph{graph} means \emph{distinguished
graph}, though this assumption is not always essential.

Note that a graph can have $V$ empty, in which case $A$ must be
empty too, and for both $W$ and $E$ we have the empty set of
ordered pairs. The graph that has $V$ empty is the \emph{empty
graph}\index{empty graph}. A graph that is not the empty graph is said to be \emph{non-empty}.\index{non-empty graph}

A graph is \emph{finite}\index{finite graph} when $A$ and $V$ are finite. In this work
we shall be concerned only with finite graphs.

A vertex $v$ of a graph is an $X$-\emph{vertex}\index{X-vertex@$X$-vertex}\index{W-vertex@$W$-vertex}
\index{E-vertex@$E$-vertex} of that graph when there is no
edge $a$ of that graph such that $\bar{X}(a)=v$. For example, in
the graph of the following picture
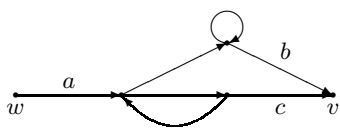
\begin{center}
\begin{picture}(120,40)(0,-5)

\put(0,10){\circle*{2}} \put(40,10){\circle*{2}}
\put(80,10){\circle*{2}} \put(120,10){\circle*{2}}
\put(80,30){\circle*{2}}

\put(80,36){\circle{12}}

\put(82.5,30.8){\vector(-2,-1){2}}

\qbezier(40,10)(60,-13)(80,10)

\put(0,10){\vector(1,0){40}} \put(40,10){\vector(1,0){40}}
\put(80,10){\vector(1,0){40}} \put(40,10){\vector(2,1){40}}

\put(80,30){\vector(2,-1){40}}

\put(42.5,7.5){\vector(-1,1){2}}

\put(0,7){\small\makebox(0,0)[t]{$w$}}
\put(20,13){\small\makebox(0,0)[b]{$a$}}
\put(100,7){\small\makebox(0,0)[t]{$c$}}
\put(120,7){\small\makebox(0,0)[t]{$v$}}
\put(102,24){\small\makebox(0,0)[b]{$b$}}

\end{picture}
\end{center}
the vertex $w$ is a $W$-vertex, while $v$ is an
$E$-vertex; the other vertices are neither $W$-vertices nor
$E$-vertices. (In the style of \cite{H69}, Chapter 16, we could
say that $W$-vertices have indegree 0, while $E$-vertices have
outdegree~0.)

The notion of $X$-vertex is given with respect to a particular
graph, and we mentioned that explicitly in the definition. We shall next define a
series of notions that should likewise be understood as given with
respect to a particular graph, but we will take this for granted,
and will not mention it explicitly.

A vertex is an \emph{inner} vertex\index{inner vertex} when it is neither a $W$-vertex
nor an $E$-vertex; $W$-vertices and $E$-vertices are accordingly called \emph{outer}\index{outer vertex}. All the vertices in our example above
except the outer vertices $w$ and $v$ are inner vertices.

An edge $a$ is an $X$-\emph{edge}\index{X-edge@$X$-edge}\index{W-edge@$W$-edge}
\index{E-edge@$E$-edge} when $X(a)$ is an $X$-vertex. An
edge is \emph{inner}\index{inner edge} when it is neither a $W$-edge nor an
$E$-edge. Alternatively, $a$ is an inner edge when $W(a)$ and
$E(a)$ are inner vertices. In our example, $a$ is a $W$-edge,
while $b$ and $c$ are $E$-edges; the remaining edges are inner.

An $X$-edge $a$ is $X$-\emph{functional}\index{X-functional edge@$X$-functional edge}\index{W-functional edge@$W$-functional edge}\index{E-functional edge@$E$-functional edge} when for every edge $b$ of
our graph different from $a$ the vertices $X(a)$ and $X(b)$ are
different. In our example, the $W$-edge $a$ is $W$-functional,
while the $E$-edges $b$ and $c$ are not $E$-functional.

A graph is $W\mbox{\rm -}E$-\emph{functional}\index{W-E-functional graph@$W\mbox{\rm -}E$-functional graph} when all its $W$-edges are $W$-functional and all its $E$-edges are $E$-functional.

We give next the definitions of a number of notions analogous to
those that may be found in \cite{H69}, and for which accordingly
we use the same terms. The reader should however keep in mind that
these are not exactly the same notions, but analogous notions
adapted to our context; the \emph{graphs} of \cite{H69} are
ordinary graphs, and not our graphs.

A \emph{semiwalk}\index{semiwalk} is either a vertex $v_0$, in which case the
semiwalk is \emph{trivial}\index{trivial semiwalk}, or for $n\geq 1$ this is a sequence
$v_0a_1v_1\ldots v_{n-1}a_nv_n$ such that for every
$i$ in $\{1,\ldots,n\}$
\begin{tabbing}
\hspace{1.6em}\=(1)\hspace{2em}\=$W(a_i)=v_{i-1}$ \=and
$E(a_i)=v_i$, or\\
\>(2)\>$W(a_i)=v_i$\>and $E(a_i)=v_{i-1}$.
\end{tabbing}
A trivial semiwalk $v_0$ is a semiwalk \emph{from} $v_0$ \emph{to}
$v_0$, while a non-trivial one is a semiwalk \emph{from} $v_0$
\emph{to} $v_n$. (Semiwalks from $v_0$ to $v_n$ correspond
bijectively to semiwalks from $v_n$ to $v_0$.) We also say that a semiwalk from $v_0$
to $v_n$ \emph{connects}\index{connect} $v_0$ \emph{with} $v_n$.  For $\sigma$ a
semiwalk and $x$ a vertex or edge, we write $x\rhd\sigma$ when $x$
occurs in~$\sigma$.

By omitting (2) from the definition of semiwalk we obtain the
definition of \emph{walk}\index{walk}.

A \emph{semipath}\index{semipath} is a semiwalk such that
\begin{tabbing}
\hspace{1.6em}\=$(*)$\hspace{2em}\=no vertex in it occurs more
than once.
\end{tabbing}
Hence all the edges in a semipath are also mutually distinct. (Examples of semipaths may be found in \S 1.9.) A
\emph{path}\index{path} is a walk such that $(*)$ holds.

A graph is \emph{weakly connected}\index{weakly connected graph} when for every two distinct
vertices $v_0$ and $v_n$ there is a semipath from $v_0$ to $v_n$
(which must be non-trivial).

A \emph{semicycle}\index{semicycle} is a non-trivial semiwalk from $v_0$ to $v_n$
such that
\begin{tabbing}
\hspace{1.6em}\=$(*)$\hspace{2em}\=no vertex in it occurs more
than once.\kill \>$(**)$\>no vertex in it occurs more than once,
except that $v_0$ is~$v_n$.
\end{tabbing}
So, in the limit case, $v_0av_0$ may be a semicycle based on a
non-trivial semiwalk. A \emph{cycle}\index{cycle} is a non-trivial walk such
that $(**)$ holds.

A graph is \emph{acyclic}\index{acyclic graph} when it has no cycles.

Now we have all we need to define one of the main kinds of graph with
which we deal in this work, and which we call \emph{D-graph}.\index{D-graph}
A D-graph is a graph that is finite, acyclic, $W\mbox{\rm -}E$-functional,
weakly connected and with an inner vertex.

A graph is \emph{incidented}\index{incidented graph} when for each of its vertices $v$ there is an edge $a$ such that $W(a)=v$ or $E(a)=v$; i.e., $v$ is
incident with $a$. It is easy to infer that every D-graph is
incidented.

A \emph{loop}\index{loop} is an edge $a$ such that $W(a)=E(a)$.
The acyclicity condition excludes loops in D-graphs.

We will draw D-graphs from left to right, and here is a picture of
one of them:
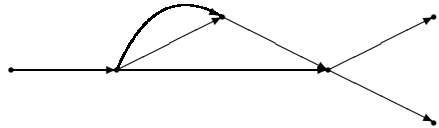
\begin{center}
\begin{picture}(160,50)(0,-10)

\put(0,10){\circle*{2}} \put(40,10){\circle*{2}}
\put(120,10){\circle*{2}} \put(80,30){\circle*{2}}
\put(160,30){\circle*{2}} \put(160,-10){\circle*{2}}

\qbezier(40,10)(55,45)(80,30) \put(77,31.5){\vector(2,-1){2}}

\put(0,10){\vector(1,0){40}} \put(40,10){\vector(1,0){80}}
\put(40,10){\vector(2,1){40}}

\put(80,30){\vector(2,-1){40}}

\put(120,10){\vector(2,1){40}} \put(120,10){\vector(2,-1){40}}

\end{picture}
\end{center}

A \emph{basic D-graph}\index{basic D-graph} is a D-graph with a single inner vertex.
Basic D-graphs are all of the form
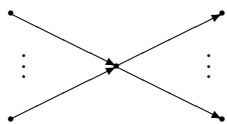
\begin{center}
\begin{picture}(80,50)(0,-5)

\put(0,0){\circle*{2}} \put(0,40){\circle*{2}}
\put(40,20){\circle*{2}} \put(80,0){\circle*{2}}
\put(80,40){\circle*{2}}

\put(0,0){\vector(2,1){40}} \put(0,40){\vector(2,-1){40}}
\put(40,20){\vector(2,1){40}} \put(40,20){\vector(2,-1){40}}

\put(5,23){\makebox(0,0){$\vdots$}}
\put(75,23){\makebox(0,0){$\vdots$}}

\end{picture}
\end{center}

\section{\large\bf Cocycles and juncture}\label{1.3}
\markright{\S 1.3. \quad Cocycles and juncture}

We say that the graph $G_1$, which is $W_1,E_1\!:A_1\rightarrow V_1$, is a \emph{subgraph}\index{subgraph} of the graph $G_2$, which is $W_2,E_2\!:A_2\rightarrow V_2$, when for $Z$ being one of $A$, $V$, $W$ and $E$ we have $Z_1\subseteq Z_2$. (It is clear that the relation of being a subgraph is a partial order.)

A \emph{component}\index{component of a graph} of a graph $G$ is a weakly connected non-empty subgraph $G'$ of $G$ such that for every weakly connected subgraph $G''$ of $G$, if $G'$ is a subgraph of $G''$, then $G'$ is~$G''$.

Consider a set $S$ of inner edges of a D-graph $D$. The
\emph{removal}\index{removal of edges} of $S$ from $D$ leaves a new graph with the same
vertices and with the edges from $S$ missing; the $W$ and $E$
functions of the new graph are obtained by restricting those of
$D$ to the new set of edges. The new graph is made of a family of components
$D_1,\ldots,D_n$ of this graph, for $n\geq 1$. Note that the graphs in this
family are not necessarily D-graphs.

When $n\geq 2$, in which case $S$ must be non-empty, we say that
$S$ is a \emph{cutset}\index{cutset} (which is a term standing for an analogous
notion of \cite{H69}, Chapter~4).

A \emph{directed graph}\index{directed graph}, in the sense of
\cite{H69} (Chapter 2; also called \emph{digraph}) is an
irreflexive binary relation on a finite set of vertices; the
ordered pairs of the binary relation are the edges. Such a graph
may be identified with a finite graph in our sense where there are
no multiple edges with the same vertices incident with them, and
no loops (see the end of \S 1.2). Various notions concerning
directed graphs, like weak connectedness and acyclicity, which we
will rely on in a moment (and other notions we have in \S 1.6),
may either be given definitions analogous to those we gave for
graphs (see \cite{H69}, Chapter 16, for weak connectedness and
acyclicity of directed graphs), or having in mind the
identification of directed graphs with a special kind of graph, we
may apply the definitions we gave for graphs.

Let $C_S(D)$\index{CSD@$C_S(D)$}, the \emph{componential graph}\index{componential graph} of $D$ with respect to
$S$, be the directed graph in the sense of \cite{H69} whose
vertices are $D_1,\ldots,D_n$, and such that for some distinct $i$ and $j$
in $\{1,\ldots,n\}$ we have that the ordered pair $(D_i,D_j)$ is
an edge of $C_S(D)$ iff there is an edge in $S$ from a vertex of
$D_i$ to a vertex of $D_j$, i.e.\ an edge $a$ in $S$ such that
\begin{tabbing}
\hspace{1.6em}\=$(*)$\hspace{2em}\=no vertex in it occurs more
than once.\kill

\>($\ddag$)\>$W(a)$ is a vertex of $D_i$ and $E(a)$
is a vertex of $D_j$.
\end{tabbing}
Less formally, we may say that the edge $a$ connects $D_i$ with~$D_j$.

It is easy to see that since the D-graph $D$ is weakly connected the directed graph $C_S(D)$ is weakly connected.

We call a cutset $S$ of $D$ \emph{strict}\index{strict cutset} when
for every $a$ in $S$ there are distinct $i$ and $j$ in
$\{1,\ldots,n\}$ such that ($\ddag$), and $C_S(D)$ is acyclic.
The acyclicity condition for the componential graph precludes $\{a,b\}$ from being a strict cutset in the D-graph of the following picture
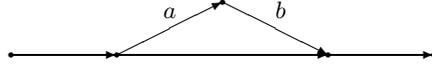
\begin{center}
\begin{picture}(160,30)(0,5)

\put(0,10){\circle*{2}} \put(40,10){\circle*{2}}
\put(120,10){\circle*{2}} \put(80,30){\circle*{2}}
\put(160,10){\circle*{2}}

\put(0,10){\vector(1,0){40}} \put(40,10){\vector(1,0){80}}
\put(40,10){\vector(2,1){40}}

\put(80,30){\vector(2,-1){40}}

\put(120,10){\vector(1,0){40}}

\put(60,24){\small\makebox(0,0)[b]{$a$}}

\put(102,24){\small\makebox(0,0)[b]{$b$}}

\end{picture}
\end{center}

A strict cutset where $n=2$ will be called a \emph{cocycle}\index{cocycle} (which
is a term standing for an analogous notion of \cite{H69},
Chapter~4).

Cocycles are related to a binary partial operation on D-graphs,
which we will call \emph{juncture}\index{juncture}, and which now we proceed to
define.

For $X$ being $W$ or $E$, let $D_X$ be $W_X,E_X\!:A_X\rightarrow V_X$,
and assume that the two graphs $D_W$ and $D_E$ are D-graphs.
Assume moreover that
\begin{tabbing}
\hspace{1.6em}\=$C=_{df} A_W\cap A_E\neq\emptyset$,\\[.5ex]
\>$(\forall a\in C)\; E_W(a)=W_E(a)$,\\[.5ex]
\>$(\forall a\in C)$ $a$ is an $\bar{X}$-edge of $D_X$.
\end{tabbing}
Let $V_C=\{v\mid (\exists a\in C)\; v=E_W(a)\}=\{v\mid (\exists
a\in C)\; v=W_E(a)\}$, and assume that
\begin{tabbing}
\hspace{1.6em}\=$V_W\cap V_E=V_C$.
\end{tabbing}
It can be inferred that every vertex in $V_C$ is an
$\bar{X}$-vertex of~$D_X$.

Then we define the D-graph $D_W\Box D_E$, which is
$W,E\!:A\rightarrow V$, in the following way:
\begin{tabbing}
\hspace{1.6em}\=$A$ \=$=A_W\cup A_E$,\\[.5ex]
\>$V$\>$=(V_W\cup V_E) - V_C$,\\[.5ex]
for $a$ in $A$,\\
\>$ X(a)=\left\{
\begin{array}{ll}
X_X(a) & \mbox{\rm{if }} a\in A_X,
\\[.5ex]
X_{\bar{X}}(a) & \mbox{\rm{if }} a\in A_{\bar{X}}-C.
\end{array}
\right.$
\end{tabbing}
This concludes the definition of the operation of juncture~$\Box$.

For example, consider the D-graphs in the following picture:
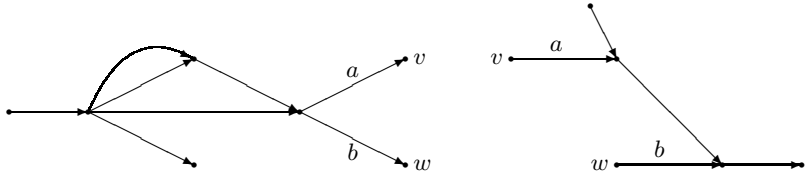
\begin{center}
\begin{picture}(300,60)(10,-10)

\put(10,10){\circle*{2}} \put(40,10){\circle*{2}}
\put(120,10){\circle*{2}} \put(80,30){\circle*{2}}
\put(80,-10){\circle*{2}} \put(160,-10){\circle*{2}}
\put(160,30){\circle*{2}}

\put(10,10){\vector(1,0){30}} \put(40,10){\vector(1,0){80}}
\put(40,10){\vector(2,1){40}}

\put(80,30){\vector(2,-1){40}}

\put(40,10){\vector(2,-1){40}} \put(120,10){\vector(2,1){40}}
\put(120,10){\vector(2,-1){40}}

\qbezier(40,10)(55,45)(80,30) \put(77,31.5){\vector(2,-1){2}}

\put(140,23){\small\makebox(0,0)[b]{$a$}}
\put(140,-2){\small\makebox(0,0)[t]{$b$}}
\put(163,30){\small\makebox(0,0)[l]{$v$}}
\put(163,-10){\small\makebox(0,0)[l]{$w$}}


\put(200,30){\circle*{2}} \put(240,-10){\circle*{2}}
\put(240,30){\circle*{2}} \put(280,-10){\circle*{2}}
\put(230,50){\circle*{2}} \put(310,-10){\circle*{2}}

\put(200,30){\vector(1,0){40}} \put(240,-10){\vector(1,0){40}}
\put(280,-10){\vector(1,0){30}} \put(240,30){\vector(1,-1){40}}
\put(230,50){\vector(1,-2){10}}

\put(217,33){\small\makebox(0,0)[b]{$a$}}
\put(256,-7){\small\makebox(0,0)[b]{$b$}}
\put(197,30){\small\makebox(0,0)[r]{$v$}}
\put(237,-10){\small\makebox(0,0)[r]{$w$}}

\end{picture}
\end{center}
The D-graph $D_W\Box D_E$ is in the picture
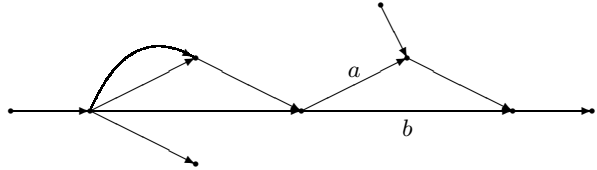
\begin{center}
\begin{picture}(220,60)(10,-10)

\put(10,10){\circle*{2}} \put(40,10){\circle*{2}}
\put(120,10){\circle*{2}} \put(80,30){\circle*{2}}
\put(80,-10){\circle*{2}} \put(160,30){\circle*{2}}
\put(200,10){\circle*{2}} \put(230,10){\circle*{2}}
\put(150,50){\circle*{2}}

\put(10,10){\vector(1,0){30}} \put(40,10){\vector(1,0){80}}
\put(40,10){\vector(2,1){40}}

\put(80,30){\vector(2,-1){40}}

\put(40,10){\vector(2,-1){40}} \put(120,10){\vector(2,1){40}}
\put(120,10){\vector(1,0){80}}

\qbezier(40,10)(55,45)(80,30) \put(77,31.5){\vector(2,-1){2}}
\put(150,50){\vector(1,-2){10}}

\put(160,30){\vector(2,-1){40}}

\put(200,10){\vector(1,0){30}}

\put(140,23){\small\makebox(0,0)[b]{$a$}}
\put(160,7){\small\makebox(0,0)[t]{$b$}}

\end{picture}
\end{center}

It is easy to check that $D_W\Box D_E$ is always a D-graph.

Note that in the resulting D-graph $D_W\Box D_E$ the set of edges
$C$ is a cocycle. By removing $C$ from $D_W\Box D_E$ we obtain the
graphs $D_W$ and $D_E$ with the edges of $C$ removed and the
isolated vertices of $V_C$ omitted.

Conversely, if we start from a D-graph $D$ and an arbitrary
cocycle $C$ of $D$, then we can construct two D-graphs $D_W$ and
$D_E$ such that $D$ is $D_W\Box D_E$ and $C$ is $A_W\cap A_E$ (see
\S 1.10 for details).

\section{\large\bf Edge-graphs}\label{1.4}
\markright{\S 1.4. \quad Edge-graphs}

In this section we consider a notion of graph without vertices, called \emph{edge-graph}, which is equivalent to the notion of incidented graph (see the end of \S 1.2). Among edge-graphs, those that correspond to D-graphs will enable us to reformulate juncture in a particularly simple way. It will boil down to union.

The notion of edge-graph shows that vertices are in principle dispensable in our exposition, but we keep them because sometimes it is more convenient to rely on them, and also because we do not want to depart too far from established terminology. We will however rely on edge-graphs in \S 6.5. Matters exposed in this section are not essential for our results later on, and this is why here we will not dwell on the details of the proofs.

The definition of edge-graph does not mention vertices, but instead it mentions the relations on edges of having a common vertex. There are three such binary relations, because the common vertex may be on the west in both edges, or on the east in both edges, or on the east in one edge and on the west in the other, in which case the first edge precedes the second; here are the three relations in pictures:
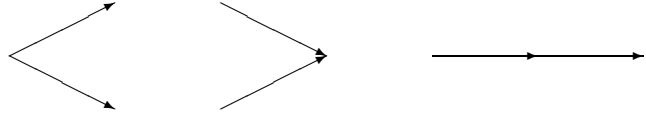
\begin{center}
\begin{picture}(240,50)(0,-5)
\put(0,20){\vector(2,1){40}} \put(0,20){\vector(2,-1){40}}

\put(80,0){\vector(2,1){40}} \put(80,40){\vector(2,-1){40}}

\put(160,20){\vector(1,0){40}} \put(200,20){\vector(1,0){40}}
\end{picture}
\end{center}

An \emph{edge-graph}\index{edge-graph} is a set $A$, whose element are called \emph{edges}\index{edge}, together with three binary relations $\textbf{W},\textbf{E},\textbf{P}\subseteq A^2$, such that $\textbf{W}$ and $\textbf{E}$ are equivalence relations, and for every $a$, $b$ and $c$ in $A$ we have
\begin{tabbing}
\hspace{1.6em}\=$a\textbf{W}b\Rightarrow(c\textbf{P}a\Rightarrow c\textbf{P}b)$,\hspace{2em}\=$a\textbf{E}b\Rightarrow(a\textbf{P}c\Rightarrow b\textbf{P}c)$,\\[.5ex]
\>$(c\textbf{P}a \;\&\; c\textbf{P}b)\Rightarrow a\textbf{W}b$,\>$(a\textbf{P}c\;\&\; b\textbf{P}c)\Rightarrow a\textbf{E}b$.
\end{tabbing}
For $\textbf{X}$ being $\textbf{W}$ or $\textbf{E}$, we read intuitively $a\textbf{X}b$ as $a$ \emph{and} $b$ \emph{have the same} \textbf{X}-\emph{end}, while $a\textbf{P}b$ is read as $a$ \emph{precedes}~$b$.

The equivalence of the notion of edge-graph with the notion of incidented graph is a result about equivalence of categories. We define first the category $\cal E$, where the objects are edge-graphs and the arrows are edge-graph morphisms, which we are now going to define.

An \emph{edge-graph morphism}\index{edge-graph morphism} from the edge-graph $\langle A_1, \textbf{W}_1, \textbf{E}_1, \textbf{P}_1\rangle$ to the edge-graph $\langle A_2, \textbf{W}_2, \textbf{E}_2, \textbf{P}_2\rangle$ is a function $F_A\!:A_1\rightarrow A_2$ such that for $\textbf{Z}$ being $\textbf{W}$, $\textbf{E}$ or $\textbf{P}$, if $a\textbf{Z}_1b$, then $F_A(a)\textbf{Z}_2F_A(b)$. The identity functions on edges serve to define the identity edge-graph morphisms, and composition in $\cal E$ is given by composition of functions.

The category $\cal I$, with which $\cal E$ is equivalent, has as its objects incidented graphs and as arrows graph morphisms (see \S 1.2). The identity graph morphisms of $\cal I$ are based on the identity functions on edges and vertices, and composition in $\cal I$ is based on composition of functions.

For a graph $G$, which is $W,E\!:A\rightarrow V$, let the edge-graph ${\cal H}(G)$, which is $\langle A, \textbf{W}, \textbf{E}, \textbf{P}\rangle$, be obtained by stipulating that for every $a$ and $b$ in $A$ we have
\begin{tabbing}
\hspace{1.6em}\=$a\textbf{W}b$ \=$\Leftrightarrow W(a)$ \=$=W(b)$,\\[.5ex]
\>$\;a\textbf{E}b$\>$\Leftrightarrow \; E(a)$\>$=E(b)$,\\[.5ex]
\>$\;a\textbf{P}b$\>$\Leftrightarrow \; E(a)$\>$=W(b)$.
\end{tabbing}

It is clear that $\textbf{W}$ and $\textbf{E}$ are equivalence relations on $A$; we also have
\begin{tabbing}
\hspace{1.6em}\=$W(a)=W(b)\Rightarrow(E(c)=W(a)\Rightarrow E(c)=W(b))$,\\[.5ex]
\>$(E(c)=W(a)\;\&\; E(c)=W(b))\Rightarrow W(a)=W(b)$,
\end{tabbing}
and analogously with $E(a)=E(b)$ instead of $W(a)=W(b)$, so that we may conclude that $\langle A, \textbf{W}, \textbf{E}, \textbf{P}\rangle$ is an edge-graph.

The empty graph $\emptyset,\emptyset\!:\emptyset\rightarrow\emptyset$ (see \S 1.2) is mapped by $\cal H$ to the \emph{empty edge-graph}\index{empty edge-graph} $\langle\emptyset,\emptyset,\emptyset,\emptyset\rangle$. A \emph{single-vertex graph}\index{single-vertex graph}, which is $\emptyset,\emptyset\!:\emptyset\rightarrow\{v\}$, is mapped by $\cal H$ to the empty edge-graph too. Note that a single-vertex graph is not an incidented graph.

For an edge-graph $H$, which is $\langle A, \textbf{W}, \textbf{E}, \textbf{P}\rangle$, we obtain as follows the incidented graph ${\cal G}(H)$, which is $W,E\!:A\rightarrow V_H$. For $a$ in $A$ let
\begin{tabbing}
\hspace{1.6em}\=$_W[a]$ \=$=\{b\in A\mid a\textbf{W}b\}$,\hspace{2em}\=$[a]_E$ \=$=\{b\in A\mid a\textbf{E}b\}$,\\[.5ex]
\>$_P[a]$\>$=\{b\in A\mid b\textbf{P}a\}$,\>$[a]_P$\>$=\{b\in A\mid a\textbf{P}b\}$,\\[.7ex]
\>$W(a)= (_P[a],_W[a])$,\>\>$E(a)=([a]_E,[a]_P)$,\\[.7ex]
\>$V_H=\{(A',A'')\mid(\exists a\in A)(W(a)=(A',A'') \;{or}\; E(a)=(A',A''))\}$.
\end{tabbing}
If $H$ is the empty edge-graph, then ${\cal G}(H)$ is the empty graph.

It is straightforward to prove the following proposition.

\prop{Proposition 1.4.1}{For every incidented graph $G$, the graph ${\cal G}({\cal H}(G))$ is isomorphic to~$G$.}

The graph isomorphism of this proposition is identity on edges,
and, for $X$ being $W$ or $E$, it maps the vertex $X(a)$ of ${\cal
G}({\cal H}(G))$ to the vertex $X(a)$ of $G$. (One must verify
that this function on vertices is well defined; part of that
consists in verifying that if
$(_P[a_1],_W[a_1])=(_P[a_2],_W[a_2])$, then $W(a_1)=W(a_2)$ in
$G$.)

The following proposition is also straightforward to prove.

\prop{Proposition 1.4.2}{For every edge-graph $H$, the edge-graph ${\cal H}({\cal G}(H))$ is~$H$.}

It is straightforward to extend $\cal H$ to a functor from the category $\cal I$ to the category $\cal E$; we just forget the $F_V$ part of a graph morphism. It is also straightforward to extend $\cal G$ to a functor from $\cal E$ to $\cal I$. (A vertex $X(a)$ of ${\cal G}(H_1)$ is mapped by $F_V$ to the vertex $X(F_A(a))$ of ${\cal G}(H_2)$.) Starting from Propositions 1.4.1 and 1.4.2, we then obtain that the categories $\cal E$ and $\cal I$ are equivalent.

Let a \emph{D-edge-graph}\index{D-edge-graph} be an edge-graph $H$
such that ${\cal G}(H)$ is a D-graph. A
$W$-\emph{edge}\index{X-edge of edge-graph@$X$-edge of
edge-graph}\index{W-edge of edge-graph@$W$-edge of edge-graph} of
an edge-graph $H$ is an edge $a$ such that there is no edge $b$ of
$H$ with $b\textbf{P}a$. An $E$-\emph{edge}\index{E-edge of
edge-graph@$E$-edge of edge-graph} is defined analogously with
$a\textbf{P}b$ replacing $b\textbf{P}a$.

Let $H_1$, which is $\langle A_1, \textbf{W}_1, \textbf{E}_1, \textbf{P}_1\rangle$, and $H_2$, which is $\langle A_2, \textbf{W}_2, \textbf{E}_2, \textbf{P}_2\rangle$, be two D-edge-graphs such that $C$, which is $A_1\cap A_2$, is non-empty, and for every edge $a$ in $C$ we have that $a$ is an $E$-edge of $H_1$ and a $W$-edge of $H_2$. Then let $H_1\cup H_2$ be $\langle A_1\cup A_2, \textbf{W}_1\cup \textbf{W}_2, \textbf{E}_1\cup \textbf{E}_2, \textbf{P}_1\cup \textbf{P}_2\rangle$. It is straightforward to verify that $H_1\cup H_2$ is a D-edge-graph.

It is also straightforward to verify that there are graphs $G_1$ and $G_2$ isomorphic to ${\cal G}(H_1)$ and ${\cal G}(H_2)$, respectively, such that ${\cal G}(H_1\cup H_2)= G_1\Box G_2$. (All that is involved in passing from ${\cal G}(H_1)$ and ${\cal G}(H_2)$ to $G_1$ and $G_2$ is the renaming of vertices incident with edges that will be in the cocycle of the juncture, in order to ensure the sharing of these vertices for $G_1$ and $G_2$.) Finally, it is straightforward to verify that ${\cal H}(D_W\Box D_E)={\cal H}(D_W)\cup {\cal H}(D_E)$.

\section{\large\bf The system {\rm S}$\Box$}\label{1.5}
\markright{\S 1.5. \quad The system {\rm S}$\Box$}

We will introduce now an equational system for juncture, called S$\Box$, which will have various associativity axioms, and which in \S 1.6 we will show sound and complete with respect to an interpretation in D-graphs.

The equations of S$\Box$ will have on their two sides terms that we will call D-terms. There will be three functions, $W$, $E$ and $A$, mapping the set of D-terms to the power set of an arbitrary infinite set so that for a D-term $\delta$ neither of $W(\delta)$ and $E(\delta)$ is empty. As before, we write $X$ for $W$ or $E$. Intuitively, $X(\delta)$ is the set of $X$-edges of a D-graph for which $\delta$ stands, while $A(\delta)$ is the set of all edges of that D-graph. Hence the sets $W(\delta)$ and $E(\delta)$ will be disjoint, and we will have $W(\delta)\cup E(\delta)\subseteq A(\delta)$.

We define \emph{D-terms}\index{D-term} inductively by starting from \emph{basic D-terms}\index{basic D-term}, which are atomic symbols. To each such symbol $\beta$ we assign three sets $W(\beta)$, $E(\beta)$ and $A(\beta)$ such that $W(\beta)$ and $E(\beta)$ are non-empty, finite and disjoint, while $A(\beta)=W(\beta)\cup E(\beta)$.

The inductive clause of our definition of D-term says that if $\delta_W$ and $\delta_E$ are D-terms such that
\begin{tabbing}
\hspace{1.6em}$C=_{df}A(\delta_W)\cap A(\delta_E)= E(\delta_W)\cap W(\delta_E)\neq\emptyset$,
\end{tabbing}
then $(\delta_W\Box\delta_E)$ is a D-term. As usual, we omit the outermost parentheses of D-terms, and take them for granted.

We define as follows the values of $W$, $E$ and $A$ for the argument $\delta_W\Box\delta_E$:
\begin{tabbing}
\hspace{1.6em}\=$X(\delta_W\Box\delta_E)$ \=$=X(\delta_X)\cup(X(\delta_{\bar{X}})-C)$,\\[.5ex]
\>$A(\delta_W\Box\delta_E)$\>$=A(\delta_W)\cup A(\delta_E)$.
\end{tabbing}
This concludes our definition of D-term.

The triple $(W(\delta),E(\delta),A(\delta))$ will be called the \emph{edge type}\index{edge type of D-term} of the D-term $\delta$. In another notation, we could have written this edge type together with $\delta$ in our language.

Note that if $(\delta_1\Box\delta_2)\Box\delta_3$ is defined---i.e., it is a D-term---then $\delta_2\Box\delta_3$ may be defined or not, but $\delta_3\Box\delta_2$ is never defined, because $E(\delta_3)\cap W(\delta_2)$ must be empty. Otherwise, $A(\delta_1\Box\delta_2)\cap A(\delta_3)$ would not be equal to $E(\delta_1\Box\delta_2)\cap W(\delta_3)$.

The equations of our system, which we call S$\Box$\index{Sbox@S$\Box$}, will be of the form $\delta=\delta'$ for $\delta$ and $\delta'$ being D-terms of the same edge type. The rules of S$\Box$ are symmetry and transitivity of $=$, and congruence with $\Box$:
\begin{tabbing}
\hspace{1.6em}if $\delta_1=\delta_2$ and $\delta_3=\delta_4$, then $\delta_1\Box\delta_3=\delta_2\Box\delta_4$,
\end{tabbing}
provided that $\delta_1\Box\delta_3$ and $\delta_2\Box\delta_4$ are defined.

The axiomatic equations of S$\Box$ are $\delta=\delta$ and the following equations:
\begin{tabbing}
\hspace{1.6em}\=(Ass~1)\hspace{2em}\=$(\delta_1\Box\delta_2)\Box\delta_3$ \=
$=\delta_1\Box(\delta_2\Box\delta_3)$,\index{Ass~1@(Ass~1)}\\[.5ex]
\>(Ass~2.1)\>$(\delta_1\Box\delta_2)\Box\delta_3$\>
$=(\delta_1\Box\delta_3)\Box\delta_2$,\index{Ass~2.1@(Ass~2.1)}\\[.5ex]
\>(Ass~2.2)\>$\delta_1\Box(\delta_2\Box\delta_3)$ \=
$=\delta_2\Box(\delta_1\Box\delta_3)$,\index{Ass~2.2@(Ass~2.2)}
\end{tabbing}
provided that for each of these equations both sides are defined. It is straightforward to verify that in all of these equations the two sides are D-terms of the same edge type.

To help intuition, for (Ass~1) we have the picture
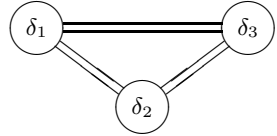
\begin{center}
\begin{picture}(80,50)

\put(0,40){\circle{20}} \put(80,40){\circle{20}}
\put(40,10){\circle{20}}

\put(10,41.5){\line(1,0){60}} \put(10,38.5){\line(1,0){60}}
\put(9,35){\line(4,-3){24}} \put(7,33){\line(4,-3){24}}
\put(71,35){\line(-4,-3){24}} \put(73,33){\line(-4,-3){24}}

\put(0,40){\small\makebox(0,0){$\delta_1$}}
\put(80,40){\small\makebox(0,0){$\delta_3$}}
\put(40,10){\small\makebox(0,0){$\delta_2$}}

\end{picture}
\end{center}
with the direct link between $\delta_1$ and $\delta_3$ at the top of the triangle perhaps missing.  On the left-hand side of (Ass~1) we have joined first $\delta_1$ and $\delta_2$, and then joined the result with $\delta_3$; on the right-hand side of (Ass~1) we have joined first $\delta_2$ and $\delta_3$, and then joined $\delta_1$ with the result.

With (Ass~2.1) and (Ass~2.2) we have the pictures
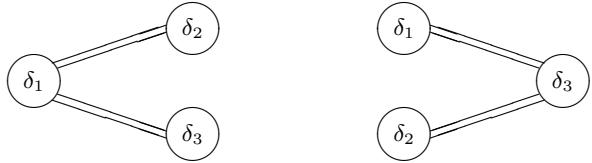
\begin{center}
\begin{picture}(200,60)

\put(0,30){\circle{20}} \put(60,50){\circle{20}}
\put(60,10){\circle{20}}

\put(9,25){\line(3,-1){41}} \put(7,23){\line(3,-1){43}}
\put(9,35){\line(3,1){41}} \put(7,37){\line(3,1){43}}

\put(0,30){\small\makebox(0,0){$\delta_1$}}
\put(60,50){\small\makebox(0,0){$\delta_2$}}
\put(60,10){\small\makebox(0,0){$\delta_3$}}


\put(200,30){\circle{20}} \put(140,50){\circle{20}}
\put(140,10){\circle{20}}

\put(191,25){\line(-3,-1){41}} \put(193,23){\line(-3,-1){43}}
\put(191,35){\line(-3,1){41}} \put(193,37){\line(-3,1){43}}

\put(200,30){\small\makebox(0,0){$\delta_3$}}
\put(140,50){\small\makebox(0,0){$\delta_1$}}
\put(140,10){\small\makebox(0,0){$\delta_2$}}

\end{picture}
\end{center}

Note that for (Ass~1) we may have $\delta_1\Box\delta_3$ defined, but this is not necessary. For (Ass~2.1) we must have that neither $\delta_2\Box\delta_3$ nor $\delta_3\Box\delta_2$ is defined, and for (Ass~2.2) we must have that neither $\delta_1\Box\delta_2$ nor $\delta_2\Box\delta_1$ is defined.

\section{\large\bf The completeness of {\rm S}$\Box$}\label{1.6}
\markright{\S 1.6. \quad The completeness of {\rm S}$\Box$}

We will now interpret the system S$\Box$ in D-graphs, and prove the completeness of S$\Box$ with respect to this interpretation. We introduce an interpretation function $\iota$\index{iota interpretation function@$\iota$ interpretation function} that assigns to a D-term a D-graph, and is defined inductively as follows.

For a basic D-term $\beta$ such that $W(\beta)$ is $\{a_1,\ldots,a_n\}$, for $n\geq 1$, and $E(\beta)$ is $\{b_1,\ldots,b_m\}$, for $m\geq 1$, let $\iota(\beta)$ be the basic D-graph of the following picture
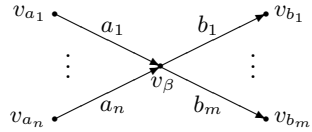
\begin{center}
\begin{picture}(80,50)(0,-5)

\put(0,0){\circle*{2}} \put(0,40){\circle*{2}}
\put(40,20){\circle*{2}} \put(80,0){\circle*{2}}
\put(80,40){\circle*{2}}

\put(0,0){\vector(2,1){40}} \put(0,40){\vector(2,-1){40}}
\put(40,20){\vector(2,1){40}} \put(40,20){\vector(2,-1){40}}

\put(5,23){\makebox(0,0){$\vdots$}}
\put(75,23){\makebox(0,0){$\vdots$}}

\put(41,16){\small\makebox(0,0)[t]{$v_\beta$}}
\put(-3,40){\small\makebox(0,0)[r]{$v_{a_1}$}}
\put(-3,0){\small\makebox(0,0)[r]{$v_{a_n}$}}
\put(83,40){\small\makebox(0,0)[l]{$v_{b_1}$}}
\put(83,0){\small\makebox(0,0)[l]{$v_{b_m}$}}

\put(22,33){\small\makebox(0,0)[b]{$a_1$}}
\put(22,7){\small\makebox(0,0)[t]{$a_n$}}
\put(58,33){\small\makebox(0,0)[b]{$b_1$}}
\put(58,9){\small\makebox(0,0)[t]{$b_m$}}

\end{picture}
\end{center}
We assume that an edge $e$ corresponds bijectively to the vertex
$v_e$, and that this bijection is the same for all basic D-graphs;
it will not vary from one basic D-graph to another. We also have a
bijection assigning the vertex $v_\beta$ to~$\beta$.

To conclude our definition of the function $\iota$ we have the inductive clause
\begin{tabbing}
\hspace{1.6em}$\iota(\delta_W\Box\delta_E)=\iota(\delta_W)\:\Box\:\iota(\delta_E)$,
\end{tabbing}
where $\Box$ on the right-hand side is juncture. As this function, analogous interpretation functions, introduced later, will be homomorphic.

It is straightforward to verify by induction on the number of occurrences of $\Box$ in the D-term $\delta$ that $A(\delta)$ is the set of edges of the D-graph $\iota(\delta)$. For $X_e(D)$\index{Xe@$X_e$} being the set of $X$-edges of the D-graph $D$,
we may also verify, by the same kind of induction, that $X(\delta)=X_e(\iota(\delta))$. (We put the subscript $e$ in $X_e$ because later we use $X(D)$ for the set of $X$-vertices of~$D$.)

With the help of that it becomes straightforward to verify that if $\delta_W\Box\delta_E$ is defined, i.e., it is a D-term, then $\iota(\delta_W)\Box\iota(\delta_E)$ is defined, i.e., it is a D-graph.

The conditions for juncture $\Box$ in D-graphs and for the
operation $\Box$ on D-terms are very similar, but not exactly the
same, because D-graphs make a somewhat more general class than the
D-graphs that are images of D-terms under $\iota$. In defining the
latter, we have introduced, for $\beta$ a basic D-term, a
bijection between the edges and the outer vertices of
$\iota(\beta)$. This bijection, together with the bijection
between basic D-terms and inner vertices, enables us not to
mention vertices when we speak of syntax, i.e.\ when we speak of
D-terms, before introducing their interpretation. This is more
economical, but it imposes a restriction on our general notion of
D-graph. (A similar restriction would be obtained with an
interpretation in edge-graphs; see \S 1.4.)

We can establish easily the following soundness proposition by induction on the length of derivation in the system S$\Box$.

\prop{Proposition 1.6.1}{If in {\rm S}$\Box$ we can derive $\delta=\delta'$, then the D-graphs $\iota(\delta)$ and $\iota(\delta')$ are the same.}\index{soundness of S$\Box$}

Our purpose next is to establish also the converse implication, i.e.\ the completeness of S$\Box$ with respect to $\iota$. For that we need a number of preliminary results, and to state them we need to introduce some terminology.

For a D-graph $D$ and a strict cutset $S$ of $D$, a vertex $D_i$ of the directed graph $C_S(D)$ (see \S 1.3) is \emph{inner}\index{inner vertex of componential graph} when, in accordance with our terminology of \S 1.2, there are two edges of $C_S(D)$ of the form $(D_j,D_i)$ and $(D_i,D_k)$. Otherwise, the vertex is \emph{outer}\index{outer vertex of componential graph}.

The \emph{removal}\index{removal of a vertex of componential graph} of a vertex $D_i$ of $C_S(D)$ leaves a directed graph with the vertices of $C_S(D)$ without $D_i$ and the edges of $C_S(D)$ without the edges in which $D_i$ occurs (cf.\ \cite{H69}, Chapter~2).

The vertex $D_i$ is a \emph{cutvertex}\index{cutvertex} of $C_S(D)$ if the removal of $D_i$ increases the number of weakly connected components of $C_S(D)$. (The notion of weak connectedness for directed graphs, which is analogous to our notion of weak connectedness for graphs of \S 1.2, may be found in \cite{H69}, Chapter~16, as well as the notion of weakly connected component, called there \emph{weak} component; cf.\ also with our notion of component of \S 1.3.)

Let a \emph{componential extreme}\index{componential extreme} of $C_S(D)$ be an outer vertex of $C_S(D)$ that is not a cutvertex. For example, in
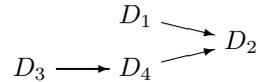
\begin{center}
\begin{picture}(80,40)

\put(10,10){\vector(1,0){20}} \put(50,12.5){\vector(4,1){20}}
\put(50,27.5){\vector(4,-1){20}}

\put(0,10){\makebox(0,0){$D_3$}} \put(40,10){\makebox(0,0){$D_4$}}
\put(80,20){\makebox(0,0){$D_2$}}
\put(40,30){\makebox(0,0){$D_1$}}

\end{picture}
\end{center}
$D_1$ and $D_3$ are componential extremes, while $D_2$, though it is an outer vertex, is not a componential extreme. The notion of componential extreme, and the notions it relies on, are not peculiar to $C_S(D)$. They could be given for arbitrary directed graphs, and the proposition on componential extremes that follows could be established for arbitrary acyclic directed graphs (which are finite by definition). We need it however for $C_S(D)$, and we formulate it accordingly.

\prop{Proposition 1.6.2}{If $S$ is a strict cutset of the D-graph $D$, then there are at least two componential extremes in $C_S(D)$.}

\dkz As a consequence of the acyclicity and finiteness of $C_S(D)$, there must be at least two outer vertices in $C_S(D)$---at least one $W$-vertex and at least one $E$-vertex. We take the outer vertices of $C_S(D)$ to be the vertices of an ordinary graph\index{ordinary graph} in the sense of \cite{H69} (Chapter~2; the edges of ordinary graphs are unordered pairs of distinct vertices), which we call $G_S(D)$; in $G_S(D)$ we have an edge $\{D_i,D_j\}$ when $i\neq j$ and there is a path from $D_i$ to $D_j$ in $C_S(D)$. The definition of path for directed graphs is analogous to the definition of path we gave for graphs in \S 1.2 (see \cite{H69}, Chapter 16). For notions concerning ordinary graphs, like the notions of path and connectedness, we rely on the definitions in \cite{H69} (Chapter 2) (but the definitions of path and connectedness for ordinary graphs is analogous to the definitions of semipath and weak connectedness we gave for graphs in \S 1.2).

That the ordinary graph $G_S(D)$ is connected is shown as follows. Since $C_S(D)$ is weakly connected, we have for every pair of distinct vertices of $G_S(D)$ a semipath of $C_S(D)$ connecting them. It is easy to pass from this semipath to a path of $G_S(D)$ connecting these two vertices, as in the following picture:
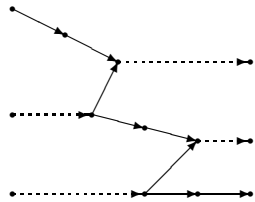
\begin{center}
\begin{picture}(90,70)(10,0)

\put(10,0){\circle*{2}} \put(10,30){\circle*{2}}
\put(10,70){\circle*{2}} \put(40,30){\circle*{2}}
\put(50,50){\circle*{2}} \put(60,0){\circle*{2}}
\put(80,20){\circle*{2}} \put(80,0){\circle*{2}}
\put(60,25){\circle*{2}} \put(30,60){\circle*{2}}
\put(100,0){\circle*{2}} \put(100,20){\circle*{2}}
\put(100,50){\circle*{2}}

\put(10,70){\vector(2,-1){20}} \put(30,60){\vector(2,-1){20}}
\put(40,30){\vector(1,2){10}} \put(40,30){\vector(4,-1){20}}
\put(60,25){\vector(4,-1){20}} \put(60,0){\vector(1,1){20}}
\put(60,0){\vector(1,0){20}} \put(80,0){\vector(1,0){20}}

\multiput(10,0)(3,0){17}{\line(1,0){1}}
\multiput(10,30)(3,0){10}{\line(1,0){1}}
\multiput(80,20)(3,0){7}{\line(1,0){1}}
\multiput(50,50)(3,0){17}{\line(1,0){1}}

\put(57,0){\vector(1,0){2}} \put(37,30){\vector(1,0){2}}
\put(97,20){\vector(1,0){2}} \put(97,50){\vector(1,0){2}}

\end{picture}
\end{center}

Since $G_S(D)$ is connected and has at least two vertices, there
are in $G_S(D)$ two distinct vertices connected by a path of
$G_S(D)$. Take two such vertices $D_i$ and $D_j$ at the greatest
possible distance from each other (this distance is the length of
the shortest path of $G_S(D)$ from $D_i$ to $D_j$; see \cite{H69},
Chapter~2.) If $D_j$ is a cutvertex of $G_S(D)$ (this notion of
cutvertex is analogous to the notion given above for directed
graphs, and may be found, under the name \emph{cutpoint}, in
\cite{H69}, Chapter~3), then its removal leaves a connected
ordinary graph $G'$ (a connected subgraph of $G_S(D)$) in which
$D_i$ is a vertex, and another connected ordinary graph $G''$ in
which we have a vertex $D_k$ different from $D_i$ and $D_j$, such
that there is a path of $G_S(D)$ connecting $D_i$ with $D_k$.
Since $D_j$ must occur in every such path, the distance between
$D_i$ and $D_k$ must be greater than the distance between $D_i$
and $D_j$, which contradicts our assumption that $D_i$ and $D_j$
are at the greatest possible distance. So $D_j$ is not a cutvertex
of $G_S(D)$, and we conclude analogously that $D_i$ is not such
(cf.\ \cite{H69}, Theorem 3.4, Chapter~3).

To conclude that $D_i$ and $D_j$ are not cutvertices of $C_S(D)$, we have the following. Suppose $D_j$ is a cutvertex of $C_S(D)$. Then the removal of $D_j$ from $C_S(D)$ would leave two weakly connected components $H'$ and $H''$ of $C_S(D)$ such that in one of them---let that be $H'$---we have the vertex $D_i$. Since $C_S(D)$ is acyclic and finite, there must be an outer vertex $D_k$ of $C_S(D)$ in $H''$, which is different from $D_i$ and $D_j$. In every semipath of $C_S(D)$ connecting $D_i$ with $D_k$ we find $D_j$. From that it is easy to conclude that in every path of $G_S(D)$ connecting $D_i$ with $D_k$ we find $D_j$. Since $G_S(D)$ is connected, we obtain that $D_j$ is a cutvertex of $G_S(D)$, which contradicts what we have established above. We conclude analogously that $D_i$ is not a cutvertex of $C_S(D)$.\qed

An inner vertex $v$ is an $X$-\emph{border}
vertex\index{X-border
vertex@$X$-border vertex}\index{W-border vertex@$W$-border
vertex}\index{E-border vertex@$E$-border vertex}\index{border
vertex} when for every edge $a$ such that $\bar{X}(a)=v$ we have
that $X(a)$ is an $X$-vertex. For example, in the D-graph of the
following picture:
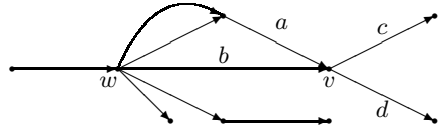
\begin{center}
\begin{picture}(160,60)(0,-18)

\put(0,10){\circle*{2}} \put(40,10){\circle*{2}}
\put(120,10){\circle*{2}} \put(80,30){\circle*{2}}
\put(80,-10){\circle*{2}} \put(160,-10){\circle*{2}}
\put(160,30){\circle*{2}} \put(60,-10){\circle*{2}}
\put(120,-10){\circle*{2}}

\put(0,10){\vector(1,0){40}} \put(40,10){\vector(1,0){80}}
\put(40,10){\vector(2,1){40}}

\put(80,30){\vector(2,-1){40}}

\put(40,10){\vector(2,-1){40}} \put(120,10){\vector(2,1){40}}
\put(120,10){\vector(2,-1){40}}

\qbezier(40,10)(55,45)(80,30) \put(77,31.5){\vector(2,-1){2}}
\put(40,10){\vector(1,-1){20}} \put(80,-10){\vector(1,0){40}}

\put(140,23){\small\makebox(0,0)[b]{$c$}}
\put(140,-2){\small\makebox(0,0)[t]{$d$}}

\put(102,25){\small\makebox(0,0)[b]{$a$}}

\put(80,12){\small\makebox(0,0)[b]{$b$}}
\put(120,7){\small\makebox(0,0)[t]{$v$}}
\put(40,7){\small\makebox(0,0)[tr]{$w$}}

\end{picture}
\end{center}
we have that $v$ is an $E$-border vertex and $w$ is a $W$-border vertex.

Let $X(v)$\index{Xv@$X(v)$}\index{Wv@$W(v)$}\index{Ev@$E(v)$} be the set of all edges such that $\bar{X}(a)=v$. In the example above, $W(v)$ is $\{a,b\}$ and $E(v)$ is $\{c,d\}$.

We say for a non-basic D-graph $D$ that it is
$n$-\emph{valent}\index{n-valent@$n$-valent}, for $n\geq 1$, with
respect to an $X$-border vertex $v$ when for the set $S$ of all
the inner edges in $\bar{X}(v)$, which is a strict cutset, we have
that $C_S(D)$ has $n\pl 1$ vertices. The D-graph in our example
above is $1$-valent with respect to $v$, with the strict cutset
$S$ having two edges, and it is $2$-valent with respect to $w$,
with $S$ now having four edges.

As usual, a subterm\index{subterm} of a D-term is a D-term that occurs in it as a part, not necessarily proper. We have the following.

\prop{Proposition 1.6.3.1}{Suppose the basic D-term $\beta$ is a subterm of the D-term $\delta$, and $v_\beta$ is a $W$-border vertex of $\iota(\delta)$. In {\rm S}$\Box$ we have an equation of the form}
\vspace{-3ex}
\begin{tabbing}
\hspace{1.6em}$\delta=(\ldots(\beta\Box\sigma_1)\Box\ldots)\Box\sigma_n$,
\end{tabbing}
\vspace{-1ex}
\noindent\emph{for $n\geq 0$, where for distinct $i$ and $j$ in $\{1,\ldots,n\}$ we have that $\sigma_i\Box\sigma_j$ is not defined. (If $n=0$, then our equation is $\delta=\beta$.)}

\vspace{2ex}

\dkz We proceed by induction on the number $k$ of occurrences of $\Box$ in $\delta$. If $k=0$, then $\delta$ is~$\beta$.

If $k>0$, then $\delta$ is of the form $\delta_1\Box\delta_2$. If $\beta$ is in $\delta_2$, then $v_\beta$ is a $W$-border vertex in $\iota(\delta_2)$, as well as in $\iota(\delta)$, and by the induction hypothesis we have in S$\Box$
\begin{tabbing}
\hspace{1.6em}$\delta_2=(\ldots(\beta\Box\tau_1)\Box\ldots)\Box\tau_m$,
\end{tabbing}
for $m\geq 1$. We cannot have $\delta_2=\beta$; otherwise, $v_\beta$ would not be a $W$-border vertex. So we have in S$\Box$ the equation $\delta_2=\tau\Box\tau_m$, with $\beta$ in $\tau$, and hence also the equation
\begin{tabbing}
\hspace{1.6em}$\delta=\delta_1\Box(\tau\Box\tau_m)$.
\end{tabbing}
If $\delta_1\Box\tau$ is defined, then, by (Ass~1), in S$\Box$ we have
\begin{tabbing}
\hspace{1.6em}$\delta=(\delta_1\Box\tau)\Box\tau_m$,
\end{tabbing}
and if $\delta_1\Box\tau$ is not defined, then, by (Ass~2.2), in S$\Box$ we have
\begin{tabbing}
\hspace{1.6em}$\delta=\tau\Box(\delta_1\Box\tau_m)$.
\end{tabbing}
So it is enough to consider the case when $\delta$ is of the form $\delta_1\Box\delta_2$ with $\beta$ in $\delta_1$.

Then, by the induction hypothesis, in S$\Box$ we have
\begin{tabbing}
\hspace{1.6em}$\delta_1=(\ldots(\beta\Box\tau_1)\Box\ldots)\Box\tau_m$,
\end{tabbing}
for $m\geq 0$. We will show that in S$\Box$ we have
\begin{tabbing}
\hspace{1.6em}$(\ast)$\hspace{2em}$\delta=((\ldots(\beta\Box\delta_2^\ast)\Box\tau_{i_1})
\Box\ldots)\Box\tau_{i_l}$,
\end{tabbing}
for some $l$ in $\{0,\ldots,m\}$ and $i_1,\ldots,i_l$ in $\{1,\ldots,m\}$, so that for every $j$ in $\{1,\ldots,l\}$ we have that $\delta_2^\ast\Box\tau_{i_j}$ is not defined. If $l=0$, then $(\ast)$ is $\delta=\beta\Box\delta_2^\ast$.

We prove the equation $(\ast)$, which suffices for our proposition, by an auxiliary induction on $m$. If $m=0$, then $\delta$ is $\beta\Box\delta_2$, and we are done. Suppose $m>0$. Then in S$\Box$ we have
\begin{tabbing}
\hspace{1.6em}$\delta=((\ldots(\beta\Box\tau_1)\Box\ldots)\Box\tau_m)\Box\delta_2$.
\end{tabbing}

If $\tau_m\Box\delta_2$ is defined, then, by (Ass~1), in S$\Box$ we have
\begin{tabbing}
\hspace{1.6em}$\delta=((\ldots(\beta\Box\tau_1)\Box\ldots)\Box\tau_{m-1})\Box\delta_2'$
\end{tabbing}
for $\delta_2'$ being $\tau_m\Box\delta_2$. We may then apply the induction hypothesis of the auxiliary induction.

If $\tau_m\Box\delta_2$ is not defined, then, by (Ass~2.1), in S$\Box$ we have
\begin{tabbing}
\hspace{1.6em}$\delta=(((\ldots(\beta\Box\tau_1)\Box\ldots)\Box\tau_{m-1})
\Box\delta_2)\Box\tau_m$,
\end{tabbing}
and we apply the induction hypothesis of the auxiliary induction~to
\begin{tabbing}
\hspace{1.6em}$((\ldots(\beta\Box\tau_1)\Box\ldots)\Box\tau_{m-1})\Box\delta_2$.\`$\dashv$
\end{tabbing}

\vspace{1ex}

We prove analogously the following dual of Proposition 1.6.3.1.

\prop{Proposition 1.6.3.2}{Suppose the basic D-term $\beta$ is a subterm of the D-term $\delta$, and $v_\beta$ is an $E$-border vertex of $\iota(\delta)$. In {\rm S}$\Box$ we have an equation of the form}
\vspace{-3ex}
\begin{tabbing}
\hspace{1.6em}$\delta=\sigma_n\Box(\ldots\Box(\sigma_1\Box\beta)\ldots)$,
\end{tabbing}
\vspace{-1ex}
\noindent\emph{for $n\geq 0$, where for distinct $i$ and $j$ in $\{1,\ldots,n\}$ we have that $\sigma_i\Box\sigma_j$ is not defined. (If $n=0$, then our equation is $\delta=\beta$.)}

\vspace{2ex}

An inner vertex $v$ of a non-basic D-graph $D$ is an $X$-\emph{extreme}\index{X-extreme@$X$-extreme}\index{W-extreme@$W$-extreme}\index{E-extreme@$E$-extreme} when it is an $X$-border vertex and $D$ is $1$-valent with respect to $v$ (see the example before Proposition 1.6.1, where $v$ is $E$-extreme). An \emph{extreme}\index{extreme of D-graph} of $D$ is a $W$-extreme or an $E$-extreme.

Proposition 1.6.2 implies that there are at least two extremes in every non-basic D-graph. For that take as $S$ in Proposition 1.6.2 the set of all inner edges of $D$. The unique inner vertex of $D$ in a componential extreme of $C_S(D)$ is an extreme of $D$. Note that when $v_\beta$ is an extreme, in Propositions 1.6.3.1 and 1.6.3.2 we have $n=1$.

We can now prove the completeness of S$\Box$ with respect to $\iota$.

\prop{Theorem 1.6.4}{In {\rm S}$\Box$ we can derive $\delta=\delta'$ iff the D-graphs $\iota(\delta)$ and $\iota(\delta')$ are the same.}\index{completeness of S$\Box$}

\dkz For the direction from left to right we have Proposition 1.6.1. For the direction from right to left we proceed by induction on the number $k$ of inner vertices in $\iota(\delta)$. If $k=1$, then $\delta$ and $\delta'$ are the same basic D-term.

If $k>1$, then $\iota(\delta)$ is not basic. Take an extreme $v$ of $\iota(\delta)$, and find the basic D-term $\beta$ that is a subterm of $\delta$ and $\delta'$ such that $v$ is $v_\beta$. Suppose $v_\beta$ is $W$-extreme. Then, by Proposition 1.6.3.1, in S$\Box$ we have $\delta=\beta\Box\sigma_1$ and $\delta'=\beta\Box\sigma_1'$. Since, by Proposition 1.6.1, we have that $\iota(\delta)$ is $\iota(\beta)\Box\iota(\sigma_1)$ and $\iota(\delta')$ is $\iota(\beta)\Box\iota(\sigma_1')$, and since $\iota(\delta)$ is $\iota(\delta')$, we must have that $\iota(\sigma_1)$ is $\iota(\sigma_1')$, and, by the induction hypothesis, in S$\Box$ we have $\sigma_1=\sigma_1'$, and hence also $\delta=\delta'$.

We proceed analogously when $v_\beta$ is an $E$-extreme, in which case we apply Proposition 1.6.3.2.\qed

Note that not every D-graph is $\iota(\delta)$ for some D-term $\delta$, but every D-graph is isomorphic to $\iota(\delta)$ for some $\delta$. This may be demonstrated by an easy argument concerning strict cutsets. A strict cutset for a D-graph that is not basic always exists (take, if nothing else, the set of all inner edges, as we did above). An arbitrary strict cutset can easily be reduced to a cocycle. (As a matter of fact, this cocycle may be made to contain an arbitrarily chosen edge of our initial cutset, but we don't need this for our results later on.) Formally, we then make an induction on the number of inner edges of our D-graph.

\section{\large\bf Compatible lists}\label{1.7}
\markright{\S 1.7. \quad Compatible lists}

Let us consider sequences of distinct elements of an arbitrary non-empty set (which later on will be mostly vertices, and sometimes edges), and let such a finite (possibly empty) sequence be called a \emph{list}.\index{list}

For $\Gamma$ a list, let $\Gamma^s$\index{s@$^s$, set of members of a list} be the set of members of $\Gamma$. We say, as expected, that $\Gamma$ is a \emph{list of}\index{list of} $\Gamma^s$.
The lists $\Gamma$ and $\Delta$ are \emph{disjoint}\index{disjoint lists} when $\Gamma^s$ and $\Delta^s$ are disjoint, and the list $\Gamma$ is \emph{empty}\index{empty list} when $\Gamma^s=\emptyset$.

Two non-empty lists are said to be \emph{compatible}\index{compatible lists} when they are either of the forms $\Phi\Xi$ and $\Xi\Psi$ or $\Phi\Xi\Psi$ and $\Xi$, for $\Phi$, $\Xi$ and $\Psi$ mutually disjoint lists, and $\Xi$ a non-empty list. As a particular case, we have that $\Xi$ is compatible with $\Xi$. (Compatibility is, of course, a symmetric relation.)

An alternative definition of compatibility is given as follows. For $\Xi$ a non-empty list, and $\Phi_1$, $\Phi_2$, $\Psi_1$, $\Psi_2$ and $\Xi$ mutually disjoint lists, the lists $\Phi_1\Xi\Psi_1$ and $\Phi_2\Xi\Psi_2$ are \emph{compatible} when at least one $\Phi_1$ and $\Phi_2$, and at least one $\Psi_1$ and $\Psi_2$, are empty lists.

The conditions we have in these definitions, and in particular the disjointness conditions, ensure that with compatible lists we have a \emph{unified} list\index{unified list} $\Phi\Xi\Psi$ with the first definition, and $\Phi_1\Phi_2\Xi\Psi_1\Psi_2$ with the second definition.

For every non-empty list $\Gamma$ of the form $x_1\ldots x_n$, with $n\geq 1$, consider the set of ordered pairs defined for $n>1$ by
\begin{tabbing}
\hspace{1.6em}$R_\Gamma=\{(x_i,x_{i+1})\mid 1\leq i\leq n\mn 1\}$,
\end{tabbing}
while for $n=1$ we have $R_\Gamma=\emptyset$. If $\Gamma$ is the empty list, then $R_\Gamma$ is again $\emptyset$.

For $\Gamma$ a list, let $\langle\Gamma^s, R_\Gamma\rangle$, i.e.\ the binary relation $R_\Gamma$ on $\Gamma^s$, be called a \emph{chain}.\index{chain} Lists correspond bijectively to chains. (The one member list $x$ corresponds to the chain $\langle\{x\}, \emptyset\rangle$, and the empty list corresponds to the chain $\langle\emptyset, \emptyset\rangle$.) Chains will serve to give another, quite natural, definition of compatibility.

We say that the chains $\langle\Gamma_1^s, R_{\Gamma_1}\rangle$ and $\langle\Gamma_2^s, R_{\Gamma_2}\rangle$ are \emph{compatible}\index{compatible chains} when there is a list $\Gamma$ such that
\begin{tabbing}
\hspace{1.6em}$\langle\Gamma^s, R_\Gamma\rangle=\langle\Gamma_1^s\cup \Gamma_2^s, R_{\Gamma_1}\cup R_{\Gamma_2}\rangle$.
\end{tabbing}
We will not go into the rather straightforward proofs of the following propositions, which show that the compatibility of lists and the compatibility of the corresponding chains are in complete agreement. These propositions are not essential for our results later~on.

\prop{Proposition 1.7.1}{If the lists $\Gamma_1$ and $\Gamma_2$ are compatible, then the chains $\langle\Gamma_1^s, R_{\Gamma_1}\rangle$ and $\langle\Gamma_2^s, R_{\Gamma_2}\rangle$ are compatible.}

\vspace{-1ex}

\prop{Proposition 1.7.2}{If the chains $\langle\Gamma_1^s, R_{\Gamma_1}\rangle$ and $\langle\Gamma_2^s, R_{\Gamma_2}\rangle$ are compatible, then}

\vspace{-3.5ex}

\begin{tabbing}
\hspace{1.6em}\=(1)\hspace{2em}\=$C=_{df}\Gamma_1^s\cap\Gamma_2^s\neq\emptyset$,\\[.5ex]
\>(2)\>$\langle C,R_{\Gamma_1}\cap C^2\rangle$ \emph{and} $\langle C,R_{\Gamma_2}\cap C^2\rangle$ \emph{are chains},\\[.5ex]
\>(3)\>$R_{\Gamma_1}\cap C^2=R_{\Gamma_2}\cap C^2$,\\[.5ex]
\>(4)\>\emph{there are no} $x,y\notin C$ \emph{and a} $z\in C$ \emph{such that}\\[.5ex] \>\>($xR_{\Gamma_1}z$ \emph{and} $yR_{\Gamma_2}z$) \emph{or} ($zR_{\Gamma_1}x$ \emph{and} $zR_{\Gamma_2}y$).
\end{tabbing}

\vspace{-1ex}

\prop{Proposition 1.7.3}{If for the chains $\langle\Gamma_1^s, R_{\Gamma_1}\rangle$ and $\langle\Gamma_2^s, R_{\Gamma_2}\rangle$ we have {\rm (1), (2), (3)} and {\rm (4)} of Proposition 1.7.2, then the lists $\Gamma_1$ and $\Gamma_2$ are compatible.}

\section{\large\bf P$'$-graphs}\label{1.8}
\markright{\S 1.8. \quad P$\,'$-graphs}

Our purpose now is to define a kind of D-graph realizable in a particular manner in the plane (see Chapter~7). First, in this section and in \S 1.9, we will have two inductive definitions, which will yield the notions of P$'$-graph and P$''$-graph. Then in \S 1.10 we will have a non-inductive definition, which will yield the notion of P$'''$-graph. All these definitions are based on juncture. In Chapters 2-5 we will show that these three notions cover the same graphs, which we will call \emph{P-graphs}\index{P-graph}.

A \emph{construction of a} \emph{P$\,'$-graph}\index{construction of a P1-graph@construction of a P$'$-graph} (for short, \emph{construction}\index{construction}) is a finite binary tree such that in each node we have a triple $(D,L_W,L_E)$ where $D$ is a D-graph and $L_X$\index{LX@$L_X$}\index{LW@$L_W$}\index{LE@$L_E$}, for $X$ being $W$ or $E$, is a list of all the $X$-vertices of $D$, the set of which is designated by $X(D)$\index{Xd@$X(D)$}\index{Wd@$W(D)$}\index{Ed@$E(D)$}. For the triple $(D,L_W,L_E)$ at the root of a construction $K$ we call $D$ the \emph{root graph}\index{root graph} of $K$, while $L_W$ and $L_E$ are the \emph{root lists}\index{root list} of $K$. We say that a construction is a construction of its root graph.

Here are the two inductive clauses of our definition of construction:
\begin{itemize}
\item[(1)]The single-node tree in whose single node we have a basic D-graph  (see the end of \S 1.2) together with an arbitrary list of all of its $W$-vertices and an arbitrary list of all of its $E$-vertices is a construction;
\item[(2)]For $X$ being $W$ or $E$, let $K_X$ be a construction that in its root has $(D_X,L_W^X,L_E^X)$ so that the lists $L_E^W$ and $L_W^E$ are compatible (see \S 1.7). Out of $K_W$ and $K_E$ we obtain a new construction $K_W\Box K_E$ by adding a new node to serve as its root, whose successors are the roots of $K_W$ and $K_E$; in the new root we have $(D_W\Box D_E,L_W,L_E)$, where if $L_E^W$ is $\Phi_E\Xi\Psi_E$ and $L_W^E$ is $\Phi_W\Xi\Psi_W$, then $L_X$ is $\Phi_XL_X^X\Psi_X$.
\end{itemize}
The compatibility of $L_E^W$ and $L_W^E$ in clause (2) implies that at least one of $\Phi_E$ and $\Phi_W$, and at least one of $\Psi_E$ and $\Psi_W$, are empty lists (see \S 1.7).

A \emph{P$\,'$-graph}\index{P1-graph@P$'$-graph} is the root graph of a construction.

Note that this definition could have relied on lists of the $X$-edges instead of the $X$-vertices of a D-graph $D$, because the $X$-edges and the $X$-vertices of $D$ correspond bijectively to each other. For some of our purposes
concentrating on the vertices seems better, and more natural, while
for other purposes it is easier to concentrate on the edges. On a
few occasions (see, for example, the proof of Proposition 2.2.1), it may seem unnecessarily tedious to the reader to pass
from one point of view to the other, but we believe that any
exposition of our subject matter would have if not this some other
kind of shortcoming.

\section{\large\bf P$''$-graphs}\label{1.9}
\markright{\S 1.9. \quad P$\,''$-graphs}

For $u$ and $v$ vertices of a D-Graph $D$, let $[u,v]$ be the set of all semipaths from $u$ to $v$. (This set is, of course, in a bijection with $[v,u]$.) Let
\begin{tabbing}
\hspace{1.6em}$[u]_X=\bigcup\,\{[u,v]\mid v\in X(D)\}$.
\end{tabbing}

For example, in
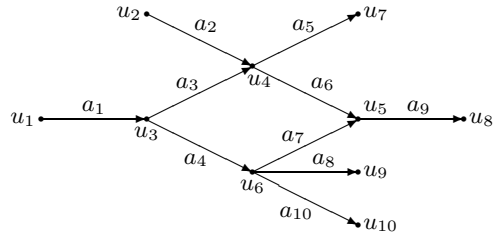
\begin{center}
\begin{picture}(160,90)(0,-5)

\put(0,40){\circle*{2}} \put(40,40){\circle*{2}}
\put(40,80){\circle*{2}} \put(80,20){\circle*{2}}
\put(80,60){\circle*{2}} \put(120,0){\circle*{2}}
\put(120,20){\circle*{2}} \put(120,40){\circle*{2}}
\put(120,80){\circle*{2}} \put(160,40){\circle*{2}}

\put(0,40){\vector(1,0){40}} \put(120,40){\vector(1,0){40}}
\put(40,40){\vector(2,1){40}} \put(40,40){\vector(2,-1){40}}

\put(40,80){\vector(2,-1){40}}

\put(80,20){\vector(2,1){40}} \put(80,20){\vector(2,-1){40}}

\put(80,60){\vector(2,-1){40}}

\put(80,20){\vector(1,0){40}} \put(80,60){\vector(2,1){40}}

\put(-2,40){\small\makebox(0,0)[r]{$u_1$}}
\put(38,80){\small\makebox(0,0)[r]{$u_2$}}
\put(40,37){\small\makebox(0,0)[t]{$u_3$}}
\put(78,57){\small\makebox(0,0)[tl]{$u_4$}}
\put(122,43){\small\makebox(0,0)[bl]{$u_5$}}
\put(80,17){\small\makebox(0,0)[t]{$u_6$}}
\put(122,80){\small\makebox(0,0)[l]{$u_7$}}
\put(162,40){\small\makebox(0,0)[l]{$u_8$}}
\put(122,20){\small\makebox(0,0)[l]{$u_9$}}
\put(122,0){\small\makebox(0,0)[l]{$u_{10}$}}

\put(20,43){\small\makebox(0,0)[b]{$a_1$}}

\put(58,74){\small\makebox(0,0)[bl]{$a_2$}}

\put(60,52){\small\makebox(0,0)[br]{$a_3$}}
\put(63,27){\small\makebox(0,0)[tr]{$a_4$}}
\put(95,74){\small\makebox(0,0)[bl]{$a_5$}}

\put(102,52){\small\makebox(0,0)[bl]{$a_6$}}

\put(100,32){\small\makebox(0,0)[br]{$a_7$}}
\put(107,23){\small\makebox(0,0)[b]{$a_8$}}
\put(143,43){\small\makebox(0,0)[b]{$a_9$}}
\put(103,7){\small\makebox(0,0)[tr]{$a_{10}$}}

\end{picture}
\end{center}
we have
\begin{tabbing}
\hspace{1.6em}$[u_7,u_9]=\{u_7a_5u_4a_6u_5a_7u_6a_8u_9, u_7a_5u_4a_3u_3a_4u_6a_8u_9\}$,\\[.55ex]
\hspace{1.6em}$[u_8]_W=\{u_8a_9u_5a_6u_4a_3u_3a_1u_1, u_8a_9u_5a_7u_6a_4u_3a_1u_1,u_8a_9u_5a_6u_4a_2u_2,$\\
\`$ u_8a_9u_5a_7u_6a_4u_3a_3u_4a_2u_2\}$.
\end{tabbing}

We say that two semipaths \emph{intersect}\index{intersecting semipaths} when they have a common vertex.

Let $u$, $v$ and $w$ be distinct $X$-vertices of a D-graph. We write $\psi_X(v,u,w)$\index{psix@$\psi_X$} when every semipath in $[v,w]$ and every semipath in $[u]_{\bar{X}}$ intersect. It is clear that $\psi_X(v,u,w)$ implies $\psi_X(w,u,v)$.

In the example above, we have $\psi_E(u_7,u_8,u_9)$ and not $\psi_E(u_7,u_9,u_8)$, because we have $u_7a_5u_4a_6u_5a_9u_8$ in $[u_7,u_8]$ and $u_9a_8u_6a_4u_3a_1u_1$ in $[u_9]_W$.

For $n\geq 3$, we write $\Gamma\!:x_1\mn x_2\mn x_3\mn\ldots\mn x_n$ to assert that in the list $\Gamma$ the distinct members $x_1,x_2,x_3,\ldots,x_n$ occur either in that order or in the order $x_n,\ldots,x_3,x_2,x_1$, where for $i$ in $\{1,2,\ldots,n\mn 1\}$ the members $x_i$ and $x_{i+1}$ are not necessarily immediate neighbours. For example, if $\Gamma$ is $75465983$, then we have $\Gamma\!:7\mn 6\mn 9\mn 8$ and $\Gamma\!:8\mn 9\mn 6\mn 7$.

For $D$ a D-graph, we say that a list $\Lambda$ of $X(D)$ is \emph{grounded}\index{grounded list} in $D$ when for every $v$, $u$ and $w$ in $\Lambda^s$ if $\Lambda\!:v\mn u\mn w$, then $\psi_X(v,u,w)$.

In our example, we have that $u_7u_8u_9u_{10}$ and $u_7u_8u_{10}u_9$ are grounded, while $u_7u_9u_8u_{10}$ and $u_7u_{10}u_9u_8$ are not grounded.

For two D-graphs $D_W$ and $D_E$ such that $D_W\Box D_E$ is defined, we say that they are \emph{P-compatible}\index{P-compatible D-graphs} when a list of $E(D_W)$ grounded in $D_W$ and a list of $W(D_E)$ grounded in $D_E$ are compatible.

A \emph{P$\,''$-graph}\index{P2-graph@P$''$-graph} is defined inductively by the following two clauses:
\begin{itemize}
\item[(1)]basic D-graphs are P$''$-graphs;\vspace{-1ex}
\item[(2)]if $D_W$ and $D_E$ are P-compatible P$''$-graphs, then $D_W\Box D_E$ is a P$''$-graph.
\end{itemize}

The remainder of this section is an appendix, which is not
essential for the later exposition, and can hence be skipped. We
defined in \S 1.6 an $X$-border vertex as an inner vertex $x$ such
that for \emph{every} edge $a$ where $\bar{X}(a)=v$ we have that
$X(a)$ is an $X$-vertex. Let an $X$-\emph{peripheral}
vertex\index{X-peripheral vertex@$X$-peripheral
vertex}\index{W-peripheral vertex@$W$-peripheral
vertex}\index{E-peripheral vertex@$E$-peripheral
vertex}\index{peripheral vertex} be an inner vertex $x$ such that for
\emph{some} edge $a$ where $\bar{X}(a)=v$ we have that $X(a)$ is
an $X$-vertex. Every $X$-border vertex is an $X$-peripheral
vertex, but not necessarily vice versa.

We show in this appendix that our notion of grounding, based on
the ternary relation $\psi_X$ on $X$-vertices, could be replaced
by an equivalent notion based on a ternary relation on
$X$-peripheral vertices. The interest of this is that it
contributes to showing that in D-graphs inner vertices are
essential. Vertices that are not inner play a secondary role.

Let $u$, $v$ and $w$ be $X$-peripheral vertices of a D-graph, not
necessarily distinct. We write
$\psi^b_X(v,u,w)$\index{psixb@$\psi^b_X$} when every semipath in
$[v,w]$ and every semipath in $[u]_{\bar{X}}$ intersect.

Take an $X$-vertex $x$, and consider the edge $a$ such that
$X(a)=x$. Then we say that the $X$-peripheral vertex $\bar{X}(a)$
is the \emph{mate}\index{mate} of $x$, which we designate by
$m(x)$.\index{m function@$m$} The function $m$ from $X$-vertices
to $X$-peripheral vertices is onto, but not one-one. We can prove
the following for every D-graph $D$ and every distinct vertices
$v$, $u$ and $w$ in $X(D)$.

\prop{Proposition 1.9.1}{We have $\psi_X(v,u,w)$ iff $\psi^b_X(m(v),m(u),m(w))$.}

\dkz For the proof from left to right, suppose we have a semipath $\sigma$ in $[m(v),m(w)]$ and a semipath $\tau$ in $[m(u)]_{\bar{X}}$. We extend $\sigma$ to $\sigma^+$ in $[v,w]$ just by adding two edges at the ends and the vertices $v$ and $w$, and we extend $\tau$ to $\tau^+$ in $[u]_{\bar{X}}$ just by adding one edge and the vertex $u$. Since $\psi_X(v,u,w)$, we have that $\sigma^+$ and $\tau^+$ intersect, but since $v$, $u$ and $w$ are distinct vertices, we obtain that $\sigma$ and $\tau$ intersect.

For the proof from right to left it is enough to remark that for every semipath $\rho$ in $[v,w]$ we have that $m(v)$ and $m(w)$ occur in $\rho$, and for every semipath $\pi$ in $[u]_{\bar{X}}$ we have that $m(u)$ occurs in $\pi$. Let $\rho^-$ in $[m(v),m(w)]$ be obtained from $\rho$ by rejecting the vertices $v$ and $w$ and two edges at the ends incident with $v$ and $w$ respectively, and let $\pi^-$ in $[m(u)]_{\bar{X}}$ be obtained for $\pi$ by rejecting $u$ and the edge incident with $u$. Since $\psi^b_X(m(v),m(u),m(w))$, we have that $\rho^-$ and $\pi^-$ intersect, and hence $\rho$ and $\pi$ intersect.\qed

\vspace{-2ex}

\section{\large\bf P$'''$-graphs}\label{1.10}
\markright{\S 1.10. \quad P$\,'''$-graphs}

Consider a cocycle $C$ of a D-graph $D$. Let the removal (see the beginning of \S 1.3) of $C$ from $D$ leave the graphs $D_1$ and $D_2$ such that for an edge $a$ in $C$ we have $W(a)$ in $D_1$ and $E(a)$ in $D_2$. Out of $D_1$ we build a D-graph $D_W$ by adding to the edges of $D_1$ all the edges in $C$, and by stipulating that for every $a$ in $C$ we have $W(a)$ equal to what it was in $D$ , while $E(a)$ is a new vertex $v_a$, which we add to the vertices of $D_1$ for every edge $a$ in $C$. We build the D-graph $D_E$ analogously out of $D_2$, by adding again the edges of $C$, and all the new vertices $v_a$ we have added to $D_1$ to obtain $D_W$; now we have $E(a)$ as in $D$, while $W(a)$ is $v_a$.

We say that $D_W$ and $D_E$ are obtained by \emph{cutting} $D$ \emph{through} $C$.\index{cutting a D-graph through a cocycle} It is obvious that $D$ is $D_W\Box D_E$.

A D-graph $D$ is a \emph{P$\,'''$-graph}\index{P3-graph@P$'''$-graph} when for every cocycle $C$ of $D$ the D-graphs $D_W$ and $D_E$ obtained by cutting $D$ through $C$ are P-compatible (see \S 1.9 for P-compatibility).

\clearpage \pagestyle{empty} \makebox[1em]{} \clearpage

\chapter{\huge\bf P$'$-Graphs and P$'''$-Graphs}\label{2}
\pagestyle{myheadings}\markboth{CHAPTER 2. \quad P$\,'$-GRAPHS AND
P$\,'''$-GRAPHS}{right-head}

\section{\large\bf Interlacing and parallelism}\label{2.1}
\markright{\S 2.1. \quad Interlacing and parallelism}

In this chapter our goal is to prove that every P$'$-graph (as defined in \S 1.8) is a P$'''$-graph (as defined in \S 1.10). Before achieving that in \S 2.3, we deal with preliminary matters. In this section we start with combinatorial matters concerning lists (see \S 1.7).

When in a list $L$ we have that $x$ and $y$ are immediate
neighbours---i.e., $x$ is the immediate predecessor or the
immediate successor of $y$---we say that $x$ and $y$ are
\emph{$L$-neighbours}\index{neighbours in a list}.

Let $A$ be the list $a_{k+1}\ldots a_{k+n}$, for $k\geq 0$ and
$n\geq 1$, and let $B$ be the list $b_{l+1}\ldots b_{l+m}$, for
$l\geq 0$ and $m\geq 1$. (We need $k$ and $l$ because we will
have lists where indexing does not start from 1; see $A'$ in the
lemmata below.) Assume the sets $\{a_{k+1},\ldots,a_{k+n}\}$ and
$\{b_{l+1},\ldots,b_{l+m}\}$ are disjoint, and let $M$ be a list
of the union of these two sets.

We have that in $A$ and $B$ respectively $a_{k+1}$ and $b_{l+1}$
are \emph{initial},\index{initial member of a list} while $a_{k+n}$ and $b_{l+m}$ are
\emph{final}.\index{final member of a list} For $F$ being $A$ or $B$, let $u^p$ be the immediate
predecessor of $u$ in $F$, provided this predecessor exists, i.e.,
$u$ is not initial in $F$, and let $u^s$ be the immediate
successor of $u$ in $F$, provided this successor exists, i.e., $u$
is not final in~$F$.

We say that two members $u$ and $v$ of $M$ are of the same
\emph{parity},\index{parity of members lists} and write $u\equiv_{\Pi}v$, when their indices are
either both even or both odd. For example, we have
$a_3\equiv_{\Pi}b_7$, as well as $a_3\equiv_{\Pi}a_{17}$.

Take two distinct members $u$ and $v$ of a list $M$, which are either one
in $A$ and other in $B$, or both in $A$, or both in $B$, and
assume that $v^s$ exists. We say that $u$ is
\emph{interlaced}\index{interlaced members of lists} in $M$ with $v$
and $v^s$, and write $M[v,u,v^s]$, when $M\!:v\mn u\mn v^s$ and
\begin{tabbing}
\hspace{1.6em}\=(1)\hspace{2em}\=if $u\equiv_\Pi v^s$ and $u^p$
exists, then not $M\!:v\mn u^p\mn v^s$, and
\\*[.5ex]
\>(2)\>if $u\equiv_\Pi v$ and $u^s$ exists, then not $M\!:v\mn
u^s\mn v^s$.
\end{tabbing}
When $M$ is clear from the context, we may omit ``in $M$'' from
``interlaced in~$M$''.

Note that clause (1) is trivially satisfied when $u^p$ does not
exist, and (2) is trivially satisfied when $u^s$ does not exist.
For example, suppose we have $M\!:a_1\mn b_4\mn a_2$; then we have
$b_4$ interlaced with $a_1$ and $a_2$ if either $M\!:b_3\mn a_1\mn
a_2$, or $M\!:a_1\mn a_2\mn b_3$, or $b_4$ is initial in $B$. We
do not have $b_4$ interlaced with $a_1$ and $a_2$ if $M\!:a_1\mn
b_3\mn a_2$.

To help the intuition, let us draw the list $M$ vertically. On the
right of the line of $M$ let us draw lines connecting the
successive members $u$ and $u^s$ of $A$ and $B$ where $u$ has an
odd index, and on the left let us draw lines connecting $u$ and
$u^s$ where $u$ has an even index. These lines make $A$ and~$B$.

For example, in the first of these two pictures, with intersecting
lines on the right, $b_4$ is interlaced with $a_1$ and $a_2$,
while in the second it is not, and lines on the right do not
intersect:
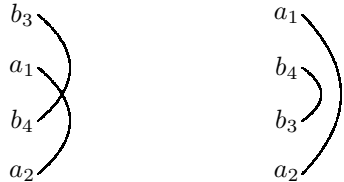
\begin{center}
\begin{picture}(100,70)(0,5)

\qbezier(6,70)(30,50)(6,30) \qbezier(6,50)(30,30)(6,10)

\put(0,10){\makebox(0,0){$a_2$}} \put(0,30){\makebox(0,0){$b_4$}}
\put(0,50){\makebox(0,0){$a_1$}} \put(0,70){\makebox(0,0){$b_3$}}

\qbezier(106,70)(135,40)(106,10) \qbezier(106,50)(120,40)(106,30)

\put(100,10){\makebox(0,0){$a_2$}}
\put(100,30){\makebox(0,0){$b_3$}}
\put(100,50){\makebox(0,0){$b_4$}}
\put(100,70){\makebox(0,0){$a_1$}}

\end{picture}
\end{center}
In the first picture, we have also $a_1$ interlaced with
$b_3$ and~$b_4$.

Nothing changes when we replace $b_3$ and $b_4$ by $a_3$ and $a_4$
respectively. When $b_3$ is final we may draw a horizontal line
from it to the right to ensure intersection, which indicates
interlacing:
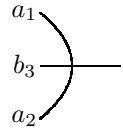
\begin{center}
\begin{picture}(0,50)(0,5)

\qbezier(6,50)(30,30)(6,10) \put(6,30){\line(1,0){30}}

\put(0,10){\makebox(0,0){$a_2$}} \put(0,30){\makebox(0,0){$b_3$}}
\put(0,50){\makebox(0,0){$a_1$}}

\end{picture}
\end{center}
For final members with an even index, the horizontal
line would go to the left, and for initial members we have dual
conventions (see the example below).

For a more involved example, let $A$ be $a_1a_2a_3$, let $B$ be
$b_5b_6b_7b_8b_9$, and let $M$ be $b_9b_8b_5a_2b_6a_3a_1b_7$:
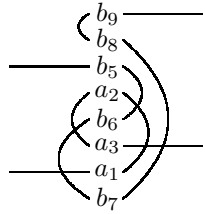
\begin{center}
\begin{picture}(0,80)(0,5)

\qbezier(-7,40)(-30,25)(-7,10) \qbezier(-7,50)(-20,40)(-7,30)
\qbezier(-7,80)(-15,75)(-7,70) \put(-7,20){\line(-1,0){30}}
\put(-7,60){\line(-1,0){30}}

\qbezier(6,70)(40,40)(6,10) \qbezier(6,60)(20,50)(6,40)
\qbezier(6,50)(25,35)(6,20) \put(6,30){\line(1,0){30}}
\put(6,80){\line(1,0){30}}

\put(0,10){\makebox(0,0){$b_7$}} \put(0,20){\makebox(0,0){$a_1$}}
\put(0,30){\makebox(0,0){$a_3$}} \put(0,40){\makebox(0,0){$b_6$}}
\put(0,50){\makebox(0,0){$a_2$}} \put(0,60){\makebox(0,0){$b_5$}}
\put(0,70){\makebox(0,0){$b_8$}} \put(0,80){\makebox(0,0){$b_9$}}

\end{picture}
\end{center}
As before, the intersections of the $A$ and $B$ lines
indicate interlacing. For example, $a_3$ is interlaced with $b_6$
and $b_7$, as well as with $a_1$ and $a_2$, and with $b_7$ and~$b_8$.

Let $F,G\in\{A,B\}$ (so $F$ may be either different from or equal
to $G$). We say that $F$ and $G$ are
\emph{parallel}\index{parallel lists@parallel lists} in $M$, and write
$F\parallel_M G$, when no member of $F$ is interlaced in $M$ with
two successive members of $G$ and no member of $G$ is interlaced
in $M$ with two successive members of $F$. (One of the two
conjuncts in this definition does not entail the other when
initial and final members of $A$ and $B$ are involved; for
example, if $A$ is $a_1a_2$ and $B$ is $b_2b_3$, with $M$ being
$b_3a_1a_2b_2$, we have that both $a_1$ and $a_2$ are interlaced
with $b_2$ and $b_3$, but neither $b_2$ nor $b_3$ is interlaced
with $a_1$ and~$a_2$.)

To obtain an example of parallelism, let $A$ be $a_1a_2a_3$, let
$B$ be $b_5b_6b_7b_8b_9$, and let $M$ be
$b_5a_3a_2b_4b_1b_2b_3a_1$:
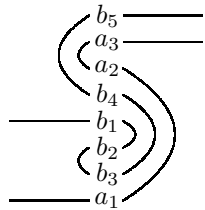
\begin{center}
\begin{picture}(0,80)(0,5)

\qbezier(-7,30)(-15,25)(-7,20) \qbezier(-7,80)(-30,65)(-7,50)
\qbezier(-7,70)(-15,65)(-7,60) \put(-7,10){\line(-1,0){30}}
\put(-7,40){\line(-1,0){30}}

\qbezier(6,60)(45,35)(6,10) \qbezier(6,50)(30,35)(6,20)
\qbezier(6,40)(15,35)(6,30) \put(6,70){\line(1,0){30}}
\put(6,80){\line(1,0){30}}

\put(0,10){\makebox(0,0){$a_1$}} \put(0,20){\makebox(0,0){$b_3$}}
\put(0,30){\makebox(0,0){$b_2$}} \put(0,40){\makebox(0,0){$b_1$}}
\put(0,50){\makebox(0,0){$b_4$}} \put(0,60){\makebox(0,0){$a_2$}}
\put(0,70){\makebox(0,0){$a_3$}} \put(0,80){\makebox(0,0){$b_5$}}

\end{picture}
\end{center}
It is straightforward to check that we have
$A\parallel_M B$, $A\parallel_M A$ and $B\parallel_M B$.

For $A$ having at least three members, assume two $A$-neighbours
$a_i$ and $a_{i+1}$ are also $M$-neighbours, and let the lists
$A'$ and $M'$ be obtained by omitting $a_i$ and $a_{i+1}$ from $A$
and $M$ respectively. We can prove the following.

\prop{Lemma 2.1.1.1}{If $A\parallel_M B$, then
$A'\parallel_{M'}B$.}

\dkz Suppose we have $a_{i-1}$ and $a_{i+2}$ in $A$, i.e., $a_i$
is not initial and $a_{i+1}$ is not final in $A$. For $w$ being
$a_i$ and $a_{i+1}$, we must have one of the following:
\begin{tabbing}
\hspace{1.6em}\=(I)\hspace{2em}\=$M\!:a_{i-1}\mn w\mn a_{i+2}$,
\\[.5ex]
\>(II)\>$M\!:w\mn a_{i-1}\mn a_{i+2}$,
\\[.5ex]
\>(III)\>$M\!:a_{i-1}\mn a_{i+2}\mn w$.
\end{tabbing}
Suppose not $A'\parallel_{M'}B$. We will infer that not
$A\parallel_M B$.

\vspace{1ex}

($B$ in $A'$) We consider first that $A'\parallel_{M'} B$ fails
because for some member $u$ of $B$ and some members $v$ and $v^s$
of $A'$ we have that $M'[v,u,v^s]$. If $v$ is different from $a_{i-1}$,
then we obtain easily that $M[v,u,v^s]$. If $v$ is $a_{i-1}$, then
$v^s$ is $a_{i+2}$, and we have to consider separately the three
cases (I)-(III) above.

In case (I) we can infer easily that either $M[a_{i-1},u,a_i]$ or
$M[a_{i+1},u,a_{i+2}]$. In case (II) we have several subcases to
consider.

\vspace{1ex}

(II.1) Suppose $u\equiv_\Pi a_{i+2}$ and $u^p$ exists. We may have
\begin{tabbing}
\hspace{1.6em}\=$M\!:u^p\mn w\mn a_{i-1}\mn u\mn a_{i+2}$, or
\\[.5ex]
\>$M\!:w\mn a_{i-1}\mn u\mn a_{i+2}\mn u^p$, or
\\[.5ex]
\>$M\!:w\mn u^p\mn a_{i-1}\mn u\mn a_{i+2}$.
\end{tabbing}
In the first two cases we conclude that $M[a_{i+1},u,a_{i+2}]$,
while in the third case we obtain that $M[a_{i-1},u^p,a_i]$.

If $u\equiv_\Pi a_{i+2}$ and $u^p$ does not exist, then
$M[a_{i+1},u,a_{i+2}]$.

\vspace{1ex}

(II.2) Suppose $u\equiv_\Pi a_{i-1}$ and $u^s$ exists. Then we
have again three cases as in (II.1), obtained by substituting $u^s$ for
$u^p$, and we continue reasoning analogously to what we had in
(II.1). In case we have (III), we reason analogously to what we had
for (II) above.

\vspace{1ex}

($A'$ in $B$) We consider now that $A'\parallel_{M'} B$ fails
because for some member $z$ of $A'$ and some members $u$ and $u^s$
of $B$ we have that $M'[u,z,u^s]$. Suppose $z$ is~$a_{i-1}$.

If $a_{i-1}\equiv_\Pi u^s$, then we conclude easily
that $M[u,a_{i-1},u^s]$. If $a_{i-1}\equiv_\Pi u$, then either we have
that $M[u,a_{i-1},u^s]$ or, in case that we have $M\!:u\mn w\mn u^s$,
we have that $M[u,a_{i+1},u^s]$.

We reason analogously if we suppose that $z$ is $a_{i+2}$. This
concludes the proof of $A'\parallel_{M'}B$, under the assumption
in the first sentence of the proof.

In case in $A$ we have $a_{i-1}$ but not $a_{i+2}$, or $a_{i+2}$
but not $a_{i-1}$, we reason by simplifying the reasoning we had
above. From ($B$ in $A'$) we keep just the easy case when $v$ is
different from $a_{i-1}$, while from ($A'$ in $B$) we keep just
the case when $z$ is $a_{i-1}$, or just the case when $z$ is
$a_{i+2}$. Since $A$ has at least three members, one of $a_{i-1}$
and $a_{i+2}$ must exist. \qed

\vspace{-2ex}

\prop{Lemma 2.1.1.2}{If $A\parallel_M A$, then $A'\parallel_{M'}
A'$.}

\dkz The proof is analogous to the proof of Lemma 2.1.1.1; we make
only some obvious adaptations. Note that in (I) we may replace $w$
by $a_i\mn a_{i+1}$, but not by $a_{i+1}\mn a_1$, because
$A\parallel_M A$; analogously, in (II) and (III) we may replace
$w$ by $a_{i+1}\mn a_i$, but not by $a_i\mn a_{i+1}$. This does
not however influence essentially the exposition of the proof.
\qed

The following holds for $A'$ being any subset of $A$.

\prop{Lemma 2.1.1.3}{If $B\parallel_M B$, then $B\parallel_{M'}
B$.}

For the lemmata that follow we assume that $A\parallel_M B$,
$A\parallel_M A$ and $B\parallel_M B$.

\prop{Lemma 2.1.2}{If $A$ or $B$ has at least two members, then
two $A$-neighbours or two $B$-neighbours are $M$-neighbours.}

\dkz For $v$ and $w$ distinct members of $M$, let $d_M(v,w)$
be the number of members of $M$ between $v$ and $w$ and let $k$ be
the minimal number in the set
\begin{tabbing}
\hspace{1.6em}$S=\{d_M(v,w)\mid v$ and $w$ are $A$-neighbours or
$B$-neighbours$\}$.
\end{tabbing}
This set is non-empty because either $A$ or $B$ has at least two
members.

When $k=0$, it is clear that the lemma holds. Next we show that
the assumption that $k\neq 0$ leads to a contradiction.

Suppose $k>0$, and suppose $v$ and $w$ are $A$-neighbours or
$B$-neighbours such that $d_M(v,w)=k$, and suppose $M\!:v\mn
u\mn w$ for $u$ a member of $F$, which is either $A$ or $B$. If
$u$ has no $F$-neighbours (so it is both initial and final), then
it is interlaced with $v$ and $w$, which contradicts our
assumptions about parallelism for $A$ and $B$. So $u$ has at least
one $F$-neighbour.

If $u$ has an $F$-neighbour $u'$, which is either $v$ or $w$, or
$M\!:v\mn u'\mn w$, then $d_M(u,u')<d_M(v,w)$, which
contradicts the assumption that $k$ is minimal in the set $S$. So
for every $F$-neighbour $u'$ of $u$ we have that $u'$ is neither
$v$ nor $w$, nor $M\!:v\mn u'\mn w$. Then we have two cases.

One case is that $u$ has two $F$-neighbours (so it is neither initial
nor final in $F$), in which case we easily obtain that $u$ is
interlaced with $v$ and $w$. The other case is that $u$ has only
one $F$-neighbour; so $u$ is either initial or final in $F$ without
being both. Suppose $u$ is initial in $F$, and $w$ is the
immediate successor of $v$ in $A$ or $B$. If $u\equiv_\Pi w$, then
(1) of the definition of interlacing is trivially satisfied, and
if $u\equiv_\Pi v$, then (2) of this definition is satisfied. So
$u$ is interlaced with $v$ and $w$. The cases when $v$ is the
immediate successor of $w$, and when $u$ is final in $F$ are
treated analogously. In any case, we contradict our assumptions
about parallelism for $A$ and~$B$.\qed

Let us write $x<_L y$ if $x$ precedes $y$ in the list $L$, not
necessarily as an immediate predecessor.

\prop{Lemma 2.1.3.1}{If $n\geq 3$, and $a_{k+1}$ and $a_{k+2}$ are
$M$-neighbours, and $a_{k+1}\equiv_\Pi b_{l+1}$, then}
\vspace{-2ex}
\begin{tabbing}
\hspace{1.6em}\=$a_{k+1}<_Mb_{l+1}$\hspace{1.1em}\=\emph{iff}\hspace{1em}\=$a_{k+3}<_M b_{l+1}$.
\end{tabbing}

\dkz It is enough to see that ${M\!:a_{k+1}\mn b_{l+1}\mn
a_{k+3}}$ entails that \linebreak $M[a_{k+2},b_{l+1},a_{k+3}]$. \qed

We prove analogously the following

\prop{Lemma 2.1.3.2}{If $n\geq 3$, and $a_{k+n}$ and $a_{k+n-1}$
are $M$-neighbours, and $a_{k+n}\equiv_\Pi b_{l+m}$, then}
\vspace{-2ex}
\begin{tabbing}
\hspace{1.6em}\=$a_{k+1}<_Mb_{l+1}$\hspace{1.1em}\=\emph{iff}\hspace{1em}\=$a_{k+3}<_M b_{l+1}$.\kill

\>$a_{k+n}<_M b_{l+m}$\>\emph{iff}\>$a_{k+n-2}<_M
b_{l+m}$.
\end{tabbing}
\vspace{-2ex}
\prop{Lemma 2.1.3.3}{If $n\geq 3$, and $a_{k+1}$ and $a_{k+2}$
are $M$-neighbours, and $a_{k+2}\equiv_\Pi b_{l+m}$, then}
\vspace{-2ex}
\begin{tabbing}
\hspace{1.6em}\=$a_{k+1}<_Mb_{l+1}$\hspace{1.1em}\=\emph{iff}\hspace{1em}\=$a_{k+3}<_M b_{l+1}$.\kill

\>$a_{k+1}<_M b_{l+m}$\>\emph{iff}\>$a_{k+3}<_M b_{l+m}$.
\end{tabbing}
\dkz It is enough to see that ${M\!:a_{k+1}\mn b_{l+m}\mn
a_{k+3}}$ entails that \linebreak $M[a_{k+2},b_{l+m},a_{k+3}]$. \qed

We prove analogously the following.

\prop{Lemma 2.1.3.4}{If $n\geq 3$, and $a_{k+n}$ and $a_{k+n-1}$
are $M$-neighbours, and $a_{k+n-1}\equiv_\Pi b_{l+1}$, then}
\vspace{-2ex}
\begin{tabbing}
\hspace{1.6em}\=$a_{k+1}<_Mb_{l+1}$\hspace{1.1em}\=\emph{iff}\hspace{1em}\=$a_{k+3}<_M b_{l+1}$.\kill

\>$a_{k+n}<_M b_{l+1}$\>\emph{iff}\>$a_{k+n-2}<_M b_{l+1}$.
\end{tabbing}
\vspace{-2ex}
\prop{Lemma 2.1.4}{If $a_{k+1}\equiv_\Pi b_{l+1}$ and $n$ and
$m$ are both odd, then}
\vspace{-2ex}
\begin{tabbing}
\hspace{1.6em}\=$a_{k+1}<_Mb_{l+1}$\hspace{1.1em}\=\emph{iff}\hspace{1em}\=$a_{k+3}<_M b_{l+1}$.\kill

\>$a_{k+1}<_M b_{l+1}$\>\emph{iff}\>$a_{k+n}<_M b_{l+m}$.
\end{tabbing}
\dkz We make an induction on $n\pl m$. The basis, when $n=m=1$, is
trivial. In the induction step, apply Lemma 2.1.2, and suppose the
members $a_i$ and $a_{i+1}$ of $A$ are $M$-neighbours. By our
assumptions on parallelism and Lemmata 2.1.1.1, 2.1.1.2 and
2.1.1.3, we have $A'\parallel_{M'} B$, $A'\parallel_{M'} A'$ and
$B\parallel_{M'} B$; so we will be able to apply the induction
hypothesis. If $i\neq k\pl 1$ and $i\pl 1\neq k\pl n$, then we are
done. If $i=k\pl 1$ or $i\pl 1=k\pl n$, then we use Lemmata
2.1.3.1 and 2.1.3.2 respectively.\qed

\vspace{-2ex}

\prop{Lemma 2.1.5}{If $a_{k+1}\equiv_\Pi b_{l+1}$ and $n$ is even
and $m$ odd, then}
\vspace{-2ex}
\begin{tabbing}
\hspace{1.6em}\=$a_{k+1}<_Mb_{l+1}$\hspace{1.1em}\=\emph{iff}\hspace{1em}\=$a_{k+3}<_M b_{l+1}$.\kill

\>$a_{k+1}<_M b_{l+1}$\>\emph{iff}\>$a_{k+n}<_M
b_{l+1}$.
\end{tabbing}
\dkz We make an induction on $n\pl m$. In the basis we have $n=2$
and $m=1$. If our equivalence did not hold, then we would have that
$M[a_{k+1},b_{l+1},a_{k+2}]$. In the induction step we proceed as
in the induction step of the proof of Lemma 2.1.4, by using
Lemmata 2.1.3.1 and 2.1.3.4 when we pass to $A'$, and by using
appropriately renamed variants of Lemmata 2.1.3.1 and 2.1.3.3 when
we pass to~$B'$. \qed

\vspace{-2ex}

\prop{Lemma 2.1.6}{If $a_{k+1}\equiv_\Pi b_{l+1}$ and $n$ and $m$
are both even, then}
\vspace{-2ex}
\begin{tabbing}
\hspace{1.6em}\=$a_{k+1}<_Mb_{l+1}$\hspace{1.1em}\=\emph{iff}\hspace{1em}\=$a_{k+3}<_M b_{l+1}$.\kill

\>$M\!:a_{k+1}\mn b_{l+1}\mn a_{k+n}\quad{\rm\emph{iff}}\quad M\!:a_{k+1}\mn b_{l+m}\mn a_{k+n}$.
\end{tabbing}
\dkz Note first that the equivalence of this lemma can be stated
equivalently as follows:
\begin{tabbing}
\hspace{1.6em}\=$(\ast)$\hspace{3em}\=$M\!:b_{l+1}\mn a_{k+1}\mn
b_{l+m}$\hspace{1em}\=\emph{iff}\hspace{1em}\=$M\!:b_{l+1}\mn a_{k+n}\mn
b_{l+m}$.\kill

\>$(\ast)$\hspace{3em}\>$M\!:b_{l+1}\mn a_{k+1}\mn
b_{l+m}$\hspace{1em}\>\emph{iff}\hspace{1em}\>$M\!:b_{l+1}\mn a_{k+n}\mn
b_{l+m}$.
\end{tabbing}
To prove that, we proceed by induction on $n\pl m$. In the basis, we have
$n=m=2$, since $n$ and $m$ are both even, but are different from
0. If our equivalence did not hold, then we would contradict our
assumptions on parallelism. For the induction step, apply first
Lemma 2.1.2, and suppose the members $a_i$ and $a_{i+1}$ of $A$
are $M$-neighbours. If $n=2$, then $(\ast)$ holds trivially.

If $n\geq 3$, and $k\pl 1<i$ and $i\pl 1<k\pl n$, then we just
apply Lemmata 2.1.1.1, 2.1.1.2 and 2.1.1.3 and the induction
hypothesis to obtain $(\ast)$. Suppose $n\geq 3$ and $k\pl 1=i$,
and suppose $b_{l+1}<_M a_{k+1}$ and $a_{k+1}<_M b_{l+m}$. Then,
since for $u$ being $a_{k+1}$ and $a_{k+3}$ we have $b_{l+1}<_M u$
iff not $u<_M b_{l+1}$, by Lemmata 2.1.3.1 and 2.1.3.3, we infer
$b_{l+1}<_M a_{k+3}$ and $a_{k+3}<_M b_{l+m}$. We can also make
the converse inference, and conclude that
\[
(b_{l+1}<_M a_{k+1}\;{\rm and}\; a_{k+1}<_M b_{l+m})\;{\rm iff}\;
(b_{l+1}<_M a_{k+3}\;{\rm and}\; a_{k+3}<_M b_{l+m}).
\]

Since we can prove analogously the equivalence obtained from this
one by interchanging $b_{l+1}$ and $b_{l+m}$, we obtain
\begin{tabbing}
\hspace{1.6em}\=$(\ast)$\hspace{3em}\=$M\!:b_{l+1}\mn a_{k+1}\mn
b_{l+m}$\hspace{1em}\=\emph{iff}\hspace{1em}\=$M\!:b_{l+1}\mn a_{k+n}\mn
b_{l+m}$.\kill

\>$(\ast\ast)$\hspace{3em}\>$M\!:b_{l+1}\mn a_{k+1}\mn
b_{l+m}$\hspace{1em}\>\emph{iff}\hspace{1em}\>$M\!:b_{l+1}\mn a_{k+3}\mn
b_{l+m}$.
\end{tabbing}
Since for $3\leq j\leq n$ we have
\begin{tabbing}
\hspace{1.6em}\=$(\ast)$\hspace{3em}\=$M\!:b_{l+1}\mn a_{k+1}\mn
b_{l+m}$\hspace{1em}\=\emph{iff}\hspace{1em}\=$M\!:b_{l+1}\mn a_{k+n}\mn
b_{l+m}$.\kill

\>${(\ast\!\ast\!\ast)}$\hspace{3em}\>$ M\!:b_{l+1}\mn a_{k+j}\mn
b_{l+m}$\hspace{1em}\>\emph{iff}\hspace{1em}\>$M'\!:b_{l+1}\mn a_{k+j}\mn
b_{l+m}$,
\end{tabbing}
and since by Lemmata 2.1.1.1, 2.1.1.2 and 2.1.1.3 and the
induction hypothesis we have
\begin{tabbing}
\hspace{1.6em}\=$(\ast)$\hspace{3em}\=$M\!:b_{l+1}\mn a_{k+1}\mn
b_{l+m}$\hspace{1em}\=\emph{iff}\hspace{1em}\=$M\!:b_{l+1}\mn a_{k+n}\mn
b_{l+m}$.\kill

\>$(\ast\!\ast\!\ast\ast)$\hspace{3em}\>$M'\!:b_{l+1}\mn a_{k+3}\mn
b_{l+m}$\hspace{1em}\>\emph{iff}\hspace{1em}\>$M'\!:b_{l+1}\mn a_{k+n}\mn
b_{l+m}$,
\end{tabbing}
we derive $(\ast)$ as follows:
\begin{tabbing}
\hspace{1.6em}$M\!:b_{l+1}\mn a_{k+1}\mn b_{l+m}$\quad\=iff\quad\=
$M\!:b_{l+1}\mn a_{k+3}\mn b_{l+m}$,\quad\=by $(\ast\ast)$,
\\[.5ex]
\>iff\>$M'\!:b_{l+1}\mn a_{k+3}\mn b_{l+m}$,\quad\>by
${(\ast\!\ast\!\ast)}$,
\\[.5ex]
\>iff\>$M'\!:b_{l+1}\mn a_{k+n}\mn b_{l+m}$,\quad\>by
$(\ast\!\ast\!\ast\ast)$,
\\[.5ex]
\>iff\>$M\!:b_{l+1}\mn a_{k+n}\mn b_{l+m}$,\quad\>by
${(\ast\!\ast\!\ast)}$.
\end{tabbing}
We reason analogously when $n\geq 3$ and $i+1=k+n$. \qed

\vspace{-2ex}

\section{\large\bf P$'$-graphs and grounding}\label{2.2}
\markright{\S 2.2. \quad P$\,'$-graphs and grounding}

Suppose we have a construction $K$ of the P$'$-graph $D$, and
suppose that in the root of $K$ we have the triple $(D,L_W,L_E)$.
Our purpose now is to prove that $L_W$ and $L_E$ are grounded in
$D$. This is (2) of the proposition we are going to prove. That
proposition asserts also something else concerning related
matters, which are involved in the proof of~(2).

\prop{Proposition 2.2.1}{{\rm (1)} If $x<_{L_W} y$ and $z<_{L_E} u$,
then every semipath in $[x,u]$ and every semipath in $[y,z]$
intersect. \\[1ex] \indent {\rm (2)} If $L_X\!:x\mn y\mn z$, then every
semipath in $[x,z]$ and every semipath in $[y]_{\bar{X}}$
intersect (i.e., we have $\psi_X(x,y,z)$). \\[1ex] \indent {\rm (3)}
If $L_X\!:x\mn y\mn u\mn z$, then every semipath in $[x,u]$ and
every semipath in $[y,z]$ intersect.}

\dkz We proceed by induction on the number $k$ of inner vertices
of $D$. If $k=1$, then $D$ is a basic D-graph, and it is easy to
convince oneself that the unique inner vertex of $D$ serves for
all the intersections we need in (1), (2) and (3).

Suppose now that $k>1$. So we have $D=D_W\Box D_E$. We prove first
(1).

\vspace{.5ex}

(1) Suppose $x<_{L_W} y$ and $z<_{L_E} u$. We have the following
cases depending on where $x$, $y$, $z$ and $u$ are.

\vspace{.5ex}

(1.1) Suppose $x,y\in W(D_W)$ and $z,u\in E(D_E)$. Out of a
semipath $\sigma$ in $[x,u]$ we construct the list $A$, which is
$a_1\ldots a_n$, with $n\geq 1$ and odd, made of all the edges of
$\sigma$ that are elements of the set of edges $C$ involved in
$D_W\Box D_E$; the edges of this list are listed in the order in
which they appear in $\sigma$. In the list $A$ we could
alternatively take that $a_i$ for $i\in\{1,\ldots, n\}$ instead of
being an edge is a vertex $v_{a_i}$ such that in $D_E$ we have
$W_E(a_i)=v_{a_i}$ and in $D_W$ we have $E_W(a_i)=v_{a_i}$. Our
lists $A$ of edges and such lists of vertices correspond
bijectively to each other. (The disadvantage of the list of
vertices $A$ is however that the vertex $v_{a_i}$ is not in $D_W\Box
D_E$, while the edge $a_i$~is.)

We construct the list $B$, which is $b_1\ldots b_m$, with $m\geq
1$ and odd, out of a semipath $\tau$ in $[y,z]$ in an analogous
manner. Since the lists $L^W_E$ and $L^E_W$ are compatible, there
is a unified list $L$ of these two lists (see \S 1.7). We
construct a list $M$ of the members of $A$ and $B$ where these
members are listed in the order in which they occur in~$L$.

The following part of our proof will be repeated in several
analogous variants later on, and this is why we mark it with~$(\dagger)$.

\vspace{.5ex}

$(\dagger)$ Suppose for $F,G\in\{A,B\}$ we do not have
$F\parallel_M G$. Suppose a member $a$ of $F$ is interlaced with
two members $b$ and $b^s$ of $G$. If $a\equiv_\Pi b^s$ and $a^p$
exists, then we have either $(L^{\bar{X}}_X)^e\!:a^p\mn b\mn a\mn
b^s$ or $(L^{\bar{X}}_X)^e\!:b\mn a\mn b^s\mn a^p$, where
$(L^{\bar{X}}_X)^e$ is the list of edges corresponding to the list
of vertices $L^{\bar{X}}_X$. If $v_e$ is the vertex of $D_X$
corresponding bijectively to $e$ (see the comment above after the
definition of $A$), then we have $L^{\bar{X}}_X\!:v_{a^p}\mn
v_b\mn v_a\mn v_{b^s}$ or $L^{\bar{X}}_X\!:v_b\mn v_a\mn
v_{b^s}\mn v_{a^p}$.

If $b$ has an odd index, then $\bar{X}$ is $E$, while if $b$ has
an even index, then $\bar{X}$ is $W$. Then we apply the induction
hypothesis (3) to $D_{\bar{X}}$, and obtain an intersection of
every semipath in $[v_a,v_{a^p}]$ with every semipath in
$[v_b,v_{b^s}]$. We cannot have that $F$ and $G$ are both $A$ or
both $B$, because then $\sigma$ or $\tau$ would not be semipaths
(a vertex cannot occur twice in a semipath). So one of $F$ and $G$
is $A$, while the other is $B$, and we may conclude that $\sigma$
and $\tau$ intersect, as required.

If $a^p$ does not exist, i.e., $a$ is initial in $A$, then we
apply the induction hypothesis (2) to $D_W$, and obtain again that
$\sigma$ and $\tau$ intersect. If $a\equiv_\Pi b$, then again we
have two cases, in one of which we apply the induction hypothesis
(3) and in the other the induction hypothesis (2). This concludes
the $(\dagger)$ part of the proof.

\vspace{.5ex}

Suppose now that $A\parallel_M B$, $A\parallel_M A$ and
$B\parallel_M B$. If we have $a_1<_M b_1$, then by Lemma 2.1.4 we
have $a_n<_M b_m$, and hence $v_{a_n}<_{L^E_W} v_{b_m}$. We may
apply the induction hypothesis (1) to $D_E$, to infer that
$\sigma$ and $\tau$ intersect. If we have $b_1<_M a_1$, then we
have $v_{b_1}<_{L^W_E} v_{a_1}$, and we apply the induction
hypothesis (1) to~$D_W$.

\vspace{.5ex}

(1.2) Suppose $x\in W(D_E)$, $y\in W(D_W)$ and $z,u\in E(D_E)$. We
construct $A$, $B$ and $M$ as in (1.1) except that $a_1$ is not an
element of $C$, but it is the edge $e$ such that $W(e)=x$. Then we
reason as in $(\dagger)$ until we reach the supposition that
$A\parallel_M B$, $A\parallel_M A$ and $B\parallel_M B$. Now we
know that we have $v_{a_1}<_{L_E} v_{b_1}$, and we apply again
Lemma 2.1.4 and the induction hypothesis (1) to $D_E$.

\vspace{.5ex}

(1.3) Suppose $x,y\in W(D_E)$ and $z,u\in E(D_E)$. We construct
$A$, $B$ and $M$ as in (1.2) except that $b_1$ too is not an element
of $C$, but it is the edge $e$ such that $W(e)=y$. In the
remainder of this case we reason as for~(1.2).

\vspace{.5ex}

(1.4) Suppose $x\in W(D_E)$, $y\in W(D_W)$, $z\in E(D_E)$ and
$u\in E(D_W)$. We construct $A$, $B$ and $M$ as in (1.2) except
that $a_n$ too is not an element of $C$, but it is the edge $e$ such
that $E(e)=u$. Since we must have $a_1<_M b_1$ and $b_m<_M a_n$,
by Lemma 2.1.4 we can conclude that $A\parallel_M B$ fails, and we
obtain that $\sigma$ intersects $\tau$ as in~$(\dagger)$.

All the other cases are treated analogously. Note that the case
$x\in W(D_E)$, $y\in W(D_W)$, $z\in E(D_W)$ and $u\in E(D_E)$ is
impossible, since it makes $L^W_E$, which is of the form
$\Phi_E\Lambda\Psi_E$, and $L^E_W$, which is of the form
$\Phi_W\Lambda\Psi_W$, not compatible; in $\Phi_E$ we have $x$ and
in $\Phi_W$ we have $z$. There are other such impossible cases,
excluded by the compatibility of $L^W_E$ and~$L^E_W$.

\vspace{.5ex}

(2) Suppose $x<_{L_X}y$ and $y<_{L_X} z$. (It will help the
intuition to suppose that $X$ is here $W$ while $\bar{X}$ is $E$;
or the other way round.) We have the following cases depending on
where $x$, $y$ and $z$ are.

\vspace{.5ex}

(2.1) Suppose $x,y,z\in X(D_X)$. Out of a semipath $\sigma$ in
$[x,z]$ we construct the list $A$, which is $a_1\ldots a_n$, with
$n\geq 0$ and even, in the same manner as in (1.1). (If $n=0$, then
$A$ i empty.) We construct out of a semipath $\tau$ in
$[y]_{\bar{X}}$ the list $B$, which is $b_1\ldots b_m$, with
$m\geq 1$ and odd, in an analogous manner. We construct $M$ as in~(1.1).

If $n=0$, i.e., the list $A$ is empty, then we apply the induction
hypothesis (2) to $D_X$ to obtain that $\sigma$ and $\tau$
intersect. If $n\geq 2$, then we continue as in $(\dagger)$ until
the supposition that $A\parallel_M B$, $A\parallel_M A$ and
$B\parallel_M B$.

If we have $a_1<_M b_1$, then by Lemma 2.1.5 we have $a_n<_M b_1$,
and hence $v_{a_n}<_{L^X_{\bar{X}}} v_{b_1}$. We may then apply
the induction hypothesis (1) to $D_X$ (with $v_{a_n}$, $v_{b_1}$,
$y$ and $z$ standing respectively for $x$, $y$, $z$ and $u$) to
infer that $\sigma$ and $\tau$ intersect.

If we have $b_1<_M a_1$, then we apply immediately the induction
hypothesis (1) to $D_X$ (with $v_{b_1}$, $v_{a_1}$, $x$ and $y$
standing respectively for $x$, $y$, $z$ and~$u$).

\vspace{.5ex}

(2.2) Suppose $x\in X(D_{\bar{X}})$ and $y,z\in X(D_X)$. We
construct $A$, $B$ and $M$ as in (2.1) except that $a_1$ is not an element of $C$, but (as in (1.2)) it is the edge $e$ such that
$\bar{X}(e)=x$. We cannot have $A$ empty now. We continue as in
$(\dagger)$ until the supposition that $A\parallel_M B$,
$A\parallel_M A$ and $B\parallel_M B$.

If we have $a_1<_M b_1$, then we continue reasoning as in (2.1) by
relying on Lemma 2.1.5. It is now excluded that $b_1<_M a_1$.

\vspace{.5ex}

(2.3) Suppose $x,y\in X(D_{\bar{X}})$ and $Z\in X(D_X)$. We
construct $A$, $B$ and $M$ as in (2.2) except that $b_1$ is not an element of $C$, but it is the edge $e$ such that $\bar{X}(e)=y$. It
is excluded that $A$ is empty. Since we must have $a_1<_M b_1$ and
$b_1<_M a_n$, because $a_n$ is in $C$, by Lemma 2.1.5 we conclude
that $A\parallel_M B$ fails, and we obtain that $\sigma$
intersects $\tau$ as in~$(\dagger)$.

In cases where $x,y,z\in X(D_{\bar{X}})$, or where $u, z\in
X(D_{\bar{X}})$ and $y\in X(D_X)$, we reason analogously to what we
had for case (2.3). We reason analogously when $z<_{L_X} y$ and
$y<_{L_X} x$.

\vspace{.5ex}

(3) Suppose $L_X\!: x\mn y\mn u\mn z$. We have the following cases
depending on where $x$, $y$, $u$ and $z$ are.

\vspace{.5ex}

(3.1) Suppose $x,y,u,z\in X(D_X)$. We construct the lists $A$, $B$
and $M$ as in (1.1) save that for $A$, which is $a_1\ldots a_n$,
and for $B$, which is $b_1\ldots b_m$, we have that $n\geq 0$ and
$m\geq 0$, and they are both even.

If $n=m=0$, i.e., both $A$ and $B$ are empty, then we apply the
induction hypothesis (3) to $D_X$ to obtain that $\sigma$ and
$\tau$ intersect.

If $n=0$ and $m>0$, then we apply the induction hypothesis (2) to
$D_X$ (with $x$, $y$ and $u$ standing respectively for $x$, $y$
and $z$), in order to infer that $\sigma$ and $\tau$ intersect. We
proceed analogously when $n>0$ and $m=0$.

If $n>0$ and $m>0$, then we continue as in $(\dagger)$ until the
supposition that $A\parallel_M B$, $A\parallel_M A$ and
$B\parallel_M B$. Then by applying Lemma 2.1.6 we obtain the
possibility to apply the induction hypothesis (1) to $D_X$, in
order to infer that $\sigma$ and $\tau$ intersect.

The cases when

\vspace{.5ex}

(3.2) $x\in X(D_{\bar{X}})$ and $y,u,z\in X(D_X)$,

\vspace{.5ex}

(3.3) $x,y\in X(D_{\bar{X}})$ and $u,z\in X(D_X)$,

\vspace{.5ex}

(3.4) $x,z\in X(D_{\bar{X}})$ and $y,u\in X(D_X)$,

\vspace{.5ex}

\noindent are treated analogously to case (3.1). In (3.2) we may
have $n>0$ and $m=0$, where we apply the induction hypothesis (2).
If $n>0$ and $m>0$, which we may have with (3.2), and which must
be the case with (3.3) and (3.4), we have $(\dagger)$ and we apply
Lemma 2.1.6 and the induction hypothesis~(1).

In the cases when

\vspace{.5ex}

(3.5) $x,y,u\in X(D_{\bar{X}})$ and $z\in X(D_X)$,

\vspace{.5ex}

(3.6) $x,y,z\in X(D_{\bar{X}})$ and $u\in X(D_X)$,

\vspace{.5ex}

(3.7) $x,y,z,u\in X(D_{\bar{X}})$,

\vspace{.5ex}

\noindent we conclude by Lemma 2.1.6 that $A\parallel_M B$ fails,
and we obtain that $\sigma$ intersects $\tau$ as in $(\dagger)$.
For example, in (3.5) we must have $L^{\bar{X}}_X\!:x\mn y\mn u\mn
v_{b_m}$, with $x$, $y$ and $u$ being respectively $v_{a_1}$,
$v_{b_1}$ and $v_{a_n}$.

All the other cases are treated analogously to these. \qed

\vspace{-2ex}

\section{\large\bf P$'$-graphs are P$'''$-graphs}\label{2.3}
\markright{\S 2.3. \quad P$\,'$-graphs are P$\,'''$-graphs}

We may enter now into the proof that every P$'$-graph is a
P$'''$-graph. The essential ingredient of that proof will be (2)
of Proposition 2.2.1, together with the following lemmata, for
which we assume that $K_1$, $K_2$ and $K_3$ are constructions of
the P$'$-graphs $D_1$, $D_2$ and $D_3$ respectively. In these
lemmata ``construction'' is short for ``construction of a
P$'$-graph'' (see \S 1.8). We have first two lemmata corresponding to (Ass~1)
(see~\S 1.5).

\prop{Lemma 2.3.1.1}{If $K_1\Box K_2$ and $(K_1\Box K_2)\Box K_3$
are constructions, and $D_2\Box D_3$ is a D-graph, then $K_2\Box
K_3$ and $K_1\Box(K_2\Box K_3)$ are constructions.}

\dkz Let us write $K\!: L_W\vdash L_E$ to indicate that in its
root the construction $K$ has $(D,L_W,L_E)$.

Suppose that $D_1\Box D_3$ is not defined. Then we have
$K_1\!:\Gamma\vdash\Delta^1_1\Theta\Delta^1_2$,
$K_2\!:\Gamma^2_1\Theta\Gamma^2_2\vdash\Delta^2_1\Xi\Delta^2_2$
and $K_3\!:\Gamma^3_1\Xi\Gamma^3_2\vdash\Delta$, and hence
\[
\f{\f{K_1\!:\Gamma\vdash\Delta^1_1\Theta\Delta^1_2\quad
K_2\!:\Gamma^2_1\Theta\Gamma^2_2\vdash
\Delta^2_1\Xi\Delta^2_2}{K_1\Box
K_2\!:\Gamma^2_1\Gamma\Gamma^2_2\vdash\Delta^1_1\Delta^2_1\Xi\Delta^2_2\Delta^1_2}\quad
\afrac{K_3\!:\Gamma^3_1\Xi\Gamma^3_2\vdash\Delta}}{(K_1\Box
K_2)\Box
K_3\!:\Gamma^3_1\Gamma^2_1\Gamma\Gamma^2_2\Gamma^3_2\vdash\Delta^1_1\Delta^2_1\Delta\Delta^2_2\Delta^1_2}
\]
with $\Theta$ and $\Xi$ non-empty, and with some requirements
concerning the emptiness of $\Delta^1_1$, $\Delta^1_2$, etc., so as
to ensure compatibility.

If we cannot obtain $K_2\Box
K_3\!:\Gamma^3_1\Gamma^2_1\Theta\Gamma^2_2\Gamma^3_2\vdash\Delta^2_1\Delta\Delta^2_2$
because $\Delta^2_1$ and $\Gamma^3_1$ are both non-empty, then we
could not obtain $(K_1\Box K_2)\Box K_3$, and if we cannot obtain
$K_2\Box K_3$ because $\Delta^2_2$ and $\Gamma^2_2$ are both
non-empty, then again we could not obtain $(K_1\Box K_2)\Box K_3$.
So $K_2\Box K_3$ is a construction.

If we cannot obtain $K_1\Box (K_2\Box
K_3)\!:\Gamma^3_1\Gamma^2_1\Gamma\Gamma^2_2\Gamma^3_2\vdash\Delta^1_1\Delta^2_1\Delta\Delta^2_2\Delta^1_2$
because $\Delta^1_1$ and $\Gamma^3_1\Gamma^2_1$ are both non-empty,
then $\Gamma^3_1$ is non-empty or $\Gamma^2_1$ is non-empty. If
$\Delta^1_1$ and $\Gamma^3_1$ are non-empty, then we could not
obtain $(K_1\Box K_2)\Box K_3$, and if $\Delta^1_1$ and
$\Gamma^2_1$ are non-empty, then we could not obtain $K_1\Box K_2$.

If we cannot obtain $K_1\Box (K_2\Box K_3)$ because $\Delta^1_2$
and $\Gamma^2_2\Gamma^3_2$ are both non-empty, then $\Gamma^2_2$ is
non-empty or $\Gamma^3_2$ is non-empty. If $\Delta^1_2$ and
$\Gamma^2_2$ are non-empty, then we could not obtain $K_1\Box K_2$,
and if $\Delta^1_2$ and $\Gamma^3_2$ are non-empty, then we could
not obtain $(K_1\Box K_2)\Box K_3$. So $K_1\Box(K_2\Box K_3)$ is a
construction.

It remains to consider the case when $D_1\Box D_3$ is defined.
Then we may have
\[
\f{\f{K_1\!:\Gamma\vdash\Delta^1_1\Phi\Theta\Delta^1_2\quad
K_2\!:\Theta\Gamma^2_2\vdash \Xi\Delta^2_2}{K_1\Box
K_2\!:\Gamma\Gamma^2_2\vdash\Delta^1_1\Phi\Xi\Delta^2_2\Delta^1_2}\quad
\afrac{K_3\!:\Gamma^3_1\Phi\Xi\Gamma^3_2\vdash\Delta}}{(K_1\Box
K_2)\Box
K_3\!:\Gamma^3_1\Gamma\Gamma^2_2\Gamma^3_2\vdash\Delta^1_1\Delta\Delta^2_2\Delta^1_2}
\]
with $\Theta$, $\Xi$ and $\Phi$ non-empty, and with some
requirements concerning the emptiness of $\Delta^1_1$,
$\Delta^1_2$, etc., so as to ensure compatibility. There is a
mirror case where $\Phi\Theta$ and $\Phi\Xi$ are replaced by
$\Theta\Phi$ and $\Xi\Phi$ (which requires that $\Gamma^2_2$ and
$\Delta^2_2$ be empty instead of $\Gamma^2_1$ and $\Delta^2_1$);
we deal with that case analogously.

If we cannot obtain $K_2\Box
K_3\!:\Gamma^3_1\Phi\Theta\Gamma^2_2\Gamma^3_2\vdash\Delta\Delta^2_2$
because $\Delta^2_2$ and $\Gamma^3_2$ are both non-empty, then we
could not obtain $(K_1\Box K_2)\Box K_3$. So $K_2\Box K_3$ is a
construction.

If we cannot obtain $K_1\Box (K_2\Box K_3)\!:
\Gamma^3_1\Gamma\Gamma^2_2\Gamma^3_2\vdash\Delta^1_1\Delta\Delta^2_2\Delta^1_2$
because $\Delta^1_1$ and $\Gamma^3_1$ are both non-empty, then we
could not obtain $(K_1\Box K_2)\Box K_3$. If we cannot obtain
$K_1\Box (K_2\Box K_3)$ because $\Delta^1_2$ and
$\Gamma^2_2\Gamma^3_2$ are both non-empty, then we continue
reasoning as in the analogous case we had above. \qed

We can prove analogously the following kind of converse of
Lemma~2.3.1.1.

\prop{Lemma 2.3.1.2}{If $K_2\Box K_3$ and $K_1\Box(K_2\Box K_3)$
are constructions, and $D_1\Box D_2$ is a D-graph, then $K_1\Box
K_2$ and $(K_1\Box K_2)\Box K_3$ are constructions.}

We also have the following two lemmata corresponding to (Ass~2.1)
and (Ass~2.2), of which we prove only the first (the proof of the
second is analogous).

\prop{Lemma 2.3.2.1}{If $K_1\Box K_2$ and $(K_1\Box K_2)\Box K_3$
are constructions, and $D_2\Box D_3$ is not defined, then $K_1\Box
K_3$ and $(K_1\Box K_3)\Box K_2$ are constructions.}

\dkz We have
\[
\f{\f{K_1\!:\Gamma\vdash\Delta^1_1\Theta\Delta^1_2\Xi\Delta^1_3\quad
K_2\!:\Gamma^2_1\Theta\vdash \Delta^2}{K_1\Box
K_2\!:\Gamma^2_1\Gamma\vdash\Delta^1_1\Delta^2\Delta^1_2\Xi\Delta^1_3}\quad
\afrac{K_3\!:\Xi\Gamma^3_2\vdash\Delta^3}}{(K_1\Box K_2)\Box
K_3\!:\Gamma^2_1\Gamma\Gamma^3_2\vdash\Delta^1_1\Delta^2\Delta^1_2\Delta^3\Delta^1_3}
\]
with $\Theta$, $\Xi$ and $\Delta^2$ non-empty, and with some
requirements concerning the emptiness of $\Delta^1_1$,
$\Delta^1_2$, etc., so as to ensure compatibility.

If we cannot obtain $K_1\Box K_3\!:
\Gamma\Gamma^3_2\vdash\Delta^1_1\Theta\Delta^1_2\Delta^3\Delta^1_3$
because $\Delta^1_3$ and $\Gamma^3_2$ are both non-empty, then we
could not obtain $(K_1\Box K_2)\Box K_3$. So $K_1\Box K_3$ is a
construction.

If we cannot obtain $(K_1\Box K_3)\Box K_2\!:
\Gamma^2_1\Gamma\Gamma^3_2\vdash\Delta^1_1\Delta^2\Delta^1_2\Delta^3\Delta^1_3$
because $\Delta^1_1$ and $\Gamma^2_1$ are both non-empty, then we
could not obtain $K_1\Box K_2$. So $(K_1\Box K_3)\Box K_2$ is a
construction. \qed

\vspace{-2ex}

\prop{Lemma 2.3.2.2}{If $K_2\Box K_3$ and $K_1\Box (K_2\Box K_3)$
are constructions, and $D_1\Box D_2$ is not defined, then $K_1\Box
K_3$ and $K_2\Box (K_1\Box K_3)$ are constructions.}

We can prove the following concerning D-terms (see \S 1.5) and their interpretations with $\iota$ (see \S 1.6) .

\prop{Lemma 2.3.3}{For every D-graph $\iota(\delta)$ for a D-term
$\delta$ and for every cocycle $C$ of $\iota(\delta)$ there are
two D-terms $\delta_W$ and $\delta_E$ such that $\iota(\delta_W)$
and $\iota(\delta_E)$ are the D-graphs obtained by cutting
$\iota(\delta)$ through~$C$.}

\dkz Let $D_W$ and $D_E$ be the D-graphs obtained by cutting
$\iota(\delta)$ through $C$. Let $\delta_W'$ and $\delta_E'$ be
D-terms such that $\iota(\delta_W')$ and $\iota(\delta_E')$ are
isomorphic respectively to $D_W$ and $D_E$ (such D-terms exist, as
noted at the end of \S 1.6). By renaming, if need there is, the
edges assigned to $\delta_W'$ and $\delta_E'$ so that they accord
with the edges assigned to $\delta$ we pass to $\delta_W$ and
$\delta_E$. \qed

Consider a map $\varphi$ that assigns to a construction $K$ of a
P$'$-graph $D$ a D-term $\varphi(K)$ such that $\iota(\varphi(K))$
is isomorphic to $D$, and which satisfies that $\varphi(K_W\Box
K_E)$ is $\varphi(K_W)\Box\varphi(K_E)$. It is clear that such
maps exist. Relying on Lemmata 2.3.1.1, 2.3.1.2, 2.3.2.1 and
2.3.2.2 it is easy to prove the following by induction on the
length of derivation.

\prop{Lemma 2.3.4}{If $\delta=\delta'$ is derivable in {\rm
S}$\Box$, then there is a construction $K$ of the P$\,'$-graph
$\iota(\delta)$ such that $\varphi(K)$ is $\delta$ iff there is a
construction $K'$ of $\iota(\delta')$ (which is equal to
$\iota(\delta))$ such that $\varphi(K')$ is~$\delta'$.}

We can then prove the following.

\prop{Lemma 2.3.5}{Every P$\,'$-graph $\iota(\delta)$ for which
there is a construction $K$ such that $\varphi(K)$ is $\delta$ is
a P$\,'''$-graph.}

\dkz If $\iota(\delta)$ is a basic D-graph, then it has no
cocycles, and it is hence trivially a P$'''$-graph. Assume
$\iota(\delta)$ is not basic, and take an arbitrary cocycle $C$ of
$\iota(\delta)$. Let $\iota(\delta_W)$ and $\iota(\delta_E)$ be
the D-graphs obtained by cutting $\iota(\delta)$ through $C$,
which we have by Lemma 2.3.3. Hence we have
\begin{tabbing}
\hspace{1.6em}$\iota(\delta)=\iota(\delta_W)\Box\iota(\delta_E)=\iota(\delta_W\Box\delta_E)$,
\end{tabbing}
and, by the completeness of S$\Box$ (see \S 1.6), we obtain that
$\delta=\delta_W\Box\delta_E$ is derivable in S$\Box$. Since there
is a construction $K$ of $\iota(\delta)$ such that $\varphi(K)$ is
$\delta$, there is, by Lemma 2.3.4, a construction $K'$ of
$\iota(\delta_W\Box\delta_E)$ such that $\varphi(K')$ is
$\delta_W\Box\delta_E$.

Since $\varphi(K_W\Box K_E)$ is $\varphi(K_W)\Box\varphi(K_E)$,
from $\varphi(K')$ being $\delta_W\Box\delta_E$ we may conclude
that $\delta_W$ and $\delta_E$ are $\varphi(K_W)$ and
$\varphi(K_E)$ respectively, and that $K'$ is $K_W\Box K_E$, for
$K_X$ a construction of $\iota(\delta_X)$. So we have that the
D-graph $\iota(\delta_X)$ is a P$'$-graph.

Take the list of $E(\iota(\delta_W))$ from the root of $K_W$ and
the list of $W(\iota(\delta_E))$ from the root of $K_E$. By (2) of
Proposition 2.2.1 these lists are grounded in $\iota(\delta_W)$
and $\iota(\delta_E)$ respectively. They are compatible because
$K'$ is a construction. So $\iota(\delta)$ is a P$'''$-graph. \qed

Every P$'$-graph is isomorphic to a P$'$-graph $\iota(\delta)$ for
which there is a construction $K$ such that $\varphi(K)$ is
$\delta$. Just take an arbitrary construction $K$ of our
P$'$-graph, and take $\delta$ to be~$\varphi(K)$.

It is also clear that every graph isomorphic to a
P$'''$-graph is a P$'''$-graph. So from Lemma 2.3.5 we
may conclude the following.

\prop{Theorem 2.3.6}{Every P$\,'$-graph is a P$\,'''$-graph.}

\clearpage \pagestyle{empty} \makebox[1em]{} \clearpage

\chapter{\huge\bf Grounding and
Pivots}\label{3} \pagestyle{myheadings}\markboth{CHAPTER 3. \quad
GROUNDING AND PIVOTS}{right-head}

\section{\large\bf Grounding and
juncture}\label{3.1} \markright{\S 3.1. \quad Grounding and
juncture}

This chapter contains preliminary results that will help us to establish in \S 4.4 and \S 5.3 that P$''$-graphs (as defined in \S 1.9) are P$'$-graphs (as defined in \S 1.8) and that P$'''$-graphs (as defined in \S 1.10) are P$''$-graphs. The main of these results is the Pivot Theorem of~\S 3.4.

In this section we prove results that will be used to establish that P$'''$-graphs are P$''$-graphs. These results are about the relationship between groundedness (see \S 1.9) in $D_W\Box D_E$ and groundedness in $D_W$ and~$D_E$.

Let $\Gamma$ and $\Delta$ be lists such that $x$, which is either
the initial or the final member of $\Gamma$, is either not a member
of $\Delta$ or it is not the only member that $\Gamma$ and
$\Delta$ share, and let $\Gamma'$ be obtained from $\Gamma$ by
removing $x$ from $\Gamma$. Then we can easily establish the
following.

\prop{Lemma 3.1.1}{If $\Gamma$ and $\Delta$ are compatible, then
$\Gamma'$ and $\Delta$ are compatible.}

We establish easily the following for D-graph $D_W\Box D_E$.

\prop{Lemma 3.1.2}{If $x,y,z\in X(D_X)$ and $\psi_X(x,y,z)$ in
$D_W\Box D_E$, then $\psi_X(x,y,z)$ in~$D_X$.}

For $\Gamma$ a list and $S$ a subset of the set of members of
$\Gamma$, let $\Gamma|_S$ be the list obtained from $\Gamma$ by
keeping the elements of $S$ and omitting the others. We have the
following as a consequence of Lemma 3.1.2.

\prop{Lemma 3.1.3}{If $\Lambda$ is a list of $X(D_W\Box D_E)$
grounded in $D_W\Box D_E$, then $\Lambda|_{ X(D_X)}$ is grounded
in~$D_X$.}

We can prove also the following for $V_C$ being the set of
vertices defined as in \S 1.3, which is equal to $E(D_W)\cap
W(D_E)$.

\prop{Lemma 3.1.4}{If $y\in \bar{X}(D_{\bar{X}})$, while
$x,z\in\bar{X}(D_X)$, $u\in V_C$ and $\psi_{\bar{X}}(y,x,z)$ in
$D_W\Box D_E$, then $\psi_{\bar{X}}(u,x,z)$ in~$D_X$.}

\dkz To help the intuition, suppose $X$ is $W$, and take a semipath
$\sigma$ of $D_W$ in $[u,z]$ and a semipath $\tau$ of $D_W$ in $[x]_W$. Since
$D_E$ is weakly connected, there is a semipath $\sigma'$ of $D_E$
in $[y,u]$. Out of $\sigma'$ and $\sigma$ we build a semipath
$\sigma''$ of $D_W\Box D_E$ in $[y,z]$ obtained by omitting the
vertex $u$ and one copy of the edge $a$ such that $E_W(a)=u$ in
$D_W$ and $W_E(a)=u$ in $D_E$. Since $\sigma''$ intersects $\tau$ in
a vertex of $D_W$, we infer that $\sigma$ intersects $\tau$ in the
same vertex. We proceed analogously when $X$ is~$E$. \qed

We prove the following in an analogous manner.

\prop{Lemma 3.1.5}{If $y\in\bar{X}(D_{\bar{X}})$, while
$x,z\in\bar{X}(D_X)$, $u\in V_C$ and $\psi_{\bar{X}}(x,y,z)$ in
$D_W\Box D_E$, then $\psi_{\bar{X}}(x,u,z)$ in $D_X$.}

\noindent The difference with the preceding proof is that we
enlarge a semipath $\tau$ of $D_W$ in $[u]_W$, instead of a semipath $\sigma$ of $D_W$ in $[x,z]$, with a semipath $\tau'$ of $D_E$ in $[y,u]$  so as to obtain a semipath $\tau''$ of $D_W\Box D_E$ in $[y]_W$.

The following lemmata are established immediately.

\prop{Lemma 3.1.6}{If $x,y,z\in\bar{X}(D_X)$ and
$\psi_{\bar{X}}(x,y,z)$ in $D_W\Box D_E$, then
$\psi_{\bar{X}}(x,y,z)$ in $D_X$.}

\vspace{-2ex}

\prop{Lemma 3.1.7}{If $x,z\in X(D_X)$ and $y\in X(D_{\bar{X}})$,
then we do not have $\psi_{\bar{X}}(x,y,z)$ in $D_W\Box D_E$.}

The following is an immediate consequence of Lemma 3.1.7.

\prop{Lemma 3.1.8}{If $\Lambda$ is a list of $X(D_W\Box D_E)$
grounded in $D_W\Box D_E$, while $x,z\in X(D_X)$ and $y\in
X(D_{\bar{X}})$, then we do not have $\Lambda\!:x\mn y\mn z$.}

\vspace{-2ex}

\section{\large\bf Pivots and their ordering}\label{3.2}
\markright{\S 3.2. \quad Pivots and their ordering}

The preliminary results of this section are about a
particular ordering of some vertices of D-graphs that we call pivots,
which we are now going to define.

For $x$ and $y$ vertices of a D-graph, we write $y\prec_\exists x$
when $x\neq y$ and $y$ occurs in some $\sigma$ in $[x]_W$, which using the abbreviated notation introduced in \S 1.2 may be written $y\rhd\sigma$.

We write $y\prec x$, and say that $y$ is a \emph{pivot}\index{pivot} of $x$, when $y\neq
x$ and $y$ occurs in every $\sigma$ in $[x]_W$. We have analogous
notions with $W$ replaced by $E$, but for the sake of definiteness
we concentrate on~$W$.

For a semipath $\sigma$ and any vertices $x$ and $y$ in $\sigma$
consider the sequence $\sigma'$ of vertices and edges of $\sigma$
that make a semipath from $x$ to $y$. We call $\sigma'$ a
\emph{subsemipath}\index{subsemipath} of $\sigma$ from $x$ to $y$, and write
$\sigma_{[x,y]}$ for $\sigma'$. (Note that in this definition the
order of vertices and edges in $\sigma'$ may either coincide or be
converse to the order of~$\sigma$.)

For every semipath $\sigma$ in $[x,y]$ and every semipath
$\tau$ in $[y,z]$ there is a vertex $v$ common to $\sigma$ and $\tau$
such that the subsemipath $\sigma_{[x,v]}$ and the subsemipath
$\tau_{[v,z]}$ have no vertex in common except $v$. (Formally,
this is shown by an induction in which the basis is the case where
$v$ is $y$; in the induction step, where $y'$ is a vertex common
to $\sigma$ and $\tau$ that differs from $y$, we apply the
induction hypothesis to $\sigma_{[x,y']}$ and $\tau_{[y',z]}$.) We
designate by $\sigma\ast\tau$ the semipath in $[x,z]$ obtained by
concatenating $\sigma_{[x,v]}$ and $\tau_{[v,z]}$, with one of the
two occurrences of $v$ deleted.

Note that $\ast$ is not an operation because $\sigma\ast\tau$ is
not uniquely determined by $\sigma$ and $\tau$. For example, in
the D-graph
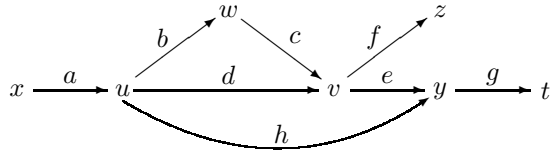
\begin{center}
\begin{picture}(200,60)(0,5)

\put(44,30){\vector(1,0){70}}

\qbezier(40,26)(100,-10)(154,26) \put(154,26){\vector(3,2){2}}

\put(6,30){\vector(1,0){28}} \put(45,35){\vector(4,3){30}}
\put(85,57){\vector(4,-3){30}} \put(125,35){\vector(4,3){30}}
\put(126,30){\vector(1,0){28}} \put(166,30){\vector(1,0){28}}

\put(0,30){\makebox(0,0){$x$}} \put(40,30){\makebox(0,0){$u$}}
\put(80,60){\makebox(0,0){$w$}} \put(120,30){\makebox(0,0){$v$}}
\put(160,30){\makebox(0,0){$y$}} \put(160,60){\makebox(0,0){$z$}}
\put(200,30){\makebox(0,0){$t$}}

\put(20,35){\makebox(0,0){$a$}} \put(55,50){\makebox(0,0){$b$}}
\put(105,50){\makebox(0,0){$c$}}

\put(80,36){\makebox(0,0){$d$}}

\put(100,13.5){\makebox(0,0){$h$}}
\put(135,50){\makebox(0,0){$f$}} \put(140,35){\makebox(0,0){$e$}}
\put(180,35){\makebox(0,0){$g$}}

\end{picture}
\end{center}
if $\sigma$ is $xaubwcvey$ and $\tau$ is $yhudvfz$, then
$\sigma\ast\tau$ can be either $xaubwcvfz$, which is
$\sigma_{[x,v]}$ concatenated with $\tau_{[v,z]}$, or $xaudvfz$,
which is $\sigma_{[x,u]}$ concatenated with $\tau_{[u,z]}$.

We can prove the following.

\prop{Lemma 3.2.1}{If $v_1\prec x$, $v_2\prec x$ and $v_1\neq
v_2$, then $v_2\prec v_1$ or $v_1\prec v_2$.}

\dkz Suppose $v_1\prec x$, $v_2\prec x$ and $v_1\neq v_2$, and
suppose there is a $\sigma_1$ in $[v_1]_W$ in which $v_2$ does not
occur and a $\sigma_2$ in $[v_2]_W$ in which $v_1$ does not occur. Let
$\sigma\in[x]_W$, and suppose $v_1$ is between $x$ and $v_2$ in
$\sigma$. Then $\sigma_{[x,v_1]}\ast\sigma_1$ is a semipath in
$[x]_W$ in which $v_2$ does not occur, which contradicts $v_2\prec
x$. We reason analogously when $v_2$ is between $x$ and $v_1$
in~$\sigma$. \qed

For $\sigma\in[x]_W$ and an arbitrary vertex $y$ in $\sigma$, let
$\sigma_{[y]_W}$ be the subsemipath of $\sigma$ obtained by
rejecting from $\sigma$ everything that occurs in the subsemipath
of $\sigma$ from $x$ to $y$ except~$y$.

\prop{Lemma 3.2.2}{If $v_2\prec v_1$, then not $v_1\prec_\exists
v_2$.}

\dkz Suppose $v_2\prec v_1$, and suppose there is a
$\sigma$ in $[v_2]_W$ in which $v_1$ occurs. Take the semipath
$\sigma_{[v_1]_W}$. Since $v_2\prec v_1$, we have that $v_2\rhd\sigma_{[v_1]_W}$, which implies that $v_2$ occurs
twice in $\sigma$. This contradicts the assumption that $\sigma$
is a semipath. \qed

As a consequence of Lemma 3.2.2, and of the nonemptiness of
$[x]_W$ for every $x$, we obtain the following.

\prop{Lemma 3.2.3}{If $v_2\prec v_1$, then not $v_1\prec v_2$.}

This means that we could replace ``$v_2\prec v_1$ or $v_1\prec
v_2$'' in Lemma 3.2.1 by ``$v_2\prec v_1$ or $v_1\prec v_2$, but
not both''. We also have that $\prec$ is transitive.

\prop{Lema 3.2.4}{If $x\prec y$ and $y\prec z$, then $x\prec z$.}

\dkz Suppose $x\prec y$ and $y\prec z$. By Lemma 3.2.3 we have
that $x\neq z$, and then it is clear that $x\prec z$. \qed

Let $y\rhd[x,z]$ mean that $y$ occurs in every semipath $\sigma$
in $[x,z]$. It is clear that $y\rhd[x,z]$ iff $y\rhd[z,x]$. Then
we can prove the following.

\prop{Lemma 3.2.5}{If $v_2\prec_\exists v_1$ and $v_1\prec x$,
then $v_1\rhd[x,v_2]$.}

\dkz Suppose $\sigma\in[x,v_2]$ and not $v_1\rhd\sigma$. Let $\tau\in[v_1]_W$ and $v_2\rhd\tau$.
There must be such a semipath because $v_2\prec_\exists v_1$. Then
$\sigma\ast\tau_{[v_2]_W}\in[x]_W$, but not $v_1\rhd\sigma\ast\tau_{[v_2]_W}$, which contradicts $v_1\prec x$. \qed

\vspace{-2ex}

\prop{Lemma 3.2.6}{If $v_2\rhd[x,z]$ and not $v_2\rhd[x,v_1]$,
then $v_2\rhd[v_1,z]$.}

\dkz Suppose $\sigma\in[v_1,z]$ and not $v_2\rhd\sigma$. Because not $v_2\rhd[x,v_1]$, for some $\tau$ in $[x,v_1]$
we have that not $v_2\rhd\tau$. Then
$\tau\ast\sigma\in[x,z]$ and not $v_2\rhd\tau\ast\sigma$, contradicting $v_2\rhd[x,z]$. \qed

\vspace{-2ex}

\prop{Lemma 3.2.7}{If $v_1\rhd[z,y]$ and $v_2\rhd[z,v_1]$, then
$v_2\rhd[z,y]$.}

\dkz Let $\sigma\in[z,y]$. Since $v_1\rhd[z,y]$, we have
$v_1\rhd\sigma$. Then
since $v_2\rhd[z,v_1]$, we have $v_2\rhd\sigma_{[z,v_1]}$, and
hence $v_2\rhd\sigma$. \qed

As a corollary of Lemma 3.2.5, we have the following.

\prop{Lemma 3.2.8.1}{If $v_1\prec v_2$ and $v_2\prec x$, then
$v_2\rhd[x,v_1]$.}

\vspace{-2ex}

\prop{Lemma 3.2.8.2}{If $v_1\prec v_2$ and $v_2\prec x$ and
$v_1\rhd[z,x]$, then $v_2\rhd[z,x]$.}

\dkz In Lemma 3.2.7 put $x$ for $z$ and $z$ for $y$, and then use
Lemma 3.2.8.1. \qed

\vspace{-2ex}

\prop{Lemma 3.2.8.3}{If $v_1\prec v_2$, $v_2\prec x$ and
$\sigma\in[x,v_2]$, then not $v_1\rhd\sigma$.}

\dkz Suppose $v_1\prec v_2$, $v_2\prec x$ and for some
$\sigma$ in $[x,v_2]$ we have $v_1\rhd\sigma$. From $v_1\prec v_2$ we
have that $v_1\neq v_2$, and, by Lemma 3.2.8.1, we have
$v_2\rhd\sigma_{[x,v_1]}$. But then we would have two occurrences
of $v_2$ in $\sigma$, which is impossible, because $\sigma$ is a
semipath. \qed

As a corollary of Lemma 3.2.8.3, we have the following.

\prop{Lemma 3.2.8.4}{If $v_1\prec v_2$ and $v_2\prec x$, then not
$v_1\rhd[x,v_2]$.}

\vspace{-2ex}

\section{\large\bf Further results for the Pivot Theorem}\label{3.3}
\markright{\S 3.3. \quad Further results for the Pivot
Theorem}

The lemmata of this section give further preliminary results for the Pivot Theorem of \S 3.4. They are grouped together because of their common combinatorial inspiration.

Let $x$, $y_1,\ldots,y_n$, $z$ and $v_1,\ldots,v_n$ be vertices of
a D-graph. Consider the following condition:
\begin{tabbing}
\hspace{1.6em}\=$(yji)$\hspace{2em}\=$v_j\rhd[z,x]\quad \&\quad
v_j\rhd[z,y_i]\quad \&\quad v_j\prec x\quad \&\quad v_j\prec y_i$.
\end{tabbing}
Then we can prove the following for $n\geq 1$.

\prop{Lemma 3.3.1}{If for every $i\in\{1,\ldots,n\}$ we have
$(yii)$, then for some $j\in\{1,\ldots,n\}$ for every
$i\in\{1,\ldots,n\}$ we have $(yji)$.}

\dkz We proceed by induction on $n$. If $n=1$, then the lemma
holds trivially. For the induction step, when $n>1$, suppose that
\vspace{.5ex}

$(\ast)$ \quad for every $i\in\{1,\ldots,n\}$ we have $(yii)$.

\vspace{.5ex}
\noindent Hence for every $k\in\{1,\ldots,n\mn 1\}$ we have
$(ykk)$, and by the induction hypothesis there is an element of
$\{1,\ldots,n\mn 1\}$, which we call $m$, such that

\vspace{.5ex}

$(\ast\ast)$ \quad for every $k\in\{1,\ldots,n\mn 1\}$ we have $(ymk)$.

\vspace{.5ex}

If $v_m=v_n$, then from $(ynn)$, which follows from $(\ast)$, and
from $(\ast\ast)$, we obtain for $j=m=n$ that

\vspace{.5ex}

${(\ast\!\ast\!\ast)}$ \quad for every $i\in\{1,\ldots,n\}$ we have $(yji)$.

\vspace{.5ex}

Suppose $v_m\neq v_n$. From $(\ast)$ we conclude that we have
$v_m\prec x$ and $v_n\prec x$. Then by Lemma 3.2.1 we have
$v_n\prec v_m$ or $v_m\prec v_n$.

Suppose that

\vspace{.5ex}

(1) \quad $v_n\prec v_m$.

\vspace{.5ex}

\noindent We want to show ${(\ast\!\ast\!\ast)}$ for $j$ being $n$, and
to achieve that, since we have $(ynn)$, it suffices to show
$(ynk)$ for every $k\in\{1,\ldots,n\mn 1\}$. We have
\begin{tabbing}
\hspace{1.6em}\=(1.1)\hspace{2em}\=$v_n\rhd[z,x]$,\quad by
$(ynn)$,
\\[.5ex]
\>(1.2)\>$v_m\prec x$, \quad by $(ymm)$, which follows from
$(\ast)$,
\\[.5ex]
\>(1.3)\>not $v_n\rhd[x,v_m]$,\quad from (1) and (1.2), by Lemma
3.2.8.4,
\\[.5ex]
\>(1.4)\>$v_n\rhd[v_m,z]$,\quad from (1.1) and (1.3), by Lemma
3.2.6,
\\[.5ex]
\>(1.5)\>$v_m\rhd[z,y_k]$,\quad by $(ymk)$, which follows from
$(\ast\ast)$,
\\[.5ex]
\>(1.6)\>$v_n\rhd[z,y_k]$,\quad from (1.5) and (1.4), by Lemma
3.2.7,
\\[.5ex]
\>(1.7)\>$v_n\prec x$,\quad by $(ynn)$,
\\[.5ex]
\>(1.8)\>$v_m\prec y_k$,\quad by $(ymk)$,
\\[.5ex]
\>(1.9)\>$v_n\prec y_k$,\quad from (1) and (1.8), by Lemma 3.2.4.
\end{tabbing}
With (1.1), (1.6), (1.7) and (1.9) we have the four conjuncts of
$(ynk)$, for every $k\in\{1,\ldots,n\mn 1\}$.

Suppose now that

\vspace{.5ex}

(2) \quad $v_m\prec v_n$.

\vspace{.5ex}

\noindent We want to show ${(\ast\!\ast\!\ast)}$ for $j$ being $m$, and
to achieve that, since we have $(\ast\ast)$, it suffices to show
$(ymn)$. We have a derivation of $(ymn)$ obtained from the
derivation of $(ynk)$ above by putting $m$ for $n$ and $n$ for $m$
and~$k$.

So, in any case, from $(\ast)$ we have inferred ${(\ast\!\ast\!\ast)}$ for
some $j$, which yields our lemma. \qed

Let $x$, $y$, $z_1,\ldots,z_n$ and $v_1,\ldots,v_n$ be vertices of
a D-graph. Consider the following condition:
\begin{tabbing}
\hspace{1.6em}$(zji)$\hspace{2em}$v_j\rhd[z_i,x]\quad \&\quad
v_j\rhd[z_i,y]\quad \&\quad v_j\prec x\quad \&\quad v_j\prec y$.
\end{tabbing}
Then we can prove the following for $n\geq 1$.

\prop{Lemma 3.3.2}{If for every $i\in\{1,\ldots,n\}$ we have
$(zii)$, then for some $j\in\{1,\ldots,n\}$ for every
$i\in\{1,\ldots,n\}$ we have $(zji)$.}

\dkz We proceed by induction on $n$. If $n=1$, then the lemma
holds trivially. For the induction step, when $n>1$, suppose that

\vspace{.5ex}

$(\ast)$ \quad for every $i\in\{1,\ldots,n\}$ we have $(zii)$.

\vspace{.5ex}

\noindent Hence for every $k\in\{1,\ldots,n\mn 1\}$ we have
$(zkk)$, and by the induction hypothesis there is an element of
$\{1,\ldots,n\mn 1\}$, which we call $m$, such that

\vspace{.5ex}

$(\ast\ast)$ \quad for every $k\in\{1,\ldots,n\mn 1\}$ we have $(zmk)$.

\vspace{.5ex}

If $v_m=v_n$, then from $(znn)$, which follows from $(\ast)$, and
from $(\ast\ast)$, we obtain for $j=m=n$ that

\vspace{.5ex}

${(\ast\!\ast\!\ast)}$ \quad for every $i\in\{1,\ldots,n\}$ we have $(zji)$.

\vspace{.5ex}

Suppose $v_m\neq v_n$. From $(\ast)$ we conclude that we have
$v_m\prec x$ and $v_n\prec x$. Then by Lemma 3.2.1 we have
$v_m\prec v_n$ or $v_n\prec v_m$.

Suppose that

\vspace{.5ex}

(1) \quad $v_m\prec v_n$.

\vspace{.5ex}

\noindent We want to show ${(\ast\!\ast\!\ast)}$ for $j$ being $n$, and
to achieve that, since we have $(znn)$, it suffices to show
$(znk)$ for every $k\in\{1,\ldots,n\mn 1\}$.

Take a semipath $\sigma$ in $[z_k,x]$. Then we have
\begin{tabbing}
\hspace{1.6em}\=(1.1)\hspace{2em}\=$v_m\rhd[z_k,x]$,\quad by
$(zmk)$, which follows from $(\ast\ast)$,
\\[.5ex]
\>(1.2)\>$v_n\prec x$, \quad by $(znn)$,
\\[.5ex]
\>(1.3)\>$v_n\rhd[v_m,x]$,\quad from (1) and (1.2), by Lemma
3.2.8.1,
\\[.5ex]
\>(1.4)\>$v_n\rhd\sigma_{[v_m,x]}$,\quad from (1.1) and (1.3),
\\[.5ex]
\>(1.5)\>$v_n\rhd\sigma$,\quad from (1.4).
\end{tabbing}
So we have $v_n\rhd[z_k,x]$. We derive analogously
$v_n\rhd[z_k,y]$, the second conjunct of $(znk)$, and we have the
third and fourth conjunct by $(znn)$.

Suppose now that

\vspace{.5ex}

(2) \quad $v_n\prec v_m$.

\vspace{.5ex}

\noindent We want to show ${(\ast\!\ast\!\ast)}$ for $j$ being $m$, and
to achieve that, since we have $(\ast\ast)$, it suffices to show
$(zmn)$. We have a derivation of $(zmn)$ obtained from the
derivation of $(znk)$ above by putting $m$ for $n$ and $n$ for $m$
and~$k$.

So, in any case, from $(\ast)$ we have inferred ${(\ast\!\ast\!\ast)}$ for
some $j$, which yields our Lemma. \qed

We can also prove the following two lemmata.

\prop{Lemma 3.3.3.1}{Let $\sigma\in[x,r]$ and $\tau\in[y,p]$. If
$v\rhd\sigma$, $v\rhd[y,r]$ and not $v\rhd\tau$, then
$\sigma_{[v,r]}$ and $\tau$ do not intersect.}

\dkz If $\sigma_{[v,r]}$ and $\tau$ intersect in $s$, then not
$v\rhd \tau_{[y,s]}\ast\sigma_{[s,r]}$, since not $v\rhd\tau$;
this contradicts $v\rhd[y,r]$. \qed

\vspace{-2ex}

\prop{Lemma 3.3.3.2}{Let $\sigma\in[x,r]$ and $\tau\in[q,o]$. If
$v\rhd\sigma$, $v\rhd[x,o]$, and not $v\rhd\tau$, then
$\sigma_{[x,v]}$ and $\tau$ do not intersect.}

\dkz If $\sigma_{[x,v]}$ and $\tau$ intersect in $s$, then not
$v\rhd \sigma_{[x,s]}\ast\tau_{[s,o]}$, since not $v\rhd\tau$;
this contradicts $v\rhd[x,o]$. \qed

As a matter of fact, one of these lemmata can be derived from the other, but it is easier to give independent proofs than to make these derivations by appropriate substitutions.

\section{\large\bf The Pivot Theorem}\label{3.4}
\markright{\S 3.4. \quad The Pivot Theorem}

In this section we give the proof of the following theorem, which is the Pivot Theorem that we have announced.

\begin{samepage}
\prop{Theorem 3.4.1}{If $\psi_E(x,y,z)$ and $\psi_E(y,x,z)$, then
there is a vertex $v$ such that}
\vspace{-2ex}
\begin{tabbing}
\hspace{1.6em}$(\ast)$\hspace{2em}$v\rhd[z,x]\quad
\&\quad v\rhd[z,y]\quad \&\quad v\prec x\quad \&\quad v\prec y$.
\end{tabbing}
\end{samepage}

\dkz Assume $\psi_E(x,y,z)$ and $\psi_E(y,x,z)$. We proceed by
induction on the number of inner vertices in our D-graph. The
basis, where it has just one inner vertex, is trivial, because $v$
is that vertex.

For the induction step, suppose our D-graph is $D_W\Box D_E$. With
respect to the distribution of $x$, $y$ and $z$ in $D_W$ and $D_E$
we have the following cases:

\begin{center}
\begin{tabular}{c|c|c}
 & $E(D_W)$ & $E(D_E)$  \\[.5ex]
\hline (I) & & $x,y,z$ \\[.5ex]
\hline (II$x$) & $x,z$ & $y$
\\[.5ex]  (II$y$) & $y,z$ & $x$ \\[.5ex]  (III) & $x,y$ &
$z$
\\[.5ex] \hline (IV$x$) & $x$ & $y,z$
\\[.5ex] (IV$y$) & $y$ & $x,z$
\\[.5ex] \hline (V) & $x,y,z$ &
\\[.5ex] \hline (VI) & $z$ & $x,y$
\end{tabular}
\end{center}
The order of these cases is dictated by the structure of our
proof.

Note first that the cases (IV$x$) and (IV$y$) are impossible
because of Lemma 3.1.7. The other cases are possible and will now
be treated.

\vspace{.5ex}

(I) Let $x$, $y$ and $z$ be distinct vertices of $E(D_E)$. (Their
distinctness is a consequence of the assumptions $\psi_E(x,y,z)$
and $\psi_E(y,x,z)$.) By Lemma 3.1.2, we have $\psi_E(x,y,z)$ and
$\psi_E(y,x,z)$ in $D_E$. So by the induction hypothesis we have a
$v$ such that $(\ast)$ for~$D_E$. We will show that for that $v$ we have $(\ast)$ for $D_W\Box D_E$.

For the first conjunct of $(\ast)$ for $D_W\Box D_E$, take a
semipath $\sigma$ of $D_W\Box D_E$ in $[z,x]$. If $\sigma$ is in
$D_E$, then $v\rhd\sigma$ by the induction hypothesis. If $\sigma$
is not in $D_E$, then by replacing in a subsemipath of $\sigma$ a single
vertex of $D_W$ by a vertex in $V_C$ we obtain a semipath
$\sigma'$ of $D_E$ in $[x]_W$. By the induction hypothesis we have
$v\rhd\sigma'$, because $v\prec x$ in $D_E$; hence $v\rhd\sigma$.
So we have the first conjunct of $(\ast)$ for $D_W\Box D_E$.

For the third conjunct of $(\ast)$ for $D_W\Box D_E$, take a
semipath $\sigma$ of $D_W\Box D_E$ in $[x]_W$. Then by replacing
in a subsemipath of $\sigma$ a single vertex of $D_W$ by a vertex in
$V_C$ we obtain a semipath $\sigma'$ of $D_E$ in $[x]_W$. By the
induction hypothesis, as above, we have $v\rhd\sigma'$; hence
$v\rhd\sigma$. So we have the third conjunct of $(\ast)$ for
$D_W\Box D_E$.

For the second and fourth conjunct we proceed analogously by
replacing $x$ by $y$. So we have $(\ast)$ for $D_W\Box D_E$.

\vspace{.5ex}

(II$x$) Let $x$ and $z$ be distinct vertices of $E(D_W)$, and let
$y$ be a vertex of $E(D_E)$. Let $V_C=\{y_1,\ldots,y_n\}$; here
$n\geq 1$. By Lemmata 3.1.4 and 3.1.5 for every
$i\in\{1,\ldots,n\}$ we obtain $\psi_E(y_i,x,z)$ and
$\psi_E(x,y_i,z)$ in $D_W$, and then by applying the induction
hypothesis $n$ times to $D_W$ we obtain $(yii)$ for every
$i\in\{1,\ldots,n\}$. By Lemma 3.3.1, for some
$j\in\{1,\ldots,n\}$ for every $i\in\{1,\ldots,n\}$ we have
$(yji)$. We will show $(\ast)$ for $D_W\Box D_E$ with $v$ being~$v_j$.

For the first conjunct, $v_j\rhd[z,x]$, suppose $\sigma$ is a
semipath of $D_W\Box D_E$ in $[z,x]$. If $\sigma$ is a semipath of
$D_W$, then we use the first conjunct of $(yji)$. If $\sigma$
passes through $D_E$, then by replacing in a subsemipath of $\sigma$ a
single vertex of $D_E$ by a vertex $y_k$ for some
$k\in\{1,\ldots,n\}$, we obtain a semipath $\sigma'$ of $D_W$ in
$[z,y_k]$. Then we use the second conjunct of $(yjk)$, and we have
that $v_j\rhd\sigma'$; hence $v_j\rhd\sigma$.

We proceed analogously for the remaining conjuncts $v_j\rhd[z,y]$,
$v_j\prec x$ and $v_j\prec y$. We look for a subsemipath which is
either a semipath of $D_W$ or which after replacement of a single
vertex by $y_k$ becomes a semipath of $D_W$. Then we apply
$(yjk)$. So we have $(\ast)$ for $D_W\Box D_E$. We proceed
analogously for~(II$y$).

\vspace{.5ex}

(III) Let $x$ and  $y$ be distinct vertices of $E(D_W)$, and let
$z$ be a vertex of $E(D_E)$. Let $V_C=\{z_1,\ldots,z_n\}$; here
$n\geq 1$. By Lemma 3.1.4 where $\psi_E(y,x,z)$ is replaced by
$\psi_E(z,x,y)$ and $\psi_E(z,y,x)$, for every
$i\in\{1,\ldots,n\}$ we obtain $\psi_E(z_i,x,y)$ and
$\psi_E(z_i,y,x)$ in $D_W$, and then by applying the induction
hypothesis $n$-times to $D_W$ we obtain $(zii)$ for every
$i\in\{1,\ldots,n\}$. By Lemma 3.3.2 , for some
$j\in\{1,\ldots,n\}$ for every $i\in\{1,\ldots,n\}$ we have
$(zji)$. We will show $(\ast)$ for $D_W\Box D_E$ with $v$ being
$v_j$. We proceed as in case (II$x$), with the help of vertex
$z_k$ from $V_C$ instead of $y_k$, if need there~is.

\vspace{.5ex}

(V) Let $x$, $y$ and $z$ be distinct vertices of $E(D_W)$. By
Lemma 3.1.6, we have that $\psi_E(x,y,z)$ and $\psi_E(y,x,z)$ in
$D_W$. Consider the set of vertices $v$ of $D_W$ such that
$(\ast)$ holds in $D_W$. By the induction hypothesis, this set is
non-empty, and let $v_1,\ldots,v_n$, for $n\geq 1$, be all its
elements. By Lemmata 3.2.1 and 3.2.4, we have that this set is
linearly ordered by the relation $\prec$. Let us assume that we
have $v_1\prec\ldots\prec v_n$.

We have two cases to consider. Assume first that in $D_W$ we have
that
\begin{tabbing}
\hspace{1.6em}(V.1)\hspace{2em}$(\forall{u}\in
V_C)(v_n\rhd[u,x]\quad\&\quad v_n\rhd[u,y])$.
\end{tabbing}
We will show $(\ast)$ for $D_W\Box D_E$ with $v$ being $v_n$.

For the first conjunct, $v_n\rhd[z,x]$, suppose $\sigma$ is a
semipath of $D_W\Box D_E$ in $[z,x]$. If $\sigma$ is a semipath of
$D_W$, then we are done. If $\sigma$ passes through $D_E$, then by
replacing in a subsemipath of $\sigma$ a single vertex of $D_E$ by a
vertex $u$ of $V_C$, we obtain a semipath $\sigma'$ of $D_W$ in
$[u,x]$. Then, by (V.1) we conclude that $v_n\rhd\sigma'$; hence
$v_n\rhd\sigma$. We prove the second conjunct analogously.

For the third conjunct, $v_n\prec x$, we take a semipath $\sigma$
of $D_W\Box D_E$ in $[x]_W$. If $\sigma$ is a semipath of $D_W$,
then we are done. If $\sigma$ passes through $D_E$, then as above
we obtain a semipath $\sigma'$ of $D_W$ in $[x,u]$ for $u\in V_C$,
and reason as above with (V.1). For the fourth conjunct we proceed
analogously, and hence we have $(\ast)$ for $D_W\Box D_E$.

If for some $u$ in $V_C$ and for $t$ being $x$ or $y$ we have in
$D_W$
\begin{tabbing}
\hspace{1.6em}\=(V.1)\hspace{2em}\= \kill
\>(V.2)\>not
$v_n\rhd[u,t]$,
\end{tabbing}
then we show $(\ast)$ for $D_W\Box D_E$ with $v$ being $v_1$.

We prove first that in $D_W$ we have
\begin{tabbing}
\hspace{1.6em}\=(V.1)\hspace{2em}\= \kill \>$(\dagger)$\>not
$(\psi_E(y,x,u)\quad\&\quad \psi_E(x,y,u))$.
\end{tabbing}
Suppose not $(\dagger)$. Then, by the induction hypothesis applied
to $D_W$, we have a vertex $w$ such that
\begin{tabbing}
\hspace{1.6em}\=(V.1)\hspace{2em}\= \kill
\>$(\dagger\dagger)$\>$w\rhd[u,x]\quad\&\quad
w\rhd[u,y]\quad\&\quad w\prec x\quad\&\quad w\prec y$.
\end{tabbing}
Since $w\rhd[u,t]$, by (V.2), we obtain $w\neq v_n$. Since
$v_n\prec t$, $w\prec t$ and $v_n\neq w$, we
obtain $v_n\prec w$ or $w\prec v_n$ by Lemma 3.2.1. We will show that not $w\prec
v_n$.

By (V.2) we have a semipath $\sigma$ of $D_W$ in $[u,t]$ such that
not $v_n\rhd\sigma$. From the conjunct $w\rhd[u,t]$ of
$(\dagger\dagger)$, we have $w\rhd\sigma$. Then we have not
$v_n\rhd\sigma_{[w,t]}$. By putting $v_n$, $t$ and $w$ for $v_2$,
$x$ and $v_1$ respectively in Lemma 3.2.8.1, we obtain not $w\prec
v_n$. Hence we have $v_n\prec w$.

We show $(\ast)$ for $D_W$ with $v$ being $w$. The third and
fourth conjunct are given by the third and fourth conjunct of
$(\dagger\dagger)$. For the first conjunct, $w\rhd[z,x]$, we apply
Lemma 3.2.8.2 with $v_1$ and $v_2$ being $v_n$ and $w$
respectively. For the second conjunct, $w\rhd[z,y]$, we proceed
analogously, and hence we have $(\ast)$ for $D_W$ with $v$ being
$w$. So $w\in\{v_1,\ldots,v_n\}$, but this is in contradiction
with $v_1\prec\ldots\prec v_n$ and $v_n\prec w$. So we can infer~$(\dagger)$.

Suppose not $\psi_E(y,x,u)$. So there is a semipath $\rho$
of $D_W$ in $[y,u]$, and a semipath $\pi$ of $D_W$ in $[x,r]$, for
$r$ in $W(D_W)$, which do not intersect. We can show
first that
\begin{tabbing}
\hspace{1.6em}\=(V.1)\hspace{2em}\= \kill \>(1)\>for every $v$ in
$W(D_E)$ we have $v\in V_C$.
\end{tabbing}

Suppose $v\in W(D_E)$ and $v\notin V_C$. Take a semipath $\sigma$
of $D_W$ in $[x,z]$, and take an $i\in\{1,\ldots, n\}$; we have
$v_i\rhd\sigma$, by $(\ast)$ for $D_W$. By putting $v_i$, $z$,
$\rho$ and $u$ for respectively $v$, $r$, $\tau$ and $p$ in Lemma
3.3.3.1 we obtain that $\sigma_{[v_i,z]}$ and $\rho$ do not
intersect. Let $u=W_E(a_u)$ in $D_E$; so $a_u\in C$. Let $\varphi$
be a semipath of $D_E$ in $[E_E(a_u),v]$, and let $\rho'$ be
obtained from $\rho$ by replacing $u$ by $E_E(a_u)$. Then
$\rho'\ast\varphi$, a semipath of $D_W\Box D_E$ in $[y]_W$, and
$\pi_{[x,v_i]}\ast\sigma_{[v_i,z]}$, a semipath of $D_W$, and
hence of $D_W\Box D_E$, in $[x,z]$, do not intersect, which
contradicts the assumption that $\psi_E(x,y,z)$ in $D_W\Box D_E$.
This proves~(1).

For every $v$ in $V_C$, we can show that $\psi_E(x,v,z)$ and
$\psi_E(v,x,z)$. We prove first that $\psi_E(x,v,z)$, which is
similar to the proof of (1) we have just given. Suppose not
$\psi_E(x,v,z)$ in $D_W$. Then there is a semipath $\sigma$ of
$D_W$ in $[x,z]$ and a semipath $\tau$ of $D_W$ in $[v]_W$ that do
not intersect. Take an $i\in\{1,\ldots,n\}$; we have
$v_i\rhd\sigma$. By putting $v_i$, $z$, $\rho$ and $u$ for $v$,
$r$, $\tau$ and $p$ in Lemma 3.3.3.1 we obtain that
$\sigma_{[v_i,z]}$ and $\rho$ do not intersect. By putting $v_i$,
$\pi$ and $\tau$ for $v$, $\sigma$ and $\tau$ in Lemma 3.3.3.2, we
obtain that $\pi_{[x,v_i]}$ and $\tau$ do not intersect.

Let $u=W_E(a_u)$ and $v=W_E(a_v)$ in $D_E$; so $a_u,a_v\in C$. Let
$\varphi$ be a semipath of $D_E$ in $[E_E(a_u),E_E(a_v)]$, and let
$\rho'$ and $\tau'$ be obtained from $\rho$ and $\tau$ by
replacing $u$ and $v$ by $E_E(a_u)$ and $E_E(a_v)$ respectively. Then
$\rho'\ast\varphi\ast\tau'$, a semipath of $D_W\Box D_E$ in
$[y]_W$, and $\pi_{[x,v_i]}\ast\sigma_{[v_i,z]}$, a semipath of
$D_W\Box D_E$ in $[x,z]$, do not intersect, which contradicts the
assumption that $\psi_E(x,y,z)$ in $D_W\Box D_E$. So we have
$\psi_E(x,v,z)$.

Suppose not $\psi_E(v,x,z)$ in $D_W$. Then there is a semipath
$\sigma$ of $D_W$ in $[v,z]$ and a semipath $\tau$ of $D_W$ in
$[x]_W$ that do not intersect. Take an $i\in\{1,\ldots,n\}$; we
have $v_i\rhd\tau$. By putting $v_i$, $\tau$, $\rho$ and $u$ for
$v$, $\sigma$, $\tau$ and $p$ in Lemma 3.3.3.1, we obtain that
$\tau_{[v_i,r]}$ and $\rho$ do not intersect. By putting $v_i$,
$\pi$ and $\sigma$ for $v$, $\sigma$ and $\tau$ in Lemma 3.3.3.2,
we obtain that $\pi_{[x,v_i]}$ and $\sigma$ do not intersect. Let
$\varphi$ and $\rho'$ be defined as above, and let $\sigma'$ be
obtained from $\sigma$ by replacing $v$ by the vertex $E_E(a_v)$ of
$D_E$. Then $\rho'\ast\varphi\ast\sigma'$, a semipath of $D_W\Box
D_E$ in $[y,z]$, and $\pi_{[x,v_i]}\ast\tau_{[v_i,r]}$, a semipath
of $D_W\Box D_E$ in $[x]_W$, do not intersect, which contradicts
the assumption that $\psi_E(y,x,z)$ in $D_W\Box D_E$. So we have
$\psi_E(v,x,z)$, and hence we have shown that $\psi_E(x,v,z)$ and
$\psi_E(v,x,z)$.

We apply then the induction hypothesis for $D_W$ and obtain for
every $v$ in $V_C$ a vertex $w$ of $D_W$ such that
\begin{tabbing}
\hspace{1.6em}\=(V.1)\hspace{2em}\= \kill \>$(\ast
vw)$\>$w\rhd[z,x]\quad\&\quad w\rhd[z,v]\quad\&\quad w\prec
x\quad\&\quad w\prec v$.
\end{tabbing}
We can then prove
\begin{tabbing}
\hspace{1.6em}\=(V.1)\hspace{2em}\= \kill \>(2)\>for every $v$ in
$V_C$ we have $v_1\prec v$ in $D_W$.
\end{tabbing}

Suppose for some $v$ in $V_C$ we do not have $v_1\prec v$ in $D_W$.
We can then infer that $w\prec v_1$. Since by $(\ast vw)$ we have
$w\prec v$, we also have $w\neq v_1$, because not $v_1\prec v$.
Since by $(\ast vw)$ we have $w\prec x$, and we have $v_1\prec x$, we obtain that $w\prec v_1$ or $v_1\prec w$ by Lemma 3.2.1. If we
had $v_1\prec w$, with $w\prec v$, we would have $v_1\prec v$,
which contradicts our assumption. Hence we have $w\prec v_1$.

We will then prove that in $D_W$ we have $(\ast yw)$, which is
$(\ast v w)$ with $v$ replaced by $y$, or $(\ast)$ with $v$
replaced by $w$. The first and third conjunct of $(\ast yw)$ are
obtained from the respective conjuncts of $(\ast vw)$. The fourth
conjunct of $(\ast yw)$, namely $w\prec y$, follows immediately
from $w\prec v_1$ and $v_1\prec y$.

For the only remaining conjunct, $w\rhd[z,y]$, we show first that
for every $i\in\{1,\ldots,n\}$ we do not have $w\rhd\pi_{[x,v_i]}$,
by Lemma 3.2.8.3; we put $w$, $v_i$ and $\pi_{[x,v_i]}$ for $v_1$,
$v_2$ and $\sigma$, and we use $w\prec v_i$, which follows from
$w\prec v_1$.

Let $\sigma\in[z,y]$. Since by $(\ast vw)$ we have $w\rhd[x,z]$,
we obtain $w\rhd\pi_{[x,v_i]}\ast\sigma_{[v_i,z]}$. Since not
$w\rhd\pi_{[x,v_i]}$, we have $w\rhd\sigma_{[v_i,z]}$, and hence
$w\rhd\sigma$. So we have $w\rhd[z,y]$, and hence $(\ast yw)$
holds.

So $w\in\{v_1,\ldots,v_n\}$, but $w\prec v_1$ contradicts the
assumption that $v_1\prec\ldots\prec v_n$. So (2) holds. We prove
also the following
\begin{tabbing}
\hspace{1.6em}\=(V.1)\hspace{2em}\= \kill \>(3)\>for every $v$ in
$V_C$ we have $v_1\rhd[v,z]$ in $D_W$.
\end{tabbing}
Suppose not $v_1\rhd[v,z]$. Then for some $\sigma$ in $[v,z]$ we have
not $v_1\rhd\sigma$. By putting $v_1$, $\pi$, $\sigma$, $v$ and
$z$ for $v$, $\sigma$, $\tau$, $y$ and $p$ in Lemma 3.3.3.1, we
obtain that $\pi_{[v_1,r]}$ and $\sigma$ do not intersect; for the
assumption $v_1\rhd[v,r]$ of Lemma 3.3.3.1 after the replacement,
we use $v_1\prec v$. Let $\pi^{-1}$ be the semipath $\pi_{[r,x]}$
in $[r,x]$, which is obtained by taking $\pi$ in reverse order,
and let $\sigma^{-1}$ be~$\sigma_{[z,v]}$.

By putting $v_1$, $\pi^{-1}$, $r$, $x$, $\sigma^{-1}$, $z$ and $v$
for $v$, $\sigma$, $x$, $r$, $\tau$, $y$ and $p$ in Lemma 3.3.3.1,
we obtain that $\pi^{-1}_{[v_1,x]}$ and $\sigma^{-1}$ do not
intersect. Hence $\pi_{[x,v_1]}$ and $\sigma$ do not intersect.
Since $\pi$ is $\pi_{[x,v_1]}\ast\pi_{[v_1,r]}$, we conclude that
$\pi$ and $\sigma$ do not intersect.

By defining $\rho'$, $\varphi$ and $\sigma'$ as before (see the
proofs of $\psi_E(x,v,z)$ and $\psi_E(v,x,z)$ in $D_W$), we have
that $\pi$ and $\rho'\ast\varphi\ast\sigma'$ do not intersect,
which contradicts the assumption that $\psi_E(y,x,z)$ in $D_W\Box
D_E$. So (3) holds.

Now we can show $(\ast)$ for $D_W\Box D_E$ with $v$ being $v_1$.
For $v_1\rhd[z,x]$, the first conjunct, take a semipath $\sigma$ of
$D_W\Box D_E$ in $[z,x]$. If $\sigma$ is a semipath of $D_W$, we
are done, by the induction hypothesis. If $\sigma$ passes through
$D_E$, then let $s$ be the first vertex in $\sigma$ not in $D_W$.
Consider $\sigma_{[z,s]}$, and let $\sigma^*$ be obtained from it
by replacing $s$ by the corresponding $v\in V_C$. (If $s=E_E(a)$ in
$D_E$, then $v=W_E(a)$ in $D_E$; alternatively $v=E_W(a)$ in $D_W$.)
By (3) we have that $v_1\rhd\sigma^*$, and since $v_1\neq v$, we
have $v_1\rhd\sigma$. For $v_1\rhd[z,y]$, the second conjunct, we
proceed analogously.

For $v_1\prec x$, the third conjunct, take a semipath $\sigma$ of
$D_W\Box D_E$ in $[x]_W$. If $\sigma$ is a semipath of $D_W$, we
are done. If $\sigma$ passes through $D_E$, we proceed in
principle as above in order to obtain a semipath $\sigma^*$ of
$D_W$. The vertex $s$ of $D_E$ is now the last vertex of $D_E$ in
$\sigma$. We use (1) to guarantee that $s$, which is not in
$W(D_E)$, has a corresponding vertex $v$ in $V_C$, as above. By
(2) we have that $v_1\rhd\sigma^*$, and hence $v_1\rhd\sigma$. For
$v_1\prec y$, the  last conjunct, we proceed analogously. This
concludes the proof of~(V).

\vspace{.5ex}

(VI) Let $z$ be a vertex of $E(D_W)$, and let $x$ and $y$ be
distinct vertices of $E(D_E)$. We proceed by an auxiliary
induction on the number $n$ of inner vertices of $D_E$. If $n=1$,
then $D_E$ is a basic D-graph, and its unique inner vertex is the
$v$ required by $(\ast)$. If $n>1$, then let $D_E$ be $D_E'\Box
D_E''$. We have the following cases.

\vspace{.5ex}

(1) Suppose $x,y\in E(D_E')$.

\vspace{.5ex}

(1.1) If $D_W\Box D_E'$ is defined, then $D_W\Box(D_E'\Box D_E'')$
is equal to $(D_W\Box D_E')\linebreak\Box D_E''$, and we are in
case (V) for $(D_W\Box D_E')\Box D_E''$, for which we have
$x,y,z\in E(D_W\Box D_E')$. We continue reasoning as for (V)
above.

\vspace{.5ex}

(1.2) If $D_W\Box D_E'$ is not defined, then $D_W\Box(D_E'\Box
D_E'')$ is equal to $D_E'\linebreak\Box(D_W\Box D_E'')$ and we are
in case (III).

\vspace{.5ex}

(2) Suppose $x,y\in E(D_E'')$.

\vspace{.5ex}

(2.1) If $D_W\Box D_E'$ is defined, then we may apply the
hypothesis of the auxiliary induction to $(D_W\Box D_E')\Box
D_E''$, since the number of inner vertices of $D_E''$ is lesser
than $n$.

\vspace{.5ex}

(2.2) If $D_W\Box D_E'$ is not defined, then we are in case (I)
for $D_E'\Box(D_W\linebreak\Box D_E'')$.

\vspace{.5ex}

(3) Suppose $x\in E(D_E')$ and $y\in E(D_E'')$.

\vspace{.5ex}

(3.1) If $D_W\Box D_E'$ is defined, then we are in case (II$x$)
for $(D_W\Box D_E')\Box D_E''$.

\vspace{.5ex}

(3.2) If $D_W\Box D_E'$ is not defined, then we are in case
(IV$x$) for $D_E'\Box(D_W\linebreak\Box D_E'')$, which is
impossible. This concludes our proof of the theorem. \qed

Let us write $y\rhd [x]_X$ when for every $\sigma$ in $[x]_X$ we
have that $y\rhd\sigma$, and let $y\prec_X x$ stand for
$y\rhd[x]_X$ and $y\neq x$. So $y\prec_W x$ iff $y\prec x$. It is
clear that for all we have proven since \S 3.2 about pivots,
$\prec_W$ and $\psi_E$ there are dual results about $\prec_E$ and~$\psi_W$.

\chapter{\huge\bf P$''$-Graphs and P$'$-Graphs}\label{4} \pagestyle{myheadings}\markboth{CHAPTER 4. \quad
P$\,''$-GRAPHS AND P$\,'$-GRAPHS}{right-head}

\section{\large\bf Petals}\label{4.1} \markright{\S 4.1. \quad Petals}

In this chapter the goal is to prove that every P$''$-graph (as defined in \S 1.9) is a P$'$-graph (as defined in \S 1.8). For that we must first deal with some preliminary matters in this and in the next two sections.

For a vertex $v$ of a D-graph $D$ let ${\cal C}(v)$, the
\emph{corolla}\index{corolla of a vertex} of $v$, be the set of all
vertices $x$ of $D$ such that $v\prec x$. For example, in
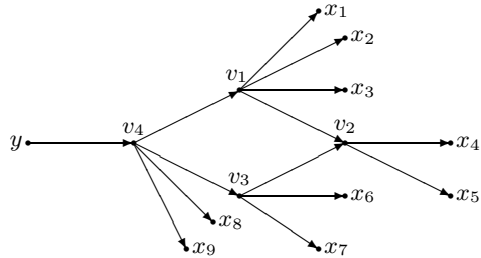
\begin{center}
\begin{picture}(160,100)(0,-5)

\put(0,40){\circle*{2}} \put(40,40){\circle*{2}}
\put(60,0){\circle*{2}} \put(70,10){\circle*{2}}
\put(80,20){\circle*{2}} \put(80,60){\circle*{2}}
\put(110,0){\circle*{2}} \put(120,20){\circle*{2}}
\put(120,40){\circle*{2}} \put(120,60){\circle*{2}}
\put(120,80){\circle*{2}} \put(110,90){\circle*{2}}
\put(160,20){\circle*{2}} \put(160,40){\circle*{2}}

\put(0,40){\vector(1,0){40}} \put(40,40){\vector(1,-2){20}}
\put(40,40){\vector(1,-1){30}} \put(40,40){\vector(2,-1){40}}
\put(40,40){\vector(2,1){40}} \put(80,20){\vector(3,-2){30}}
\put(80,20){\vector(1,0){40}} \put(80,20){\vector(2,1){40}}
\put(80,60){\vector(2,-1){40}} \put(80,60){\vector(1,0){40}}
\put(80,60){\vector(2,1){40}} \put(80,60){\vector(1,1){30}}
\put(120,40){\vector(2,-1){40}} \put(120,40){\vector(1,0){40}}

\put(-2,40){\small\makebox(0,0)[r]{$y$}}
\put(40,44){\small\makebox(0,0)[b]{$v_4$}}
\put(79,64){\small\makebox(0,0)[b]{$v_1$}}
\put(80,25){\small\makebox(0,0)[b]{$v_3$}}
\put(120,45){\small\makebox(0,0)[b]{$v_2$}}
\put(112,90){\small\makebox(0,0)[l]{$x_1$}}
\put(122,80){\small\makebox(0,0)[l]{$x_2$}}
\put(122,60){\small\makebox(0,0)[l]{$x_3$}}
\put(162,40){\small\makebox(0,0)[l]{$x_4$}}
\put(162,20){\small\makebox(0,0)[l]{$x_5$}}
\put(122,20){\small\makebox(0,0)[l]{$x_6$}}
\put(112,0){\small\makebox(0,0)[l]{$x_7$}}
\put(72,10){\small\makebox(0,0)[l]{$x_8$}}
\put(62,0){\small\makebox(0,0)[l]{$x_9$}}

\end{picture}
\end{center}
${\cal C}(v_1)=\{x_1,x_2,x_3\}$, while ${\cal C}(v_4)$
is made of all the vertices except $v_4$ and~$y$.

The binary relation that holds between the elements $x$ and $y$ of
${\cal C}(v)$ whenever \emph{not} $v\rhd[x,y]$ is an equivalence
relation on ${\cal C}(v)$. For reflexivity, we have that
\emph{not} $v\rhd[x,x]$ because $v\neq x$ (which is assumed with
$v\prec x$), and $v$ does not belong to the one-vertex trivial semipath
from $x$ to $x$. Symmetry is trivial, because we can always read a
semipath in reverse order. For transitivity, assume we have
the semipaths $\sigma$ in $[x,y]$ and $\tau$ in $[y,z]$ such that not
$v\rhd\sigma$ and not $v\rhd\tau$. Then we have that
$\sigma\ast\tau\in[x,z]$ and not $v\rhd\sigma\ast\tau$.

For $x$ in ${\cal C}(v)$, let the equivalence class
\begin{tabbing}
\hspace{1.6em}$|\![x]\!|_v=_{df}\{y\in{\cal C}(v)\mid {\rm\emph{not }}
v\rhd[x,y]\}$
\end{tabbing}
be called a \emph{petal}\index{petal}. In the example above
we have $|\![x_1]\!|_{v_1}=\{x_1\}$, while $|\![x_1]\!|_{v_4}$ is
made of all the vertices except $x_8$, $x_9$, $v_4$ and $y$.

\prop{Lemma 4.1.1}{If $v\prec x$ and not $v\rhd[x,u]$, then
$v\prec u$.}

\dkz Suppose $v\prec x$ and not $v\prec u$. Hence for some
$\sigma\in[u]_W$ we have not $v\rhd\sigma$. Take a $\tau$ in $[x,u]$.
Then for $\tau\ast\sigma\in[x]_W$ we have $v\rhd\tau\ast\sigma$.
Since not $v\rhd\sigma$, we must have $v\rhd\tau$. \qed

As a corollary we have the following.

\prop{Lemma 4.1.2}{If $v\prec x$ and not $v\rhd[x,u]$, then
$u\in|\![x]\!|_v$.}

\vspace{-2ex}

\prop{Lemma 4.1.3}{If $v\prec x$, $v\rhd[x,y]$, $\tau\in[y]_W$
and $u\rhd\tau$, then $v\rhd[x,u]$.}

\dkz Suppose $v\prec x$, $v\rhd[x,y]$, $\tau\in[y,z]$, $z\in
W(D)$, and $u\rhd\tau$. Suppose not $v\rhd[x,u]$. Hence for some
$\sigma$ in $[x,u]$ we have not $v\rhd\sigma$. Then
$\sigma\ast\tau_{[u,z]}\in[x]_W$, and since $v\prec x$, we have
$v\rhd\sigma\ast\tau_{[u,z]}$. Since not $v\rhd\sigma$, we have
$v\rhd\tau_{[u,z]}$. Hence not $v\rhd\tau_{[u,y]}$; otherwise,
$\tau$ would not be a semipath. So not
$v\rhd\sigma\ast\tau_{[u,y]}$, which together with
$\sigma\ast\tau_{[u,y]}\in[x,y]$ contradicts our assumption that
$v\rhd[x,y]$. \qed

\vspace{-2ex}

\prop{Lemma 4.1.4}{Let $x$, $x'$ and $y$ be distinct vertices of
$E(D)$. If $v\prec x$, $v\prec x'$, $v\rhd[x,y]$ and not
$v\rhd[x,x']$, then not $\psi_E(x,y,x')$.}

\dkz Suppose $v\prec x$, $v\prec x'$, $v\rhd[x,y]$,
$\sigma\in[x,x']$ and not $v\rhd\sigma$. Suppose for
$\tau$ in $[y]_W$ we have a vertex $u$ such that $u\rhd\sigma$ and
$u\rhd\tau$. By Lemma 4.1.3, we have $v\rhd[x,u]$, and hence
$v\rhd\sigma_{[x,u]}$, which contradicts not $v\rhd\sigma$. Hence
not $\psi_E(x,y,x')$. \qed

As a corollary we have the following.

\prop{Lemma 4.1.5}{Let $x_1$, $x_2$ and $y$ be distinct vertices
of $E(D)$. If $v\prec x$, $y\notin |\![x]\!|_v$ and $x_1,x_2\in
|\![x]\!|_v$, then not $\psi_E(x_1,y,x_2)$. }

\dkz From $v\prec x$, $y\notin |\![x]\!|_v$ and $x_1,x_2\in
|\![x]\!|_v$, we conclude $v\prec x_1$, $v\prec x_2$,
$v\rhd[x_1,y]$ and not $v\rhd[x_1,x_2]$, with the help of Lemma
4.1.2; then we apply Lemma 4.1.4. \qed

\vspace{-2ex}

\section{\large\bf P-moves}\label{4.2}
\markright{\S 4.2. \quad P-moves}

For any list $A$, let $\bar{A}$ be the converse list, i.e.\ $A$
read in reverse order. For the lists $A$ and $A'$, let $P_{A,A'}$
be the set of all ordered pairs $(a,b)$ such that $a$ precedes $b$
in $A$ and $b$ precedes $a$ in $A'$.

Let $|\![x]\!|^E_v$ be $|\![x]\!|_v\cap E(D)$, and let
$\Lambda_v(x)$ stand for a list of $|\![x]\!|^E_v$. Since the
equivalence class $|\![x]\!|_v$ is non-empty, we conclude that
$|\![x]\!|^E_v$ and hence also $\Lambda_v(x)$ are always non-empty.

For $x$ and $y$ distinct $E$-vertices, let $V(x,y)$ be $\{v\mid
v\prec x\quad\&\quad v\prec y\}$, i.e.\ the set of common pivots
of $x$ and $y$. We say that $v$ is the \emph{closest common
pivot}\index{closest common pivot} of $x$ and
$y$, and write $v{\rm CCP}(x,y)$,\index{CCP} when $v\in V(x,y)$ and for every
$w$ in $V(x,y)$ either $w\prec v$ or $w=v$.

Let $\Pi$ and $\Theta$ be two lists of $E(D)$. Then we call
P$_{\Pi,\Theta}$-\emph{moves}\index{P-move}, or sometimes
P-\emph{moves} for short, the following rewrite rules for lists of
$E(D)$; we read these rules as stating that we can pass from the
list of $E(D)$ above the horizontal line, which we call $\Pi$, to the
list of $E(D)$ below, which we call $\Pi'$, provided $(x,y)\in
P_{\Pi,\Theta}$:
\[
\begin{array}{l}
\lpravilo{Tr-$(x,y)$}\f{\Gamma\Lambda_v(x)\Lambda_v(y)\Delta}{\Gamma\Lambda_v(y)\Lambda_v(x)\Delta}
\\[1ex]
\lpravilo{Sf-$(x,y)$}\f{\Gamma\Lambda_v(z)\Delta}{\Gamma\overline{\Lambda_v(z)}\Delta}
\end{array}
\]
provided that in Sf-$(x,y)$ we have that $x$ precedes $y$ in $\Lambda_v(z)$ and $v{\rm CCP}
(x,y)$,
\[
\lpravilo{Bf-$(x,y)$}\f{\;\Pi\;}{\overline{\Pi}}\mbox{\hspace{5em}}
\]
provided that in Bf-$(x,y)$ we have that $x$ precedes $y$ in $\Pi$ and $V(x,y)=\emptyset$. In the
names of these rules, Tr\index{Tr} stands for \emph{transposition}, Sf\index{Sf}
stands for \emph{small flip}, and Bf\index{Bf} stands for \emph{big flip}.

Let $D$ be the D-graph form the beginning of \S 4.1. As an example, consider the
following lists of~$E(D)$:
\[
\begin{array}{l}
\Pi\!: x_1 x_2 x_3 x_4 x_5 x_6 x_7 x_8 x_9
\\[1ex]
\Theta\!: x_8 x_6 x_7 x_5 x_4 x_1 x_2 x_3 x_9
\end{array}
\]
The P-move Tr-$(x_5,x_8)$ is
\[
\f{x_1 x_2 x_3 x_4 x_5 x_6 x_7 x_8 x_9}{x_8 x_1 x_2 x_3 x_4 x_5
x_6 x_7 x_9}
\]
with $\Gamma$ being empty, $\Lambda_{v_4}(x_5)\!:x_1 x_2 x_3 x_4
x_5 x_6 x_7$, $\Lambda_{v_4}(x_8)\!: x_8$ and $\Delta\!:x_9$.

\noindent The P-move Sf-$(x_5,x_6)$ is
\[
\f{x_1 x_2 x_3 x_4 x_5 x_6 x_7 x_8 x_9}{x_7 x_6 x_5 x_4 x_3 x_2
x_1 x_8 x_9}
\]
with $\Gamma$ being empty, $\Lambda_{v_4}(x_5)$ as above and
$\Delta$ being $x_8 x_9$.

Note that for Tr-$(x,y)$ we can infer that $v{\rm CCP} (x,y)$, as
in the proviso for Sf-$(x,y)$. Otherwise, for some $w$ we would
have $v\prec w$, $w\prec x$ and $w\prec y$. From that by Lemma
3.2.8.4 we would obtain that $w\in |\![x]\!|_v$ and $w\in
|\![y]\!|_v$, and hence $|\![x]\!|_v=|\![y]\!|_v$, which is
contradictory to our assumptions (lists are without repetitions,
and hence $\Lambda_v(x)$ and $\Lambda_v(y)$ are lists with
different members).

\prop{Lemma 4.2.1}{For every P$_{\Pi,\Theta}$-move, if $\Pi$ is
grounded in $D$, then $\Pi'$ is grounded in $D$.}

\dkz Suppose the vertices $r$,$s$ and $t$ occur in that order in
$\Pi$, and we have $\psi_E(r,s,t)$. Suppose first our
P$_{\Pi,\Theta}$-move is Tr-$(x,y)$. All the cases where not more
than one of $r$, $s$ and $t$ occur in $\Lambda_v(x)\Lambda_v(y)$
are settled in the obvious manner, as well as the two cases where
they are all in $\Lambda_v(x)$ or all in $\Lambda_v(y)$, and the
cases where two of $r$, $s$ and $t$ are in one of $\Lambda_v(x)$
and $\Lambda_v(y)$, while the remaining vertex is in $\Gamma$ or
$\Delta$. In all these cases, $r$, $s$ and $t$ occur in the same
order in $\Pi$ and~$\Pi'$.

Note also that the case where $r$ and $t$ are in one of
$\Lambda_v(x)$ and $\Lambda_v(y)$, while $s$ is in the other, is
impossible, because we have $\Pi\!: r\mn s\mn t$. As interesting
cases, only the following remain.

\vspace{.5ex}

(1) Suppose $r$ and $s$ are in $\Lambda_v(x)$, and $t$ is in
$\Lambda_v(y)$. We need to show that $\psi_E(t,r,s)$. Take a
$\sigma$ in $[t,s]$. We have $v\prec s$ and $t\notin|\![s]\!|_v$. By
Lemma 4.1.2, we obtain $v\rhd\sigma$. Since $v\prec r$, we obtain
$\psi_E(t,r,s)$. We reason analogously when $r$ is in
$\Lambda_v(x)$, while $s$ and $t$ are in~$\Lambda_v(y)$.

\vspace{.5ex}

(2) Suppose $r$ is in $\Gamma$, $s$ is in $\Lambda_v(x)$ and $t$
is in $\Lambda_v(y)$. We need to show that $\psi_E(r,t,s)$. Take a
$\sigma$ in $[r,s]$. Then we have $v\prec s$ and
$r\notin|\![s]\!|_v$. By Lemma 4.1.2, we obtain $v\rhd \sigma$.
Since $v\prec t$, we obtain $\psi_E(r,t,s)$. We reason analogously
when $r$ is in $\Lambda_v(x)$, $s$ is in $\Lambda_v(y)$ and $t$ is
in~$\Delta$.

\vspace{.5ex}

Suppose next our P$_{\Pi,\Theta}$-move is Sf-$(x,y)$. Excluding
obvious cases, like those mentioned above, we have as interesting
cases only the following. Suppose $r$ is in $\Gamma$, while $s$
and $t$ are in $\Lambda_v(z)$. Then we reason as in (2), by using
Lemma 4.1.2. We reason analogously when $r$ and $s$ are in
$\Lambda_v(z)$ and $t$ is in $\Delta$. The case where our
P$_{\Pi,\Theta}$-move is Bf-$(x,y)$ is settled in the obvious way.
\qed

Note that this Lemma holds without taking account of the provisos
for Sf-$(x,y)$ and Bf-$(x,y)$. As a matter of fact, our proof of
the lemma shows that for every P-move, respecting the provisos or
not, $\Pi$ is grounded in $D$ iff $\Pi'$ is grounded in $D$.

The following proposition shows that making a Tr-$(x,y)$ move
brings us closer to $\Theta$, ``closer'' in a sense that will be
made precise later (see Proposition 4.3.2 below). For this
proposition we assume that $\Pi$ is
$\Gamma\Lambda_v(x)\Lambda_v(y)\Delta$, as for Tr-$(x,y)$.

\prop{Proposition 4.2.2}{Suppose $\Pi$ and $\Theta$ are grounded
in $D$, and $r,r'\in|\![x]\!|^E_v$ and $s,s'\in|\![y]\!|^E_v$.
Then $(r,s)\in P_{\Pi,\Theta}$ implies $(r',s')\in
P_{\Pi,\Theta}$.}

\dkz Suppose $(r,s)\in P_{\Pi,\Theta}$ and $(r',s')\notin
P_{\Pi,\Theta}$. Hence $r$ precedes $s$ in $\Pi$ and $s$ precedes
$r$ in $\Theta$. Suppose $s'$ precedes $r'$ in both $\Pi$ and~$\Theta$.

\vspace{.5ex}

(1) Suppose $r'$ precedes $r$ in $\Pi$. Then we have $\Pi\!:s'\mn
r'\mn r\mn s$, and by using Lemma 4.1.5 with $x_1$, $x_2$ and $y$
replaced by $s$, $s'$ and $r$ or $r'$ we obtain that not
$\psi_E(s,r,s')$ or not $\psi_E(s,r',s')$, which contradicts our
assumption.

\vspace{.5ex}

(2) Suppose $r$ precedes $r'$ in~$\Pi$.

\vspace{.5ex}

(2.1) Suppose $s$ precedes $s'$ in $\Pi$. Then we have $\Pi\!:r\mn
s\mn s'\mn r'$, and by using Lemma 4.1.5 with $x_1$, $x_2$ and $y$
replaced by $r$, $r'$ and $s$ or $s'$ we obtain that not
$\psi_E(r,s,r')$ or not $\psi_E(r,s',r')$, which contradicts our
assumption.

\vspace{.5ex}

(2.2) Suppose $s'$ precedes $s$ in~$\Pi$.

\vspace{.5ex}

(2.21) Suppose $r$ precedes $s'$ in $\Pi$. Then we have
$\Pi\!:r\mn s'\mn r'$, and by using Lemma 4.1.5 as in (2.1) we
obtain a contradiction.

\vspace{.5ex}

(2.22) Suppose $s'$ precedes $r$ in $\Pi$. Then we have
$\Pi\!:s'\mn r\mn s$, and by using Lemma 4.1.5 as in (1) we obtain
a contradiction.

When $r'$ precedes $s'$ in both $\Pi$ and $\Theta$ we proceed in
the same manner after replacing $\Pi$, $r$, $r'$, $s$, and $s'$ by
respectively $\Theta$, $r'$, $r$, $s'$ and $s$. We obtain in that
case contradictions with our assumption that $\Theta$ is grounded
in~$D$. \qed

We can prove the following.

\prop{Lemma 4.2.3}{If $u_1, u_2, u_3\in |\![z]\!|^E_v$,
$\psi_E(u_1,u_2,u_3)$ and $\psi_E(u_2,u_1,u_3)$, then not $v{\rm
CCP} (u_1,u_2)$.}

\dkz Since $\psi_E(u_1,u_2,u_3)$ and $\psi_E(u_2,u_1,u_3)$, we have, by
Theorem 3.4.1, a vertex $w$ such that $w\rhd[u_3,u_1]$ and
$w\in V(u_1,u_2)$. Since $u_1,u_3\in |\![z]\!|^E_v$, there is a
$\sigma$ in $[u_3,u_1]$ such that not $v\rhd\sigma$. Since
$w\rhd[u_3,u_1]$, we have $w\rhd\sigma$, and hence not
$v\rhd\sigma_{[u_1,w]}$. So not $v\rhd[u_1,w]$, and since $v\prec
u_1$, we obtain $v\prec w$ by Lemma 4.1.1. \qed

For the following two lemmata we assume that $v{\rm CCP} (x,y)$,
and that $\Pi$ and $\Theta$ are grounded in $D$. First we have a
lemma that is a direct corollary of Lemma 4.2.3.

\prop{Lemma 4.2.4}{If $x,y,u\in|\![z]\!|^E_v$, then it is
impossible that $\Pi\!: x\mn y\mn u$ and $\Theta\!: y\mn x\mn u$.}

Next we have the following.

\prop{Lemma 4.2.5}{If $x,y,x',y'\in |\![z]\!|^E_v$, then it is
impossible that}

\vspace{-2ex}

\begin{tabbing}
\hspace{1.6em}\=(1)\hspace{2em}\=$\Pi\!:y'\mn x\mn y\mn x'$
\=\emph{and} \= $\Theta\!:y\mn y'\mn x'\mn x$,
\\[.5ex]
\>(2)\>$\Pi\!:x\mn y'\mn x'\mn y$ \>\emph{and} \> $\Theta\!:y'\mn
y\mn x\mn x'$.
\end{tabbing}

\dkz Suppose we have (1). Then since $\Pi\!:y'\mn x\mn y$ and
$\Theta\!:x\mn y'\mn y$, we have, by Lemma 4.2.3, a vertex $w_1$ in
$V(x,y')$ such that $v\prec w_1$. And since $\Pi\!:y'\mn y\mn x'$
and $\Theta\!:y\mn y'\mn x'$, we have, by Lemma 4.2.3, a vertex
$w_2$ in $V(y',y)$ such that $v\prec w_2$.

It is impossible that $w_1=w_2$, because otherwise not $v{\rm CCP}
(x,y)$. Then since $w_1\prec y'$ and $w_2\prec y'$,
we have, by Lemma 3.2.1 that either $w_1\prec w_2$ or $w_2\prec w_1$. If $w_1\prec
w_2$, then $w_1\in V(x,y)$ and since $v\prec w_1$, we obtain a
contradiction with $v{\rm CCP} (x,y)$, , and we reason analogously if $w_2\prec w_1$. This proves that (1) is
impossible.

An alternative proof that (1) is impossible is obtained by showing
that we have a vertex $w_1$ in $V(y,x')$ such that $v\prec w_1$, and
a vertex $w_2$ in $V(x',x)$ such that $v\prec w_2$. We prove (2)
analogously with two applications of Lemma 4.2.3. \qed

\vspace{-2ex}

\prop{Lemma 4.2.6}{If $V(x,y)\neq\emptyset$, then there is a $w$
in $V(x,y)$ such that $w{\rm CCP} (x,y)$.}

\dkz It follows from Lemmata 3.2.1 and 3.2.4 that $V(x,y)$ is
linearly ordered by $\prec$. Then $w$ is the greatest element of
this linear order. \qed

Besides the assumptions that $v{\rm CCP} (x,y)$, and that $\Pi$
and $\Theta$ are grounded in $D$, which we made before Lemma
4.2.4, we assume also that $(x,y)\in P_{\Pi,\Theta}$, that $x,y\in
|\![z]\!|^E_v$, that $\Pi$ is $\Gamma\Lambda_v(z)\Delta$, and that $\Pi'$ is $\Gamma\overline{\Lambda_v(z)}\Delta$, as for
Sf-$(x,y)$. Then we have the following.

\prop{Proposition 4.2.7}{If $(x',y')\in
P_{\Pi',\Theta}-P_{\Pi,\Theta}$, then not $v{\rm CCP} (x',y')$.}

\dkz Let $P_\Pi$ and $P_{\Pi'}$ abbreviate $P_{\Pi,\Theta}$ and
$P_{\Pi',\Theta}$ respectively. Suppose $(x',y')\in
P_{\Pi'}-P_\Pi$; i.e., $(x',y')\notin P_\Pi$ and $(x',y')\in
P_{\Pi'}$. From that we derive that $y'$ precedes $x'$ in both
$\Theta$ and $\Pi$, and that $x',y'\in |\![z]\!|^E_v$. It is
excluded that $x=y$ because $(x,y)\in P_\Pi$, and that $x'=y'$
because $(x',y')\in P_{\Pi'}$. It is excluded that $x'=x$ and
$y'=y$ because $(x,y)\in P_\Pi$ and $(x',y')\notin P_\Pi$, and
that $x'=y$ and $y'=x$ because $(x',y')\in P_{\Pi'}$ and
$(y,x)\notin P_{\Pi'}$.

\vspace{.5ex}

(1) Suppose $x'=x$ and $y'\neq y$. Then, since $y'$ precedes $x'$
in $\Pi$, we have $\Pi\!: y'\mn x\mn y$.

It is excluded that $\Theta\!: y'\mn y\mn x$, because of Lemma
4.2.4, and that $\Theta\!: y\mn x\mn y'$, because $y$ precedes $x$
and $y'$ precedes $x'$ in $\Theta$; hence $\Theta\!: y\mn y'\mn
x$. From $\Pi\!: y'\mn x\mn y$ and $\Theta\!: y\mn y'\mn x$ we
conclude that $\psi_E(y',x',y)$ and $\psi_E(x',y',y)$, and by,
Lemma 4.2.3, not $v{\rm CCP} (x',y')$. The case where $x'\neq x$
and $y'=y$ is treated analogously.

\vspace{.5ex}

(2) Suppose $x'=y$ and $y'\neq x$. Then, since $y'$ precedes $x'$
in $\Theta$, we have $\Theta\!: y'\mn y\mn x$.

It is excluded that $\Pi\!: y'\mn x\mn y$, because of Lemma 4.2.4,
and that $\Pi\!: x\mn y\mn y'$, because $x$ precedes $y$ and $y'$
precedes $x'$ in $\Pi$. Hence $\Pi\!: x\mn y'\mn y$, and by reasoning
as in (1), we obtain, by Lemma 4.2.3, not $v{\rm CCP} (x',y')$. The
case where $x'\neq y$ and $y'=x$ is treated analogously. The only
remaining case is the following.

\vspace{.5ex}

(3) Suppose $x$, $y$, $x'$ and $y'$ are all mutually distinct.

\vspace{.5ex}

(3.1) Suppose $\Pi\!: y'\mn x\mn y$. It is excluded that
$\Theta\!: y'\mn y\mn x$, because of Lemma 4.2.4.

\vspace{.5ex}

(3.11) Suppose $\Theta\!: y\mn y'\mn x$. It is excluded that
$\Pi\!: x'\mn y'\mn x\mn y$, because $x$ precedes $y$ and $y'$
precedes $x'$ in~$\Pi$.

\vspace{.5ex}

(3.111) Suppose $\Pi\!: y'\mn x'\mn x\mn y$ or $\Pi\!: y'\mn x\mn
x'\mn y$. It is excluded that $\Theta\!: x'\mn y\mn y'\mn x$ and
$\Theta\!: y\mn x'\mn y'\mn x$, because $y$ precedes $x$ and $y'$
precedes $x'$ in $\Theta$. If $\Theta\!: y\mn y'\mn x'\mn x$ or
$\Theta\!: y\mn y'\mn x\mn x'$, then since $\Pi\!: y\mn x'\mn y'$
and $\Theta\!: y\mn y'\mn x'$, we obtain, by Lemma 4.2.3, not
$v{\rm CCP} (x',y')$.

\vspace{.5ex}

(3.112) Suppose $\Pi\!: y'\mn x\mn y\mn x'$. It is excluded that
$\Theta\!: x'\mn y\mn y'\mn x$ and $\Theta\!: y\mn x'\mn y'\mn x$,
for the reasons given in (3.111). It is excluded also that
$\Theta\!: y\mn y'\mn x'\mn x$, because of (1) of Lemma 4.2.5, and
that $\Theta\!: y\mn y'\mn x\mn x'$, because of Lemma 4.2.4. So
(3.112) is impossible.

\vspace{.5ex}

(3.12) Suppose $\Theta\!: y\mn x\mn y'$. Hence $\Theta\!: y\mn
x\mn y'\mn x'$, because $y$ precedes $x$ and $y'$ precedes $x'$ in
$\Theta$. From now on we will take for granted this sort of
justification based on precedence. We may have that either $\Pi\!:
y'\mn x'\mn x\mn y$ or $\Pi\!: y'\mn x\mn x\mn y$, as in (3.111),
and we reason as for (3.111). If we suppose $\Pi\!: y'\mn x\mn
y\mn x'$ (as in (3.112)), then, since $\Pi\!: x\mn y\mn x'$ and
$\Theta\!: y\mn x\mn x'$, we obtain a contradiction with Lemma 4.2.4.

\vspace{.5ex}

(3.2) Suppose $\Pi\!: x\mn y'\mn y$.

\vspace{.5ex}

(3.21) Suppose $\Theta\!: y'\mn y\mn x$.

\vspace{.5ex}

(3.211) Suppose $\Pi\!: x\mn y'\mn x'\mn y$. Then with either
$\Theta\!: y'\mn x'\mn y\mn x$ or $\Theta\!: y'\mn y\mn x'\mn x$
we apply Lemma 4.2.3 to obtain not $v{\rm CCP} (x',y')$. It is
excluded that $\Theta\!: y'\mn y\mn x\mn x'$, because of (2) of
Lemma 4.2.5.

\vspace{.5ex}

(3.212) Suppose $\Pi\!: x\mn y'\mn y\mn x'$. Then with either
$\Theta\!: y'\mn x'\mn y\mn x$ or $\Theta\!: y'\mn y\mn x'\mn x$
we apply Lemma 4.2.3 to obtain not $v{\rm CCP} (x',y')$. It is
excluded that $\Theta\!: y'\mn y\mn x\mn x'$, because of Lemma
4.2.4.

\vspace{.5ex}

(3.22) Suppose $\Theta\!: y\mn y'\mn x$.

\vspace{.5ex}

(3.221) Suppose $\Pi\!: x\mn y'\mn x'\mn y$. Then with either
$\Theta\!: y\mn y'\mn x'\mn x$ or $\Theta\!: y\mn y'\mn x\mn x'$
we apply Lemma 4.2.3 to obtain not $v{\rm CCP} (x',y')$.

\vspace{.5ex}

(3.222) Suppose $\Pi\!: x\mn y'\mn y\mn x'$. Then with $\Theta\!:
y\mn y'\mn x'\mn x$ we apply Lemma 4.2.3 to obtain not $v{\rm CCP}
(x',y')$. It is excluded that $\Theta\!: y\mn y'\mn x\mn x'$,
because of Lemma 4.2.4.

\vspace{.5ex}

(3.23) Suppose $\Theta\!: y\mn x\mn y'$. Hence $\Theta\!: y\mn
x\mn y'\mn x'$. Then with $\Pi\!: x\mn y'\mn x'\mn y$ we apply
Lemma 4.2.3 to obtain not $v{\rm CCP} (x',y')$. It is excluded
that $\Pi\!: x\mn y'\mn y\mn x'$, because of Lemma 4.2.4.

\vspace{.5ex}

(3.3) Suppose $\Pi\!: x\mn y\mn y'$. Hence $\Pi\!: x\mn y\mn y'\mn
x'$.

\vspace{.5ex}

(3.31) Suppose $\Theta\!: y'\mn y\mn x$. With either $\Theta\!:
y'\mn x'\mn y\mn x$ or $\Theta\!: y'\mn y\mn x'\mn x$, we obtain, by Lemma
4.2.3, not $v{\rm CCP} (x',y')$. It is excluded that
$\Theta\!: y'\mn y\mn x\mn x'$, because of Lemma 4.2.4.

\vspace{.5ex}

(3.32) Suppose $\Theta\!: y\mn y'\mn x$. With $\Theta\!: y\mn
y'\mn x'\mn x$, we obtain, by Lemma 4.2.3, not $v{\rm CCP} (x',y')$.
It is excluded that $\Theta\!: y\mn y'\mn x\mn x'$, because of
Lemma 4.2.4.

\vspace{.5ex}

(3.33) Suppose $\Theta\!: y\mn x\mn y'$. Hence $\Theta\!: y\mn
x\mn y'\mn x'$, which is excluded because of Lemma 4.2.4. \qed

Under the assumptions that $x$ precedes $y$ in $\Pi$, that
$V(x,y)=\emptyset$, and that $\Pi'$ is $\overline{\Pi}$, as for
Bf-$(x,y)$, we can prove the following version of Proposition
4.2.7.

\prop{Proposition 4.2.8}{If $(x',y')\in
P_{\Pi',\Theta}-P_{\Pi,\Theta}$, then $V(x',y')\neq\emptyset$.}

\dkz Let the $W$-edges of $D$ be $a_1,\ldots,a_n$ for $n\geq 2$.
(If $n$ where 1, then $V(x,y)$ would not be $\emptyset$.) Let $B$
be the basic D-graph whose $E$-edges are $a_1,\ldots,a_n$, whose
inner vertex is $v$, and whose $W$-edges are arbitrarily chosen so
that $B\Box D$ is defined. Then a P-move Bf-$(x,y)$ of $D$ becomes
Sf-$(x,y)$ of $B\Box D$. Proposition 4.2.7 for that P-move
Sf-$(x,y)$ of $B\Box D$ yields Proposition 4.2.8 for~$D$. \qed

\vspace{-2ex}

\section{\large\bf Completeness of P-moves}\label{4.3}
\markright{\S 4.3. \quad Completeness of P-moves}

For a vertex $v$ of $D$, let ${\cal C}_E(v)={\cal C}(v)\cap E(D)$,
where ${\cal C}(v)$ is the corolla of $v$ defined at the beginning
of \S 4.1. In the D-graph from the beginning of \S 4.1 we have, for example, ${\cal
C}_E(v_1)={\cal C}(v_1)=\{x_1,x_2,x_3\}$ and ${\cal
C}_E(v_4)=\{x_1,\ldots,x_9\}$. Let $U_E$ be the set of all the
vertices $v$ of $D$ such that ${\cal C}_E(v)$ has at least two
elements.

In the set-theoretic sense, a \emph{tree}\index{tree, set-theoretic} is a partially ordered
set in which for every element the set of predecessors is
well-ordered. (Such a tree need not have a single root.) It
follows from Lemmata 3.2.1 and 3.2.4 that $\langle
U_E,\prec\rangle$ is a finite tree. (This finite tree corresponds
to a forest in the graph-theoretic sense; see \cite{H69},
Chapter~4.)

Let $\langle U^+_E,\prec\rangle$ be the single-rooted tree
obtained by adding to $U_E$ a new element $\ast$ and by assuming
that for every $v$ in $U_E$ we have $\ast\prec v$. The new element
$\ast$ has the same function as the inner vertex of $B$ in the
proof of Proposition 4.2.8. The tree $\langle U^+_E,\prec\rangle$
is interesting when  $\langle U_E,\prec\rangle$ is not
single-rooted. In our example above, the tree $\langle
U_E,\prec\rangle$ is pictured by
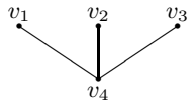
\begin{center}
\begin{picture}(60,40)

\put(0,30){\circle*{2}} \put(30,10){\circle*{2}}
\put(30,30){\circle*{2}} \put(60,30){\circle*{2}}

\put(30,10){\line(-3,2){30}} \put(30,10){\line(0,1){20}}
\put(30,10){\line(3,2){30}}

\put(30,7){\small\makebox(0,0)[t]{$v_4$}}
\put(0,33){\small\makebox(0,0)[b]{$v_1$}}
\put(30,33){\small\makebox(0,0)[b]{$v_2$}}
\put(60,33){\small\makebox(0,0)[b]{$v_3$}}

\end{picture}
\end{center}
In the D-graph
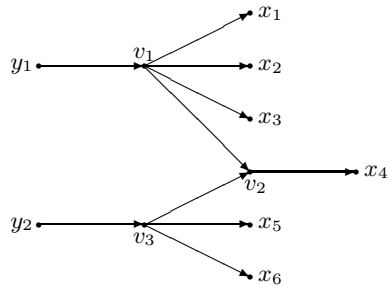
\begin{center}
\begin{picture}(120,100)(0,-5)

\put(0,20){\circle*{2}} \put(0,80){\circle*{2}}
\put(40,20){\circle*{2}} \put(40,80){\circle*{2}}
\put(80,0){\circle*{2}} \put(80,20){\circle*{2}}
\put(80,40){\circle*{2}} \put(80,60){\circle*{2}}
\put(80,80){\circle*{2}} \put(80,100){\circle*{2}}
\put(120,40){\circle*{2}}

\put(0,20){\vector(1,0){40}} \put(0,80){\vector(1,0){40}}
\put(40,20){\vector(2,-1){40}} \put(40,20){\vector(1,0){40}}
\put(40,20){\vector(2,1){40}} \put(40,80){\vector(1,-1){40}}
\put(40,80){\vector(2,-1){40}} \put(40,80){\vector(1,0){40}}
\put(40,80){\vector(2,1){40}} \put(80,40){\vector(1,0){40}}

\put(-2,20){\small\makebox(0,0)[r]{$y_2$}}
\put(-2,80){\small\makebox(0,0)[r]{$y_1$}}
\put(40,17){\small\makebox(0,0)[t]{$v_3$}}
\put(40,83){\small\makebox(0,0)[b]{$v_1$}}
\put(82,37){\small\makebox(0,0)[t]{$v_2$}}
\put(83,0){\small\makebox(0,0)[l]{$x_6$}}
\put(83,20){\small\makebox(0,0)[l]{$x_5$}}
\put(123,40){\small\makebox(0,0)[l]{$x_4$}}
\put(83,60){\small\makebox(0,0)[l]{$x_3$}}
\put(83,80){\small\makebox(0,0)[l]{$x_2$}}
\put(83,100){\small\makebox(0,0)[l]{$x_1$}}

\end{picture}
\end{center}
the tree $\langle U_E,\prec\rangle$ would be made just
of two roots $v_1$ and $v_3$, and $\langle U^+_E,\prec\rangle$
would be
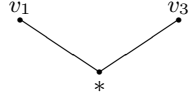
\begin{center}
\begin{picture}(60,40)

\put(0,30){\circle*{2}} \put(30,10){\circle*{2}}
\put(60,30){\circle*{2}}

\put(30,10){\line(-3,2){30}} \put(30,10){\line(3,2){30}}

\put(30,7){\small\makebox(0,0)[t]{$\ast$}}
\put(0,33){\small\makebox(0,0)[b]{$v_1$}}
\put(60,33){\small\makebox(0,0)[b]{$v_3$}}

\end{picture}
\end{center}

For every element $v$ of $U^+_E$ let $S(v)=\{w\in U^+_E\mid v\prec
w\}$; i.e., $S(v)$\index{Sv@$S(v)$} is the set of successors of $v$ in the tree,
not necessarily immediate successors; $v$ is a \emph{leaf}\index{leaf} when $S(v)$ is
empty. In the tree $\langle U^+_E,\prec\rangle$ obtained from the
first $\langle U_E,\prec\rangle$ tree above
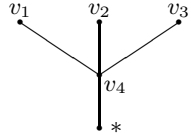
\begin{center}
\begin{picture}(60,50)(0,5)

\put(0,50){\circle*{2}} \put(30,30){\circle*{2}}
\put(30,50){\circle*{2}} \put(60,50){\circle*{2}}
\put(30,10){\circle*{2}}

\put(30,30){\line(-3,2){30}} \put(30,10){\line(0,1){20}}
\put(30,30){\line(0,1){20}} \put(30,30){\line(3,2){30}}

\put(36,28){\small\makebox(0,0)[t]{$v_4$}}
\put(0,53){\small\makebox(0,0)[b]{$v_1$}}
\put(30,53){\small\makebox(0,0)[b]{$v_2$}}
\put(60,53){\small\makebox(0,0)[b]{$v_3$}}
\put(34,10){\small\makebox(0,0)[l]{$\ast$}}

\end{picture}
\end{center}
$S(\ast)=\{v_1,v_2,v_3,v_4\}$, $S(v_4)=\{v_1,v_2,v_3\}$,
and $v_1$, $v_2$ and $v_3$ are leaves.

Let $k_v$ be the cardinality of ${\cal C}_E(v)$. We assign
inductively to every element $v$ of $U^+_E$ a natural number
$m(v)\geq 1$. For a leaf $v$, let $m(v)$ be 1, and, for other
elements $v$ of $U^+_E$, let $m(v)$ be $(\sum_{w\in
S(v)}{{k_w}\choose 2}\cdot m(w))+1$. Since $\sum_{w\in
S(v)}{{k_w}\choose 2}\cdot m(w)$ is 0 when $S(v)=\emptyset$, we
compute $m(v)$ for a leaf $v$ in the same way. The number
${k_w}\choose 2$ is the number of pairs of distinct $E$-vertices
in~${\cal C}(w)$.

For the last example for $\langle U^+_E,\prec\rangle$ above we
have
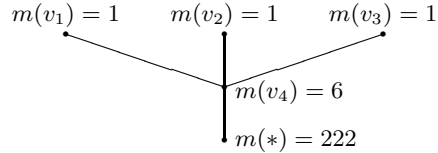
\begin{center}
\begin{picture}(120,60)(0,5)

\put(0,50){\circle*{2}} \put(60,30){\circle*{2}}
\put(60,50){\circle*{2}} \put(120,50){\circle*{2}}
\put(60,10){\circle*{2}}

\put(60,30){\line(-3,1){60}} \put(60,10){\line(0,1){20}}
\put(60,30){\line(0,1){20}} \put(60,30){\line(3,1){60}}

\put(64,28){\small\makebox(0,0)[l]{$m(v_4)=6$}}
\put(0,55){\small\makebox(0,0)[b]{$m(v_1)=1$}}
\put(60,55){\small\makebox(0,0)[b]{$m(v_2)=1$}}
\put(120,55){\small\makebox(0,0)[b]{$m(v_3)=1$}}
\put(64,10){\small\makebox(0,0)[l]{$m(\ast)=222$}}

\end{picture}
\end{center}
\begin{tabbing}
\hspace{1.6em}\=$m(v_4)$ \=$={{k_{v_1}}\choose 2} \cdot m(v_1)+
{{k_{v_2}}\choose 2} \cdot m(v_2)+ {{k_{v_3}}\choose 2} \cdot
m(v_3)+1$,
\\[1ex]
\>\>$={3\choose 2}\cdot 1+{2\choose 2}\cdot 1+{2\choose 2}\cdot
1+1=5+1=6$,
\\[1.5ex]
\>$m(\ast)$\>$={{k_{v_1}}\choose 2}\cdot m(v_1)+\ldots +
{{k_{v_4}}\choose 2}\cdot m(v_4) +1$,
\\[1ex]
\>\>$= 5+{9\choose 2}\cdot 6+1=222$.
\end{tabbing}

For $(x,y)$ a pair of distinct elements of $E(D)$, let
\[
M(x,y)=\left\{
\begin{array}{ll}
\min\{m(v)\mid v\in V(x,y)\}, & {\rm when }\; V(x,y)\neq\emptyset
\\[1ex]
m(\ast), & {\rm when }\; V(x,y)=\emptyset.
\end{array}
\right.
\]
An alternative way to define $M(x,y)$ when $V(x,y)\neq\emptyset$
is to say that it is $m(v)$ for $v$ such that $v{\rm CCP} (x,y)$;
it yields the same number as the definition above.

Finally, we define $\mu_\Theta(\Pi)$ as $\sum_{(x,y)\in
P_{\Pi,\Theta}}M(x,y)$. For the D-graph from the beginning of \S 4.1 and the lists
$\Pi$ and $\Theta$ given as an example in \S 4.2 after introducing
the P-moves, we have
\begin{tabbing}
\hspace{1.6em}$P_{\Pi,\Theta}=\{$\=$(x_1,x_8),\ldots,(x_7,x_8),$
\\[.5ex]
\>$(x_1,x_6),\ldots,(x_5,x_6),$
\\[.5ex]
\>$(x_1,x_7),\ldots,(x_5,x_7),$
\\[.5ex]
\>$(x_1,x_5),\ldots,(x_4,x_5),$
\\[.5ex]
\>$(x_1,x_4),\ldots,(x_3,x_4)\}$.
\end{tabbing}
With $\Pi'$ of Tr-$(x_5,x_8)$, which we call $\Pi'_{\scriptsize\rm
Tr}$, we have
\begin{tabbing}
\hspace{1.6em}\=$P_{\Pi'_{\tiny\rm
Tr},\Theta}=P_{\Pi,\Theta}-\{(x_1,x_8),\ldots,(x_7,x_8)\}$.
\end{tabbing}
With $\Pi'$ of Sf-$(x_5,x_6)$, which we call $\Pi'_{\scriptsize\rm
Sf}$, we have
\begin{tabbing}
\hspace{1.6em}\=$P_{\Pi'_{\tiny\rm
Sf},\Theta}=\{(x_1,x_8),\ldots,(x_7,x_8)\}\cup\{(x_7,x_6),(x_3,x_2),(x_3,x_1),(x_2,x_1)\}$.
\end{tabbing}
We then have
\begin{tabbing}
\hspace{1.6em}\=$\mu_\Theta(\Pi_{\scriptsize\rm Tr}\;$ \= \kill

\>$\mu_\Theta(\Pi)$ \> $=7\cdot 6+5\cdot 6+5\cdot 6+(3\cdot
6+1)+3\cdot 6=139$,
\\[.5ex]
\>$\mu_\Theta(\Pi'_{\scriptsize\rm Tr})$ \> $=139-42=97$,
\\[.5ex]
\>$\mu_\Theta(\Pi'_{\scriptsize\rm Sf})$ \> $=42+4=46$.
\end{tabbing}

Note that if in our example we replace
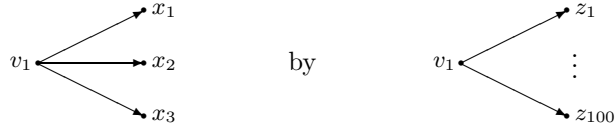
\begin{center}
\begin{picture}(200,50)(0,-5)

\put(0,20){\circle*{2}} \put(40,0){\circle*{2}}
\put(40,20){\circle*{2}} \put(40,40){\circle*{2}}
\put(160,20){\circle*{2}} \put(200,0){\circle*{2}}
\put(200,40){\circle*{2}}

\put(0,20){\vector(1,0){40}} \put(0,20){\vector(2,-1){40}}
\put(0,20){\vector(2,1){40}} \put(160,20){\vector(2,-1){40}}
\put(160,20){\vector(2,1){40}}

\put(100,20){\makebox(0,0){by}}

\put(-2,20){\small\makebox(0,0)[r]{$v_1$}}
\put(43,0){\small\makebox(0,0)[l]{$x_3$}}
\put(43,20){\small\makebox(0,0)[l]{$x_2$}}
\put(43,40){\small\makebox(0,0)[l]{$x_1$}}

\put(203,22){\makebox(0,0){$\vdots$}}

\put(158,20){\small\makebox(0,0)[r]{$v_1$}}
\put(203,0){\small\makebox(0,0)[l]{$z_{100}$}}
\put(203,40){\small\makebox(0,0)[l]{$z_1$}}

\end{picture}
\end{center}
and in $\Pi$ and $\Theta$ we replace $x_1x_2x_3$ by
$z_1\ldots z_{100}$, then $P_{\Pi,\Theta}$ has 509 elements, while
$P_{\Pi'_{\tiny\rm Sf},\Theta}$ has 5055 elements, but
$\mu_\Theta(\Pi)=2\, 516\, 124$, while
$\mu_\Theta(\Pi'_{\scriptsize\rm Sf})=520\, 063$.

We make the same assumptions as for Proposition 4.2.7, and we
prove the following.

\prop{Proposition 4.3.1}{We have $\sum_{(x',y')\in
P_{\Pi',\Theta}-P_{\Pi,\Theta}} M(x',y')<M(x,y)$.}

\dkz We have that $v$, for which we have $v{\rm CCP} (x,y)$, is in
$U^+_E$, because $x$ and $y$ are distinct members of a list. If
$P_{\Pi',\Theta}-P_{\Pi,\Theta}$ is empty, then the sum on the
left is 0; this is lesser than $M(x,y)$, which is $m(v)$, and is
at least 1.

If $P_{\Pi',\Theta}-P_{\Pi,\Theta}$ is non-empty, then, by using
Proposition 4.2.7 and Lemma 4.2.6, we conclude that for every pair
$(x',y')$ in it there is a $w$ such that $w{\rm CCP} (x',y')$ and
$w\in U^+_E$, because $x'$ and $y'$ are distinct members of a
list, and $w\in S(v)$. We have
\begin{tabbing}
\hspace{1.6em}$\sum_{(x',y')\in
P_{\Pi',\Theta}-P_{\Pi,\Theta}}M(x',y')$
\\[1ex]
\hspace{5em}\=$=\sum_{w\in S(v)}(\sum_{(x',y')\in
P_{\Pi',\Theta}-P_{\Pi,\Theta}\;\&\;w{\scriptsize\rm CCP}(x',y')}
M(x',y'))$
\\[1ex]
\>$\leq\sum_{w\in S(v)}{{k_w}\choose 2}\cdot m(w)$, \=since the
number of pairs $(x',y')$ in
\\
\>\>$P_{\Pi',\Theta}-P_{\Pi,\Theta}$ such that $w{\rm CCP}
(x',y')$
\\
\>\>is lesser than or equal to ${k_w}\choose 2$,
\\[1ex]
\>$<(\sum_{w\in S(v)}{{k_w}\choose 2}\cdot m(w))+1$
\\[1ex]
\>$=M(x,y)$. \`$\dashv$
\end{tabbing}

We also have Proposition 4.3.1 under the same assumptions as for
Proposition 4.2.8. The proof is very much analogous, with $\ast$
functioning as $v$, and Proposition 4.2.8 replacing Proposition
4.2.7. If we have a P-move Tr-$(x,y)$, then Proposition 4.2.2
guarantees that $P_{\Pi',\Theta}$ is a proper subset
of~$P_{\Pi,\Theta}$.

Then Proposition 4.2.2 and Proposition 4.3.1 in both versions, the
Sf and Bf versions, yield the following proposition where $\Pi$
and $\Pi'$ are from any P-move.

\prop{Proposition 4.3.2}{We have
$\mu_\Theta(\Pi')<\mu_\Theta(\Pi)$.}

Then we can prove the following proposition, which says that
P-moves are complete, in the sense that they enable us to pass
from a grounded list to any other grounded list.

\prop{Proposition 4.3.3}{If $\Pi$ and $\Theta$ are grounded in
$D$, then they are either the same or there is a finite sequence
of P-moves $P_1,\ldots,P_n$, with $n\geq 1$, such that $\Pi$ is
the upper list of $P_1$, while $\Theta$ is the lower list of
$P_n$, and for every $P_i$ such that $1\leq i<n$ we have that the
lower list of $P_i$ is the upper list of~$P_{i+1}$.}

\dkz If $\mu_\Theta(\Pi)=0$, then $P_{\Pi,\Theta}=\emptyset$, and
$\Pi$ and $\Theta$ are the same. If $\mu_\Theta(\Pi)>0$, then
$P_{\Pi,\Theta}\neq\emptyset$, and $\Pi$ and $\Theta$ are
distinct.

The \emph{distance}\index{distance between members of a list} $d(a,b)$ between the distinct members $a$ and
$b$ of a list $A$ is the number of members of $A$ between $a$ and
$b$. Among all the pairs in $P_{\Pi,\Theta}$ take a pair $(x,y)$
with a minimal distance $d(x,y)$. We have the following
possibilities.

Suppose $V(x,y)\neq\emptyset$ and let $v{\rm CCP}(x,y)$. This $v$
exists by Lemma 4.2.6. We have two subcases. Suppose
$|\![x]\!|_v\neq |\![y]\!|_v$. It follows from Lemma 4.1.5 that
$\Pi$ must be of the form
$\Gamma\Lambda_v(x)\Xi\Lambda_v(y)\Delta$. That $\Xi$ must be the
empty list follows from our assumption about the minimality of
$d(x,y)$. (If there were a member $z$ in $\Xi$, then by the
minimality of $d(x,y)$ we would have that $(x,z)$ and $(z,y)$ are
not in $P_{\Pi,\Theta}$, from which it would follow that $(x,y)$
is not in $P_{\Pi,\Theta}$.) Then we may apply Tr-$(x,y)$.

If $|\![x]\!|_v = |\![y]\!|_v$, then we appeal again to Lemma
4.1.5, and we apply Sf-$(x,y)$. If $V(x,y)=\emptyset$, then we
apply Bf-$(x,y)$. Our proof is formally an induction on
$\mu_\Theta(\Pi)$, which relies on Proposition 4.3.2. \qed

\vspace{-2ex}

\section{\large\bf P$''$-graphs are P$'$-graphs}\label{4.4}
\markright{\S 4.4. \quad P$\,''$-graphs are P$\,'$-graphs}

In this section we will prove that every P$''$-graph is a
P$'$-graph. For that we need some more preliminary results. (Here we use
the notation $V_W$, $V_E$ and $V_C$ introduced in \S 1.3.)

\begin{samepage}

\prop{Lemma 4.4.1}{Suppose $v$ and $x$ are vertices of $D_W\Box
D_E$ and $v\prec x$. If $v$ is from $D_W$, then either $(1)$
$V_E\cap |\![x]\!|_v=\emptyset$ or $(2)$ $V_E-V_C\subseteq
|\![x]\!|_v$.}

\end{samepage}

\dkz Suppose $y\in V_E\cap|\![x]\!|_v$, and let $z\in V_E-V_C$.
Then for some semipath $\sigma$ of $D_W\Box D_E$ in $[y,z]$ we
have not $v\rhd\sigma$, because $D_E$ is weakly connected. Hence
$z\in|\![x]\!|_v$. \qed

\vspace{-2ex}

\prop{Lemma 4.4.2}{If there is a construction $K$ of a $P'$-graph with the root
list $L_E$, which is $\Gamma\Lambda_v(x)\Lambda_v(y)\Delta$, then
there is a construction $K'$ of the same $P'$-graph with the root list $L_E'$ being
$\Gamma\Lambda_v(y)\Lambda_v(x)\Delta$, while the root lists $L_W$
and $L_W'$ of $K$ and $K'$ respectively are the same.}

\dkz We proceed by induction on the number of nodes in $K$. In the
basis, when $K$ has a single node, then in this node, which is the
root of $K$, we have $(B,L_W,L_E)$, and $v$ is the inner vertex of
the basic D-graph $B$. All petals with respect to $v$ are
singletons, and we pass from $L_E$ to $L_E'$ by one transposition
of $x$ and $y$. The construction $K'$ has $(B,L_W,L_E')$ in its
root.

For the induction step, we have that $K$ is $K_W\Box K_E$, and the
root graphs of $K_W$ and $K_E$ are $D_W$ and $D_E$ respectively.

\vspace{.5ex}

(1) If the vertex $v$ of the petals $|\![x]\!|_v$ and
$|\![y]\!|_v$ is in $D_E$, then we just apply the induction
hypothesis to $K_E$ to obtain $K_E'$, and $K'$ is $K_W\Box K_E'$.
(Note that $|\![x]\!|_v$ in $D_W\Box D_E$ and $|\![x]\!|_v$ in
$D_E$ are here the same sets of vertices.)

If $v$ is in $D_W$, then, according to Lemma 4.4.1, we have three
possibilities.

\vspace{.5ex}

(2) Suppose $V_E\cap |\![x]\!|_v=V_E\cap |\![y]\!|_v=\emptyset$.
Then we apply the induction hypothesis to $K_W$ to obtain $K_W'$,
and $K'$ is $K_W'\Box K_E$.

\vspace{.5ex}

(3) Suppose $V_E-V_C\subseteq |\![x]\!|_v$. Then $V_E\cap
|\![y]\!|_v=\emptyset$. Replace in $\Lambda_v(x)$ the subset
$L^E_E$ by the list $\Xi$, which is the common sublist of $L^W_E$
and $L^E_W$, made of the vertices of $V_C$, and let the result of
this replacement be $\Lambda_v(x_1)\ldots\Lambda_v(x_k)$, where
$k\geq 1$, with $\Lambda_v(x_i)$ for $i\in\{1,\ldots,k\}$ being a
list of an $|\![x_i]\!|^E_v$ for $|\![x_i]\!|_v$ a petal of $D_W$.
The induction hypothesis allows us to make $k$ applications of
moves like P-moves of the Tr-$(x_i,y)$ type to obtain $K_W'$, and
$K'$ will be $K_W'\Box K_E$. Note that $\Xi$ coincides with
$L^E_W$, which follows from $V_E-V_C\subseteq |\![x]\!|_v$.

This is important to ascertain that $L_W'$ is $L_W$ and that
$(L^W_E)'$ and $L^E_W$ are compatible. The case when
$V_E-V_C\subseteq |\![y]\!|_v$ is treated analogously. \qed

\vspace{-2ex}

\prop{Lemma 4.4.3}{If there is a construction $K$ of a $P'$-graph with the root
lists $L_W$ and $L_E$, then there is a construction $K'$ of the same $P'$-graph with the
root lists $L_W'$ and $L_E'$ being respectively $\overline{L_W}$
and $\overline{L_E}$.}

\dkz In every leaf of $K$ replace $(B,L_W,L_E)$ by
$(B,\overline{L_W},\overline{L_E})$. Formally, we have again an
induction on the number of nodes in $K$, with the induction step
trivial. \qed

\vspace{-2ex}

\prop{Lemma 4.4.4}{If there is a construction $K$ of a $P'$-graph with the root
list $L_E$, which is $\Gamma\Lambda_v(z)\Delta$, then there is a
construction $K'$ of the same $P'$-graph with the root list $L_E'$ being
$\Gamma\overline{\Lambda_v(z)}\Delta$, while the root lists $L_W$
and $L_W'$ of $K$ and $K'$ respectively are the same.}

\dkz We proceed by induction on the number of nodes in $K$. In the
basis, when $K$ has a single node, we take that $K$ is just $K'$.
This is because all petals are singletons, as in the proof of
Lemma 4.4.2.

For the induction step, if $v$ is in $D_E$, we proceed as in (1)
of the proof of Lemma 4.4.2. Suppose $v$ is in $D_W$. If we have
$V_E\cap |\![x]\!|_v=\emptyset$, then we proceed as for (2) of the
proof of Lemma 4.4.2.

Suppose $V_E-V_C\subseteq |\![x]\!|_v$. Replace in
$\Lambda_v(x)$ the sublist $L^E_E$ by the list $\Xi$, as in (3) of
the proof of Lemma 4.4.2, and let the result be
$\Lambda_v(x_1)\ldots\Lambda_v(x_k)$, where $k\geq 1$. By the
induction hypothesis, we obtain a construction with
$\overline{\Lambda_v(x_1)}\ldots\overline{\Lambda_v(x_k)}$, and
then by Lemma 4.4.2 we have a construction $K_W'$ with
$\overline{\Lambda_v(x_k)}\ldots\overline{\Lambda_v(x_1)}$ in its
root.

By Lemma 4.4.3, we have a construction $K_E'$ obtained by
replacing the lists $L^E_W$ and $L^E_E$ of $K_E$ by
$\overline{L^E_W}$ and $\overline{L^E_E}$ respectively. Note that
$\Xi$ coincides with $L^E_W$, which follows from
$V_E-V_C\subseteq |\![x]\!|_v$, and is important for the reasons
mentioned in the proof of Lemma 4.4.2. Then the construction $K'$
will be $K_W'\Box K_E'$. \qed

It is clear that for all the results based on $\prec$, which is
$\prec_W$, as at the end of \S 3.4, and obtained starting from \S
4.1 up to now, we have analogous results based on $\prec_E$, with
completely analogous proofs. Then we can prove the following.

\prop{Theorem 4.4.5}{For every P$\,''$-graph $D$ there is a
construction of a P$\,'$-graph with root graph $D$, which means
that $D$ is a P$\,'$-graph.}

\dkz We proceed by induction on the number $k$ of inner vertices
of $D$. When $k$ is 1, the theorem is trivial.

For the induction step, suppose $D$ is $D_W\Box D_E$. So we have a
list $\Theta_W$ of $E(D_W)$ grounded in $D_W$ and a list
$\Theta_E$ of $W(D_E)$ grounded in $D_E$, which are compatible. By
the induction hypothesis we have the constructions $K_W$ and $K_E$
with root graphs $D_W$ and $D_E$ respectively.

Let the root lists $L_W$ and $L_E$ of $K_X$ be respectively
$\Pi^X_W$ and $\Pi^X_E$. By Proposition 4.3.3 and Lemmata 4.4.2,
4.4.3 and 4.4.4, there is a finite sequence of constructions of the
P$'$-graph $D_W$ starting with $K_W$, for which $L_E$ is $\Pi^W_E$, and
ending with $K_W'$, for which the root list $L_E'$ is $\Theta_W$.

By analogous results, in an analogous manner, we obtain out of $K_E$ a construction $K_E'$ of the P$'$-graph $D_E$, for which the root list $L_W'$ is $\Theta_E$. Since
$\Theta_W$ and $\Theta_E$ are compatible, we have that $K_W'\Box
K_E'$ is a construction with root graph~$D$. \qed

\clearpage \pagestyle{empty} \makebox[1em]{} \clearpage

\chapter{\huge\bf P$'''$-Graphs and P$''$-Graphs}\label{5} \pagestyle{myheadings}\markboth{CHAPTER 5. \quad
P$\,'''$-GRAPHS AND P$\,''$-GRAPHS}{right-head}

\section{\large\bf B$_m$-moves}\label{5.1} \markright{\S 5.1. \quad B$_m$-moves}

In this chapter we will finish establishing that the three definitions of P-graph are equivalent by proving that every P$'''$-graph (as defined in \S 1.10) is a P$''$-graph (as defined in \S 1.9). For that we must first deal with some preliminary matters in this and in the next section. The present section is based on matters introduced in~\S 4.2.

Let $\Pi$, which is $\Theta\Psi\Theta'$, be a list of $E(D)$. Let
$T$ be the set of members of $\Theta$, while $T'$ is the set of
members of $\Theta'$. Consider a set $F\subseteq T\cup T'$, and
let $F'$ be the relative complement of $F$ with respect to $T\cup
T'$. Let $B$ be the set $(T\cap F')\cup(F\cap T')$, which amounts
to the symmetric difference of the sets $T$ and $F$.

For $m$ a member of $\Psi$, we call B$_m$-moves the following
rewrite rules from $\Pi$ to $\Pi'$, provided $x\in B$:
\[
\begin{array}{l}
\lpravilo{Tr-$(x,m)$}\f{\Gamma\Lambda_v(x)\Phi\Lambda_v(m)\Delta}{\Gamma\Phi\Lambda_v(m)\Lambda_v(x)\Delta}
\\[1ex]
\lpravilo{Tr-$(m,x)$}\f{\Gamma\Lambda_v(m)\Phi\Lambda_v(x)\Delta}{\Gamma\Lambda_v(x)\Lambda_v(m)\Phi\Delta}
\\[1ex]
\lpravilo{Sf-$(x,m)$}\f{\Gamma\Lambda_v(m)\Delta}{\Gamma\overline{\Lambda_v(m)}\Delta}
\end{array}
\]
{\samepage provided that in Sf-$(x,m)$ we have that $x$ is a member of $\Lambda_v(m)$ and $v{\rm CCP} (x,m)$,}
\[
\lpravilo{Bf-$(x,m)$}\f{\;\Pi\;}{\overline{\Pi}}\mbox{\hspace{6em}}
\]
provided that in Bf-$(x,m)$ we have that $V(x,m)=\emptyset$.
Note that, as in \S 4.2, we can infer
for Tr-$(x,m)$ and Tr-$(m,x)$ that $v{\rm CCP} (x,m)$.

From now on
we assume that $D$ is $D_W\Box D_E$, and let $\Psi$ of
$\Theta\Psi\Theta'$ be a list of $E(D_E)$. Hence $x$ is an
$E$-vertex of $D_W$. We may infer that $v$ of Tr-$(x,m)$, Tr-$(m,x)$ and Sf-$(x,m)$ is a vertex of $D_W$;
otherwise, if $v$ were a vertex of $D_E$, then, since $D_W$ is
weakly connected, for some $\sigma$ in $[x]_W$ we would have not
$v\rhd\sigma$, and so we would not have $v\prec x$, which is
presupposed by $\Lambda_v(x)$. We may also infer that $\Psi$ is a
sublist of $\Lambda_v(m)$, by Lemma 4.4.1 for $x$ being $m$. We
may further infer that in Tr-$(x,m)$ we have that
$\Gamma\Lambda_v(x)\Phi$ is a sublist of $\Theta$ and $\Delta$
a sublist of $\Theta'$; that in Tr-$(m,x)$ we have that $\Gamma$
is a sublist of $\Theta$ and $\Phi\Lambda_v(x)\Delta$ a sublist
of $\Theta'$; and that in Sf-$(x,m)$ we have that $\Gamma$ is a
sublist of $\Theta$ and $\Delta$ a sublist of $\Theta'$. With a
proof analogous to the proof of Lemma 4.2.1, we establish that for
every B$_m$-move, if $\Pi$ is grounded in $D$, then $\Pi'$ is
grounded in~$D$.

Let $V_C$ be as in \S 1.3 for $D_W\Box D_E$, and let $k\in V_C$.
Then we have the following lemma, which will help us to prove Proposition 5.1.4, a proposition that will play a similar role
to Proposition 4.2.2.

\prop{Lemma 5.1.1}{For $\Pi$ being the upper list of Tr-$(x,m)$ or
Tr-$(m,x)$, let $x'$ be a member of $\Lambda_v(x)$. Then not
$\psi_E(x,k,x')$ in~$D_W$.}

\dkz It is easy to infer that $v\prec x$ in $D_W$. Next we show
that $k\notin |\![x]\!|^E_v$ in $D_W$. Otherwise, not
$v\rhd[x,k]$ in $D_W$, which would yield not $v\rhd[x,m]$ in
$D_W\Box D_E$, and this contradicts $m\notin |\![x]\!|^E_v$ in
$D_W\Box D_E$. We show also that $x'\in |\![x]\!|^E_v$ in $D_W$.
Otherwise, $x$ and $x'$ would be connected by a semipath with
vertices of $D_E$ in which $v$ does not occur, and this would
again contradict $m\notin |\![x]\!|^E_v$ in $D_W\Box D_E$. Hence
by Lemma 4.1.5 we have not $\psi_E(x,k,x')$ in~$D_W$. \qed

\vspace{-2ex}

\prop{Lemma 5.1.2}{For $\Pi$ being the upper list of Sf-$(x,m)$,
let $x'$ be, as $x$, an element of $T\cup T'$ and a member of
$\Lambda_v(m)$, and assume that for every $k$ in $V_C$ we have
$\psi_E(x,k,x')$ and $\psi_E(k,x,x')$ in $D_W$. Then not $v{\rm
CCP} (x,m)$ in $D_W\Box D_E$.}

\dkz By Theorem 3.4.1 and Lemma 3.3.1 we obtain a vertex $w$ of
$D_W$ such that for every $k$ in~$V_C$
\[
w\rhd[x',x]\quad\&\quad w\rhd[x',k]\quad\&\quad w\prec
x\quad\&\quad w\prec k
\]
in $D_W$. Since for every $k$ in $V_C$ we have $w\prec k$ in
$D_W$, we infer $w\in V(x,m)$ in $D_W\Box D_E$. Since $v\in
V(x,m)$ in $D_W\Box D_E$, we must have by Lemma 3.2.1 either
$v=w$, or $w\prec v$, or $v\prec w$ in $D_W\Box D_E$. We show next
that $v=w$ or $w\prec v$ implies a contradiction.

Since there is a semipath $\sigma$ of $D_W\Box D_E$ in $[x',m]$ such that
not $v\rhd\sigma$, because $x'$ is a member of $\Lambda_v(m)$,
there is a $k'$ in $V_C$ and a semipath $\sigma'$ of $D_W$ in
$[x',k']$ such that not $v\rhd\sigma'$. Since for every $k$ in
$V_C$ we have $w\rhd[x',k]$, we have that $w\rhd\sigma'$.

Since not $v\rhd\sigma'_{[x',w]}$, we have not $v\rhd[x',w]$, and
since $v\prec x'$ in $D_W$, because $v\prec x'$ in $D_W\Box D_E$,
which follows from $x'$ being a member of $\Lambda_v(m)$, we conclude, by Lemma
4.1.1, that $v\prec w$ in $D_W$. This contradicts $v=w$
immediately, and it contradicts also $w\prec v$ in $D_W\Box D_E$,
which implies $w\prec v$ in $D_W$; we rely on Lemma 3.2.2. Hence
we must have $v\prec w$ in $D_W\Box D_E$, which implies not $v{\rm
CCP} (x,m)$. \qed

The following lemma is the analogue of Lemma 5.1.2 for Bf-$(x,m)$.

\prop{Lemma 5.1.3}{For $\Pi$ being the upper list of Bf-$(x,m)$,
let $x'$ be, as $x$, a member of $\Pi$, and assume that for every
$k$ in $V_C$ we have $\psi_E(x,k,x')$ and $\psi_E(k,x,x')$ in
$D_W$. Then $V(x,m)\neq\emptyset$ in $D_W\Box D_E$.}

\noindent The proof is as the proof of Lemma 5.1.2 until we reach
the conclusion that $w\in V(x,m)$ in $D_W\Box D_E$.

Suppose we have the D-graphs $D_1$, $D_2$ and $D_3$ such that
$D_1\Box D_2$ and $D_1\Box D_3$ are defined, but neither $D_2\Box
D_3$ nor $D_3\Box D_2$ is defined. This is the situation analogous
to what we had with (Ass~2.1) in \S 1.5. Alternatively, it is
equivalent to suppose that $(D_1\Box D_2)\Box D_3$ and $(D_1\Box
D_3)\Box D_2$ are defined.

Let $\Pi$, which is a list of $E(D_1\Box D_2)$, be of the form
$\Theta\Psi\Theta'$ for $\Psi$ a list of $E(D_2)$; here $D_1$ and
$D_2$ correspond respectively to what was above $D_W$ and $D_E$.
As before, the sets of vertices $T$ and $T'$ of $E(D_1)$ are
respectively the sets of members of $\Theta$ and~$\Theta'$.

Let $\Sigma$ be a list of $E(D_1\Box D_3)$, which is of the form
$\Omega\Xi\Omega'$ for $\Xi$ a list of $V_C$, which is $E(D_1)\cap
W(D_2)$, and for all the members of the list $\Omega$ being
elements of $E(D_1)$. If the members of $\Omega$ are not all in
$E(D_1)$, but some are elements of $E(D_3)$, then all the members
of the list $\Omega'$ are elements of $E(D_1)$, and all that we
do up to the end of \S 5.2 would be done in a dual manner,
involving $\Theta'$ and $\Omega'$ instead of $\Theta$ and
$\Omega$. Let $F$ be the set of members of $\Omega$, and let $T'$,
$F'$ and $B$ be defined with respect to $T$ and $F$ as at the
beginning of this section.

Our purpose now is to show that B$_m$-moves are complete, in the
sense that they enable us to pass from any list $\Pi$ of
$E(D_1\Box D_2)$ grounded in $D_1\Box D_2$ to a list
$\Theta\Psi\Theta'$ of $E(D_1\Box D_2)$ grounded in $D_1\Box D_2$
such that $\Theta$ is a list of $F$ and $\Psi$ is a list of
$E(D_2)$. For that we assume that $\Sigma$ is grounded in $D_1\Box
D_3$.

For the propositions that follow we assume that $\Pi$ is grounded
in $D_1\Box D_2$ and that $\Sigma$ is grounded in $D_1\Box D_3$.
First, we have a proposition analogous to Proposition 4.2.2, which
says that with B$_m$-moves of the kind Tr-$(x,m)$ or Tr-$(m,x)$
the set $B$ diminishes, in a sense which will be made precise
later (see Proposition 5.2.2).

\prop{Proposition 5.1.4}{For $\Pi$ being the upper list of
Tr-$(x,m)$ or Tr-$(m,x)$, let $x'$ be a member of $\Lambda_v(x)$.
Then $x'\in B$.}

\dkz By our assumption for B-moves, we have that $x\in B$. Suppose
$\Pi$ is the upper list of Tr-$(x,m)$, and suppose $x'\notin B$.
Suppose $x\in T\cap F'$. We have that $x$ is a member of $\Theta$,
i.e., $x$ precedes $m$ in $\Pi$. Hence $x'$ precedes $m$ in $\Pi$ by Lemma 4.1.5 and by the groundedness of $\Pi$. For every $k$ in $V_C$ we have $\Sigma\!:
x'\mn k\mn x$, since $x\in B$ and $x'\notin B$. Since $\Sigma$ is
grounded in $D_1\Box D_3$, we have $\psi(x',k,x)$ in $D_1\Box
D_3$. By Lemma 3.1.6 we have $\psi_E(x',k,x)$ in $D_1$, which
contradicts Lemma 5.1.1.

If $x\in F\cap T'$, then we proceed analogously, and obtain again a
contradiction with Lemma 5.1.1. Hence $x'\in B$. We proceed
analogously when $\Pi$ is the upper list of Tr-$(m,x)$. \qed

Next we have a proposition related to Proposition 4.2.7.

\prop{Proposition 5.1.5}{For $\Pi$ being the upper list of
Sf-$(x,m)$, let $x'$ be, as $x$, an element of $T\cup T'$ and a
member of $\Lambda_v(m)$, and suppose $x'\notin B$. If
\\[.5ex] \indent $(1)$ $\Pi\!: x'\mn x\mn m$, or
\\[.5ex] \indent $(2)$ $\Pi\!: x\mn m\mn x'$,
\\ then not $v{\rm CCP}(x,m)$ in $D_1\Box D_2$, and if
\\[.5ex] \indent $(3)$ $\Pi\!: x\mn x'\mn m$,
\\ then not $v{\rm CCP}(x',m)$ in $D_1\Box D_2$.
}

\dkz We prove first the implication from (1) or
(2) to not $v{\rm CCP}(x,m)$ in $D_1\Box D_2$.

\vspace{.5ex}

(I) Suppose $x\notin E(D_1)\cap W(D_3)$, and suppose we have (1).
Then for every $k$ in $V_C$ we have $\Sigma\!: x'\mn k\mn x$,
since $x\in B$ and $x'\notin B$. Since $\Sigma$ is grounded in
$D_1\Box D_3$, we have $\psi_E(x',k,x)$ in $D_1\Box D_3$. By Lemma
3.1.6 we obtain $\psi_E(x',k,x)$ in $D_1$. From (1), and the
groundedness of $\Pi$ in $D_1\Box D_2$, we infer $\psi_E(x',x,k)$
in $D_1$ by Lemma 3.1.4. Then by Lemma 5.1.2 we obtain that not
$v{\rm CCP}(x,m)$ in $D_1\Box D_2$. If we have (2), then we
proceed analogously, with Lemma 3.1.4 replaced by Lemma 3.1.5.

\vspace{.5ex}

(II) Suppose $x'\in E(D_1)\cap W(D_3)$. Suppose we have (1). Then
for every $k$ in $V_C$ and for some $w$ in $E(D_3)$ we have
$\Sigma\!: w\mn k\mn x$, since $x\in B$ and $x'\notin B$. Since
$\Sigma$ is grounded in $D_1\Box D_3$, we have $\psi_E(w,k,x)$ in
$D_1\Box D_3$. By Lemma 3.1.4 we obtain $\psi_E(x',k,x)$ in $D_1$.
After that we proceed as in (I) to show that not $v{\rm CCP}(x,m)$
in $D_1\Box D_2$. If we have (2), then we proceed analogously,
with Lemma 3.1.4 replaced by Lemma 3.1.5.

To prove the implication from (3) to not $v{\rm CCP}(x',m)$ in
$D_1\Box D_2$, we proceed analogously to what we had with (1)
above. Instead of Lemma 5.1.2, we now apply the lemma obtained
from Lemma 5.1.2 by interchanging $x$ and $x'$. \qed

Finally, we have a proposition related to Proposition 4.2.8.

\prop{Proposition 5.1.6}{For $\Pi$ being the upper list of
Bf-$(x,m)$, let $x'$ be, as $x$, a member of $\Pi$, and suppose
$x'\notin B$. If $(1)$ or $(2)$ of Proposition 5.1.5, then
$V(x,m)\neq\emptyset$ in $D_1\Box D_2$, and if $(3)$ of Proposition
5.1.5, then $V(x',m)\neq\emptyset$ in $D_1\Box D_2$.}

\noindent The proof is as for Proposition 5.1.5 by relying on
Lemma 5.1.3 instead of Lemma 5.1.2.

\section{\large\bf Completeness of B$_m$-moves}\label{5.2} \markright{\S 5.2. \quad Completeness of B$_m$-moves}

Let $D$ be a D-graph $D_1\Box D_2$. With the tree $\langle
U^+_E,\prec\rangle$ of this D-graph, we define $m(v)$ and $M(x,y)$
as in \S 4.3. As a matter of fact, we could now modify the
definition of $m(v)$ by replacing $k_w\choose 2$ in it by $k_w$,
or by $k_w$ diminished by the number of vertices in $E(D_2)$. We
write $M(x)$ as an abbreviation for $M(x,m)$, where $m$ is the
vertex involved in our B$_m$-moves. We define $\mu_B$, which is
analogous to $\mu_\Theta(\Pi)$, as $\sum_{x\in B}M(x)$.

With $D$ being $D_1\Box D_2$, for every B$_m$-move, $\Pi'$ is the
lower list, which, as $\Pi$, may be conceived as being of the form
$\Theta\Psi\Theta'$ for $\Psi$ a list of $E(D_2)$. This is
because, as we remarked after introducing the B$_m$-moves, the
sublist $\Psi$ of $\Pi$ is a sublist of $\Lambda_v(m)$. Let $B'$
be defined for this $\Pi'$, as $B$ was defined for $\Pi$, namely
as the symmetric difference of $T$ and $F$, with $T$ being the set
of members of the sublist $\Theta$ of $\Pi'$, and $F$ being as for
$B$ the set of members of $\Omega$ (see the assumptions concerning
$\Sigma$ before Proposition 5.1.4).

We can prove the following proposition analogous to Proposition
4.3.1. Assume for that proposition that $\Pi$ and $\Pi'$ are as
for an Sf-$(x,m)$ move.

\prop{Proposition 5.2.1}{We have $\sum_{x'\in B'-B}M(x')<M(x)$.}

\dkz As a corollary of Proposition 5.1.5 we may ascertain that if
$x'\in B'-B$, then not $v{\rm CCP}(x',m)$, a proposition analogous
to Proposition 4.2.7. Note first that $x'\in B'-B$ implies that
$x'$ is an element of $T\cup T'$ and a member of $\Lambda_v(m)$.
We also have, of course, $x'\notin B$. Since we have a Sf-$(x,m)$ move,
we have $v{\rm CCP}(x,m)$, and hence (1) and (2) of Proposition
5.1.5 are impossible. The only remaining possibility is (3), which
yields not $v{\rm CCP}(x',m)$.

We may then continue reasoning as in the proof of Proposition
4.2.7. (Now, $y'$ is either omitted or replaced by $m$.) \qed

By relying on Proposition 5.1.4 and Proposition 5.1.6, and by
reasoning in a manner analogous to what we had before Proposition
4.3.2, we obtain the following for any B$_m$-move.

\prop{Proposition 5.2.2}{We have $\mu(\Pi')<\mu(\Pi)$.}

Then we can prove that B$_m$-moves are complete as explained
before Proposition 5.1.4. The proof of this completeness proceeds
as the proof of Proposition 4.3.3. (We replace $y$ by $m$, and
disregard matters concerning the distance $d(x,m)$.) So we may
assume that in the list $\Pi$, which is $\Theta\Psi\Theta'$, of
$E(D_1\Box D_2)$ grounded in $D_1\Box D_2$, the members of
$\Theta$ make $F$; i.e., the members of $\Theta$ and $\Omega$ are
the same.

Out of this list $\Pi$ we make the list $\Pi_\Xi$ of $E(D_1)$ by
replacing $\Psi$ by the $\Xi$ of $\Sigma$; the list $\Sigma$,
which is $\Omega\Xi\Omega'$, is a list of $E(D_1\Box D_3)$
grounded in $D_1\Box D_3$. So $\Pi_\Xi$ is $\Theta\Xi\Theta'$.

For $\Pi_\Xi$, with the assumptions that the members of $\Theta$
make $F$, as the members of $\Omega$ do, we can prove the
following.

\prop{Proposition 5.2.3}{The list $\Pi_\Xi$ is grounded in $D_1$.}

\dkz Suppose $\Pi_\Xi\!: x\mn y\mn z$. We have the following
cases.

If $x,y,z\notin V_C$, where $V_C=E(D_1)\cap W(D_2)$, then we
appeal to the groundedness of $\Pi$ in $D_1\Box D_2$ and to Lemma
3.1.6 to obtain that $\psi_E(x,y,z)$ in~$D_1$.

If $x\in V_C$ and $y,z\notin V_C$, then for any $m$ in $E(D_E)$
we have that $\Pi\!: m\mn y\mn z$, since $\Pi_\Xi\!: x\mn y\mn z$.
Hence $\psi_E(m,y,z)$ in $D_1\Box D_2$, and, by Lemma 3.1.4, we
have that $\psi_E(x,y,z)$ in $D_1$.

If $y\in V_C$ and $x,z\notin V_C$, then we proceed analogously by
applying Lemma 3.1.5 to obtain $\psi_E(x,y,z)$ in~$D_1$.

If $z\in V_C$ and $x,y\notin V_C$, then we reason as when $x\in
V_C$ and $y,z\notin V_C$.

If $y,z\in V_C$ and $x\notin V_C$, with $x\in F$, then $\Sigma\!:
x\mn y\mn z$, since $\Pi_\Xi\!: x\mn y\mn z$. Since $\Sigma$ is
grounded in $D_1\Box D_3$, by Lemma 3.1.6, we obtain that
$\psi_E(x,y,z)$ in~$D_1$.

If $y,z\in V_C$ and $x\notin V_C$, with $x\notin F$, then $z$
precedes $y$ in $\Sigma$, because $\Pi_\Xi\!: z\mn y\mn x$. The
vertex $x$, which is not in $F$, is a member of $\Theta'$. We have
two subcases.

If $x\notin E(D_1)\cap W(D_3)$, then $x$ is a member of
$\Omega'$, and then $y$ precedes $x$ in $\Sigma$. So we have
$\Sigma\!: z \mn y\mn x$, and, by the groundedness of $\Sigma$ in
$D_1\Box D_3$ and by Lemma 3.1.6, we obtain $\psi_E(z,y,x)$ in~$D_1$.

If $x\in E(D_1)\cap W(D_3)$, then there is a vertex $w$ in
$E(D_3)$ such that $y$ precedes $w$ in $\Sigma$; since $\Theta$
and $\Omega$ have the same members, $w$ cannot be a member of
$\Omega$, and is hence a member of $\Omega'$. So we have
$\Sigma\!: z\mn y\mn w$ and, by the groundedness of $\Sigma$ in
$D_1\Box D_3$ and by Lemma 3.1.4, we obtain $\psi_E(z,y,x)$ in
$D_1$. This concludes the case when $y,z\in V_C$ and $x\notin
V_C$. The case when $x,y\in V_C$ and $z\notin V_C$ is treated
analogously.

The final case is when $x,y,z\in V_C$. Then we rely on the
groundedness of $\Sigma$ in $D_1\Box D_3$ and on Lemma 3.1.6 to
obtain that $\psi_E(x,y,z)$ in~$D_1$. \qed

\vspace{-2ex}

\section{\large\bf P$'''$-graphs are P$''$-graphs}\label{5.3}
\markright{\S 5.3. \quad P$\,'''$-graphs are P$\,''$-graphs}

In this section we will prove the assertion that is in its title.

Suppose both $D_1\Box(D_2\Box D_3)$ and $(D_1\Box D_2)\Box D_3$
are defined, i.e., stand for a D-graph. This is analogous to what
we have with (Ass~1) (see \S 1.5). Then we can prove the
following.

\prop{Proposition 5.3.1.1}{If $D_1$ and $D_2\Box D_3$ are
P-compatible, then $D_1$ and $D_2$ are P-compatible.}

\dkz Suppose a list $\Lambda_W$ of $E(D_1)$ grounded in $D_1$ and a list $\Lambda_E$ of $W(D_2\Box D_3)$ grounded in
$D_2\Box D_3$ are compatible. By Lemma 3.1.8, we conclude that it is impossible
that for some $x$ and $z$ in $W(D_2)$ and some $y$ in $W(D_3)$ we have
$\Lambda_E\!: x\mn y\mn z$.

Remove from $\Lambda_E$ all the $W$-vertices of $D_2\Box D_3$ that
belong to $W(D_3)$. The resulting list $\Lambda_E'$ is a list of
$W(D_2)$ compatible with $\Lambda_W$ by Lemma 3.1.1. (Note that
since $D_1\Box D_2$ is defined, there must be in $\Lambda_E'$ a
member of $\Lambda_W$.) By Lemma 3.1.3, this list is grounded in
$D_2$, because $\Lambda_E$ was grounded in $D_2\Box D_3$. \qed

With the same assumptions as above Proposition 5.3.1.1, we establish the following in an analogous manner by using Lemmata 3.1.8, 3.1.1 and 3.1.3.

\prop{Proposition 5.3.1.2}{If $D_1\Box D_2$ and $D_3$ are
P-compatible, then $D_2$ and $D_3$ are P-compatible.}

Suppose both $(D_1\Box D_2)\Box D_3$ and $(D_1\Box D_3)\Box D_2$
are defined. This is analogous to what we have with (Ass~2.1) (see
\S 1.5). Then we can prove the following.

\prop{Proposition 5.3.2.1}{If $D_1\Box D_2$ and $D_3$ are
P-compatible, and if $D_1\Box D_3$ and $D_2$ are P-compatible,
then $D_1$ and $D_2$, as well as $D_1$ and $D_3$, are
P-compatible.}

\dkz Assume for $\Pi$ and $\Sigma$ all that was assumed for them
before Proposition 5.1.4. Assume moreover that there is a list
$\Phi$ of $W(D_2)$ grounded in $D_2$ such that $\Sigma$ and $\Phi$
are compatible. Because we have assumed that the members of
$\Omega$ are elements of $E(D_1)$, we have that $\Phi$ is of the
form $\Phi'\Xi$, since $E(D_3)\neq\emptyset$.

According to what we concluded after Proposition 5.2.2, we may
assume that the members of $\Theta$ and $\Omega$ are the same.
Then by applying Proposition 5.2.3 we obtain that the list
$\Pi_\Xi$ of $E(D_1)$ is grounded in $D_1$.

If $\Phi'$ is the empty list, then $\Pi_\Xi$, which is
$\Theta\Xi\Theta'$, and $\Phi$, which is $\Xi$, are compatible. If
$\Phi'$ is not empty, then $\Omega$ must be empty, and hence
$\Theta$ is empty. It follows that $\Pi_\Xi$, which is
$\Xi\Theta'$, and $\Phi$, which is $\Phi'\Xi$, are compatible. So
$\Pi_\Xi$ and $\Phi$ are compatible.

If we assume that the members of $\Omega'$, instead of those of
$\Omega$, are elements of $E(D_1)$, we proceed in a dual manner,
making the members of $\Theta'$ and $\Omega'$ coincide, instead of
those of $\Theta$ and $\Omega$. So $D_1$ and $D_2$ are
P-compatible. The proof that $D_1$ and $D_3$ are P-compatible is
obtained by renaming. \qed

By a general dualizing of all that we had done to prove Proposition
5.3.2.1 we may prove the following. Suppose $D_1\Box(D_2\Box D_3)$
and $D_2\Box(D_1\Box D_3)$ are defined. This is analogous to what
we have with (Ass~2.2) (see \S 1.5).

\prop{Proposition 5.3.2.2}{If $D_1$ and $D_2\Box D_3$ are
P-compatible, and if $D_2$ and $D_1\Box D_3$ are P-compatible,
then $D_1$ and $D_3$, as well as $D_2$ and $D_3$, are
P-compatible.}

We can now prove the following.

\prop{Theorem 5.3.3}{Every P$\,'''$-graph is a P$\,''$-graph.}

\dkz We proceed by induction on the number of inner vertices of a
P$'''$-graph $D$. In the basis, when $D$ has a single inner vertex,
it is a basic D-graph, and we are done. For the induction step,
suppose $D$ is $D_W\Box D_E$. We will prove that $D_W$ and $D_E$
are P$'''$-graphs.

Take $D_W$. If $D_W$ has no cocycles, then it is trivially a
P$'''$-graph. If it has a cocycle, then take an arbitrary cocycle
of $D_W$, and assume $D_W$ is $D_W'\Box D_W''$. Then we have two
cases.

If $D_W''\Box D_E$ is defined, then, since $D$ is a P$'''$-graph,
$D_W'$ and $D_W''\Box D_E$ are P-compatible. By Proposition
5.3.1.1 we obtain that $D_W'$ and $D_W''$ are P-compatible.

If $D_W''\Box D_E$ is not defined, then $D_W'\Box D_E$ is defined,
and, since $D$ is a P$'''$-graph, $D_W'\Box D_E$ and $D_W''$ are
P-compatible. Since $D_W'\Box D_W''$ and $D_E$ are P-compatible,
for the same reason, we conclude by Proposition 5.3.2.1 that
$D_W'$ and $D_W''$ are P-compatible. So $D_W$ is a P$'''$-graph,
and by the induction hypothesis it is a P$''$-graph. We conclude
analogously by relying on Propositions 5.3.1.2 and 5.3.2.2 that
$D_E$ is a P$''$-graph, and by the inductive clause of the
definition of P$''$-graphs, since $D_W$ and $D_E$ are P-compatible by the assumption that $D$ is a P$'''$-graph, we obtain that $D_W\Box D_E$ is a
P$''$-graph. \qed

By Theorem 2.3.6, Theorem 4.4.5 and the theorem we have just
proven, we conclude that the notions of P$'$-graph, P$''$-graph
and P$'''$-graph define the same class of D-graphs, which we call simply P-graphs.

\clearpage \pagestyle{empty} \makebox[1em]{} \clearpage

\chapter{\huge\bf The Systems {\rm S1} and {\rm S2}}\label{6}
\pagestyle{myheadings}\markboth{CHAPTER 6. \quad
THE SYSTEMS {\rm S1} AND {\rm S2}}{right-head}

\section{\large\bf The system {\rm S}$\Box_P$}\label{6.1}
\markright{\S 6.1. \quad The system {\rm S}$\Box_P$}

In this chapter we put juncture into a wider context, which from the point of view of 2-categories involves, besides vertical composition, the remaining operations on 2-cells---horizontal composition and identity 2-cells. As a preliminary, we show in this section that the equations of S$\Box$ are complete not only with respect to D-graphs, as we proved in \S 1.6, but also with respect to P-graphs. We introduce next the system S1, which is an extension of S$\Box$ with unit terms and appropriate axiomatic equations. This system is proven equivalent, i.e.\ homomorphically intertranslatable, with the system S2, which has operations corresponding to the standard operations on 2-cells---viz., vertical composition, horizontal composition and identity 2-cells---and as axiomatic equations the standard assumptions for 2-categories. The systems S1 and S2 are then proven complete with respect to interpretations in appropriate kinds of graphs. For S1 these are graphs based on P-graphs, while for S2 these are graphs, dual in a certain sense, which correspond to the usual diagrams of category theory, and which will be called M-graphs. The duality in question will be investigated in more detail in \S 7.6, but it is already described in this chapter by our completeness results for S1 and S2, and the equivalence of these two systems.

The system S$\Box_P$ will not differ essentially from S$\Box$. Its equations will be of the same form, but, instead of equations between D-terms, they will now be equations between what we will call P-terms. These P-terms are also of the same form as D-terms, but they will be interpretable in P-graphs, and not in any D-graph.

A P-term will be a D-term (see \S 1.5) $\delta$ for which, in addition to the functions $W$, $E$ and $A$, we provide two functions ${\cal L}_W$\index{Lw@${\cal L}_W$} and ${\cal L}_E$\index{Le@${\cal L}_E$} such that ${\cal L}_X$\index{Lx@${\cal L}_X$} is a list (see \S 1.7) of $X(\delta)$, for $X$ being $W$ or $E$. The ordered pair (${\cal L}_W(\delta),{\cal L}_E(\delta))$ is the \emph{sequential type}\index{sequential type of P-term} of the P-term $\delta$. The equations of the system S$\Box_P$ will be equations between P-terms of the same edge type (see \S 1.5) and the same sequential type.

We define \emph{P-terms}\index{P-term} inductively by starting from a set of \emph{basic P-terms}\index{basic P-term}, which are basic D-terms $\beta$ (these are atomic symbols; see \S 1.5), and we have that ${\cal L}_X(\beta)$ is an arbitrary list of $X(\beta)$. Next we have the following inductive clause:
\begin{itemize}
\item[]if $\delta_W$ of sequential type $(\Lambda_W^W,\Phi_E\Xi\Psi_E)$ and $\delta_E$ of sequential type $(\Phi_W\Xi\Psi_W,\Lambda_E^E)$ are P-terms, then $\delta_W\Box\delta_E$ is a P-term of sequential type $(\Phi_W\Lambda_W^W\Psi_W, \Phi_E\Lambda_E^E\Psi_E)$, provided that $\Phi_E\Xi\Psi_E$ and $\Phi_W\Xi\Psi_W$ are compatible (see \S 1.7) and $\Xi$ is a list of $C=_{df}A(\delta_W)\cap A(\delta_E)=E(\delta_W)\cap W(\delta_E)\neq\emptyset$.
\end{itemize}
The condition concerning $C$ ensures that $\delta_W\Box\delta_E$ is a D-term. We define the values of $W$, $E$ and $A$ for the argument $\delta_W\Box\delta_E$ as we did for the definition of D-term in \S 1.5. (As before, we take the outermost parentheses of $\delta_W\Box\delta_E$ for granted.)

The system S$\Box_P$\index{SboxP@S$\Box_P$} is defined as S$\Box$ (see \S 1.5) with ``P-term'' substituted for ``D-term''. We can prove the following.

\prop{Proposition 6.1.1}{If $\delta=\delta'$ is derivable in {\rm S}$\Box$, then $\delta$ is a P-term of sequential type $(\Gamma,\Delta)$ iff $\delta'$ is a P-term of sequential type $(\Gamma,\Delta)$.}

\dkz We proceed by induction on the length of derivation of $\delta=\delta'$ in S$\Box$. If $\delta'$ is $\delta$, then we are done. If $\delta=\delta'$ is an instance of (Ass~1), (Ass~2.1) or (Ass~2.2), then we rely on lemmata analogous to Lemmata 2.3.1.1, 2.3.1.2, 2.3.2.1 and  2.3.2.2. For example, the lemma analogous to Lemma 2.3.1.1 says that if $\delta_1\Box\delta_2$ and $(\delta_1\Box\delta_2)\Box\delta_3$ are P-terms, and $\delta_2\Box\delta_3$ is a D-term, then $\delta_2\Box\delta_3$ and $\delta_1\Box(\delta_2\Box\delta_3)$ are P-terms. The proof of that is analogous to the proof of Lemma 2.3.1.1, and likewise for the analogues for the other lemmata.

With that we have proven the basis of the induction. The induction step, which involves the symmetry and transitivity of $=$, and congruence with $\Box$, is straightforward.\qed

\vspace{-2ex}

\prop{Proposition 6.1.2}{For $\delta$ and $\delta'$ being P-terms, $\delta=\delta'$ is derivable in {\rm S}$\Box$ iff $\delta=\delta'$ is derivable in {\rm S}$\Box_P$.}

\dkz The implication from right to left is trivial. From left to right we proceed by induction on the length of derivation of $\delta=\delta'$ in S$\Box$. If $\delta=\delta'$ is an axiomatic equation of S$\Box$, then it is an axiomatic equation of S$\Box_P$ as well. If $\delta=\delta'$ is derived in S$\Box$ by the symmetry of $=$, or by the congruence with $\Box$, then we proceed easily by applying the induction hypothesis. The only more difficult case is when $\delta=\delta'$ is derived in S$\Box$ by the transitivity of $=$ from $\delta=\delta''$ and $\delta''=\delta'$. Then, by Proposition 6.1.1, we have that $\delta''$ is a P-term, and, by the induction hypothesis, we obtain that $\delta=\delta''$ and $\delta''=\delta'$ are derivable in S$\Box_P$. Hence $\delta=\delta'$ is derivable in {\rm S}$\Box$.\qed

The following proposition is proven by a straightforward induction on the number of occurrences of $\Box$ in~$\delta$.

\prop{Proposition 6.1.3}{If $\delta$ is a P-term, then there is a construction with $(\iota(\delta),{\cal L}_W(\delta),{\cal L}_E(\delta))$ in its root.}

As a corollary of Proposition 6.1.3, and of the fact that P-graphs may be defined as P$'$-graphs, we have the following.

\prop{Proposition 6.1.4}{If $\delta$ is a P-term, then $\iota(\delta)$ is a P-graph.}

From Theorem 1.6.4, the completeness theorem for S$\Box$, with the help of Propositions 6.1.2 and 6.1.4, we obtain the following completeness theorem.

\prop{Theorem 6.1.5}{In {\rm S}$\Box_P$ we can derive $\delta=\delta'$ iff the P-graphs $\iota(\delta)$ and $\iota(\delta')$ are the same.}\index{completeness of S$\Box_P$}

\vspace{-2ex}

\section{\large\bf The system {\rm S1}}\label{6.2}
\markright{\S 6.2. \quad The system {\rm S1}}

The functions $W$, $E$ and $A$ associated with D-terms, map D-terms into the power set of an infinite set, which we will now call $\cal A$\index{A, set of all edges@$\cal A$, set of all edges}. (Intuitively, $\cal A$ is the set of all possible edges.)

Let a \emph{unit term}\index{unit term} be $\mj_\Gamma$ where $\Gamma$ is a list (see \S 1.7) of some elements of $\cal A$. We stipulate that $W(\mj_\Gamma)=E(\mj_\Gamma)=A(\mj_\Gamma)=\Gamma^s$, which is the finite (possibly empty) set of the members of $\Gamma$, and we stipulate that ${\cal L}_W(\mj_\Gamma)={\cal L}_E(\mj_\Gamma)=\Gamma$.

The P1-terms we are now going to define have, as P-terms (see \S 6.1), the functions $W$, $E$, $A$, ${\cal L}_W$ and ${\cal L}_E$ associated with them, subject to the same conditions, and their sequential types are defined analogously. The equations of the system S1 will be equations between P1-terms of the same edge type (see \S 1.5) and the same sequential type (see~\S 6.1).

We define \emph{P1-terms}\index{P1-term} inductively by starting from a set of \emph{atomic P1-terms},\index{atomic P1-term} which are either basic P-terms (see \S 6.1) or unit terms. We have an inductive clause for P1-terms involving $\Box$, which is obtained by substituting ``P1-term'' for ``P-term'' in the inductive clause of the definition of P-term (see \S 6.1), and with this clause we define the values of $W$, $E$ and $A$ for the argument $\delta_W\Box\delta_E$ as we did for the definition of D-term in \S 1.5. We also have one more inductive clause involving~$\Box$:
\begin{itemize}
\item[]if $\delta_W$ of sequential type $(\Gamma_W,\Delta_W)$ and $\delta_E$ of sequential type $(\Gamma_E,\Delta_E)$ are P1-terms, and $A(\delta_W)$ or $A(\delta_E)$ is empty, then $\delta_W\Box\delta_E$ is a P1-term of sequential type $(\Gamma_W\Gamma_E,\Delta_W\Delta_E)$.
\end{itemize}
In this case, for $Z$ being one of $W$, $E$ and $A$, we define $Z(\delta_W\Box\delta_E)$ as $Z(\delta_W)\cup Z(\delta_E)$. This concludes the definition of P1-term.

Note that in the second inductive clause above we must have that $(\Gamma_W\Gamma_E,\Delta_W\Delta_E)$ is either $(\Gamma_W,\Delta_W)$ or $(\Gamma_E,\Delta_E)$.

The system S1\index{S1} is defined as S$\Box$ (see \S 1.5) with ``P1-term'' substituted for ``D-term'' and with the following additional axiomatic equations for $\delta$ of sequential type $(\Gamma_1\Phi\Gamma_2,\Delta_1\Psi\Delta_2)$:
\begin{tabbing}
\hspace{1.6em}\=(\mj 1)\hspace{3em}\=$\mj_\Phi\Box\delta=\delta=
\delta\Box\mj_\Psi$,\index{one~1@(\mj 1)}\\[1ex]
for $\Gamma_1$ and $\Delta_1$ empty,\\[.5ex]
\>(\mj 2L)\>$\mj_{\Theta\Phi}\Box\delta=
\delta\Box\mj_{\Theta\Psi}$,\index{one 2L@(\mj 2L)}\\[1ex]
for $\Gamma_2$ and $\Delta_2$ empty,\\[.5ex]
\>(\mj 2R)\>$\mj_{\Phi\Theta}\Box\delta=
\delta\Box\mj_{\Psi\Theta}$.\index{one 2R@(\mj 2R)}
\end{tabbing}
This defines the system~S1.

The equations (\mj 2L) and (\mj 2R) could be replaced by two of their instances: the equation (\mj 2L$\Phi$), which is (\mj 2L) with $\Gamma_2$ empty, and the equation(\mj 2R$\Phi$), which is (\mj 2R) with $\Gamma_1$ empty. We can write down these new equations as follows, for $\delta$ of sequential type $(\Gamma_1\Phi\Gamma_2,\Delta_1\Psi\Delta_2)$:
\begin{tabbing}
\hspace{1.6em}\=(\mj 1)\hspace{3em}\=$\mj_\Phi\Box\delta=\delta=
\delta\Box\mj_\Psi$,\index{one~1@(\mj 1)}\kill

for $\Gamma_1$ and $\Delta_1$ empty,\\[.5ex]
\>(\mj 2L$\Phi$)\>$\mj_{\Theta\Phi\Gamma_2}\Box\delta=
\delta\Box\mj_{\Theta\Psi}$,\index{one 2LPhi@(\mj 2L$\Phi$)}\\[1ex]
for $\Gamma_2$ and $\Delta_2$ empty,\\[.5ex]
\>(\mj 2R$\Phi$)\>$\mj_{\Gamma_1\Phi\Theta}\Box\delta=
\delta\Box\mj_{\Psi\Theta}$.\index{one 2RPhi@(\mj 2R$\Phi$)}
\end{tabbing}

By taking $\delta$ in (\mj 2L$\Phi$) to be $\mj_{\Phi\Gamma_2}$, and by using (\mj 1), we obtain
\begin{tabbing}
\hspace{1.6em}\=(\mj 1)\hspace{3em}\=$\mj_\Phi\Box\delta=\delta=
\delta\Box\mj_\Psi$,\index{one~1@(\mj 1)}\kill

\>\>$\mj_{\Theta\Phi\Gamma_2}=\mj_{\Phi\Gamma_2}\Box\mj_{\Theta\Phi}$.
\end{tabbing}
From this equation, by using (\mj 1) and (Ass~1), we obtain
\begin{tabbing}
\hspace{1.6em}\=(\mj 1)\hspace{3em}\=$\mj_\Phi\Box\delta=\delta=
\delta\Box\mj_\Psi$,\index{one~1@(\mj 1)}\kill

\>\>$\mj_{\Theta\Phi\Gamma_2}=\mj_{\Theta\Phi}\Box\mj_{\Phi\Gamma_2}$.
\end{tabbing}
Hence, after renaming $\Gamma_2$ into $\Gamma$, we have
\begin{tabbing}
\hspace{1.6em}\=(\mj 1)\hspace{3em}\=$\mj_\Phi\Box\delta=\delta=
\delta\Box\mj_\Psi$,\index{one~1@(\mj 1)}\kill

\>(\mj$\Box$\mj)\>$\mj_{\Theta\Phi}\Box\mj_{\Phi\Gamma}=\mj_{\Theta\Phi\Gamma}
=\mj_{\Phi\Gamma}\Box\mj_{\Theta\Phi}$.
\end{tabbing}

Then we derive (\mj 2L) as follows:
\begin{tabbing}
\hspace{1.6em}\=(\mj 1)\hspace{3em}\=$\mj_\Phi\Box\delta=\delta=
\delta\Box\mj_\Psi$,\index{one~1@(\mj 1)}\kill

\>\>$\mj_{\Theta\Phi}\Box\delta$ \=$=\mj_{\Theta\Phi}\Box(\mj_{\Phi\Gamma_2}
\Box\delta)$,\hspace{1em}\=by (\mj1),\\[.5ex]
\>\>\>$=\mj_{\Theta\Phi\Gamma_2}\Box\delta$,\>by (Ass~1) and (\mj$\Box$\mj),\\[.5ex]
\>\>\>$=\delta\Box\mj_{\Theta\Psi}$,\>by (\mj 2L$\Phi$).
\end{tabbing}
We derive analogously (\mj 2R) by using (\mj 2R$\Phi$). Hence (\mj 2L$\Phi$) and (\mj 2R$\Phi$) can replace (\mj 2L) and (\mj 2R).

An alternative is to replace (\mj 2L) and (\mj 2R) by their instances (\mj 2L$\Psi$) and (\mj 2R$\Psi$), in which we have, respectively, $\Delta_2$ and $\Delta_1$ empty. We could write down these new equations as follows, for $\delta$, as before, of sequential type $(\Gamma_1\Phi\Gamma_2,\Delta_1\Psi\Delta_2)$:
\begin{tabbing}
\hspace{1.6em}\=(\mj 1)\hspace{3em}\=$\mj_\Phi\Box\delta=\delta=
\delta\Box\mj_\Psi$,\index{one~1@(\mj 1)}\kill

for $\Gamma_1$ and $\Delta_1$ empty,\\[.5ex]
\>(\mj 2L$\Psi$)\>$\mj_{\Theta\Phi}\Box\delta=
\delta\Box\mj_{\Theta\Psi\Delta_2}$,\index{one 2LPsi@(\mj 2L$\Psi$)}\\[1ex]
for $\Gamma_2$ and $\Delta_2$ empty,\\[.5ex]
\>(\mj 2R$\Psi$)\>$\mj_{\Phi\Theta}\Box\delta=
\delta\Box\mj_{\Delta_1\Psi\Theta}$.\index{one 2RPsi@(\mj 2R$\Psi$)}
\end{tabbing}
Still other alternatives are to replace (\mj 2L) and (\mj 2R) by (\mj 2L$\Phi$) and (\mj 2R$\Psi$), or by (\mj 2L$\Psi$) and (\mj 2R$\Phi$).

In S1 one of (Ass~2.1) and (Ass~2.2) is superfluous as an axiom; it is derivable in the presence of the other. Here is a derivation of (Ass~2.2):
\begin{tabbing}
\hspace{1.6em}\=$\delta_1\Box(\delta_2\Box\delta_3)$ \=
$=((\mj_{{\cal L}_W(\delta_1\Box(\delta_2\Box\delta_3))}\Box\delta_1)\Box\delta_2)\Box\delta_3$,
by (\mj 1) and (Ass~1),\\[.5ex]
\>\>$=((\mj_{{\cal L}_W(\delta_1\Box(\delta_2\Box\delta_3))}\Box\delta_2)\Box\delta_1)\Box\delta_3$,
by (Ass~2.1),\\[.5ex]
\>\>$=\delta_2\Box(\delta_1\Box\delta_3)$, by (Ass~1) and (\mj 1).
\end{tabbing}

\vspace{-2ex}

\section{\large\bf The system {\rm S2}}\label{6.3}
\markright{\S 6.3. \quad The system {\rm S2}}

The P2-terms we are now going to define have, as P-terms and P1-terms (see \S 6.1 and \S 6.2), the functions $W$, $E$, $A$, ${\cal L}_W$ and ${\cal L}_E$ associated with them, subject to the same conditions, and their sequential types are defined analogously. The equations of the system S2 will be equations between P2-terms of the same edge type  (see \S 1.5) and the same sequential type (see~\S 6.1).

We define \emph{P2-terms}\index{P2-term} inductively by starting from a set of \emph{atomic} P2-terms\index{atomic P2-term}, which are the same as the atomic P1-terms; i.e., they are either basic P-terms (see \S 6.1) or unit terms (see \S 6.2). We have the following two inductive clauses:
\begin{itemize}
\item[]if $\delta_W$ of sequential type $(\Lambda_W^W,\Xi)$ and $\delta_E$ of sequential type $(\Xi,\Lambda_E^E)$ are P2-terms, then $\delta_W\cirk\delta_E$ is a P2-term of sequential type $(\Lambda_W^W,\Lambda_E^E)$, provided $\Xi$ is a list of $C=_{df}A(\delta_W)\cap A(\delta_E)=E(\delta_W)=W(\delta_E)$;
\end{itemize}
for $X$ being $W$ or $E$, we have $X(\delta_W\cirk\delta_E)=X(\delta_X)$, and $A(\delta_W\cirk\delta_E)=A(\delta_W)\cup A(\delta_E)$;
\begin{itemize}
\item[]if $\delta_N$ of sequential type $(\Gamma_N,\Delta_N)$ and $\delta_S$ of sequential type $(\Gamma_S,\Delta_S)$ are P2-terms, then $\delta_N\otimes\delta_S$ is a P2-term of sequential type $(\Gamma_N\Gamma_S,\linebreak\Delta_N\Delta_S)$, provided $A(\delta_N)$ and $A(\delta_S)$ are disjoint;
\end{itemize}
for $Z$ being one of $W$, $E$ and $A$, we define $Z(\delta_N\otimes\delta_S)$ as $Z(\delta_N)\cup Z(\delta_S)$. This concludes the definition of P2-term. (The reason for using in the second clause of this definition the indices $N$ and $S$, rather than $1$ and $2$, will become apparent in clause (2$\otimes$) of the definition of M-graph in \S 6.6, an later~on.)

Note that in the first inductive clause above, for $\cirk$, we may have $C$ also empty, but, with the atomic P2-terms at our disposal, this will happen only if $A(\delta_W)$ and $A(\delta_E)$ are both empty. With that, we obtain P2-terms like $(\mj_\Lambda\cirk\mj_\Lambda)\cirk\mj_\Lambda$ for $\Lambda$ the empty list. With our atomic P2-terms, we cannot have $\delta_W\cirk\delta_E$ defined when one of $A(\delta_W)$ and $A(\delta_E)$ is empty and the other is not.

It is easy to establish that for every P2-term $\delta$ we have $A(\delta)$ empty iff all the atomic P2-terms occurring in $\delta$ are $\mj_\Lambda$ for $\Lambda$ the empty list.

The rules of the system S2\index{S2} are symmetry and transitivity of $=$, and congruence with $\cirk$ and $\otimes$ (these two congruence rules are obtained from the congruence with $\Box$ of \S 1.5 by substituting $\cirk$ and $\otimes$ respectively for $\Box$). The axiomatic equations of S2 are $\delta=\delta$ and the following equations:
\begin{tabbing}
\hspace{1.6em}\=(Ass~$\cirk$)\hspace{2em}\=$(\delta_1\cirk\delta_2)\cirk\delta_3$ \=
$=\delta_1\cirk(\delta_2\cirk\delta_3)$,\index{Ass comp@(Ass~$\cirk$)}\\[.5ex]
\>(\mj $\cirk$)\>$\mj_{{\cal L}_W(\delta)}\cirk\delta=\delta=
\delta\cirk\mj_{{\cal L}_E(\delta)}$,\index{one comp@(\mj $\cirk$)}\\[1.5ex]
\>(Ass~$\otimes$)\>$(\delta_1\otimes\delta_2)\otimes\delta_3
=\delta_1\otimes(\delta_2\otimes\delta_3)$,\index{Ass tensor@(Ass~$\otimes$)}\\[.5ex]
for $\Lambda$ the empty list,\\*[.5ex]
\>(\mj$\otimes$)\>$\mj_\Lambda\otimes\delta=\delta=
\delta\otimes\mj_\Lambda$,\index{one tensor@($\mj\,\otimes$)}\\[1.5ex]
\>($\otimes\cirk$)\>$(\delta_1\cirk\delta_2)\otimes(\delta_3\cirk\delta_4)=
(\delta_1\otimes\delta_3)\cirk(\delta_2\otimes\delta_4)$,\index{tensor comp@($\otimes\cirk$)}\\[.5ex]
\>($\otimes\,$\mj)\>$\mj_\Gamma\otimes\mj_\Delta=\mj_{\Gamma\Delta}$,\index{tensor one@($\otimes\,$\mj)}
\end{tabbing}
provided that for each of these equations both sides are defined, i.e., they are P2-terms. It is straightforward to verify that in all of these equations the two sides are P2-terms of the same edge type and the same sequential type.

The axiomatic equations of S2 are like the assumptions for 2-cells in 2-categories, where $\cirk$ is interpreted as vertical composition, $\otimes$ as horizontal composition and unit terms as identity 2-cells (see \cite{KS74} and \cite{ML98}, Sections XII.3 and XII.6). Note however that in category theory the notation is usually different (and so it is in the references we gave), not only because it uses different symbols, but also because, contrary to what we do here, the terms composed are written from right to left. The P2-terms stand for 2-cells, while the elements of $\cal A$, i.e.\ the edges, stand for 1-cells. Nothing is provided in this syntax for 0-cells, i.e.\ vertices. It would be more in the spirit of this reading of S2, but not very perspicuous, to write $\cirk_2$ for $\cirk$ and $\cirk_1$ for $\otimes$.

From the point of view of ordinary category theory, in S2 we
assume that we have lists as objects and arrows between these
lists. With the axiomatic equations (Ass~$\cirk$) and (\mj
$\cirk$) of  S2 we assume that we have a category with composition
$\cirk$ and identity arrows $\mj_\Gamma$. We have moreover a
strict monoidal structure with a bifunctor $\otimes$ and unit
object $\mj_\Lambda$ for $\Lambda$ the empty list (see
\cite{ML98}, Sections VII.1 and XI.3); on the objects, $\otimes$
is concatenation. The axiomatic equations (Ass~$\otimes$) and
(\mj$\otimes$) tell that this monoidal structure is strict
(associativity isomorphisms and isomorphisms involving the unit
are identity arrows), while ($\otimes\cirk$) and ($\otimes\,$\mj)
are the assumptions of bifunctoriality. This reading of the
axiomatic equations of S2 explains our notation.

By the equation (\mj$\otimes$), the P2-term $\mj_\Lambda$ with the empty list $\Lambda$ behaves like the unit for horizontal composition. Having this P1-term and P2-term is helpful, from a notational, computational and aesthetic point of view (like having zero), but it is not essential. Every P2-term that is not equal in S2 to $\mj_\Lambda$ is equal to a P2-term in which $\mj_\Lambda$ does not occur. In the graphs corresponding to diagrams of 2-cells, which we will call M-graphs (see \S 6.6), we have allowed the empty graph, because we will interpret $\mj_\Lambda$ by the empty graph (see \S 6.7). Had we however omitted $\mj_\Lambda$ from the language of S2, the empty graph would be excluded from M-graphs, and nothing would change essentially. In the notion of pasting scheme (see \S 7.3), which is a planar realization of an M-graph, the empty graph is not taken into account.

Omitting $\mj_\Lambda$ from P2-terms would not make it difficult to axiomatize the remaining complete fragment of S2. From the axiomatic equations of S2 we would just omit (\mj$\otimes$). For the equivalent fragment of S1, in the axiomatic equations of S1 in \S 6.2 we would require that the lists $\Phi$ and $\Psi$, as well as $\Theta$, are not empty.

\section{\large\bf The equivalence of {\rm S1} and {\rm S2}}\label{6.4}
\markright{\S 6.4. \quad The equivalence of {\rm S1} and {\rm S2}}

We show in this section that there are two translations, i.e.\ homomorphic maps, one from P1-terms to P2-terms and the other from P2-terms to P1-terms, which are inverse to each other up to derivable equality in S2 and S1 (see Propositions 6.4.1 and 6.4.2). These translations preserve derivability of equality in S1 and S2 (see Propositions 6.4.3 and 6.4.4).

We define first inductively a map $t_2$ from P1-terms to P2-terms:
\begin{tabbing}
\hspace{1.6em}$t_2(\delta)=\delta$,\hspace{1em}when $\delta$ is atomic;
\end{tabbing}
for $\delta_W$ a P1-term of sequential type $(\Lambda_W^W,\Phi_E\Xi\Psi_E)$ and $\delta_E$ a P1-term of sequential type $(\Phi_W\Xi\Psi_W,\Lambda_E^E)$, where $\Phi_E\Xi\Psi_E$ and $\Phi_W\Xi\Psi_W$ are compatible lists (which means that
at least one of $\Phi_W$ and $\Phi_E$, and at least one of $\Psi_W$ and $\Psi_E$, are the empty list; see \S 1.7),
\begin{tabbing}
\hspace{1.6em}$t_2(\delta_W\Box\delta_E)=(\mj_{\Phi_W}\otimes t_2(\delta_W)\otimes\mj_{\Psi_W})\cirk
(\mj_{\Phi_E}\otimes t_2(\delta_E)\otimes\mj_{\Psi_E})$,
\end{tabbing}
where the P2-term on the right-hand side is of sequential type $(\Phi_W\Lambda_W^W\Psi_W,$ $\Phi_E\Lambda_E^E\Psi_E)$. Since we have (Ass~$\otimes$) in S2, we may restore the missing parentheses involving $\otimes$ in this P2-term as we wish.

Next define inductively a map $t_1$ from P2-terms to P1-terms:
\begin{tabbing}
\hspace{1.6em}\=$t_1(\delta)=\delta$,\hspace{1em}when $\delta$ is atomic,\\[.5ex]
\>$t_1(\delta_W\cirk\delta_E)$ \=$=t_1(\delta_W)\Box t_1(\delta_E)$,\\[.5ex]
\>$t_1(\delta_N\otimes\delta_S)$\>$=(\mj_{{\cal L}_W(\delta_N){\cal L}_W(\delta_S)}
\Box t_1(\delta_N))\Box t_1(\delta_S)$.
\end{tabbing}

We can prove the following.

\prop{Proposition 6.4.1}{In {\rm S1} we can derive $t_1(t_2(\delta))=\delta$.}

\dkz We proceed by induction on the number $k$ of occurrences of $\Box$ in the P1-term $\delta$. If $k=0$, then $\delta$ is atomic, and $t_1(t_2(\delta))$ is~$\delta$.

If $\delta$ is $\delta_W\Box\delta_E$, then in S1 we have
\begin{tabbing}
\hspace{1.1em}\=$t_1(t_2(\delta_W\Box\delta_E))$ \=$=
t_1(((\mj_{\Phi_W}\otimes t_2(\delta_W))\otimes\mj_{\Psi_W})\cirk
((\mj_{\Phi_E}\otimes t_2(\delta_E))\otimes\mj_{\Psi_E}))$,\\[.5ex]
$=(\mj_{\Phi_W{\cal L}_W(\delta_W)\Psi_W}
\Box(\mj_{\Phi_W{\cal L}_W(\delta_W)}\Box \delta_W))\Box(\mj_{\Phi_E{\cal L}_W(\delta_E)\Psi_E}
\Box(\mj_{\Phi_E{\cal L}_W(\delta_E)}\Box \delta_E))$,\\[.5ex]
\`by the induction hypothesis and (\mj 1),\\[.5ex]
\>\>$=(\mj_{\Phi_W\Lambda_W^W\Psi_W}\Box \delta_W)\Box
(\mj_{\Phi_E\Phi_W\Xi\Psi_W\Psi_E}\Box \delta_E)$,\\[.5ex]
\` by (Ass~1) and (\mj 1),\\[.5ex]
\>\>$=(\delta_W\Box\mj_{\Phi_E\Phi_W\Xi\Psi_W\Psi_E})\Box\delta_E$,\hspace{1em} by (Ass~1) and (\mj 1),\\[.5ex]
\>\>$=\delta_W\Box\delta_E$, \hspace{1.2em} by (\mj 2L) or (\mj 2R), (Ass~1) and
(\mj 1).\`$\dashv$
\end{tabbing}

\prop{Proposition 6.4.2}{In {\rm S2} we can derive $t_2(t_1(\delta))=\delta$.}

\dkz We proceed again by induction on the number $k$ of occurrences of $\cirk$ or $\otimes$ in the P2-term $\delta$. If $k=0$, then $\delta$ is atomic, and $t_2(t_1(\delta))$ is~$\delta$.

If $\delta$ is $\delta_W\cirk\delta_E$, then in S2 we have
\begin{tabbing}
\hspace{1.6em}\=$t_2(t_1(\delta_W\cirk\delta_E))$ \=$=t_2(t_1(\delta_W)\Box t_1(\delta_E))$,\\[.5ex]
\>\>$=(\mj_\Lambda\otimes\delta_W\otimes\mj_\Lambda)\cirk
(\mj_\Lambda\otimes\delta_E\otimes\mj_\Lambda)$,\\[.5ex]
\`by the induction hypothesis, for $\Lambda$ the empty list,\\[.5ex]
\>\>$=\delta_W\cirk\delta_E$,\hspace{1em} by (\mj$\otimes$).
\end{tabbing}

If $\delta$ is $\delta_N\otimes\delta_S$, then in S2 we have
\begin{tabbing}
\hspace{1.6em}\=$t_2(t_1(\delta_N\otimes\delta_S))$ \=$=t_2((\mj_{{\cal L}_W(\delta_N){\cal L}_W(\delta_S)}\Box t_1(\delta_N))\Box t_1(\delta_S))$,\\[.5ex]
\>\>$=(\delta_N\otimes\mj_{{\cal L}_W(\delta_S)})\cirk(\mj_{{\cal L}_E(\delta_N)}\otimes\delta_S)$,\\*[.5ex]
\`by the induction hypothesis, (\mj$\otimes$) and (\mj $\cirk$),\\[.5ex]
\>\>$=\delta_N\otimes\delta_S$,\hspace{1em} by ($\otimes\cirk$) and (\mj $\cirk$).\`$\dashv$
\end{tabbing}

\prop{Proposition 6.4.3}{If in {\rm S1} we can derive $\delta=\delta'$, then in {\rm S2} we can derive $t_2(\delta)=t_2(\delta')$.}

To prove this proposition we proceed by induction on the length of the derivation of $\delta=\delta'$ in S1. Most of the work is in the basis, where for the axiomatic equations of S1 we have lengthy, but straightforward, derivations in S2 of their $t_2$ variants. Finally, we have the following proposition.

\prop{Proposition 6.4.4}{If in {\rm S2} we can derive $\delta=\delta'$, then in {\rm S1} we can derive $t_1(\delta)=t_1(\delta')$.}

To prove this proposition, we proceed as for Proposition 6.4.3.
Note that we use (Ass~2.1) only to derive in S1 the $t_1$ variant of~($\otimes\cirk$).

\section{\large\bf The completeness of {\rm S1}}\label{6.5}
\markright{\S 6.5. \quad The completeness of {\rm S1}}

In this section we show that S1 is complete with respect to an interpretation in graphs of a particular kind, with a juncture operation and unit graphs. For convenience, we rely here on edge-graphs (see \S 1.4). By the equivalence of this notion with incidented graphs, we obtain a completeness result with respect to a notion of graph based on the notion of P-graph. This requires however a modification of our notion of juncture in the presence of units, a modification involving the vertices incident with the edges of the cocycle of the juncture. To disregard vertices, as we do by relying on edge-graphs, seemed to us the best way to get around the trivial, but annoying, difficulties involved in these modifications. (The kind of difficulty we avoid this way may be sensed in the definition of $\mu$ in~\S 6.7.)

Let a \emph{D1-graph}\index{D1-graph} be a graph that is finite, acyclic, $W\mbox{\rm -}E$-functional and incidented (see \S 1.2 for these notions, and for the related notion of D-graph). Every D-graph is a D1-graph. We have that D1-graphs differ from D-graphs by possibly not being weakly connected and by possibly lacking inner vertices. (Note that D-graphs are incidented.)

Let a \emph{D1-edge-graph}\index{D1-edge-graph} be an edge-graph $H$ such that ${\cal G}(H)$ is a D1-graph (see \S 1.4 for~$\cal G$).

The empty graph (see \S 1.2), which is not a D-graph, is a D1-graph, but a single-vertex graph (see \S 1.4) is neither a D-graph nor a D1-graph, because it is not incidented.

In a \emph{straight single-edge graph}\index{straight single-edge graph} $W,E\!:A\rightarrow V$ we have $A=\{a\}$ and $V=\{W(a),E(a)\}$ with $W(a)\neq E(a)$ (see also \S 6.6; a single-edge graph is not straight when $W(a)=E(a)$).

It is easy to infer that an equivalent alternative definition of D1-graph is that these are graphs where every component is either a D-graph or a straight single-edge graph (see the beginning of \S 1.3 for the notion of component). In the empty graph the set of components is empty, and hence every component is trivially what is required to make the empty graph a D1-graph.

Joining two components of a D1-graph into a single D1-graph could be conceived as the result of a new kind of juncture, via an empty set of edges $C$ (see \S 1.3). Such a juncture, which we disallowed before, is an operation related to the $\otimes$ of~\S 6.3.

A \emph{basic D-edge-graph}\index{basic D-edge-graph} is $\langle A, \textbf{W}, \textbf{E}, \textbf{P}\rangle$ where $A=A_W\cup A_E$, $A_W\neq\emptyset$, $A_E\neq\emptyset$ and $A_W\cap A_E=\emptyset$, and, moreover,
\begin{tabbing}
\hspace{1.6em}\=$a_1\neq a_2$ \=$\Rightarrow(a_1\textbf{E}a_2\;\;$\=$\Leftrightarrow a_1,a_2\in A_W)$,\\[.5ex]
\>$a_1\neq a_2$\>$\Rightarrow(a_1\textbf{W}a_2$\>$\Leftrightarrow a_1,a_2\in A_E)$,\\[.5ex]
\>$\;a_1\textbf{P}a_2$\>$\Leftrightarrow(a_1\in A_W\;\&\; a_2\in A_E)$.
\end{tabbing}
It is straightforward to verify that $H$ is a basic D-edge-graph iff ${\cal G}(H)$ is a basic D-graph.

A \emph{unit D1-edge-graph}\index{unit D1-edge-graph} is $\mj_A=\langle A, I_A,I_A,\emptyset\rangle$, where $I_A=\{(a,a)\mid a\in A\}$ is the identity relation on $A$. Here $A$ can be the empty set, in which case $I_A$ is the empty set too. If $H$ is a unit D1-edge-graph, then in ${\cal G}(H)$ every component is a straight single-edge graph.

We define now inductively for every P1-term (see \S 6.2) $\delta$ a D1-edge-graph $\eta(\delta)$.\index{eta interpretation function@$\eta$ interpretation function}

If $\beta$ is a basic P-term (see \S 6.1 and \S 1.5), then let $\eta(\beta)$ be the basic D-edge-graph such that $A=A(\beta)$ and $A_X=X(\beta)$.

If $\mj_\Gamma$ is a unit term (see the beginning of \S 6.2), then let $\eta(\mj_\Gamma)$ be the unit D1-edge graph $\mj_{\Gamma^s}$, where $\Gamma^s$ is the set, possibly empty, of members of the list~$\Gamma$.

If $\delta_W\Box\delta_E$ is a P1-term, then let $\eta(\delta_W\Box\delta_E)$ be the D1-edge-graph $\eta(\delta_W)\cup\eta(\delta_E)$, where this union of D1-edge graphs is defined as the union of D-edge graphs that corresponds to juncture (we take the union of the two sets of edges, and the unions of the two functions \textbf{W}, \textbf{E} and \textbf{P}; see the end of~\S 1.4).

We also have the following definitions for every P1-term~$\delta$:
\begin{tabbing}
\hspace{1.6em}\=$\rho(\delta)$ \=$=\{\beta\mid\beta$ is a basic P-term occurring in $\delta\}$,\\[.5ex]\index{rho@$\rho$}
\>$\eta^*(\delta)$\>$=\langle\eta(\delta),\rho(\delta),{\cal L}_W(\delta),{\cal L}_E(\delta)\rangle$.\index{eta star interpretation function@$\eta^*$ interpretation function}
\end{tabbing}

We can prove the following lemma.

\prop{Lemma 6.5.1}{For P-terms $\delta$ and $\delta'$ we have $\iota(\delta)=\iota(\delta')$ iff $\eta^*(\delta)=\eta^*(\delta')$.}

\dkz From left to right we rely on the right-to-left direction of Theorem 6.1.5, and on the implication from $\delta=\delta'$ in S$\Box_P$ to $\eta^*(\delta)=\eta^*(\delta')$. This implication is easy to establish because the $\Box$ of S$\Box_P$ is interpreted in terms of union, and because in all the axiomatic equations the two sides have the same basic P-terms. We rely also on Proposition 6.1.1.

From right to left we pass from $\eta^*(\delta)$ to a unique $\iota(\delta)$ by taking
for each basic P-term, i.e.\ basic D-term, $\beta$ occurring in $\delta$ the inner vertex $v_\beta$ from $\iota(\beta)$ (see the beginning of \S 1.6), and adding this vertex to $\eta(\delta)$. Remember that for every pair $(A',A'')$ of sets of edges of $\eta(\delta)$ such that $(A',A'')$ belongs to $V_{\eta(\delta)}$ (see the definition of $V_H$ in \S 1.4) we have a unique $\beta$ in $\rho(\delta)$ such that
$W(\beta)=A'$ and $E(\beta)=A''$.  The $X$-vertices of $\iota(\delta)$ are induced by the edges (see the beginning of~\S 1.6).\qed

The following proposition is analogous to Proposition 1.6.3.1.

\prop{Proposition 6.5.2.1}{For every P1-term $\delta$, in {\rm S1} we have an equation of the form}
\vspace{-3ex}
\begin{tabbing}
\hspace{1.6em}$\delta=(\ldots(\mj_{{\cal L}_W(\delta)}\Box\sigma_1)\Box\ldots)\Box\sigma_n$,
\end{tabbing}
\vspace{-1ex}
\noindent\emph{for $n\geq 0$, where for distinct $i$ and $j$ in $\{1,\ldots,n\}$ we have that $\sigma_i\Box\sigma_j$ is not defined, and for every $i$ in $\{1,\ldots,n\}$ we have that $\sigma_i$ is a P-term. (If $n=0$, then our equation is} $\delta=\mj_{{\cal L}_W(\delta)}$\emph{.)}

\vspace{2ex}

\dkz We proceed by induction on the number $k$ of occurrences of $\Box$ in $\delta$. If $k=0$, then $\delta$ is either basic, in which case we have $\mj_{{\cal L}_W(\delta)}\Box\delta$ by (\mj 1), or $\delta$ is $\mj_{{\cal L}_W(\delta)}$.

If $k>0$, then $\delta$ is of the form $\delta_W\Box\delta_E$, and by the induction hypothesis, for $X$ being $W$ or $E$, in S1 we have
\begin{tabbing}
\hspace{1.6em}$\delta_X=(\ldots(\mj_{\Gamma_X}\Box\sigma_1^X)\Box\ldots)
\Box\sigma_{n_X}^X$,
\end{tabbing}
where the right-hand side is abbreviated by $\mj_{\Gamma_X}\Box\vec{\sigma}_{n_X}$. By applying (Ass~1), in S1 we obtain
\begin{tabbing}
\hspace{1.6em}$\delta_W\Box\delta_E=(\delta_W\Box\mj_{\Gamma_E})\Box\vec{\sigma}_{n_E}$,
\end{tabbing}
with an abbreviated notation of the same kind.

Then we make an auxiliary induction on $n_W$ to prove that in S1 we have
\begin{tabbing}
\hspace{1.6em}$\delta_W\Box\mj_{\Gamma_E}=\mj_{{\cal L}_W(\delta)}
\Box\vec{\sigma}_{n_W}$.
\end{tabbing}
If $n_W=0$, then $\delta_W$ is $\mj_{\Gamma_W}$, and in S1 we have $\mj_{\Gamma_W}\Box\mj_{\Gamma_E}=\mj_{{\cal L}_W(\delta)}$ by (\mj$\Box$\mj) (see~\S 6.2).

If $n_W>0$, then $\delta_W$ is $\delta_W'\Box\sigma_{n_W}^W$. Then in S1 we have either
\begin{tabbing}
\hspace{1.6em}$(\delta_W'\Box\sigma_{n_W}^W)\Box\mj_{\Gamma_E}=\delta_W'\Box\sigma_{n_W}^W$,
\end{tabbing}
by (\mj 1), or we have
\begin{tabbing}
\hspace{1.6em}$(\delta_W'\Box\sigma_{n_W}^W)\Box\mj_{\Gamma_E}=
(\delta_W'\Box\mj_{\Gamma_E'})\Box\sigma_{n_W}^W$,
\end{tabbing}
by using either (Ass~2.1), in which case $\Gamma_E'$ is $\Gamma_E$, or (Ass~1) together with (\mj 2L$\Phi$) or (\mj 2R$\Phi$) (see \S 6.2). Then we apply the induction hypothesis to $\delta_W'\Box\mj_{\Gamma_E'}$. This concludes the auxiliary induction.
So in S1 we have
\begin{tabbing}
\hspace{1.6em}$\delta=(\mj_{{\cal L}_W(\delta)}\Box\vec{\sigma}_{n_W})\Box\vec{\sigma}_{n_E}$.
\end{tabbing}

For the remainder of the proof we proceed as in the proof of Proposition 1.6.3.1, by applying only (Ass~1) and (Ass~2.1). Formally, we need auxiliary inductions on $n_W$ and $n_E$ to show that in S1 we have $\delta=\mj_{{\cal L}_W(\delta)}\Box\vec{\sigma}_n$, with the conditions of the proposition satisfied.

If $\sigma_{n_W}^W\Box\sigma_1^E$ is defined, then, by (Ass~1), in S1 we have
\begin{tabbing}
\hspace{1.6em}$((\mj_{{\cal L}_W(\delta)}\Box\vec{\sigma}_{n_W-1})
\Box\sigma_{n_W}^W)\Box\sigma_1^E=
(\mj_{{\cal L}_W(\delta)}\Box\vec{\sigma}_{n_W-1})\Box(\sigma_{n_W}^W\Box\sigma_1^E)$,
\end{tabbing}
and if $\sigma_{n_W}^W\Box\sigma_1^E$ is not defined, then, by (Ass~2.1), in S1 we have
\begin{tabbing}
\hspace{1.6em}$((\mj_{{\cal L}_W(\delta)}\Box\vec{\sigma}_{n_W-1})
\Box\sigma_{n_W}^W)\Box\sigma_1^E=
((\mj_{{\cal L}_W(\delta)}\Box\vec{\sigma}_{n_W-1})\Box\sigma_1^E)\Box\sigma_{n_W}^W$.\`$\dashv$
\end{tabbing}

\vspace{1ex}

\noindent An example of a P1-term in the form on the right-hand side of the equation of Proposition 6.5.2.1 may be found at the end of~\S 6.7.

By relying on (\mj 2L$\Psi$) or (\mj 2R$\Psi$) (see \S 6.2) instead of (\mj 2L$\Phi$) or (\mj 2R$\Phi$), we prove analogously the following proposition, which corresponds to Proposition 1.6.3.2.

\prop{Proposition 6.5.2.2}{For every P1-term $\delta$, in {\rm S1} we have an equation of the form}
\vspace{-3ex}
\begin{tabbing}
\hspace{1.6em}$\delta=\sigma_n\Box(\ldots\Box(\sigma_1\Box\mj_{{\cal L}_E(\delta)})\ldots)$,
\end{tabbing}
\vspace{-1ex}
\noindent\emph{for $n\geq 0$, where for distinct $i$ and $j$ in $\{1,\ldots,n\}$ we have that $\sigma_i\Box\sigma_j$ is not defined, and for every $i$ in $\{1,\ldots,n\}$ we have that $\sigma_i$ is a P-term. (If $n=0$, then our equation is} $\delta=\mj_{{\cal L}_W(\delta)}$\emph{.)}

\vspace{2ex}

Then we can prove the completeness of S1 with respect to~$\eta^*$.

\prop{Theorem 6.5.3}{In S1 we can derive $\delta=\delta'$ iff $\eta^*(\delta)=\eta^*(\delta')$.}\index{completeness of S1}

\dkz From left to right we proceed by an easy induction on the length of the derivation of $\delta=\delta'$ in S1. For the axiomatic equations (Ass~1), (Ass~2.1) and (Ass~2.2) we rely on the interpretation of $\Box$ in terms of union in $\eta(\delta)$ (see the proof of the left-to-right direction of Lemma 6.5.1). We rely moreover on a variant of Proposition 6.1.1, which says that if $\delta=\delta'$ is derivable in S1, then the sequential types of $\delta$ and $\delta'$ are the same. The proof of that is straightforward.

For the proof from right to left, suppose $\eta^*(\delta)=\eta^*(\delta')$. By Proposition 6.5.2.1 we have that equations of the form $\delta=\mj_\Gamma\Box\vec{\sigma}_{n}$ and $\delta'=\mj_\Gamma\Box\vec{\tau}_{n'}$ are derivable in S1 (see the proof of that proposition for the vector abbreviations). By the direction from left to right, which we have just proven, $\eta^*(\mj_\Gamma\Box\vec{\sigma}_{n})=\eta^*(\mj_\Gamma\Box\vec{\tau}_{n'})$. From that we infer that $n=n'$, because $n$ is the number of components in $\eta(\delta)$ and $\eta(\delta')$, components for edge-graphs being defined as for the corresponding incidented graphs. We infer also that there is a bijection $\pi$ of $\{1,\ldots,n\}$ to itself such that for every $i$ in $\{1,\ldots,n\}$ we have $\eta^*(\sigma_i)=\eta^*(\tau_{\pi(i)})$. By Lemma 6.5.1, we infer that $\iota(\sigma_i)=\iota(\tau_{\pi(i)})$, which, by Theorem 6.1.5, yields that $\sigma_i=\tau_{\pi(i)}$ is derivable in S$\Box_P$, and hence also in S1. We apply then (Ass~2.1), if needed, to derive $\mj_\Gamma\Box\vec{\sigma}_{n}=\mj_\Gamma\Box\vec{\tau}_{n'}$ in~S1.\qed

Let a \emph{P1-graph}\index{P1-graph} be a D1-graph where every component is either a P-graph or a straight single-edge graph. Note that the empty graph is trivially a P1-graph.

Let a \emph{P1-edge-graph}\index{P1-edge-graph} be a D1-edge-graph $H$ such that ${\cal G}(H)$ is a P1-graph (see \S 1.4 for $\cal G$, and \S 6.7 for an example). Analogously, let a \emph{P-edge-graph}\index{P-edge-graph} be a D-edge-graph $H$ such that ${\cal G}(H)$ is a P-graph.

The following proposition shows that the interpretation function $\eta^*$ is based on P1-graphs.

\prop{Proposition 6.5.4}{For every P1-term $\delta$ we have that $\eta(\delta)$ is a P1-edge-graph.}

\dkz By Proposition 6.5.2.1, in S1 we have
\begin{tabbing}
\hspace{1.6em}$\delta=(\ldots(\mj_{{\cal L}_W(\delta)}\Box\sigma_1)\Box\ldots)\Box\sigma_n$,
\end{tabbing}
where for every $i$ in $\{1,\ldots,n\}$ we have that $\sigma_i$ is a P-term. By Proposition 6.1.4, we have that $\iota(\sigma_i)$ is a P-graph.

We can then verify easily by induction on the number of occurrences of $\Box$ in the P-term $\tau$ that ${\cal H}(\iota(\tau))$ is $\eta(\tau)$ (see \S 1.4 for $\cal H$). For that we need that ${\cal H}(D_W\Box D_E)={\cal H}(D_W)\cup {\cal H}(D_E)$.

By Proposition 1.4.1, we know that ${\cal G}({\cal H}(\iota(\sigma_i)))$ is isomorphic to $\iota(\sigma_i)$, and since ${\cal H}(\iota(\sigma_i))$ is $\eta(\sigma_i)$, we have that $\eta(\sigma_i)$ is a P-edge-graph. It is then easy to conclude that $\eta(\delta)$, which~is
\begin{tabbing}
\hspace{1.6em}$(\ldots(\eta(\mj_{{\cal L}_W(\delta)})\Box\eta(\sigma_1))\Box\ldots)\Box\eta(\sigma_n)$,
\end{tabbing}
is a P1-edge-graph.\qed

\vspace{-2ex}

\section{\large\bf M-graphs}\label{6.6}
\markright{\S 6.6. \quad M-graphs}

The M-graphs (M comes from \emph{mandorla}; see the picture below), which we will introduce in this section, correspond to diagrams of 2-cells in 2-categories. This correspondence will be made manifest in the pictures of \S 6.7, and in \S 7.3 through the notion of pasting scheme---a notion of plane graph from \cite{P90} (see \S 7.3). Every non-empty M-graph is realizable in the plane as a pasting scheme, and every pasting scheme is an M-graph. In \S 6.7 we will prove the completeness of S2 with respect to an interpretation in M-graphs.

Except for the empty graph (see \S 1.2), which is an M-graph too, every M-graph $M$ will have two special distinct vertices $N(M)$ and $S(M)$\index{SM@$S(M)$} (which are respectively the source and sink of $M$; see \S 7.3). (The names of the functions $N$ and $S$ come from \emph{North} and \emph{South}.) When $M$ is not the empty graph, we define also two distinct paths $W(M)$ and $E(M)$ from $N(M)$ to $S(M)$, which are the \emph{domain}\index{domain} and \emph{codomain}\index{codomain} of $M$; since $N(M)$ and $S(M)$ are distinct, the paths $W(M)$ and $E(M)$ are non-trivial (see~\S 1.2).

Note that diagrams of 2-cells in 2-categories are usually drawn so that the domain and codomain of a 2-cell are not in the West and in the East, respectively, but in the North and in the South, i.e.\ above and below, as in the most common maps. The terminology of \emph{vertical} and \emph{horizontal} composition is suggested by this way of drawing. The usual drawings are reflected with respect to the axis $y=-x$, as well as the dual diagrams of 2-cells (see the pictures of \S 6.7), to connect them with our way of drawing P-graphs and P1-graphs. These graphs are the main subject of our work, and to draw them as we did seems more practical. This is however done at the price of having vertical composition going from West to East, and horizontal composition from North to South.

For $v$ distinct from $w$ and $n,m\geq 1$, consider a graph $B$ that corresponds to the following picture (in the shape of a mandorla)
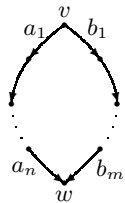
\begin{center}
\begin{picture}(60,80)(0,-10)

\put(30,0){\circle*{2}} \put(30,60){\circle*{2}}
\put(16.5,47){\circle*{2}} \put(43.5,47){\circle*{2}}
\put(16.5,13){\circle*{2}} \put(43.5,13){\circle*{2}}
\put(10,30){\circle*{2}} \put(50,30){\circle*{2}}

\qbezier(30,0)(23,5)(16.5,13) \qbezier(10,30)(11,46)(30,60)
\qbezier[4](16.5,13)(10,21)(10,30)

\qbezier(30,0)(36,5)(43.5,13) \qbezier(50,30)(49,46)(30,60)
\qbezier[4](43.5,13)(50,21)(50,30)

\put(18.5,49.5){\vector(-1,-1){2}}
\put(11.7,38.5){\vector(-1,-4){2}}

\put(28,1.7){\vector(4,-3){2}}

\put(41.5,49.5){\vector(1,-1){2}}
\put(48.3,38.5){\vector(1,-4){2}}

\put(32,1.7){\vector(-4,-3){2}}

\put(30,-3){\small\makebox(0,0)[t]{$w$}}
\put(30,63){\small\makebox(0,0)[b]{$v$}}
\put(23.5,56){\small\makebox(0,0)[br]{$a_1$}}
\put(38.5,56){\small\makebox(0,0)[bl]{$b_1$}}
\put(20,8){\small\makebox(0,0)[tr]{$a_n$}}
\put(42,9){\small\makebox(0,0)[tl]{$b_m$}}

\end{picture}
\end{center}
Let $N(B)$ be the vertex $v$, let $S(B)$ be the vertex $w$, let $W(B)$ be the path from $v$ to $w$ with the edges $a_1,\ldots,a_n$, and let $E(M)$ be the path from $v$ to $w$ with the edges $b_1,\ldots,b_m$.

We say that $B$ is a \emph{basic M-graph}\index{basic M-graph} when $W(B)$ and $E(B)$ have no common edge and no common vertex except for $v$ and~$w$.

A straight single-edge graph\index{straight single-edge graph} (as defined in \S 6.5) is $B$ where $n=m=1$ and $a_1=b_1$.

We can now give the clauses of our inductive definition of \emph{M-graph}\index{M-graph}.
\begin{itemize}
\item[(1)]Every basic M-graph, every straight single-edge graph and the empty graph are M-graphs.
\item[(2$\cirk$)]For $X$ being $W$ or $E$, let $M_X$, which is $W_X,E_X\!:A_X\rightarrow V_X$, be an M-graph. If $M_W$ and $M_E$ are not the empty graph, and they have in common as vertices and edges just the vertices and edges of $E(M_W)$, which is the same as $W(M_E)$, then $M_W\cirk M_E$, which is the graph $W,E\!:A_W\cup A_E\rightarrow V_W\cup V_E$ such that for every $a$ in $A_W\cup A_E$
\begin{tabbing}
\hspace{1.6em}$ X(a)=\left\{
\begin{array}{ll}
X_W(a) & \mbox{\rm{if }} a\in A_W,
\\[.5ex]
X_E(a) & \mbox{\rm{if }} a\in A_E
\end{array}
\right.$
\end{tabbing}
is an M-graph. (Note that for an edge $a$ in $A_W\cap A_E$, i.e.\ an edge in the path $E(M_W)$, which coincides with the path $W(M_E)$, we have $X_W(a)=X_E(a)$.) For $Y$ being $N$ or $S$, let $Y(M_W\cirk M_E)=Y(M_W)=Y(M_E)$, let $W(M_W\cirk M_E)=W(M_W)$, and let $E(M_W\cirk M_E)$ $=E(M_E)$.

\vspace{-1ex}

\hspace{1.6em}If one of $M_W$ and $M_E$ is the empty graph, then $M_W\cirk M_E$ is defined only if the other is the empty graph too, and $M_W\cirk M_E$ for both $M_W$ and $M_E$ being the empty graph is the empty graph.
\item[(2$\otimes$)]For $Y$ being $N$ or $S$, let $M_Y$, which is $W_Y,E_Y\!:A_Y\rightarrow V_Y$, be an M-graph. If $M_N$ and $M_S$ are not the empty graph, and they have in common as vertices and edges just the vertex $S(M_N)$, which is the same vertex as $N(M_S)$, then $M_N\otimes M_S$, which is the graph $W,E\!:A_N\cup A_S\rightarrow V_N\cup V_S$ where $X(a)$ is defined as in clause (2$\cirk$) above, with $N$ and $S$ substituted respectively for $W$ and $E$, is an M-graph. For $Y$ being $N$ or $S$, let $Y(M_N\otimes M_S)=Y(M_Y)$, while $X(M_N\otimes M_S)$ is the path from $N(M_N)$ to $S(M_S)$ obtained by concatenating the paths $X(M_N)$ and $X(M_S)$ with one of the two occurrences of $S(M_N)$, which is equal to $N(M_S)$, deleted. (This joining of paths is analogous to what we had with $\ast$ and semipaths in~\S 3.2.)

\vspace{-1ex}

\hspace{1.6em}If one of $M_N$ and $M_S$ is the empty graph, then $M_N\otimes M_S$ is the other graph of these two graphs, from which, if this other graph is not the empty graph, it inherits the functions $N$, $S$, $W$ and~$E$.
\end{itemize}
This concludes our definition of M-graph. Examples of M-graphs, with pictures, may be found in~\S 6.7.

An \emph{M-edge-graph}\index{M-edge-graph} is an edge-graph $H$
such that ${\cal G}(H)$ is an M-graph. When ${\cal G}(H)$ is a
basic M-graph, $H$ is a \emph{basic M-edge-graph}\index{basic
M-edge-graph}, and when ${\cal G}(H)$ is a straight single-edge
graph, $H$ is the \emph{straight single-edge
edge-graph}\index{straight single-edge edge-graph}
$\langle\{a\},I_{\{a\}},I_{\{a\}},\emptyset\rangle$, which is the
unit D1-edge graph $\mj_{\{a\}}$ (see~\S 6.5).

\section{\large\bf The completeness of {\rm S2}}\label{6.7}
\markright{\S 6.7. \quad The completeness of {\rm S2}}

We will now interpret the system S2 in M-graphs, and prove the completeness of S2 with respect to this interpretation. We introduce an interpretation function $\mu$\index{mu interpretation function@$\mu$ interpretation function} that assigns to a P2-term an M-graph, and is defined inductively as follows. As an auxiliary for this definition, we have a function $\alpha$\index{alpha@$\alpha$} that assigns to an atomic P2-term an M-edge-graph.

For atomic P2-terms, which are P1-terms, we have first that for a basic P-term $\beta$ such that ${\cal L}_W(\beta)$ is $a_1\ldots a_n$, for $n\geq 1$, and ${\cal L}_E(\beta)$ is $b_1\ldots b_m$, for $m\geq 1$, the M-edge-graph $\alpha(\beta)$ is the basic M-edge-graph that corresponds to the following picture:
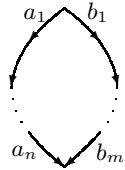
\begin{center}
\begin{picture}(60,80)(0,-10)

\qbezier(30,0)(23,5)(16.5,13) \qbezier(10,30)(11,46)(30,60)
\qbezier[4](16.5,13)(10,21)(10,30)

\qbezier(30,0)(36,5)(43.5,13) \qbezier(50,30)(49,46)(30,60)
\qbezier[4](43.5,13)(50,21)(50,30)

\put(18.5,49.5){\vector(-1,-1){2}}
\put(11.7,38.5){\vector(-1,-4){2}}

\put(28,1.7){\vector(4,-3){2}}

\put(41.5,49.5){\vector(1,-1){2}}
\put(48.3,38.5){\vector(1,-4){2}}

\put(32,1.7){\vector(-4,-3){2}}

\put(23.5,56){\small\makebox(0,0)[br]{$a_1$}}
\put(38.5,56){\small\makebox(0,0)[bl]{$b_1$}}
\put(20,8){\small\makebox(0,0)[tr]{$a_n$}}
\put(42,9){\small\makebox(0,0)[tl]{$b_m$}}

\end{picture}
\end{center}

For the unit term $\mj_\Gamma$ where $\Gamma$ is $a_1\ldots a_n$ for $n\geq 1$ we have that $\alpha(\mj_\Gamma)$ is the M-edge-graph that corresponds to the following picture:
\begin{center}
\begin{picture}(40,80)(0,-20)

\put(20,60){\vector(0,-1){20}}
\put(20,40){\vector(0,-1){20}}

\put(20,0){\vector(0,-1){20}}

\multiput(20,14.5)(0,-4){3}{\line(0,-1){.5}}

\put(17,50){\small\makebox(0,0)[br]{$a_1$}}

\put(17,-10){\small\makebox(0,0)[br]{$a_n$}}
\end{picture}
\end{center}

If $\Lambda$ is the empty list, then $\alpha(\mj_\Lambda)$ is the empty edge-graph $\langle\emptyset,\emptyset,\emptyset,\emptyset\rangle$ (see~\S 1.4).

For every atomic P1-term $\delta$ we have that $\mu(\delta)$ is ${\cal G}(\alpha(\delta))$. Note that $\mu(\mj_\Lambda)$ for $\Lambda$ the empty list is the empty graph.

Suppose we have the P2-term $\delta_W\cirk\delta_E$, and we are given the M-graphs $\mu(\delta_W)$ and $\mu(\delta_E)$. Consider the paths $E(\mu(\delta_W))$ and $W(\mu(\delta_E))$, which are made of the same edges in the same order. Let the $k$-th vertex in the first path be $(A_W',A_W'')$, and let the $k$-th vertex in the second path be $(A_E',A_E'')$; here $k\geq 1$. Let the M-graph $M_X$ isomorphic to $\mu(\delta_X)$ be obtained from $\mu(\delta_X)$ by replacing every such vertex $(A_X',A_X'')$ by $(A_W'\cup A_E',A_W''\cup A_E'')$. We have that the paths $E(M_W)$ and $W(M_E)$ coincide. If $\mu(\delta_X)$ is the empty graph, then the M-graph $M_X$ is also the empty graph. We take $\mu(\delta_W\cirk\delta_E)$ to be $M_W\cirk M_E$ (see (2$\cirk$) in~\S 6.6).

Suppose we have the P2-term $\delta_N\otimes\delta_S$, and we are given the M-graphs $\mu(\delta_N)$ and $\mu(\delta_S)$. If neither $\mu(\delta_N)$ nor $\mu(\delta_S)$ is the empty graph, then we have that $S(\mu(\delta_N))$ is a vertex of the form $(A',\emptyset)$, and $N(\mu(\delta_S))$ is a vertex of the form $(\emptyset, A'')$. Let $Y$ be $N$ or $S$, and let $\bar{N}$ be $S$, and $\bar{S}$ be $N$. Let the M-graph $M_Y$ isomorphic to $\mu(\delta_Y)$ be obtained from $\mu(\delta_Y)$ by replacing the vertex $\bar{Y}(\mu(\delta_Y))$ by $(A',A'')$.
If $\mu(\delta_Y)$ is the empty graph, then the M-graph $M_Y$ is also the empty graph. We take $\mu(\delta_N\otimes\delta_S)$ to be $M_N\otimes M_S$ (see (2$\otimes$) in \S 6.6), and this concludes our definition of~$\mu$.

For the clause concerning $\mu(\delta_W\cirk\delta_E)$ in this definition to be correct, i.e.\ for $\mu(\delta_W\cirk\delta_E)$ to be defined, we cannot have that one of $M_W$ and $M_E$ is the empty graph and the other is not, because this is required by the definition of $M_W\cirk M_E$. Suppose $\mu(\delta_W)$ and $\mu(\delta_E)$ are both defined. It is easy to see that the edges of $\mu(\delta_X)$ are the elements of $A(\delta_X)$. Since the edges of $\mu(\delta_X)$ and $M_X$ are the same, if we had that one of $M_W$ and $M_E$ is the empty graph and the other is not, then we would have that one of $A(\delta_W)$ and $A(\delta_E)$ is empty and the other is not. In that case, however, as we noted after the definition of P2-graph in \S 6.3, we would not have that $\delta_W\cirk\delta_E$ is a P2-term.

Let $\mu^*(\delta)=\langle\mu(\delta),\rho(\delta)\rangle$,\index{mu star interpretation function@$\mu^*$ interpretation function} where $\rho$ is defined as for $\eta^*$ in \S 6.5. We can establish in a straightforward manner the following soundness proposition by induction on the length of derivation in the system~S2.

\prop{Proposition 6.7.1}{If in {\rm S2} we can derive $\delta=\delta'$, then $\mu^*(\delta)=\mu^*(\delta')$.}\index{soundness of S2}

\noindent To establish also the converse implication, i.e.\ the completeness of S2 with respect to $\mu^*$, we consider first some preliminary matters.

We say that a P2-term is \emph{developed}\index{developed term} when it is of the form
$\delta_0\cirk\delta_1\cirk\ldots\cirk\delta_n$, where $n\geq 0$, parentheses tied to $\cirk$ are associated arbitrarily, $\delta_0$ is a unit term $\mj_\Gamma$, and if $n>0$, then for each $i$ in $\{1,\ldots,n\}$ we have that $\delta_i$ is of the form $(\mj_{\Gamma_i'}\otimes\beta_i\otimes\mj_{\Gamma_i''})$, with parentheses tied to the two occurrences of $\otimes$ associated arbitrarily, and $\beta_i$ a basic P-term. Here $\beta$ is the \emph{core}\index{core} of $\delta_i$. If $n=0$, then $\delta_0\cirk\delta_1\cirk\ldots\cirk\delta_n$ is just $\delta_0$, which is of the form $\mj_\Gamma$. An example of a developed P2-term may be found in $\gamma$, at the end of this section.

We can prove the following development lemma.

\prop{Lemma 6.7.2}{For every P2-term $\delta$ there is a developed P2-term $\delta^\dagger$ such that $\delta=\delta^\dagger$ is derivable in~S2.}

\dkz We proceed by induction on the complexity of $\delta$. If $\delta$ is an atomic P2-term, the lemma is established easily by using, if need there is, (\mj $\cirk$) or~(\mj$\otimes$).

If $\delta$ is of the form $\delta'\cirk\delta''$, then we apply the induction hypothesis to $\delta'$ and $\delta''$, and in S2 we have
\begin{tabbing}
\hspace{.5em}$(\mj_{\Gamma'}\cirk\delta_1'\cirk\ldots\cirk\delta_{n'}')\cirk
(\mj_{\Gamma''}\cirk\delta_1''\cirk\ldots\cirk\delta_{n''}'')=
\mj_{\Gamma'}\cirk\delta_1'\cirk\ldots\cirk\delta_{n'}'\cirk
\delta_1''\cirk\ldots\cirk\delta_{n''}''$,
\end{tabbing}
by (\mj $\cirk$) and (Ass~$\cirk$).

If $\delta$ is of the form $\delta'\otimes\delta''$, then we apply again the induction hypothesis to $\delta'$ and $\delta''$, and we make an auxiliary induction on $n'\pl n''$. If $n'\pl n''=0$, then, by ($\otimes\,$\mj), in S2 we have
\begin{tabbing}
\hspace{1.6em}$\mj_{\Gamma'}\otimes\mj_{\Gamma''}=\mj_{\Gamma'\Gamma''}$.
\end{tabbing}
If $n'>0$, then in S2 we have
\begin{tabbing}
\hspace{.5em}$(\mj_{\Gamma'}\cirk\delta_1'\cirk\ldots\cirk\delta_{n'}')\otimes
(\mj_{\Gamma''}\cirk\delta_1''\cirk\ldots\cirk\delta_{n''}'')$\\*[.5ex]
\hspace{2.2em}\=$=
(\mj_{\Gamma'}\cirk\delta_1'\cirk\ldots\cirk\delta_{n'}')\otimes
(\mj_{\Gamma''}\cirk\delta_1''\cirk\ldots\cirk\delta_{n''}''\cirk\mj_\Delta)$,
\hspace{1em}by (\mj $\cirk$),\\*[.5ex]
\>$=
((\mj_{\Gamma'}\cirk\delta_1'\cirk\ldots\cirk\delta_{n'-1}')\otimes
(\mj_{\Gamma''}\cirk\delta_1''\cirk\ldots\cirk\delta_{n''}''))\cirk
(\delta_{n'}'\otimes\mj_\Delta)$,
by ($\otimes\cirk$).
\end{tabbing}
Then we apply the induction hypothesis of the auxiliary induction and ($\otimes\,$\mj). We proceed analogously if $n''>0$.

In all that we rely on (Ass~$\cirk$) and (Ass~$\otimes$) to associate parentheses as we wish.\qed

Then we can prove the completeness of S2 with respect to~$\mu^*$.

\prop{Theorem 6.7.3}{In S2 we can derive $\delta=\delta'$ iff $\mu^*(\delta)=\mu^*(\delta')$.}\index{completeness of S2}

\dkz From left to right we have Proposition 6.7.1. For the other direction we proceed as follows.

By Lemma 6.7.2, in S2 we have
\begin{tabbing}
\hspace{1.6em}$\delta=\delta_0\cirk\delta_1\cirk\ldots\cirk\delta_n$ and $\delta'=\delta_0'\cirk\delta_1'\cirk\ldots\cirk\delta_{n'}'$,
\end{tabbing}
for the right-hand sides developed. Since $\mu(\delta)=\mu(\delta')$, by Proposition 6.7.1, we infer that in S2 we have
\begin{tabbing}
\hspace{1.6em}$\mu^*(\delta_0\cirk\delta_1\cirk\ldots\cirk\delta_n)=
\mu^*(\delta_0'\cirk\delta_1'\cirk\ldots\cirk\delta_{n'}')$.
\end{tabbing}
From $\rho(\delta_0\cirk\delta_1\cirk\ldots\cirk\delta_n)=
\rho(\delta_0'\cirk\delta_1'\cirk\ldots\cirk\delta_{n'}')$ we infer that $n=n'$. We proceed then by induction on~$n$.

If $n=0$, then $\delta$ and $\delta'$, which are respectively $\delta_0$ and $\delta_0'$, must both be the same unit term $\mj_\Gamma$. If $n>0$, then since
\begin{tabbing}
\hspace{1.6em}$\rho(\delta_0\cirk\delta_1\cirk\ldots\cirk\delta_n)=
\rho(\delta_0'\cirk\delta_1'\cirk\ldots\cirk\delta_{n'}')$,
\end{tabbing}
there must be an $i$ in $\{1,\ldots,n\}$ such that $\delta_n$ and $\delta_i'$ have the same core $\beta$. If $i\neq n$, then by using equations of the form
\begin{tabbing}
$(\mj_{\Gamma'}\otimes\beta_1\otimes\mj_{\Gamma''})\cirk
(\mj_{\Delta'}\otimes\beta_2\otimes\mj_{\Delta''})=
(\mj_{\Xi'}\otimes\beta_2\otimes\mj_{\Xi''})\cirk
(\mj_{\Pi'}\otimes\beta_1\otimes\mj_{\Pi''})$,
\end{tabbing}
with the proviso that $E(\beta_1)\cap W(\beta_2)$ and $E(\beta_2)\cap W(\beta_1)$ are empty, which are derivable in S2 with the help of ($\otimes\,$\mj) and (Ass~$\otimes$), and the essential use of (\mj $\cirk$) and two applications of ($\otimes\cirk$), we obtain that, for $\delta''$ being
\begin{tabbing}
\hspace{1.6em}$\delta_0\cirk\delta_1'\cirk\ldots\cirk\delta_{i-1}'\cirk
\delta_{i+1}''\cirk\ldots\cirk\delta_n''$,
\end{tabbing}
in S2 we can derive $\delta'=\delta''\cirk\delta_n$. If $i=n$, then, for $\delta''$ being $\delta_0\cirk\delta_1'\cirk\ldots\cirk\delta_{n-1}'$,
in S2 we can derive $\delta'=\delta''\cirk\delta_n$.

For $\delta''$ being either of these two, we infer that $\mu(\delta_0\cirk\delta_1\cirk\ldots\cirk\delta_{n-1})=\mu(\delta'')$, and, by the induction hypothesis, we obtain that $\delta_0\cirk\delta_1\cirk\ldots\cirk\delta_{n-1}=\delta''$ is derivable in S2. From that we infer that $\delta_0\cirk\delta_1\cirk\ldots\cirk\delta_n=\delta''\cirk\delta_n$, and hence also $\delta=\delta'$, are derivable in~S2.\qed

Note that $\rho$ plays an essential role in this completeness proof. Without involving $\rho$ in $\mu^*$, and by having just the interpretation function $\mu$, completeness for S2 would fail for the simple reason that there may be two different basic P-terms $\beta_1$ and $\beta_2$ of the same sequential type $(a,b)$; we have $\mu(\beta_1)=\mu(\beta_2)$, but $\beta_1=\beta_2$ is not derivable in S2. However, even if we secured that there are no different basic P-terms of the same sequential type, we would still need $\rho$, as the following example shows. For
\begin{tabbing}
\hspace{1.6em}\=$\beta_1$ of sequential type $(a,b)$,
\hspace{2em}\=$\beta_1'$ of sequential type $(a,c)$,\\[.5ex]
\>$\beta_2$ of sequential type $(b,c)$,\>$\beta_2'$ of sequential type $(c,b)$,\\[.5ex]
\>$\beta_3$ of sequential type $(c,d)$,\>$\beta_3'$ of sequential type $(b,d)$
\end{tabbing}
we have $\mu(\beta_1\cirk\beta_2\cirk\beta_3)=\mu(\beta_1'\cirk\beta_2'\cirk\beta_3')$, but $\beta_1\cirk\beta_2\cirk\beta_3=\beta_1'\cirk\beta_2'\cirk\beta_3'$ is not derivable in~S2.

The completeness of S2 of Theorem 6.7.3 corresponds to the unicity part of Theorem 3.3 of \cite{P90}, and the proof just given provides details for the sketch of the proof in the last paragraph of \cite{P90}. The remaining part of Theorem 3.3 of \cite{P90} is tied to our Lemma 6.7.2. One may understand P2-terms as formalizing what is called there ``ways to obtain composites'', and M-graphs correspond, as we already said at the beginning of \S 6.6, to what is called there ``pasting schemes''. Pasting schemes are not defined inductively as M-graphs are, and an essential ingredient of their definition, which may be found in \S 7.3, is planarity. One of the main purposes of this work is a combinatorial analysis of this planarity in terms of the notion of P-graph.

The connection between M-graphs and P-graphs may be derived from Theorem 6.5.3, the completeness of S1 with respect to $\eta^*$, which is based on P1-graphs and P-graphs, next from Theorem 6.7.3 above, the completeness of S2 with respect to $\mu^*$, which is based on M-graphs, and finally from the translations that establish the equivalence of S1 and S2 in \S 6.4. In P1-graphs one forgets about the lists of edges, which are incorporated in M-graphs in the paths made of the duals of theses edges. This duality, which is treated more precisely for planar realizations in \S 7.6, will here be only illustrated by some pictures, and the accompanying comments.

We have here on the left a picture of a basic M-edge-graph (see the end of \S 6.6) and on the right a picture of the corresponding D-edge-graph (see~\S 1.4):
\begin{center}
\begin{picture}(220,80)(0,-10)

\qbezier(30,0)(-10,30)(30,60) \qbezier(30,0)(70,30)(30,60)

\put(14,43.5){\vector(-1,-2){2}} \put(12.3,20){\vector(1,-2){2}}
\put(28,1.7){\vector(4,-3){2}}

\put(41.5,49.5){\vector(1,-1){2}}
\put(48.3,38.5){\vector(1,-4){2}}
\put(45.6,16.4){\vector(-2,-3){2}} \put(32,1.7){\vector(-4,-3){2}}

\put(23.5,56){\small\makebox(0,0)[br]{$a_1$}}
\put(8,30){\small\makebox(0,0)[r]{$a_2$}}
\put(20,8){\small\makebox(0,0)[tr]{$a_3$}}
\put(38.5,56){\small\makebox(0,0)[bl]{$b_1$}}
\put(49,42){\small\makebox(0,0)[bl]{$b_2$}}
\put(52,26){\small\makebox(0,0)[tl]{$b_3$}}
\put(43,9){\small\makebox(0,0)[tl]{$b_4$}}


\put(150,60){\vector(1,-1){30}}

\put(150,30){\vector(1,0){30}}

\put(150,0){\vector(1,1){30}}

\put(180,30){\vector(1,1){30}} \put(180,30){\vector(3,1){30}}
\put(180,30){\vector(3,-1){30}} \put(180,30){\vector(1,-1){30}}

\put(167,46){\small\makebox(0,0)[bl]{$a_1$}}
\put(165,32){\small\makebox(0,0)[b]{$a_2$}}
\put(155,19){\small\makebox(0,0)[tl]{$a_3$}}

\put(195,51){\small\makebox(0,0)[b]{$b_1$}}
\put(200,39){\small\makebox(0,0)[b]{$b_2$}}
\put(200,26){\small\makebox(0,0)[b]{$b_3$}}
\put(209,14){\small\makebox(0,0)[t]{$b_4$}}

\end{picture}
\end{center}
The region between the two paths in the left picture is replaced
by a vertex in the right picture. The order of the edges in the
paths in the left picture is replaced by their lists in the right
picture. When in the right picture we forget about this order, and
deal not with a given order, but with \emph{orderability}, then we
reach the level at which we have dealt with P-graphs.

Here is next on the left a picture for the M-edge-graph (see the end of \S 6.6) corresponding to $\mu(\mj_{a_1a_2a_3})$, and on the right a picture of the P1-edge-graph (see \S 6.5, before Proposition 6.5.4) $\eta(\mj_{a_1a_2a_3})$:
\begin{center}
\begin{picture}(220,70)(0,-5)

\put(30,60){\vector(0,-1){20}} \put(30,40){\vector(0,-1){20}}
\put(30,20){\vector(0,-1){20}}

\put(33,51){\small\makebox(0,0)[l]{$a_1$}}
\put(33,31){\small\makebox(0,0)[l]{$a_2$}}
\put(33,11){\small\makebox(0,0)[l]{$a_3$}}


\put(165,50){\vector(1,0){30}} \put(165,30){\vector(1,0){30}}
\put(165,10){\vector(1,0){30}}

\put(180,52){\small\makebox(0,0)[b]{$a_1$}}
\put(180,32){\small\makebox(0,0)[b]{$a_2$}}
\put(180,12){\small\makebox(0,0)[b]{$a_3$}}

\end{picture}
\end{center}

Here is finally on the left a picture for a more complex M-edge-graph, and on the right a picture for the corresponding P1-edge-graph:
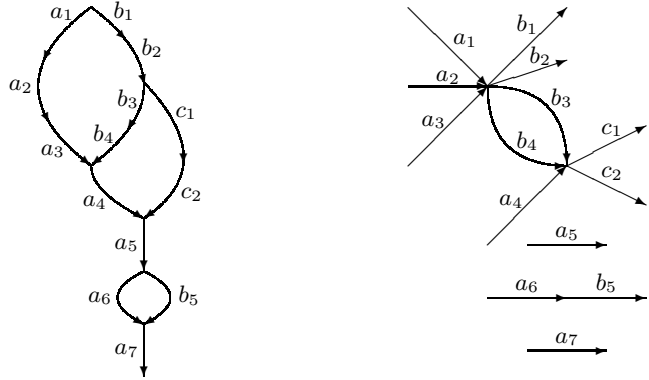
\begin{center}
\begin{picture}(220,145)(0,-80)

\qbezier(30,0)(-10,30)(30,60)

\qbezier(30,0)(70,30)(30,60)

\put(14,43.5){\vector(-1,-2){2}}

\put(12.3,20){\vector(1,-2){2}}
\put(28,1.7){\vector(4,-3){2}}

\put(41.5,49.5){\vector(1,-1){2}}
\put(48.3,38.5){\vector(1,-4){2}}
\put(45.6,16.4){\vector(-2,-3){2}} \put(32,1.7){\vector(-4,-3){2}}

\put(23.5,56){\small\makebox(0,0)[br]{$a_1$}}
\put(8,30){\small\makebox(0,0)[r]{$a_2$}}
\put(20,8){\small\makebox(0,0)[tr]{$a_3$}}

\put(38.5,56){\small\makebox(0,0)[bl]{$b_1$}}
\put(49,42){\small\makebox(0,0)[bl]{$b_2$}}
\put(48,29){\small\makebox(0,0)[tr]{$b_3$}}
\put(39,12){\small\makebox(0,0)[r]{$b_4$}}
\put(62,20){\small\makebox(0,0)[l]{$c_1$}}
\put(64,-11){\small\makebox(0,0)[l]{$c_2$}}

\put(36,-14){\small\makebox(0,0)[r]{$a_4$}}
\put(48,-30){\small\makebox(0,0)[r]{$a_5$}}
\put(38,-50){\small\makebox(0,0)[r]{$a_6$}}
\put(48,-70){\small\makebox(0,0)[r]{$a_7$}}
\put(63,-49){\small\makebox(0,0)[l]{$b_5$}}

\qbezier(50,31.5)(80,0)(50,-20) \qbezier(30,0)(30,-10)(50,-20)
\put(50,-20){\vector(0,-1){20}} \put(50,-60){\vector(0,-1){20}}
\qbezier(50,-40)(30,-50)(50,-60) \qbezier(50,-40)(70,-50)(50,-60)

\put(65,2.5){\vector(0,-1){2}} \put(48,-19){\vector(3,-2){2}}
\put(52,-18.5){\vector(-3,-2){2}} \put(48,-59){\vector(3,-2){2}}
\put(52,-59){\vector(-3,-2){2}}


\put(150,60){\vector(1,-1){30}}

\put(150,30){\vector(1,0){30}}

\put(150,0){\vector(1,1){30}}

\put(180,30){\vector(1,1){30}} \put(180,30){\vector(3,1){30}}

\qbezier(180,30)(180,0)(210,0) \qbezier(180,30)(210,30)(210,0)

\put(180,-30){\vector(1,1){30}}

\put(210,2.5){\vector(0,-1){2}}
\put(207.5,0){\vector(1,0){2}}

\put(210,0){\vector(2,1){30}} \put(210,0){\vector(2,-1){30}}
\put(195,-30){\vector(1,0){30}} \put(180,-50){\vector(1,0){30}}
\put(210,-50){\vector(1,0){30}} \put(195,-70){\vector(1,0){30}}

\put(167,46){\small\makebox(0,0)[bl]{$a_1$}}
\put(165,32){\small\makebox(0,0)[b]{$a_2$}}
\put(155,19){\small\makebox(0,0)[tl]{$a_3$}}

\put(189,-15){\small\makebox(0,0)[b]{$a_4$}}
\put(210,-27){\small\makebox(0,0)[b]{$a_5$}}
\put(195,-47){\small\makebox(0,0)[b]{$a_6$}}
\put(210,-67){\small\makebox(0,0)[b]{$a_7$}}

\put(195,51){\small\makebox(0,0)[b]{$b_1$}}
\put(200,39){\small\makebox(0,0)[b]{$b_2$}}
\put(208,23){\small\makebox(0,0)[b]{$b_3$}}
\put(195,13){\small\makebox(0,0)[t]{$b_4$}}
\put(223,13.5){\small\makebox(0,0)[l]{$c_1$}}
\put(223,-3){\small\makebox(0,0)[l]{$c_2$}}

\put(225,-47){\small\makebox(0,0)[b]{$b_5$}}

\end{picture}
\end{center}
With $\beta_1$ of sequential type $(a_1a_2a_3,\:b_1b_2b_3b_4)$, $\beta_2$ of sequential type $(b_3b_4a_4,$ $c_1c_2)$ and $\beta_3$ of sequential type $(a_6,b_5)$,
and with $\gamma$ being
\begin{tabbing}
\hspace{1.6em}$\mj_{a_1a_2a_3a_4a_5a_6a_7}\cirk(\mj_{\,}\otimes\beta_1\otimes\mj_{a_4a_5a_6a_7})
\cirk(\mj_{b_1b_2}\otimes\beta_2\otimes\mj_{a_5a_6a_7})$\\*[.5ex]
\`\cirk$(\mj_{b_1b_2c_1c_2a_5}\otimes\beta_3\otimes\mj_{a_7})$,
\end{tabbing}
for the picture on the left we have $\mu^*(\gamma)=\langle\mu(\gamma),\{\beta_1,\beta_2,\beta_3\}\rangle$, and with $\gamma'$ being
\begin{tabbing}
\hspace{1.6em}$(\mj_{a_1a_2a_3a_4a_5a_6a_7}\Box(\beta_1\Box\beta_2))\Box\beta_3$,
\end{tabbing}
for the picture on the right we have
\begin{tabbing}
\hspace{1.6em}$\eta^*(\gamma')=\langle\eta(\gamma'),
\{\beta_1,\beta_2,\beta_3\},\:a_1a_2a_3a_4a_5a_6a_7,\:b_1b_2c_1c_2a_5b_5a_7\rangle$.
\end{tabbing}
The P-term $\gamma$ is developed, while the P1-term $\gamma'$ is in the form on the right-hand side of the equation of Proposition 6.5.2.1.

\chapter{\huge\bf Disk D-Graphs and P-Graphs}\label{7}
\pagestyle{myheadings}\markboth{CHAPTER 7. \quad
DISK D-GRAPHS AND P-GRAPHS}{right-head}

\section{\large\bf Disk D-graphs}\label{7.1}
\markright{\S 7.1. \quad Disk D-graphs}

In this chapter we deal with geometrical matters concerning our graphs. We deal in particular with a special kind of realization of P-graphs in the plane. Such a realization is a plane graph situated within a disk with the boundary divided into two halves, one for the $W$-vertices and the other for the $E$-vertices. The plane graphs in question are called disk D-graphs, and the D1-graphs (see \S 6.5) based on disk D-graphs are called disk D1-graphs. We prove that every P-graph is isomorphic to a disk D-graph, and that, conversely, every disk D-graph is a P-graph. It follows that a graph is a P-graph iff it is isomorphic to a disk D-graph. This entails an analogous relationship between P1-graphs and disk D1-graphs.

We introduce in this chapter what we will call D1$'$-graphs, which are obtained from D1-graphs by adding a source and sink; namely, a single $W$-vertex and a single $E$-vertex. The disk D1$'$-graphs, i.e.\ the disk realizations of D1$'$-graphs, are known in the literature as pasting schemes, and we provide here two presumably new definitions of this notion. The equivalence of various definitions of P-graph, which we established in Chapters 2-5, enables us to obtain through the notion of P$'''$-graph a viable criterion for testing whether a D1$'$-graph is isomorphic to a pasting scheme. In the last section we state precisely the duality illustrated at the end \S 6.7. This is, namely, the particular relationship that exists between plane graphs that correspond to diagrams of 2-cells and disk D1-graphs, which we mentioned already in~\S 1.1.

In this section, before defining disk D-graphs, to fix terminology, we introduce as preliminary notions the notion of plane graph and a few associated notions. Our terminology and these notions are pretty standard, but they should be adapted to the notion of graph of~\S 1.2.

A \emph{plane graph}\index{plane graph} is a graph $W,E\!:A\rightarrow V$ where $A$ is a set of simple, open or closed, Jordan curves in ${\mathbf R}^2$ and $V$ is a set of points in ${\mathbf R}^2$ such that
\begin{itemize}
\item[(1)]for every open $a$ in $A$ the points $W(a)$ and $E(a)$ are the two distinct end points of $a$, and for every closed $a$ in $A$ we have that $W(a)$ and $E(a)$ are the same point of $a$,
\item[(2)]for every distinct $a$ and $b$ in $A$, if $v\in a\cap b$, then $v=W(a)$ or $v=E(a)$.
\end{itemize}
It follows immediately that we have also (2) with the consequent replaced by ``$v=W(b)$ or $v=E(b)$''.

An unessentially different notion of plane graph is obtained by requiring further that
\begin{itemize}
\item[]for every $a$ in $A$ and every $v$ in $V$, if $v\in a$, then $v=W(a)$ or $v=E(a)$
\end{itemize}
(cf.\ \cite{WKW01}, Section 2.2, Definition 2.1). If the graph is incidented (see \S 1.2), then this additional requirement is met anyway.

For every plane graph $G$, which is $W,E\!:A\rightarrow V$, let the \emph{point set}\index{point set of a plane graph} ${\cal U}(G)$\index{UG@${\cal U}(G)$}\index{UD@${\cal U}(D)$} of $G$ be the set of points of ${\mathbf R}^2$ that belong either to an edge in $A$ or are elements of~$V$.

When a graph $G$ is isomorphic to a plane graph $G'$ we say that $G'$ is a \emph{realization}\index{realization of a graph in the plane} of $G$. The graph $G$ is \emph{planar}\index{planar graph}, or \emph{realizable in the plane},\index{realizable in the plane}\index{graph realizable in the plane} when there is such a~$G'$.

(We work all the time with the assumption that our graphs are distinguished; see \S 1.2. Relinquishing this assumption for a moment, note that non-distinguished plane graphs do not exist, although non-distinguished planar graphs would be possible.)

A \emph{topological disk}\index{topological disk}\index{disk} in ${\mathbf R}^2$ is a closed subset of ${\mathbf R}^2$ homeomorphic to the unit disk $\{(x,y)\mid x^2\pl y^2\leq 1\}$. A \emph{compass disk}\index{compass disk} $\kappa$ is a topological disk in ${\mathbf R}^2$ with two distinct points on its boundary, called the \emph{north pole}\index{north pole} and the \emph{south pole}\index{south pole}\index{pole} of $\kappa$.
The north pole and the south pole of $\kappa$ determine within the boundary of $\kappa$ two disjoint subsets not including the two poles, called the \emph{$W$-meridian}\index{W-meridian@$W$-meridian} and the \emph{$E$-meridian}.\index{E-meridian@$E$-meridian}\index{meridian} Which of these two subsets is the $W$-meridian and which is the $E$-meridian is not arbitrary. We suppose from now on that we have fixed an orientation of ${\mathbf R}^2$, and it is with respect to this orientation that the sequence \emph{north pole, $E$-meridian, south pole, $W$-meridian} proceeds clockwise.

Let a \emph{disk D-graph}\index{disk D-graph} be a plane graph $D$ that is a D-graph such that, for a compass disk $\kappa$, all the $W$-vertices of $D$ are in the $W$-meridian of $\kappa$, all the $E$-vertices of $D$ are in the $E$-meridian of $\kappa$, and the remaining points in ${\cal U}(D)$ are in the interior of~$\kappa$.

We say that the compass disk $\kappa$ of this definition is \emph{associated}\index{associated compass disk}\index{compass disk associated with a disk graph} with the disk D-graph $D$. Although this associated compass disk is not uniquely determined, it is unique up to homeomorphisms that are identity maps on $D$. Examples of disk D-graphs with the associated compass disks may be found in the pictures of the proof of Proposition 7.2.1.

For a disk D-graph $D$, and for $X$ being $W$ or $E$, let $L_X(D)$ be the list of the $X$-vertices of $D$ obtained by going along the $X$-meridian from the north pole to the south pole.

For every edge $a$ of a plane graph, since $a$ is a Jordan curve, we have a one-one continuous map $f_a$ from the interval $[0,1]$ onto $a$. We say that a plane graph is \emph{eastward-growing}\index{eastward-growing plane graph} when for $r_W$ and $r_E$ in $[0,1]$, and $f_a(r_X)=(x_X,y_X)$, if $r_W<r_E$, then $x_W<x_E$. (An analogous notion is called \emph{upward planarity}\index{upward planarity} in the literature; see~\cite{GT95}.)

\section{\large\bf P-graphs are realizable as disk D-graphs}\label{7.2}
\markright{\S 7.2. \quad P-graphs are realizable as disk D-graphs}

In this section we prove what is announced in its title, which follows from the following proposition.

\prop{Proposition 7.2.1}{Every P-graph is isomorphic to an eastward-growing disk D-graph.}

\dkz Take a P-graph $D$ conceived as a P$'$-graph (see \S 1.8). So there is a construction $K$ with $(D,L_W,L_E)$ in its root. We will show by induction on the number $k$ of nodes in the tree of $K$ that there is a graph isomorphism from $D$ to an eastward-growing disk D-graph $R$ such that $L_X$ is $L_X(R)$. If $k=1$, then $D$ is a basic D-graph, for which the proposition is obvious (see the pictures at the end of \S 1.2 and at the beginning of~\S 1.6).

Suppose $K$ is $K_W\Box K_E$. For $X$ being $W$ or $E$, let $D_X$ be the root graph of $K_X$. By applying the induction hypothesis to $D_X$, we obtain the eastward-growing disk D-graph $R_X$, which is a realization of $D_X$ such that $L_E^W$ is $L_E(R_W)$ and $L_W^E$ is $L_W(R_E)$. If $\kappa_X$ is a compass disk associated with $R_X$, then by appealing to the compatibility of the lists $L_E^W$ and $L_W^E$, we may assume that part of the $E$-meridian of $\kappa_W$ coincides with part of the $W$-meridian of $\kappa_E$, so that the vertices that $D_W$ and $D_E$ share are realized by the same points of these two meridians. For example, we may have $R_W$ and $R_E$, with $\kappa_W$ and $\kappa_E$ drawn with dotted lines, as in the following pictures (which corresponds to the first two pictures in the example at the end of~\S 1.3):
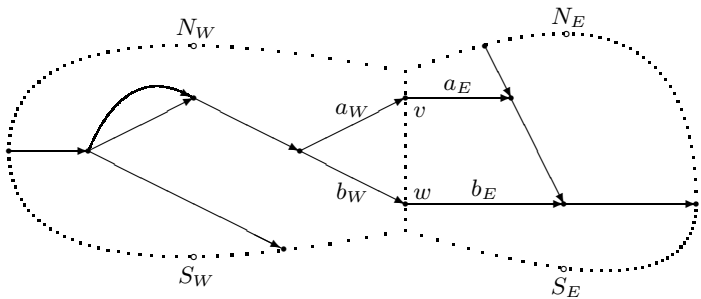
\begin{center}
\begin{picture}(260,110)(10,-45)

\put(10,10){\circle*{2}} \put(40,10){\circle*{2}}
\put(120,10){\circle*{2}} \put(80,30){\circle*{2}}
\put(114,-27){\circle*{2}} \put(160,-10){\circle*{2}}
\put(160,30){\circle*{2}}

\put(10,10){\vector(1,0){30}} \put(40,10){\vector(2,1){40}}

\put(80,30){\vector(2,-1){40}}

\put(40,10){\vector(2,-1){73}} \put(120,10){\vector(2,1){40}}
\put(120,10){\vector(2,-1){40}}

\qbezier(40,10)(55,45)(80,30) \put(77,31.5){\vector(2,-1){2}}

\put(140,24){\small\makebox(0,0)[b]{$a_W$}}
\put(140,-2){\small\makebox(0,0)[t]{$b_W$}}
\put(163,25){\small\makebox(0,0)[l]{$v$}}
\put(163,-6){\small\makebox(0,0)[l]{$w$}}

{\thicklines \qbezier[40](10,10)(10,70)(160,40)
\qbezier[40](10,10)(10,-50)(160,-20)}

\put(80,50){\circle{2}} \put(80,-30){\circle{2}}

\put(81,54){\small\makebox(0,0)[b]{$N_W$}}
\put(81,-34){\small\makebox(0,0)[t]{$S_W$}}


\put(200,30){\circle*{2}} \put(220,-10){\circle*{2}}
\put(190,50){\circle*{2}} \put(270,-10){\circle*{2}}

\put(160,30){\vector(1,0){40}} \put(160,-10){\vector(1,0){60}}
\put(220,-10){\vector(1,0){50}} \put(200,30){\vector(1,-2){20}}
\put(190,50){\vector(1,-2){10}}

\put(180,33){\small\makebox(0,0)[b]{$a_E$}}
\put(190,-7){\small\makebox(0,0)[b]{$b_E$}}

{\thicklines \qbezier[15](160,40)(160,10)(160,-20)

\qbezier[40](160,40)(270,85)(270,-10)
\qbezier[35](160,-20)(270,-55)(270,-10)}

\put(220,-34.5){\circle{2}} \put(221,54.5){\circle{2}}

\put(222,58){\small\makebox(0,0)[b]{$N_E$}}
\put(221,-38){\small\makebox(0,0)[t]{$S_E$}}

\end{picture}
\end{center}
The eastward-growing disk D-graph $R$ that is a realization of $D$ such that $L_X$ is $L_X(R)$ is obtained by removing the vertices $v$ and $w$ that were common to $R_W$ and $R_E$, and by gluing into one edge $a$ the two edges $a_W$ and $a_E$, and into one edge $b$ the two edges $b_W$ and $b_E$, as in the following picture (which corresponds to the last pictures in~\S 1.3):
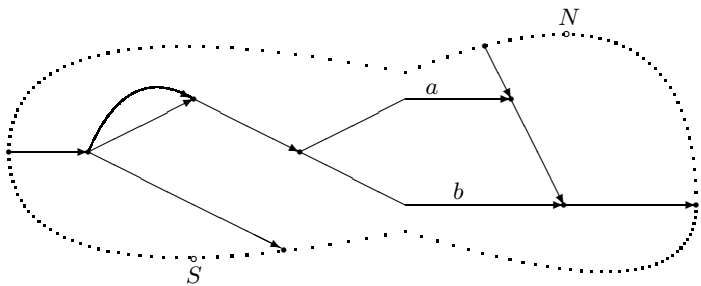
\begin{center}
\begin{picture}(260,110)(10,-45)

\put(10,10){\circle*{2}} \put(40,10){\circle*{2}}
\put(120,10){\circle*{2}} \put(80,30){\circle*{2}}
\put(114,-27){\circle*{2}}

\put(10,10){\vector(1,0){30}} \put(40,10){\vector(2,1){40}}

\qbezier(40,10)(55,45)(80,30) \put(77,31.5){\vector(2,-1){2}}

\put(80,30){\vector(2,-1){40}}

\put(40,10){\vector(2,-1){73}} \put(120,10){\line(2,1){40}}
\put(120,10){\line(2,-1){40}}

{\thicklines \qbezier[40](10,10)(10,70)(160,40)
\qbezier[40](10,10)(10,-50)(160,-20)}

\put(80,-30.5){\circle{2}} \put(80,-33.5){\small\makebox(0,0)[t]{$S$}}


\put(200,30){\circle*{2}} \put(220,-10){\circle*{2}}
\put(190,50){\circle*{2}} \put(270,-10){\circle*{2}}

\put(160,30){\vector(1,0){40}} \put(160,-10){\vector(1,0){60}}
\put(220,-10){\vector(1,0){50}} \put(200,30){\vector(1,-2){20}}
\put(190,50){\vector(1,-2){10}}

\put(170,32){\small\makebox(0,0)[b]{$a$}}
\put(180,-8){\small\makebox(0,0)[b]{$b$}}

{\thicklines \qbezier[40](160,40)(270,85)(270,-10)
\qbezier[35](160,-20)(270,-55)(270,-10)}

\put(221,54.5){\circle{2}}
\put(222,58){\small\makebox(0,0)[b]{$N$}}

\end{picture}
\end{center}

The boundary of the compass disk $\kappa$ associated here with $R$ is drawn with dotted lines in this picture of $R$. This disk is obtained from $\kappa_W$ and $\kappa_E$ by omitting the part of their boundaries that they share, and by taking as the north pole $N$ of $\kappa$ the north pole $N_E$ of $\kappa_E$, while the south pole $S$ of $\kappa$ will be the south pole $S_W$ of~$\kappa_W$.

The rules for choosing the poles of $\kappa$ are the following. We have that the list $L_E^W$ is $\Phi_E\Xi\Psi_E$, while $L_W^E$ is $\Phi_W\Xi\Psi_W$. If $\Phi_E$ and $\Psi_W$ are empty, as in our example above, then we take $N=N_E$ and $S=S_W$.  If $\Phi_W$ and $\Psi_E$ are empty, then we take $N=N_W$ and $S=S_E$, and, for $X$ being $W$ or $E$, if $\Phi_X$ and $\Psi_X$ are empty, then we take $N=N_X$ and $S=S_X$. If more than two of the four lists $\Phi_W$, $\Phi_E$, $\Psi_W$ and $\Psi_E$ are empty, then we have more than one of these rules for poles applying. The results are not the same, but the differences are not important.\qed

\vspace{-2ex}

\section{\large\bf Disk D-graphs are P-graphs}\label{7.3}
\markright{\S 7.3. \quad Disk D-graphs are P-graphs}

In this section we are going to prove that every disk D-graph is a P-graph. This implies the converse of the proposition that every P-graph is isomorphic to a disk D-graph, which follows from Proposition 7.2.1. More precisely, we will show that every disk D-graph is a P$'''$-graph. For Proposition 7.2.1 we relied, on the other hand, on the notion of P$'$-graph. Here is where our proof of the equivalence of the notions of P$'$-graph and P$'''$-graph helps us. We must first deal however with a number of preliminary matters.

The following lemma is a variant of Lemma~2 of \cite{M84}, which is proven there with the help of Brouwer's Fixed Point Theorem.

\prop{Lemma 7.3.1}{For four distinct points $v_1$, $v_2$, $v_3$ and $v_4$ occurring in that order on the boundary of a topological disk $\kappa$ in ${\mathbf R}^2$, every open Jordan curve joining $v_1$ and $v_3$ and every open Jordan curve joining $v_2$ and $v_4$, which are both included in $\kappa$, must intersect.}

We can then prove the following, for $X$ being $W$ or $E$.

\prop{Lemma 7.3.2}{For every disk D-graph $D$ the list $L_X(D)$ is grounded in~$D$.}

\dkz In this proof we use the notation introduced in \S 1.9. Suppose we have $L_X(D)\!:v\mn u\mn w$. Take a semipath $\sigma$ in $[v,w]$ and a semipath $\tau$ in $[u,t]$ such that $t$ is an $\bar{X}$-vertex of $D$. So we have $t,v,u$ and $w$ occurring in that order on the boundary of a compass disk associated with $D$. It follows from Lemma 7.3.1 that $\sigma$ and $\tau$ must intersect, and hence $\psi_X(v,u,w)$, since $D$ is a plane graph.\qed

Our purpose next is to show that the edges of every cocycle can be linearly ordered, so as to make a list. A theorem in \cite{B73} (Theorem~3, Section 2.2) asserts and proves that to a cycle of a plane pseudograph in the sense of \cite{H69} (Chapter~2, with multiple edges and loops) there corresponds a cocycle in the dual graph.  This theorem asserts also the converse---namely, that to a cocycle there corresponds a cycle in the dual graph---but without proof. This converse assertion is close to what we need for the linear orderability of the edges of every cocycle, but since it is not exactly the same (we pass from a cocycle not to a cycle in the dual graph, but to a list, which is not cyclic), since the notions of graph in question are not exactly the same, and since \cite{B73} does not provide a proof, we give an independent proof of what we need.

For a plane graph $G$, a \emph{face}\index{face of a plane graph}
of $G$ is a connected component (in the topological sense) of
${\mathbf R}^2-{\cal U}(G)$. A face of $G$ is an open subset of
${\mathbf R}^2$, in which we do not find the edges and vertices of
$G$. Assuming that ${\cal U}(G)$ is included in some sufficiently
large disk $\cal D$, which we may do when $G$ is finite, exactly
one of the faces of $G$ is unbounded---namely, the face in which
${\mathbf R}^2-{\cal D}$ is included. This face is the \emph{outer
face}\index{outer face of a plane graph} of $G$; the other faces
are the \emph{inner faces}\index{inner face of a plane graph} of
$G$. (The terminology of this paragraph agrees with that of
\cite{D10}, Section 4.2, and is close to that of \cite{B73},
Section 2.2, and \cite{H69}, Chapter 11.) The
\emph{boundary}\index{boundary of a face of a plane graph} of a
face $f$ is the closure of $f$ minus~$f$.

An inner face $f$ in a plane graph will be called \emph{bipolar}\index{bipolar face of a plane graph} when its boundary is made of two paths from a vertex $w$ to a vertex $v$. The two distinct paths must be non-trivial and $w$ must be distinct from $v$. We call the path with the face on the right-hand side the \emph{north path}\index{north path of a bipolar face}, and the path with the face on the left-hand side the \emph{south path},\index{south path of a bipolar face} as in the following picture:
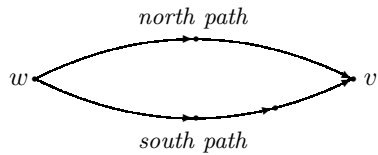
\begin{center}
\begin{picture}(120,50)(0,20)

\put(0,45){\circle*{2}} \put(61,30){\circle*{2}}
\put(61,60){\circle*{2}} \put(91,34){\circle*{2}}
\put(120.5,45){\circle*{2}}

\qbezier(0,45)(60,15)(120,45)
\qbezier(0,45)(60,75)(120,45)

\put(58,30){\vector(1,0){2}} \put(88,33.5){\vector(4,1){2}}
\put(118,44){\vector(2,1){2}}

\put(118,46){\vector(2,-1){2}}
\put(58,60){\vector(1,0){2}}

\put(-6,45){\makebox(0,0){$w$}}

\put(127,45){\makebox(0,0){$v$}}

\put(60,68){\makebox(0,0){\small \emph{north path}}}

\put(60,21){\makebox(0,0){\small \emph{south path}}}

\end{picture}
\end{center}
The left and right position of the face is here determined by the orientation we have assumed for ${\mathbf R}^2$ (see \S 7.1, where we have decided upon the $W$-meridian and $E$-meridian). The notion of bipolar face (not under that name) may be found in Proposition 2.6 of \cite{P90}, which we will state after introducing some other notions.

When a graph $G$ has a single $W$-vertex that vertex is the \emph{source}\index{source of a graph} of $G$, and when $G$ has a single $E$-vertex that vertex is the \emph{sink}\index{sink of a graph} of~$G$.

Let a \emph{pasting scheme}\index{pasting scheme} be a finite plane graph with source and sink, which are distinct, which are both on the boundary of the outer face, and which are such that for every vertex $v$ there is a path from the source to $v$ and a path from $v$ to the sink; moreover, every inner face is bipolar. In \cite{P90}, where one may find this definition, it is shown in Proposition 2.6 that an equivalent alternative definition of pasting scheme is obtained by replacing the requirement of bipolarity for inner faces by the requirement of acyclicity for the graph.

It is not difficult to show by induction that every M-graph (see \S 6.6) that is not the empty graph is isomorphic to a pasting scheme. Conversely, one can show that every pasting scheme is an M-graph. The proof of that would proceed by induction on the number of inner faces in the pasting scheme (cf.\ Proposition 2.10 of \cite{P90}). It is easy to see that every graph isomorphic to an M-graph is an M-graph, and so we may conclude that a non-empty graph is an M-graph iff it is isomorphic to a pasting scheme. (We have found it more convenient in Chapter~6 to allow the empty graph as an M-graph for the reasons mentioned at the end of~\S 6.3.)

For $D$ a disk D-graph, consider a plane graph $D'$ obtained by adding two new vertices $s$ and $t$, and new edges from $s$ to every $W$-vertex of $D$, and from every $E$-vertex of $D$ to $t$, such that $s$ and $t$ are on the boundary of the outer face of $D'$. It is easy to conclude that $D'$ is a pasting scheme by relying on the alternative definition mentioned above. We call $D'$ a \emph{source-sink closure}\index{source-sink closure of a disk D-graph} of~$D$.

Then, by the right-to-left direction of Proposition 2.6 of \cite{P90}, which says that acyclicity implies the bipolarity of inner faces, we have the following.

\prop{Lemma 7.3.3}{Every inner face of a disk D-graph is bipolar.}

\noindent It is enough to note that, for a disk D-graph $D$, a source-sink closure $D'$ of $D$ is acyclic, and if every inner face of $D'$ is bipolar, so is every inner face of~$D$.

When two distinct bipolar inner faces $f$ and $g$ of a plane graph share an edge $a$ so that $a$ is in the south path of $f$ and in the north path of $g$, we will say that $f$ \emph{precedes}\index{precedes, for faces} $g$. We can prove the following.

\prop{Lemma 7.3.4}{There is no sequence $f_1,\ldots,f_n$, with $n\geq 2$, of inner faces of a disk D-graph such that for every $i$ in $\{1,\ldots,n\mn 1\}$ we have that $f$ precedes $f_{i+1}$ and $f_n$ precedes $f_1$.}

\noindent Otherwise, we would not have a disk D-graph, because we would have either something like
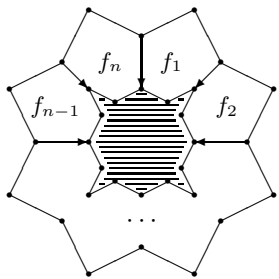
\begin{center}
\begin{picture}(100,110)(0,-5)

\put(0,30){\circle*{1.5}} \put(0,70){\circle*{1.5}}
\put(10,50){\circle*{1.5}} \put(20,20){\circle*{1.5}}
\put(20,80){\circle*{1.5}} \put(30,0){\circle*{1.5}}
\put(30,30){\circle*{1.5}} \put(30,50){\circle*{1.5}}
\put(30,70){\circle*{1.5}} \put(30,100){\circle*{1.5}}
\put(35,40){\circle*{1.5}} \put(35,60){\circle*{1.5}}
\put(40,35){\circle*{1.5}} \put(40,65){\circle*{1.5}}
\put(50,10){\circle*{1.5}} \put(50,30){\circle*{1.5}}
\put(50,70){\circle*{1.5}} \put(50,90){\circle*{1.5}}
\put(60,35){\circle*{1.5}} \put(60,65){\circle*{1.5}}
\put(65,40){\circle*{1.5}} \put(65,60){\circle*{1.5}}
\put(70,0){\circle*{1.5}} \put(70,30){\circle*{1.5}}
\put(70,50){\circle*{1.5}} \put(70,70){\circle*{1.5}}
\put(70,100){\circle*{1.5}} \put(80,20){\circle*{1.5}}
\put(80,80){\circle*{1.5}} \put(90,50){\circle*{1.5}}
\put(100,30){\circle*{1.5}} \put(100,70){\circle*{1.5}}

\put(50,90){\vector(0,-1){20}} \put(90,50){\vector(-1,0){20}}
\put(10,50){\vector(1,0){20}}

\put(80,80){\vector(-1,-1){10}} \put(20,80){\vector(1,-1){10}}

\put(50,90){\line(2,1){20}} \put(70,100){\line(1,-2){10}}
\put(80,80){\line(2,-1){20}} \put(100,70){\line(-1,-2){10}}
\put(90,50){\line(1,-2){10}} \put(100,30){\line(-2,-1){20}}
\put(80,20){\line(-1,-2){10}} \put(70,0){\line(-2,1){20}}
\put(50,10){\line(-2,-1){20}} \put(30,0){\line(-1,2){10}}
\put(20,20){\line(-2,1){20}} \put(0,30){\line(1,2){10}}
\put(10,50){\line(-1,2){10}} \put(0,70){\line(2,1){20}}
\put(20,80){\line(1,2){10}} \put(30,100){\line(2,-1){20}}

\put(50,70){\line(2,-1){10}} \put(60,65){\line(2,1){10}}
\put(70,70){\line(-1,-2){5}} \put(65,60){\line(1,-2){5}}
\put(70,50){\line(-1,-2){5}} \put(65,40){\line(1,-2){5}}
\put(70,30){\line(-2,1){10}} \put(60,35){\line(-2,-1){10}}

\put(50,70){\line(-2,-1){10}} \put(40,65){\line(-2,1){10}}
\put(30,70){\line(1,-2){5}} \put(35,60){\line(-1,-2){5}}
\put(30,50){\line(1,-2){5}} \put(35,40){\line(-1,-2){5}}
\put(30,30){\line(2,1){10}} \put(40,35){\line(2,-1){10}}

\put(57,80){\small\makebox(0,0)[l]{$f_1$}}
\put(78,63){\small\makebox(0,0)[l]{$f_2$}}
\put(27,63){\small\makebox(0,0)[r]{$f_{n-1}$}}
\put(43,80){\small\makebox(0,0)[r]{$f_n$}}

{\thinlines \put(33,50){\line(1,0){34}}
\put(34,52){\line(1,0){32}} \put(35,54){\line(1,0){30}}
\put(36,56){\line(1,0){28}} \put(37,58){\line(1,0){26}}
\put(38,60){\line(1,0){24}} \put(37,62){\line(1,0){26}}
\put(36,64){\line(1,0){28}} \put(34,66){\line(1,0){2}}
\put(44,66){\line(1,0){12}} \put(64,66){\line(1,0){2}}
\put(48,68){\line(1,0){4}}

\put(34,48){\line(1,0){32}} \put(35,46){\line(1,0){30}}
\put(36,44){\line(1,0){28}} \put(37,42){\line(1,0){26}}
\put(38,40){\line(1,0){24}} \put(37,38){\line(1,0){26}}
\put(36,36){\line(1,0){28}} \put(34,34){\line(1,0){2}}
\put(44,34){\line(1,0){12}} \put(64,34){\line(1,0){2}}
\put(48,32){\line(1,0){4}}}

\put(50,20){\makebox(0,0){$\cdots$}}

\end{picture}
\end{center}
where at least one $E$-vertex would be in the shaded area, and hence it would not be on the boundary of a compass disk associated with our disk D-graph, or we would have something like the graph in the dual picture with inverted arrows, where the same thing holds for at least one $W$-vertex. The $E$-vertex or $W$-vertex in question must be in the shaded area because all the inner faces $f_1,\ldots,f_n$ are bipolar.

For a cocycle $C$ of a disk D-graph $D$, a face of $D$ is
\emph{$C$-cocyclic}\index{C-cocyclic face of a disk D-graph@$C$-cocyclic face of a disk D-graph}
when an edge of $C$ belongs to the boundary of that face. We can
prove the following.

\prop{Lemma 7.3.5}{Every $C$-cocyclic inner face of a disk D-graph $D$ contains exactly two edges of the cocycle $C$ of $D$, one of which is in the north path and the other in the south path.}

\dkz If either in the north or in the south path we had more than one edge from $C$, then the componential graph $C_C(D)$ (see \S 1.3) would not be acyclic. The first edge from $C$ in the north or south path would connect $D_1$ with $D_2$, for $D_1$ and $D_2$ vertices of $C_C(D)$, while from the second edge from $C$ in that path we would have that it must connect $D_2$ with~$D_1$.

There cannot be a single edge from $C$ in a $C$-cocyclic face of $D$; otherwise, ($\ddag$) of \S 1.3 would not hold.\qed

For a cocycle $C$ of a disk D-graph $D$, and for $f$ and $g$ being $C$-cocyclic inner faces of $D$, let us write  $fP_Cg$,\index{PC@$P_C$} and say that $f$ is a \emph{$C$-predecessor}\index{C-predecessor@$C$-predecessor} of $g$, while $g$ is a \emph{$C$-successor}\index{C-successor@$C$-successor} of $f$, when there is an edge in $C$ that is in the south path of $f$ and in the north path of $g$. It is clear that if $fP_Cg$, then $f$ precedes $g$, according to the definition before Lemma 7.3.4. By Lemma 7.3.5, the relation $P_C$ is linear in the following sense: if $f_1P_Cg$ and $f_2P_Cg$, then $f_1=f_2$, and if $fP_Cg_1$ and $fP_Cg_2$, then $g_1=g_2$.

From Lemmata 7.3.4 and 7.3.5, and from $C_C(D)$ having exactly two vertices, we may conclude that if we have $C$-cocyclic inner faces in our disk D-graph $D$, then they make a list $f_1\ldots f_n$, with $n\geq 1$, such that $f_1$ has no $C$-predecessor, $f_n$ has no $C$-successor, and if $n\geq 2$, then for every $i$ in $\{1,\ldots,n\mn 1\}$ we have $fP_Cf_{i+1}$. The cocycle $C$ has a single edge iff there are no $C$-cocyclic inner faces of $D$, and our list is empty.

Out of such a non-empty list we make a list $L(C)$ of the edges of $C$ by starting with the edge in the north path of $f_1$, and by passing to the edge in the south path of $f_1$.
If $n\geq 2$, and we have reached the edge in the south path of $f_i$ for $i$ in $\{1,\ldots,n\mn 1\}$, then that edge is the edge in the north path of $f_{i+1}$, and we pass to the edge in the south path of $f_{i+1}$. We proceed in that manner until we reach the edge in the south path of $f_n$ (for an example, see the next picture). If $C$ has a single edge, then the list made of that edge is $L(C)$.

We will show next how to make out of a disk D-graph $D$ two disk D-graphs $D_W'$ and $D_E'$ closely related to the D-graphs $D_W$ and $D_E$ obtained by cutting $D$ through a cocycle $C$ (see \S 1.10). How we obtain $D_W'$ and $D_E'$ should be clear from the following picture, and the explanations we give:
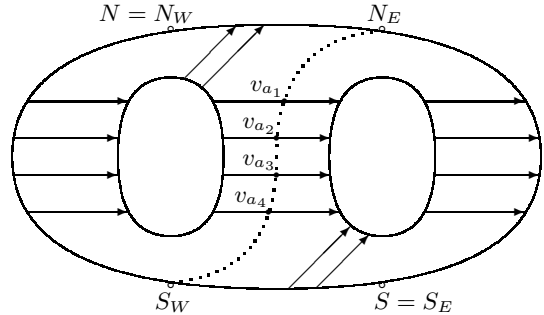
\begin{center}
\begin{picture}(200,120)(0,-10)

\qbezier(0,50)(0,0)(100,0) \qbezier(100,0)(200,0)(200,50)
\qbezier(200,50)(200,100)(100,100) \qbezier(100,100)(0,100)(0,50)

{\thicklines \qbezier[20](60,2.5)(100,5)(100,50)
\qbezier[20](140,97.5)(100,95)(100,50)}

\qbezier(40,50)(40,20)(60,20) \qbezier(60,20)(80,20)(80,50)
\qbezier(80,50)(80,80)(60,80) \qbezier(60,80)(40,80)(40,50)

\put(97,29){\circle*{2}} \put(103,71){\circle*{2}}
\put(100,43){\circle*{2}} \put(100,57){\circle*{2}}

\put(103,75){\small\makebox(0,0)[br]{$v_{a_1}$}}
\put(100,61){\small\makebox(0,0)[br]{$v_{a_2}$}}
\put(100,47){\small\makebox(0,0)[br]{$v_{a_3}$}}
\put(97,33){\small\makebox(0,0)[br]{$v_{a_4}$}}

\put(6,29){\vector(1,0){38}} \put(6,71){\vector(1,0){38}}
\put(.5,43){\vector(1,0){39.5}} \put(.5,57){\vector(1,0){39.5}}

\put(60,98.5){\circle{2}} \put(60,1.5){\circle{2}}

\put(51,101.5){\small\makebox(0,0)[b]{$N=N_W$}}
\put(61,-1.5){\small\makebox(0,0)[t]{$S_W$}}

\put(76,29){\vector(1,0){48}} \put(76,71){\vector(1,0){48}}
\put(80,43){\vector(1,0){40.5}} \put(80,57){\vector(1,0){40.5}}

\put(156,29){\vector(1,0){38}} \put(156,71){\vector(1,0){38}}
\put(160,43){\vector(1,0){39.5}} \put(160,57){\vector(1,0){39.5}}

\put(65,79.5){\vector(1,1){20}} \put(72,76.5){\vector(1,1){23.5}}


\qbezier(120,50)(120,20)(140,20) \qbezier(140,20)(160,20)(160,50)
\qbezier(160,50)(160,80)(140,80) \qbezier(140,80)(120,80)(120,50)

\put(140,98.5){\circle{2}} \put(140,1.5){\circle{2}}

\put(141,101.5){\small\makebox(0,0)[b]{$N_E$}}
\put(152,-1.5){\small\makebox(0,0)[t]{$S=S_E$}}

\put(115,.5){\vector(1,1){20}} \put(104.5,0){\vector(1,1){23.5}}

\end{picture}
\end{center}
Let $C=\{a_1,a_2,a_3,a_4\}$ be our cocycle, which in the picture is represented by the four edges on which we have chosen the points $v_{a_1}$, $v_{a_2}$, $v_{a_3}$ and $v_{a_4}$; these points are not end points. Here $L(C)$ is $a_1a_2a_3a_4$.

The points $N$ and $S$ are respectively the north pole and the south pole of the compass disk $\kappa$ associated with $D$, whose boundary is the outermost circle in the picture. The point $N$ becomes the north pole $N_W$ of the compass disk $\kappa_W$ associated with $D_W'$, and $N_E$ is a point on the boundary of $\kappa$ that may be joined with $v_{a_1}$ by a Jordan curve---a dotted line in our picture---which besides $v_{a_1}$ does not contain any point from ${\cal U}(D)$ (see \S 7.1). The point $N_E$ is the north pole of the compass disk $\kappa_E$ associated with $D_E'$. We connect analogously by Jordan curves, represented by dotted lines, $v_{a_1}$ with $v_{a_2}$, the point $v_{a_2}$ with $v_{a_3}$, the point $v_{a_3}$ with $v_{a_4}$, and finally $v_{a_4}$ with a point $S_W$ on the boundary of $\kappa$. This last point is the south pole of $\kappa_W$, while $S$, which is the south pole of $\kappa$, is also the south pole of $\kappa_E$. The boundary of $\kappa_W$ is made of the west side of the boundary of $\kappa$ from $N_E$ to $S_W$ together with the Jordan curve represented by the dotted line, which is the union of all the dotted lines introduced above. The boundary of $\kappa_E$ is made analogously with the east side. For $X$ being $W$ or $E$, the disk D-graph $D_X'$ is that part of $D$ within $\kappa_X$, with $v_{a_1}$, $v_{a_2}$, $v_{a_3}$ and $v_{a_4}$ as new $\bar{X}$-vertices.

Another possible situation is
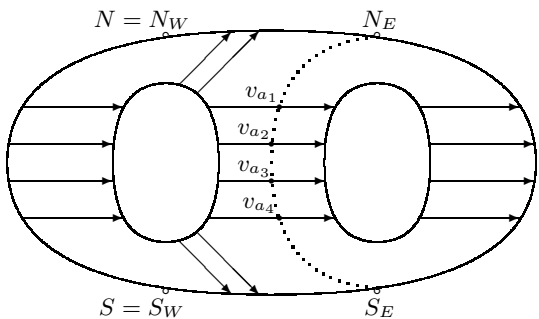
\begin{center}
\begin{picture}(200,120)(0,-10)

\qbezier(0,50)(0,0)(100,0) \qbezier(100,0)(200,0)(200,50)
\qbezier(200,50)(200,100)(100,100) \qbezier(100,100)(0,100)(0,50)

{\thicklines \qbezier[20](140,2.5)(100,5)(100,50)
\qbezier[20](140,97.5)(100,95)(100,50)}

\qbezier(40,50)(40,20)(60,20) \qbezier(60,20)(80,20)(80,50)
\qbezier(80,50)(80,80)(60,80) \qbezier(60,80)(40,80)(40,50)

\put(103,29){\circle*{2}} \put(103,71){\circle*{2}}
\put(100,43){\circle*{2}} \put(100,57){\circle*{2}}

\put(103,75){\small\makebox(0,0)[br]{$v_{a_1}$}}
\put(100,61){\small\makebox(0,0)[br]{$v_{a_2}$}}
\put(100,47){\small\makebox(0,0)[br]{$v_{a_3}$}}
\put(102,33){\small\makebox(0,0)[br]{$v_{a_4}$}}

\put(6,29){\vector(1,0){38}} \put(6,71){\vector(1,0){38}}
\put(.5,43){\vector(1,0){39.5}} \put(.5,57){\vector(1,0){39.5}}

\put(60,98.5){\circle{2}} \put(60,1.5){\circle{2}}

\put(51,101.5){\small\makebox(0,0)[b]{$N=N_W$}}
\put(51,-1.5){\small\makebox(0,0)[t]{$S=S_W$}}

\put(76,29){\vector(1,0){48}} \put(76,71){\vector(1,0){48}}
\put(80,43){\vector(1,0){40.5}} \put(80,57){\vector(1,0){40.5}}

\put(156,29){\vector(1,0){38}} \put(156,71){\vector(1,0){38}}
\put(160,43){\vector(1,0){39.5}} \put(160,57){\vector(1,0){39.5}}

\put(65,79.5){\vector(1,1){20}} \put(72,76.5){\vector(1,1){23.5}}
\put(65,20.5){\vector(1,-1){20}}
\put(72,23.5){\vector(1,-1){23.5}}


\qbezier(120,50)(120,20)(140,20) \qbezier(140,20)(160,20)(160,50)
\qbezier(160,50)(160,80)(140,80) \qbezier(140,80)(120,80)(120,50)

\put(140,98.5){\circle{2}} \put(140,1.5){\circle{2}}

\put(141,101.5){\small\makebox(0,0)[b]{$N_E$}}
\put(141,-1.5){\small\makebox(0,0)[t]{$S_E$}}

\end{picture}
\end{center}
and other possibilities are treated analogously.

The disk D-graph $D_W'$ differs from the D-graph $D_W$ obtained by cutting $D$ through $C$ by having, for every $i$ in $\{1,2,3,4\}$, the edge $a_i$ replaced by ``the west half'' of $a_i$; the vertices are the same. Note that with our choice of $v_{a_i}$, the D-graph $D_W$ is not a plane graph: $v_{a_i}$ is not an end point of $a_i$. Analogously, in $D_E'$ we have ``the east half'' of $a_i$ instead $a_i$, which is in~$D_E$.

By Lemma 7.3.2, the list $L_{\bar{X}}(D_X')$ is grounded in $D_X'$, and it is easy to conclude that the same list is grounded in $D_X$. It can be verified that the lists $L_E(D_W)$ and $L_W(D_E)$ are compatible. In our example, both lists have $v_{a_1}v_{a_2}v_{a_3}v_{a_4}$ as the common sublist, and it is easy to see that the requirements of compatibility are met. From that we may infer that $D_W$ and $D_E$ are P-compatible. Hence we have that every disk D-graph is a P$'''$-graph, and hence we have the following.

\prop{Proposition 7.3.6}{Every disk D-graph is a P-graph.}

\noindent It is easy to see that every graph isomorphic to a P-graph is a P-graph (as we remarked already at the end of \S 2.3), and so this proposition and Proposition 7.2.1 yield the following.

\prop{Proposition 7.3.7}{A graph is a P-graph iff it is isomorphic to a disk D-graph.}

As other corollaries, we have that a graph is a P-graph iff it is isomorphic to an eastward-growing disk D-graph, and that every disk D-graph is isomorphic to an eastward-growing disk D-graph. A survey of criteria for eastward growing in the plane (i.e.\ for upward planarity) and of related questions may be found in \cite{GT95}.

\section{\large\bf D1$'$-graphs}\label{7.4}
\markright{\S 7.4. \quad D1$'$-graphs}

Let a \emph{disk D1-graph}\index{disk D1-graph} be defined as a disk D-graph in \S 7.1 by substituting ``D1-graph'' for ``D-graph'' (for the notion of D1-graph see \S 6.5). (Note that the empty graph is trivially a disk D1-graph.)

It is easy to derive from Propositions 7.2.1, 7.3.6 and 7.3.7 that the assertions of these propositions hold when ``D1-graph'' is substituted for ``D-graph'' and ``P1-graph'' is substituted for ``P-graph''. (Related matters are considered in~\S 7.5.)

Let a \emph{D1$'$-graph}\index{D1prime-graph@D1$'$-graph} be a graph that is finite, acyclic, incidented and has a source and sink (see \S 1.2 and \S 7.3 for these notions). This definition differs from the definition of D1-graph by replacing the requirement of $W\mbox{\rm -}E$-functionality by the requirement of possessing a source and sink. While D1-graphs need not be weakly connected, D1$'$-graphs are always such.

Let a \emph{disk D1$'$-graph}\index{disk D1prime-graph@disk D1$'$-graph} be defined as a disk D-graph in \S 7.1 by substituting ``D1$'$-graph'' for ``D-graph'', as we did for the notion of disk D1-graph above.

A \emph{source-sink closure}\index{source-sink closure of a disk D1-graph} of a non-empty disk D1-graph is defined in the same manner as a source-sink closure of a disk D-graph (see \S 7.3, before Lemma 7.3.3); just substitute ``D1-graph'' for ``D-graph''. We can prove the following.

\prop{Proposition 7.4.1}{A graph is a source-sink closure of a disk D1-graph iff it is a disk D1$'$-graph.}

\dkz The direction from left to right is obtained immediately by a stretching of the compass disk towards West and East. The other direction is also straightforward. After removing the source and sink together with their small neighbourhoods, and replacing them by new distinct end points on the remainder of the incident edges so as to ensure $W\mbox{\rm -}E$-functionality, we shrink the compass disk as in the following picture:
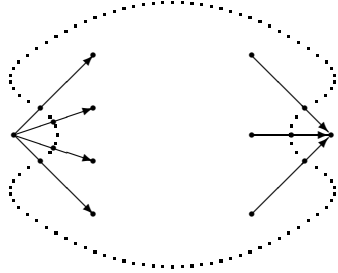
\begin{center}
\begin{picture}(120,100)(0,5)

{\thicklines \qbezier[6](10,60)(23,50)(10,40)
\qbezier[6](110,60)(99,50)(110,40)
\qbezier[30](5,80)(60,120)(115,80)
\qbezier[30](5,20)(60,-20)(115,20)
\qbezier[8](5,80)(-10,70)(10,60)
\qbezier[8](115,80)(130,70)(110,60)
\qbezier[8](5,20)(-10,30)(10,40)
\qbezier[8](115,20)(130,30)(110,40)}

\put(0,50){\circle*{2}} \put(120,50){\circle*{2}}
\put(30,20){\circle*{2}} \put(30,40){\circle*{2}}
\put(30,60){\circle*{2}} \put(30,80){\circle*{2}}
\put(90,20){\circle*{2}} \put(90,50){\circle*{2}}
\put(90,80){\circle*{2}} \put(10,40){\circle*{2}}
\put(15,45){\circle*{2}} \put(15,55){\circle*{2}}
\put(10,60){\circle*{2}} \put(110,40){\circle*{2}}
\put(105,50){\circle*{2}} \put(110,60){\circle*{2}}

\put(0,50){\vector(1,1){30}} \put(0,50){\vector(3,1){30}}
\put(0,50){\vector(3,-1){30}} \put(0,50){\vector(1,-1){30}}

\put(90,50){\vector(1,0){29}} \put(90,80){\vector(1,-1){29}}
\put(90,20){\vector(1,1){29}}

\end{picture}
\end{center}

\qed

\vspace{-2ex}

\prop{Proposition 7.4.2}{A graph is a disk D1$'$-graph iff it is a pasting scheme.}

\dkz From left to right it is enough to remark that the source and sink of a D1$'$-graph must be on the boundary of the outer face because they are on the boundary of an associated compass disk. From finiteness, acyclicity and incidentedness we infer that for every vertex $v$ there is a path from the source to $v$ and a path from $v$ to the sink. The other direction is even easier to prove.\qed

With Propositions 7.4.1 and  7.4.2 we have obtained two other alternative definitions of the notion of pasting scheme.

\section{\large\bf Realizing D1$'$-graphs}\label{7.5}
\markright{\S 7.5. \quad Realizing D1$'$-graphs}

We will consider in this section the question when a D1$'$-graph is isomorphic to a disk D1$'$-graph, i.e., when it is realizable in the plane as a pasting scheme.

We define first for a D1-graph $D$ the D1$'$-graph $D'$, which is the \emph{abstract source-sink closure}\index{abstract source-sink closure of a D1-graph} of $D$: the graph $D'$ differs from $D$ by replacing all its $W$-vertices by a new vertex $u_W$, all its $E$-vertices by a new vertex $u_E$, and by assuming for every $X$-edge $a$, where $X$ is $W$ or $E$, that $X(a)=u_X$.

Next we define for a D1$'$-graph $D$ the D1-graph $D^-$, which is the \emph{D1-interior}\index{D1-interior of a D1prime-graph@D1-interior of a D1$'$-graph}\index{D minus@$D^-$} of $D$: the graph $D^-$ is obtained from $D$ by rejecting its source and sink, and by assuming that for every $X$-edge $a$ of $D$ the vertex $X(a)$ is a new vertex $v_a^X$. We have that if $a_1\neq a_2$, then $v_{a_1}^X\neq v_{a_2}^X$.

It is clear that for every D1-graph $D$ we have that $D'^-$ is isomorphic to $D$, and that for every D1$'$-graph $D$ we have that ${D^-}'$ is isomorphic to~$D$.

We say that a D1$'$-graph $D$ is \emph{disk realizable}\index{disk realizable D1prime-graph@disk realizable D1$'$-graph} when there is a disk D1$'$-graph isomorphic to $D$. We have two more analogous definitions obtained by substituting ``D1-graph'' and ``D-graph'' respectively for ``D1$'$-graph''\index{disk realizable D1-graph}\index{disk realizable D-graph}. It is clear that we have the following two lemmata.

\prop{Lemma 7.5.1}{If the D1$'$-graph $D$ is disk realizable, then the D1-graph $D^-$ is disk realizable.}

\vspace{-2ex}

\prop{Lemma 7.5.2}{If the D1-graph $D$ is disk realizable, then the D1$'$-graph $D'$ is disk realizable.}

From these lemmata we infer the following proposition.

\prop{Proposition 7.5.3}{For every D1$'$-graph $D$ we have that $D$ is disk realizable iff the D1-graph $D^-$ is disk realizable.}

We can also prove the following.

\prop{Proposition 7.5.4}{For every D1-graph $D$ we have that $D$ is disk realizable iff every component of $D$ is disk realizable.}

\noindent Every component of a D1-graph $D$ is either a D-graph or a straight single-edge graph (see \S 6.5). Straight single-edge graphs are of course always disk realizable (as disk D1-graphs of a particularly simple kind), and a D-graph is disk realizable iff it is a P-graph, by Proposition 7.3.7. So we have reduced the question of disk realizability of D1$'$-graphs to the question whether a graph is a P-graph, and to answer this last question the notion of P$'''$-graph suggests the most viable criterion, among those we have considered.

Another, easier, way of reducing the question of disk realizability of D1$'$-graphs to the notion of P-graph is to pass from a D1$'$-graph $D$ to a D-graph $D^{\dag}$, which is obtained from $D$ by adding a new source $s^{\dag}$ and a new sink $t^{\dag}$, and two new edges, from $s^{\dag}$ to the source $s$ of $D$, and from the sink $t$ of $D$ to $t^{\dag}$. It is easy to see that the D1$'$-graph $D$ is disk realizable iff the D-graph $D^{\dag}$ is disk realizable, i.e., iff $D^{\dag}$ is a P-graph.

It is clear that the added edges from $s^{\dag}$ to $s$ and from $t$ to $t^{\dag}$ do not play an essential role here. We have them only to conform to our definition of D-graph. With a more general notion, we could dispense with this addition.

This more general notion could be based on a wider class of basic D-graphs, which could also look like
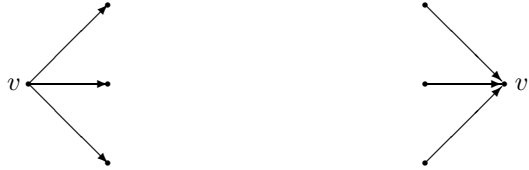
\begin{center}
\begin{picture}(180,70)(0,-5)

\put(0,30){\circle*{2}} \put(30,30){\circle*{2}}
\put(30,0){\circle*{2}} \put(30,60){\circle*{2}}

\put(150,30){\circle*{2}} \put(150,0){\circle*{2}}
\put(150,60){\circle*{2}} \put(180,30){\circle*{2}}

\put(0,30){\vector(1,0){30}} \put(0,30){\vector(1,1){30}}
\put(0,30){\vector(1,-1){30}}

\put(150,30){\vector(1,0){29}} \put(150,60){\vector(1,-1){29}}
\put(150,0){\vector(1,1){29}}

\put(-3,30){\makebox(0,0)[r]{$v$}}
\put(184,30){\makebox(0,0)[l]{$v$}}

\end{picture}
\end{center}
with the vertex $v$ in both pictures considered as an \emph{inner} vertex. This presupposes a new notion of inner vertex. Out of these enlarged basic graphs, we would obtain the other graphs in our enlarged family of P-graphs with the operation of juncture.

One could envisage further generalizations, and investigate juncture in these wider contexts as an operation for building graphs that correspond to diagrams of ordinary categories that are not only commuting diagrams of arrows (the diagrams of the last two pictures are not commuting diagrams). Juncture, which consists in identifying the tokens of the same edges in different diagrams, would in this perspective replace ordinary composition of arrows in categories, which consists in such an identifying of tokens of the same object.

This was not our point of view in this work. For us, juncture was an operation on graphs that correspond to diagrams of 2-cells, or from the point of view of ordinary categories, just to commuting diagrams of arrows. The operation of juncture was applied not to graphs that correspond directly to the diagrams of 2-cells, but to graphs that are a kind of dual of these graphs. In \S 7.6, the last section of this work, we deal with this duality.

\section{\large\bf Duality}\label{7.6}
\markright{\S 7.6. \quad Duality}

We need the following notions for the definition of dual of a disk D1$'$-graph. The edges of a disk D1$'$-graph $D$ may be of four kinds.
\begin{itemize}
\item[(1)]An edge $a$ may separate two inner faces $f$ and $g$ of $D$, in which case when $f$ precedes $g$ (see \S 7.3 before Lemma 7.3.4) we say that $a$ is an \emph{interior edge on the way from $f$ to~$g$}\index{interior edge}.

\vspace{-1ex}

\item[(2)]An edge $a$ may separate the outer face of $D$ from an inner face $f$ so that $a$ is in the north path of $f$; in that case we say that $a$ is a \emph{northern outer edge of~$f$}\index{northern outer edge}\index{outer edge}.

    \vspace{-1ex}

\item[(3)]An edge $a$ may separate an inner face $f$ from the outer face of $D$ so that $a$ is in the south path of $f$; in that case we say that $a$ is a \emph{southern outer edge of~$f$}\index{southern outer edge}.

    \vspace{-1ex}

\item[(4)]An edge $a$ may be such that there is no inner face with $a$ belonging to its  boundary; in that case we say that $a$ is a \emph{totally outer edge}\index{totally outer edge}.
\end{itemize}

When from a compass disk $\kappa$ associated with a disk D1$'$-graph $D$ we reject all the points of ${\cal U}(D)$ (see \S 7.1) and all the inner faces of $D$, we are left with two disjoint sets of points of $\kappa$, which we call $\kappa_N$,\index{kappaN@$\kappa_N$} in which we find the north pole of $\kappa$, and $\kappa_S$,\index{kappaS@$\kappa_S$} in which we find the south pole of $\kappa$. We can then pass to our definition of dual.

A \emph{dual}\index{dual of a disk D1prime-graph@dual of a disk D1$'$-graph} of a disk D1$'$-graph $D$ is a plane graph $D^*$ obtained as follows. For every inner face $f$ of $D$ a point $f^*$ from $f$ will be a vertex of $D^*$, and, moreover, we have as additional vertices of $D^*$ the north pole $N$ and the south pole $S$ of a compass disk $\kappa$ associated with $D$. For every interior edge $a$ on the way from the inner face $f$ to the inner face $g$ we have as an edge of $D^*$ a Jordan curve $a^*$ joining $f^*$ with $g^*$, such that $a^*\subseteq f\cup g\cup a$; we take that $W(a^*)=f^*$ and $E(a^*)=g^*$. For every northern outer edge $a$ of an inner face $f$ we have as an edge of $D^*$ a Jordan curve $a^*$ joining the vertex $N$ with $f^*$, such that $a^*\subseteq\kappa_N\cup f\cup a$; we take that $W(a^*)=N$ and $E(a^*)=f^*$. For every southern outer edge $a$ of an inner face $f$ we have as an edge of $D^*$ a Jordan curve $a^*$ joining $f^*$ with the vertex $S$, such that $a^*\subseteq f\cup\kappa_S\cup a$; we take that $W(a^*)=f^*$ and $E(a^*)=S$. For every totally outer edge $a$ we have as an edge of $D^*$ a Jordan curve joining the vertex $N$ with the vertex $S$, such that $a^*\subseteq\kappa_N\cup\kappa_S\cup a$; we take that $W(a^*)=N$ and $E(a^*)=S$. This concludes our definition of~$D^*$.

Note that we have required that $D^*$ be a plane graph. So we must ensure that the Jordan curves that make its edges intersect only in the vertices of $D^*$ that are the end points of these edges, as in condition (2) of the definition of plane graph (see~\S 7.1).

We can now prove the following.

\prop{Proposition 7.6.1}{For every disk D1$'$-graph $D$ the graph $D^*$ is a disk D1$'$-graph.}

\dkz We check first that $D^*$ is a D1$'$-graph. It is clear that it is finite, and, by a lemma for D1$'$-graphs analogous to Lemma 7.3.4, we obtain acyclicity. It is clear that $D^*$ is incidented, and finally the vertex $N$ is the source, while the vertex $S$ is the sink of $D^*$. The compass disk associated with $D$ will also be associated with $D^*$, with the new north pole being the sink of $D$, and the new south pole being the source of~$D$.\qed

For a graph $G$, which is $W,E\!:A\rightarrow V$, let $G^{\rm op}$\index{op@$^{\rm op}$} be the
graph $W^{\rm op},E^{\rm op}\!:A\rightarrow V$ such that $W^{\rm op}=E$ and $E^{\rm op}=W$. It is possible to show that for a D1$'$-graph $D$ the graph $D^{**}$ is isomorphic to $D^{\rm op}$, but we will not go into the proof of that. (Analogous facts in graph theory are usually skipped over, as in \cite{B73}, Section 2.2, or left as exercises, as in \cite{BM76},  Exercise 9.2.4, Section~9.2.)

Pasting schemes may be combined one with another with two operations that correspond to vertical and horizontal composition in 2-categories (cf.\ \S 6.6). Our goal was to study the operation definable in terms of these two operations that consists in gluing two pasting schemes along a common path on the boundaries, as in the first picture of \S 1.1. When applied to 2-cells, we called this operation \emph{juncture} in \S 1.1, but when applied to pasting schemes, we better find now another name for it, not to create confusion. We could  call it \emph{gluing}\index{gluing}.

Instead of dealing with pasting schemes, i.e.\ disk D1$'$-graphs, we pass to modified duals of these graphs. For every disk D1$'$-graph $D$ we take the D1-interior $D^{*-}$ of $D^*$ (see \S 7.5), which is easily seen to be isomorphic to a disk D1-graph (by proceeding as in  the proof of Proposition 7.4.1; for examples of passing from $D$ to $D^{*-}$ see the end of \S 6.7). For the modified duals $D^{*-}$, gluing becomes juncture, and the passing from $D$ to $D^{*-}$ was made to obtain the operation of juncture, more manageable in wider classes of graph, which need not be plane. We did not stop at the dual $D^*$, which is a D1$'$-graph, but passed further to the D1-graph $D^{*-}$, because the analogue of juncture for D1$'$-graphs would be less manageable. The analogue of juncture for two D1$'$-graphs $D_W$ and $D_E$ is best defined as corresponding to $(D_W^-\Box D_E^-)'$ (with $'$ being the abstract source-sink closure of~\S 7.5).

We forgot about disk realizability, to obtain more general notions, and we ended up with the notion of D1-graph and the essential ingredient of that notion, which is the notion of D-graph. Neither of these two notions has a natural dual correlate in the world of pasting schemes. (These would be, roughly, pasting schemes with vertices removed.) For the notion of D1-graph we may then ask when it is disk realizable, and this disk realizability reduces to the disk realizability of D-graphs (see Proposition 7.5.4). This last question is answered by Proposition 7.3.7, and our work serves to explain the notion of P-graph of that proposition.

\begin{theindex}\label{Index}
\pagestyle{myheadings}\markboth{Index}{Index}

  \item $\cal A$, set of all edges, 87
  \item abstract source-sink closure of a D1-graph, 118
  \item acyclic graph, 8
  \item $\alpha$, 101
  \item (Ass~$\cirk$), 90
  \item (Ass~$\otimes$), 90
  \item associated compass disk, 109
  \item (Ass~1), 15
  \item (Ass~2.1), 15
  \item (Ass~2.2), 15
  \item atomic P1-term, 88
  \item atomic P2-term, 90

  \indexspace

  \item basic D-edge-graph, 94
  \item basic D-graph, 8
  \item basic D-term, 14
  \item basic M-edge-graph, 101
  \item basic M-graph, 100
  \item basic P-term, 86
  \item Bf, 61
  \item bipolar face of a plane graph, 112
  \item border vertex, 19
  \item boundary of a face of a plane graph, 112

  \indexspace

  \item $C$-cocyclic face of a disk D-graph, 114
  \item $C$-predecessor, 114
  \item $C$-successor, 114
  \item CCP, 61
  \item chain, 23
  \item closest common pivot, 61
  \item cocycle, 10
  \item codomain, 99
  \item compass disk, 108
  \item compass disk associated with a disk graph, 109
  \item compatible chains, 23
  \item compatible lists, 22
  \item completeness of S$\Box$, 21
  \item completeness of S$\Box_P$, 87
  \item completeness of S1, 97
  \item completeness of S2, 103
  \item component of a graph, 9
  \item componential extreme, 17
  \item componential graph, 9
  \item connect, 7
  \item construction, 24
  \item construction of a P$'$-graph, 24
  \item core, 103
  \item corolla of a vertex, 59
  \item $C_S(D)$, 9
  \item cutset, 9
  \item cutting a D-graph through a cocycle, 27
  \item cutvertex, 17
  \item cycle, 8

  \indexspace

  \item $D^-$, 118
  \item D-edge-graph, 13
  \item D-graph, 8
  \item D-term, 14
  \item D1-edge-graph, 94
  \item D1-graph, 94
  \item D1-interior of a D1$'$-graph, 118
  \item D1$'$-graph, 117
  \item developed term, 102
  \item directed graph, 9
  \item disjoint lists, 22
  \item disk, 108
  \item disk D-graph, 109
  \item disk D1-graph, 117
  \item disk D1$'$-graph, 117
  \item disk planarity, 3
  \item disk realizable D-graph, 119
  \item disk realizable D1-graph, 119
  \item disk realizable D1$'$-graph, 118
  \item distance between members of a list, 71
  \item distinguished graph, 6
  \item domain, 99
  \item dual of a disk D1$'$-graph, 121

  \indexspace

  \item $E$, 5
  \item $\bar{E}$, 5
  \item $E$-border vertex, 19
  \item $E$-edge, 7
  \item $E$-edge of edge-graph, 13
  \item $E$-extreme, 21
  \item $E$-functional edge, 7
  \item $E$-meridian, 108
  \item $E$-peripheral vertex, 26
  \item $E$-vertex, 6
  \item eastward-growing plane graph, 109
  \item $E(D)$, 24
  \item edge, 5, 12
  \item edge type of D-term, 14
  \item edge-graph, 12
  \item edge-graph morphism, 12
  \item empty edge-graph, 12
  \item empty graph, 6
  \item empty list, 22
  \item $\eta$ interpretation function, 95
  \item $\eta^*$ interpretation function, 95
  \item $E(v)$, 19
  \item extreme of D-graph, 21

  \indexspace

  \item face of a plane graph, 112
  \item final member of a list, 29
  \item finite graph, 6

  \indexspace

  \item gluing, 122
  \item graph, 5
  \item graph isomorphism, 6
  \item graph morphism, 5
  \item graph realizable in the plane, 108
  \item grounded list, 25

  \indexspace

  \item incident, 5
  \item incidented graph, 8
  \item initial member of a list, 29
  \item inner edge, 7
  \item inner face of a plane graph, 112
  \item inner vertex, 7
  \item inner vertex of componential graph, 17
  \item interior edge, 120
  \item interlaced members of lists, 29
  \item intersecting semipaths, 25
  \item $\iota$ interpretation function, 16
  \item isomorphism, 6

  \indexspace

  \item juncture, 1, 10

  \indexspace

  \item $\kappa_N$, 121
  \item $\kappa_S$, 121

  \indexspace

  \item $L_E$, 24
  \item ${\cal L}_E$, 86
  \item leaf, 68
  \item list, 22
  \item list of, 22
  \item loop, 8
  \item $L_W$, 24
  \item ${\cal L}_W$, 86
  \item $L_X$, 24
  \item ${\cal L}_X$, 86

  \indexspace

  \item $m$, 26
  \item M-edge-graph, 101
  \item M-graph, 100
  \item mate, 26
  \item meridian, 108
  \item $\mu$ interpretation function, 101
  \item $\mu^*$ interpretation function, 102

  \indexspace

  \item $n$-valent, 19
  \item neighbours in a list, 29
  \item non-empty graph, 6
  \item north path of a bipolar face, 112
  \item north pole, 108
  \item northern outer edge, 120

  \indexspace

  \item (\mj 2L), 88
  \item (\mj 2L$\Phi$), 88
  \item (\mj 2L$\Psi$), 89
  \item (\mj 2R), 88
  \item (\mj 2R$\Phi$), 88
  \item (\mj 2R$\Psi$), 89
  \item (\mj $\cirk$), 90
  \item ($\mj\,\otimes$), 90
  \item (\mj 1), 88
  \item $^{\rm op}$, 122
  \item ordinary graph, 18
  \item outer edge, 120
  \item outer face of a plane graph, 112
  \item outer vertex, 7
  \item outer vertex of componential graph, 17

  \indexspace

  \item P-compatible D-graphs, 25
  \item P-edge-graph, 98
  \item P-graph, 3, 23
  \item P-move, 61
  \item P-term, 86
  \item P1-edge-graph, 98
  \item P1-graph, 98
  \item P$'$-graph, 24
  \item P1-term, 88
  \item P$''$-graph, 25
  \item P2-term, 90
  \item P$'''$-graph, 27
  \item parallel lists, 31
  \item parity of members lists, 29
  \item pasting scheme, 112
  \item path, 8
  \item $P_C$, 114
  \item peripheral vertex, 26
  \item petal, 60
  \item pivot, 47
  \item planar graph, 108
  \item plane graph, 108
  \item point set of a plane graph, 108
  \item pole, 108
  \item precedes, for faces, 113
  \item $\psi_X$, 25
  \item $\psi^b_X$, 26

  \indexspace

  \item realizable in the plane, 108
  \item realization of a graph in the plane, 108
  \item removal of a vertex of componential graph, 17
  \item removal of edges, 9
  \item $\rho$, 95
  \item root graph, 24
  \item root list, 24

  \indexspace

  \item $^s$, set of members of a list, 22
  \item S1, 88
  \item S2, 90
  \item S$\Box$, 15
  \item S$\Box_P$, 86
  \item semicycle, 8
  \item semipath, 8
  \item semiwalk, 7
  \item sequential type of P-term, 86
  \item Sf, 61
  \item single-vertex graph, 12
  \item sink of a graph, 112
  \item $S(M)$, 99
  \item soundness of S$\Box$, 17
  \item soundness of S2, 102
  \item source of a graph, 112
  \item source-sink closure of a disk D-graph, 113
  \item source-sink closure of a disk D1-graph, 117
  \item south path of a bipolar face, 112
  \item south pole, 108
  \item southern outer edge, 120
  \item straight single-edge edge-graph, 101
  \item straight single-edge graph, 94, 100
  \item strict cutset, 10
  \item subgraph, 9
  \item subsemipath, 47
  \item subterm, 19
  \item $S(v)$, 68

  \indexspace

  \item ($\otimes\cirk$), 91
  \item ($\otimes\,$\mj), 91
  \item topological disk, 108
  \item totally outer edge, 121
  \item Tr, 61
  \item tree, set-theoretic, 67
  \item trivial semiwalk, 7

  \indexspace

  \item ${\cal U}(D)$, 108
  \item ${\cal U}(G)$, 108
  \item unified list, 22
  \item unit D1-edge-graph, 95
  \item unit term, 87
  \item upward planarity, 109

  \indexspace

  \item vertex, 5

  \indexspace

  \item $W$, 5
  \item $\bar{W}$, 5
  \item $W$-border vertex, 19
  \item $W\mbox{\rm -}E$-functional graph, 7
  \item $W$-edge, 7
  \item $W$-edge of edge-graph, 13
  \item $W$-extreme, 21
  \item $W$-functional edge, 7
  \item $W$-meridian, 108
  \item $W$-peripheral vertex, 26
  \item $W$-vertex, 6
  \item walk, 7
  \item $W(D)$, 24
  \item weakly connected graph, 8
  \item $W(v)$, 19

  \indexspace

  \item $X$, 5
  \item $\bar{X}$, 5
  \item $X$-border vertex, 19
  \item $X$-edge, 7
  \item $X$-edge of edge-graph, 13
  \item $X$-extreme, 21
  \item $X$-functional edge, 7
  \item $X$-peripheral vertex, 26
  \item $X$-vertex, 6
  \item $X(D)$, 24
  \item $X_e$, 16
  \item $X(v)$, 19

\end{theindex}


\begin{thebibliography}{99}
\bibitem[\S...]{} The sections where a reference is mentioned are
listed at the end in \mbox{\hspace{5pt}square} brackets.
\end{thebibliography}
\end{document}